\documentclass[a4paper, 12pt,oneside,reqno]{amsart}
\usepackage[a4paper]{geometry}
\usepackage{tabularx}
\usepackage{graphicx}
\usepackage[matrix,arrow,curve,cmtip]{xy}
\usepackage{amsfonts}
\usepackage{hyperref}
\usepackage{txfonts,xcolor,ulem,comment}
\usepackage{amsmath,amssymb,amsthm}

\def\op{\textrm{op}}
\def\id{\textrm{id}}
\def\As{\mathrm{As}}

\numberwithin{equation}{section}

\newtheorem{Theorem}{Theorem}
\newtheorem{theorem}{Theorem}[section]
\newtheorem{lemma}[theorem]{Lemma}

\newtheorem{corollary}[theorem]{Corollary}

\theoremstyle{definition}
\newtheorem{definition}[theorem]{Definition}
\newtheorem{example}[theorem]{Example}

\title{Gr\"{o}bner---Shirshov bases of Rota---Baxter algebra of weight~$\lambda$ with spectrum lying in $\{0,-\lambda\}$}
\author{H. Alhussein$^{1),2),3)}$}

\address{$^{1)}$Novosibirsk State University, Novosibirsk,  Russia.}
\address{$^{2)}$Siberian State University of Telecommunication and Informatics, Novosibirsk, Russia.}
\address{$^{3)}$Novosibirsk State University of Economics and Management, Novosibirsk,  Russia.}

\begin{document}

\maketitle

\begin{abstract}
It is known that if $A$ is a finite-dimensional unital algebra equipped with a Rota---Baxter operator~$R$ of weight~$\lambda$, then spectrum of $R$ is a subset of~$\{0,-\lambda\}$. We are interested in finding all consequences of the Rota---Baxter relation and the relation of the form~$R^k(R+\lambda\id)^l = 0$.
In 2024, H.~Qiu, S. Zheng, Y. Dan solved this problem for $k = l = 1$ and $\lambda\neq0$.
We find a~Gr\"obner---Shirshov basis of the ideal generated by these two relations—the Rota--Baxter relation and the cancellation relation $R(R + \lambda\id)$—in the general case.
\end{abstract}

\section{Introduction}
The theory of Gröbner bases, introduced by B. Buchberger in his 1965 thesis \cite{Buchberger65}, and its noncommutative generalization developed independently by A. Shirshov \cite{Shirshov62}, has become one of the most fundamental tools in computational algebra. These methods provide algorithmic solutions to fundamental problems such as the ideal membership problem, computation of normal forms, and the study of algebraic relations in polynomial and free algebras. The extension of this theory to various classes of operator algebras has been a subject of intense research in recent decades, particularly in the context of Rota---Baxter algebras.

Rota---Baxter algebras were defined by G. Baxter~\cite{Baxter60} in 1960, they are associative algebras equipped with a linear operator~$R$ satisfying the identity 
\begin{equation} \label{RB}
R(x)R(y) = R(xR(y)) + R(R(x)y) + \lambda R(xy),
\end{equation}
where $\lambda$ is a fixed scalar from the ground field.  

These algebras have emerged as important structures in both pure and applied mathematics, with applications ranging from combinatorics and number theory to quantum field theory and operad theory \cite{ConnesKreimer00,Guo12}. 

When a finite-dimensional associative unital algebra~$A$ is equipped with a Rota---Baxter operator~$R$ of weight~$\lambda$, then of spectrum of $R$ is a subset of~$\{0,-\lambda\}$~\cite{Miller66,Miller69,Spectrum}.
The same result holds true~\cite{Spectrum} over a field of characteristic zero, when $A$ is algebraic, it means that each of its finite subsets generates a~finite-dimensional subalgebra. 
There are many examples of nilpotent RB-operators of weight zero defined on infinite-dimensional algebras.

Thus, it is important to find all consequences of the Rota---Baxter relation and the relation of the form~$R^k(R+\lambda\id)^l = 0$.
The current work is devoted to this problem.
Formally speaking, we find the Gr\"{o}bner---Shirshov basis of the ideal generated by~\eqref{RB} and the relation~$R^k(R+\lambda\id)^l = 0$ in the free Rota---Baxter algebra.
In~\cite{Qiu}, this problem was solved for~$k = l = 1$ and $\lambda\neq0$. 
In this case, it occurs that these two relations already form a~Gr\"{o}bner---Shirshov basis
in the corresponding order. 



In 2010, L. Bokut, Y. Chen and J. Qiu~\cite{Bokut} established the foundation for Gr\"{o}bner---Shirshov bases in free associative Rota---Baxter algebras. 
In~2013, L. Guo, W. Sit, and R. Zhang~\cite{Guo2013} suggested another order to develop Gr\"{o}bner---Shirshov bases technique in free associative Rota---Baxter algebras. 
We will apply the order from the last work.

In~\S2, we recall the foundations of the Gr\"{o}bner---Shirshov bases for associative RB-algebras. In~\S3, we find the Gr\"{o}bner---Shirshov base for the ideal generated by~\eqref{RB} and~the relation~$R^k(R+\lambda\id)^l = 0$ in the free associative RB-algebra $R\As\langle X\rangle$. 
In Example~1, we consider the case~$\lambda = 0$, and in Example~2 we write down the 
Gr\"{o}bner---Shirshov base for the ideal generated by~\eqref{RB} and~the relation~$R^2(R+\lambda\id) = 0$.

Throughout the work, we assume that the characteristic of the ground field~$\Bbbk$ equals zero.

\section{Gr\"{o}bner---Shirshov bases for associative RB-algebras}

By $R\As\langle X\rangle$ we denote the free associative algebra generated by a set $X$ with a linear map $R$ in the signature.
One can construct a linear basis of $R\As\langle X\rangle$ (see, e.g., \cite{EbrahimiGuo08})
by induction. At first, all elements from $S(X)$, the free semigroup generated by~$X$,
lie in the basis. At second, if we have basic elements $a_1,a_2,\ldots,a_k$, $k\geq1$,
then the word $w_1R(a_1)w_2\ldots w_kR(a_k)w_{k+1}$ lies in the basis
of~$R\As\langle X\rangle$.
Here $w_2,\ldots,w_k\in S(X)$ and $w_1,w_{k+1}\in S(X)\cup \{\emptyset\}$,
where $\emptyset$ denotes the empty word.
Let us denote the basis obtained as $RS(X)$.
Given a word~$u$ from $RS(X)$, the number of appearances of the symbol~$R$
in~$u$ is denoted by $\deg_R(u)$, the $R$-degree of~$u$.
We call an element from $RS(X)$ of the form $R(w)$ as $R$-letter.
By~$X_\infty$ we denote the union of $X$ and the set of all $R$-letters.
Given $u\in RS(X)$, define $\deg u$ (degree of $u$) as the length of $u$
in the alphabet~$X_\infty$.

Suppose that $X$ is a well-ordered set with respect to $<$.
Let us introduce by induction the deg-lex order on $S(X)$.
At first, we compare two words $u$ and $v$ by the length:
$u < v$ if $|u|<|v|$. At second, when $|u| = |v|$,
$u = x_i u'$, $v = x_j v'$, $x_i, x_j\in X$,
we have $u < v$ if either $x_i < x_j$ or $x_i = x_j$, $u'< v'$.
We compare two words $u$ and $v$ from $RS(X)$
by $R$-degree: $u<v$ if $\deg_R(u)<\deg_R(v)$.
If $\deg_R(u) = \deg_R(v)$, we compare $u$ and $v$ in deg-lex order as words
in the alphabet $X_\infty$. Here we define each $x$ from $X$ to
be less than all $R$-letters and $R(a)<R(b)$ if and only if $a<b$.

Let $*$ be a symbol not containing in $X$.
By a $*$-bracketed word on $X$, we mean a basic word
from $R\As\langle X\cup\{*\}\rangle$ with exactly one occurrence of $*$.
The set of all $*$-bracketed words on $X$ is denoted by $RS^*(X)$.
For $q\in RS^*(X)$ and $u\in R\As\langle X\rangle$, we define
$q|_u$ as the bracketed word obtained by replacing the letter $*$ in $q$ by $u$.

The order defined above is monomial: $u < v$ implies that $q|_u < q|_v$ for all $u,v\in RS(X)$ and $q\in RS^*(X)$.
Given $f\in R\As\langle X\rangle$, by $\bar{f}$ we mean the leading word in~$f$.
We call $f$ monic if the coefficient of $\bar{f}$ in $f$ equals~1.

\begin{definition}[\cite{Guo2013}]
Let $f,g\in R\As\langle X\rangle$.
If there exist $\mu,\nu,w\in RS(X)$ such that
$w = \bar{f} \mu = \nu \bar{g}$ with $\deg w<\deg(\bar{f})+\deg(\bar{g})$,
then we define $(f,g)_w$ as $f\mu - \nu g$ and call it
the composition of intersection of $f$ and $g$ with respect to $(\mu,\nu)$.

If there exist $q\in RS^*(X)$ and $w\in RS(X)$
such that $w = \bar{f} = q|_{\bar{g}}$, then we define
$(f,g)^q_w$ as $f-q|_g$ and call it the composition of inclusion
of $f$ and $g$ with respect to $q$.
\end{definition}
\begin{definition}[\cite{Guo2013}]
Let $S$ be a subset of monic elements from $R\As\langle X\rangle$
and $w\in RS(X)$.

(1) For $u,v\in R\As\langle X\rangle$, we call $u$ and $v$
congruent modulo $(S, w)$ and denote this by $u \equiv v \mod (S, w)$
if $u - v = \sum c_i q_i|_{s_i}$
with $c_i \in \Bbbk$, $q_i\in RS^*(X)$, $s_i\in S$ and $q_i|_{\overline{s_i}} < w$.

(2) For $f,g\in R\As\langle X\rangle$ and suitable
$w,\mu,\nu$ or $q$ that give a composition of intersection
$(f,g)_w$ or a~composition of inclusion $(f,g)^q_w$,
the composition is called trivial modulo $(S, w)$ if
$(f,g)_w$ or $(f,g)^q_w \equiv 0 \mod (S, w)$.

(3) The set $S\subset R\As\langle X\rangle$ is
called a Gr\"{o}bner---Shirshov basis if, for all $f,g\in S$,
all compositions of intersection $(f,g)_w$
and all compositions of inclusion $(f,g)^q_w$ are trivial modulo $(S, w)$.

Now, we formulate what is known as composition-diamond lemma for associative operated algebras. 
\end{definition}
\begin{corollary}[\cite{Guo2013}] Let $I$ be an  ideal of $R\As\langle X\rangle$. If $I$ has a generating set $S$ that is a~Gr\"{o}bner---Shirshov basis, then $Irr(S)$ is a $\Bbbk-$ basis of $R\As\langle X\rangle/I$.
\end{corollary}

\section{Gr\"{o}bner---Shirshov bases for nilpotent Rota---Baxter algebra}

Let $A$ be an associative algebra over $\Bbbk$. 
Let a $\Bbbk$-linear operator $R\colon A \to A$ satisfy
$$
R(x)R(y)=R(xR(y)+R(x)y  + \lambda xy ),\quad R^k(R+\lambda id)^l=0,\ k,l\geq 1,
$$
where $\lambda \in \Bbbk$ is a fixed scalar.
We may rewrite the second relation as follows, 
$$
R^k(R+\lambda id)^l
 = \sum\limits^l_{i=0}\binom{l}{i}\lambda ^{l-i} R^{k+i}
 = \sum\limits^l_{i=0}\lambda^{(i)} R^{k+i},
$$
where for fixed $l$, $\lambda^{(i)}  := \binom{l}{i} \lambda^{l-i}$. 

Let $n,m\in \mathbb{N}$ and $n,m\geq2$.
Denote $I_n=\{1,2,\ldots,n-1\}$.
Given $A^{n}_m=(i_1,\ldots,i_m)\in I_n^m$
and $B^n_{m}=(j_1,\ldots,j_{m-1})\in I_n^{m-1}$, 
we say that $A^n_m\in \mathcal{A}^n_m$ and $B^n_m\in \mathcal{B}^n_m$, if the conditions

$(1)$ $i_m>1$, 

$(2)$ $\sum\limits^m_{t=1}i_t=n$,

$(3)$ $j_{m-1}>1$, 

$(4)$ $\sum\limits^{m-1}_{t=1}j_t=n$

\noindent hold, respectively.
We say that $A^n_m\in \mathcal{C}^n_m$ and $B^n_m\in \mathcal{D}^n_m$, if only the condition~$(2)$ and $(4)$ holds, respectively.
If \( A_m^n = (i_1, \dotsc, i_m) \in I_m^n \), then we introduce the following notations for the tuples:
\begin{align*}
A_m^{n - q\mathbf{1}_s} &:= (i_1, \dotsc, i_{s-1}, i_s - q, i_{s+1}, \dotsc, i_m), \\
A_m^n|_{s} &:= (i_s, i_{s+1}, \dotsc, i_m), \\
A_m^n \downarrow_s &:= (i_1, \dotsc, i_s).
\end{align*}
Analogous denotations are used for \( B_m^n \in I_{m-1}^n \).
For $A^n_m\in \mathcal{A}^n_m$ and $B^n_m\in \mathcal{B}^n_m$, we define
\begin{gather*}
L_{A^n_m}[\vec{x}]=L_{A^n_m}[x_1,x_2,\ldots,x_m]=R^{i_1}\Big(R^{i_2}(\ldots(R^{i_{m-1}}(R^{i_m}(x_1)x_2)x_3\ldots)x_m\Big), \\
L_{B^n_m}[\vec{x}]=L_{B^n_{m}}[x_1,x_2,\ldots,x_m]=R^{j_1}\Big(\ldots(R^{j_{m-2}}(R^{j_{m-1}}(x_1)x_2)x_3)\ldots\Big)x_m, \\
L_{A^n_m}^{\op}[\vec{x}]=L_{A^n_m}^{\op}[x_1,x_2,\ldots,x_m]=R^{i_1}\Big(x_1R^{i_2}(\ldots(x_{m-2}R^{i_{m-1}}(x_{m-1}R^{i_m}(x_m))\ldots)\Big), \\
L_{B^n_m}^{\op}[\vec{x}]=L_{B^n_m}^\op[x_1,x_2,\ldots,x_m]=x_1R^{j_1}\Big(\ldots(x_{m-2}R^{j_{m-2}}(x_{m-1}R^{j_{m-1}}(x_m)))\ldots\Big), 
\end{gather*}

\vspace{-0.6cm}

\begin{multline*}
L^*_{A_m^n}[x,\vec{x}]
 = \sum \limits_{s=1}^m\sum^{i_s}_{t=1}
 \bigg( L_{(A_m^n\downarrow_{s-1},t)}[xL_{(i_s-t+1,{A_m^n}|_{s+1})}[x_1,\ldots,x_{m+1-s}],x_{m+2-s},\ldots,x_m] \\
 + \lambda L_{(A_m^n\downarrow_{s-1},t)}[xL_{(i_s-t,{A_m^n}|_{s+1})}[x_1,\ldots,x_{m+1-s}],x_{m+2-s},\ldots,x_m] \bigg);
\end{multline*}

\vspace{-0.6cm}

\begin{multline*}
 L^*_{B^n_m}[x,\vec{x}] 
  = \sum \limits_{s=1}^m\sum^{i_s}_{t=1}
 \bigg( L_{(B_m^n\downarrow_{s-1},t)}[xL_{(i_s-t+1,{B_m^n}|_{s+1})}[x_1,\ldots,x_{m-s}],x_{m+1-s},\ldots,x_m] \\
  + \lambda L_{(B_m^n\downarrow_{s-1},t)}[xL_{(i_s-t,{B_m^n}|_{s+1})}[x_1,\ldots,x_{m-s}],x_{m+1-s},\ldots,x_m] \bigg); 
\end{multline*}

\vspace{-0.6cm}

\begin{multline*}
{L^*}^{\op}_{A^n_m}[x,\vec{x}]
 = \sum\limits_{s=1}^m \sum^{i_s}_{t=1} \bigg(
 L^\op_{(A_m^n\downarrow_{s-1},t)}[x_1,\ldots,x_s L^\op_{(i_s-t+1,A_m^n|_{s+1})}[x_{s+1},\ldots,x_m]x] \\
  + \lambda L^\op_{(A_m^n\downarrow_{s-1},t)}[x_1,\ldots,x_s L^\op_{(i_s-t,A_m^n|_{s+1})}[x_{s+1},\ldots,x_m]x] \bigg);
\end{multline*}

\vspace{-0.6cm}

\begin{multline*}
 {L^*}^{\op}_{B^n_m}[x,\vec{x}]
  = \sum\limits_{s=1}^m \sum^{i_s}_{t=1} \bigg(
 L^\op_{(B_m^n\downarrow_{s-1},t)}[x_1,\ldots,R^{i_s-t+1}(x_{s+1} L^\op_{B_m^n|_{s+1}}[x_{s+2},\ldots,x_m])x] \\
  + \lambda L^\op_{(B_m^n\downarrow_{s-1},t)}[x_1,\ldots,R^{i_s-t}(x_{s+1} L^\op_{B_m^n|_{s+1}}[x_{s+2},\ldots,x_m])x] \bigg).
\end{multline*}

We apply the same formulas to define 
$L_{C^n_m}$, $L_{D^n_m}$, 
$L^\op_{C^n_m}$, $L^\op_{D^n_m}$, 
${L^*}^\op_{C^n_m}$, ${L^*}^\op_{D^n_m}$
for $C^n_m \in \mathcal{C}^n_m$ and $D^n_m \in \mathcal{D}^n_m$. 
Throughout the paper we assume that $L_{\mathcal{A}^n_m}[\vec{x}]=\sum_{A^n_m\in \mathcal{A}^n_m}L_{A^{n}_m}[\vec{x}]$. Analogously we deal with  ${B^{n}_m\in \mathcal{B}^n_m}$, ${C^{n}_m\in \mathcal{C}^n_m}$ and $D^n_m \in \mathcal{D}^n_m$.
\begin{lemma}\label{lem:operator-identities}
The operators $L_{{A}^n_m}$, $L_{{B}^n_m}$, $L_{{C}^n_m}$, and $L_{{D}^n_m}$ satisfy the following identities:
\begin{gather}
 L_{B_{m-1}^{n}} [x_1,\ldots,x_{m-1}]R(x_m)
 =  L_{B_{m-1}^{n}} [x_1,\ldots,x_{m-1}R(x_m)]; \label{eq1} \\
 L_{B_{m-1}^{n}} [x_1,\ldots,x_{m-1}]x_m= L_{B_{m-1}^{n}} [x_1,\ldots,x_{m-1}x_m]; \label{eq2} \\
R\bigg( L_{B_{m-1}^{n}} [\vec{x}]+L_{A_{m-1}^{n}} [\vec{x}]\bigg)x_m 
 =  L_{B_{m}^{n+1}} [\vec{x},x_m]; \quad |\vec{x}| = m-1; \label{eq3}
\end{gather} 

\vspace{-0.5cm}

\begin{multline} \label{eq4}
 L_{A_{m-1}^{n}} [x_1,\ldots,x_{m-1}]R(x_m)
 =  L_{A_{m-1}^{n}} [x_1,\ldots,x_{m-1}R(x_m)] \\
 + L_{A_{m}^{n+1}} [x_1,\ldots,x_{m-1},x_m]+\lambda  L_{A_{m-1}^{n}} [x_1,\ldots,x_{m-1}x_m]+\lambda  L_{A_{m}^{n}} [x_1,\ldots,x_{m-1},x_m]; 
\end{multline} 

\vspace{-0.5cm}
\begin{equation} \label{eq5}
  L_{A_{m-1}^{n}} [\vec{x}]x_m
 =  L_{B_{m}^{n}} [\vec{x},x_m], \quad |\vec{x}| = m-1.
\end{equation}

\begin{gather}\label{3.7}
R(x) L_{\mathcal{A}^n_{m}}[\vec{x}]
 =  L^*_{\mathcal{A}^n_{m}}[x,\vec{x}]
 +  L_{\mathcal{A}^n_{m}}[R(x)x_1+xR(x_1),x_2,\ldots,x_{m}]
 +\lambda  L_{\mathcal{A}^n_{m}}[xx_1,x_2,\ldots,x_{m}]; \\
 L_{\mathcal{A}^n_{m}}[\vec{x}] R(x_{m+1})
 =  L_{\mathcal{A}^{n+1}_{m+1}}[\vec{x},x_{m+1}]
 +  L_{\mathcal{A}^n_{m}}[x_1,x_2,\ldots,x_{m}R(x_{m+1})]
 + \lambda  L_{\mathcal{A}^n_{m}}[x_1,x_2,\ldots,x_{m}x_{m+1}];
\end{gather}
\begin{multline}\label{3.8}
L^*_{\mathcal{A}^{n+1}_{m+1}}[x,\vec{x}] R(x_{m+1})
 = L^*_{\mathcal{A}^{n+1}_{m+1}}[x,\vec{x},x_{m+1}]
 + L_{\mathcal{A}^{n+1}_{m+1}}[R(x)x_1+xR(x_1),x_2,\ldots,x_{m+1}] \\
 +\lambda  L_{\mathcal{A}^{n+1}_{m+1}}[xx_1,x_2,\ldots,x_{m+1}]
 + L^*_{\mathcal{A}^{n}_{m}}[x,x_1,x_2,\ldots,x_mR(x_{m+1})] \\
 + L_{\mathcal{A}^{n+1}_{m+1}}[R(x)x_1+xR(x_1),x_2,\ldots,x_mR(x_{m+1})]
 +\lambda  L_{\mathcal{A}^{n}_{m}}[xx_1,x_2,\ldots,x_mR(x_{m+1})]\\
 + L^*_{\mathcal{A}^{n}_{m}}[x,x_1,x_2,\ldots,x_mx_{m+1}]
 +\lambda  L_{\mathcal{A}^{n+1}_{m+1}}[R(x)x_1+xR(x_1),x_2,\ldots,x_mx_{m+1}]\\
 +\lambda^2  L_{\mathcal{A}^{n}_{m}}[xx_1,x_2,\ldots,x_mx_{m+1})]\\
 - L_{\mathcal{A}^{n+1}_{m+1}}[R(x)x_1+xR(x_1),x_2,\ldots,x_{m+1}]
 -  L_{\mathcal{A}^{n}_{m}}[R(x)x_1+xR(x_1),x_2,\ldots,x_mR(x_{m+1})]\\-\lambda  L_{\mathcal{A}^{n}_{m}}[R(x)x_1+xR(x_1),x_2,\ldots,x_mx_{m+1}]\\
 -  L_{\mathcal{A}^{n+1}_{m+1}}[xx_1,x_2,\ldots,x_{m+1}]
 - \lambda L_{\mathcal{A}^{n}_{m}}[xx_1,x_2,\ldots,x_mR(x_{m+1})]
 -\lambda^2  L_{\mathcal{A}^{n}_{m}}[xx_1,x_2,\ldots,x_mx_{m+1}];
\end{multline}

\vspace{-0.8cm}

\begin{multline}\label{3.9}
R(y) L^*_{\mathcal{A}^{n}_{m}}[R(y)x+yR(x)+\lambda xy,z,\vec{x}]
 =  L^*_{\mathcal{A}^n_{m}}[R(y)x+yR(x)+\lambda xy,z,\vec{x}]\\
 + L_{\mathcal{A}^n_{m}}[R(R(y)x+yR(x)+\lambda xy)z+(R(y)x+yR(x)+\lambda xy)R(z),\vec{x}]\\
 +\lambda  L_{\mathcal{A}^n_{m}}[R(y)xz+yR(x)z+\lambda xyz,\vec{x}]
 - L^*_{\mathcal{A}^n_{m}}[y,R(x)z+xR(z),\vec{x}] \\ \allowdisplaybreaks
 - L_{\mathcal{A}^n_{m}}[R(y)(R(x)z+xR(z))+yR(R(x)z+xR(z)),\vec{x}]
 -\lambda  L_{\mathcal{A}^n_{m}}[y(R(x)z+xR(z)),\vec{x}] \\
 -\lambda  L^*_{\mathcal{A}^n_{m}}[y,xz,\vec{x}] 
 -\lambda L_{\mathcal{A}^n_{m}}[R(y)xz+yR(xz),\vec{x}]
 -\lambda^2  L_{\mathcal{A}^n_{m}}[yxz,\vec{x}].
\end{multline}

\begin{multline}\label{3.10}
R(x)R^{i_1}\Big(R^{i_2}(\ldots(R^{i_{m-1}}(R^{i_m}(x_1)x_2)x_3)\ldots)x_m\Big) \\
 = \sum\limits^{i_1}_{t=1}R^{t}\bigg(xR^{i_1-t+1}\Big(R^{i_2}(\ldots(R^{i_{m-1}}(R^{i_m}(x_1)x_2)x_3)\ldots)x_m\Big)\bigg)
 +R^{i_1}\bigg(\sum\limits^{i_2}_{t=1}R^{t}\Big(xR^{i_2-t+1}(\ldots(R^{i_{m-1}}(R^{i_m}(x_1)x_2)x_3)\ldots\Big)x_m\bigg) \\
 +\ldots+R^{i_1}\Big(R^{i_2}(\ldots(R^{i_{m-1}}(\sum\limits^{i_m-1}_{t=1}R^{t}\Big(xR^{i_m-t+1}(x_1)\Big)x_2)x_3)\ldots)x_m\Big)
+\lambda\sum\limits^{i_1-1}_{t=1}R^{t}\bigg(xR^{i_1-t}\Big(R^{i_2}(\ldots(R^{i_{m-1}}(R^{i_m}(x_1)x_2)x_3)\ldots)x_m\Big)\bigg)\\
+\lambda R^{i_1}\bigg(\sum\limits^{i_2-1}_{t=1}R^{t}\Big(xR^{i_2-t}(\ldots(R^{i_{m-1}}(R^{i_m}(x_1)x_2)x_3)\ldots\Big)x_m\bigg)
+\ldots+\lambda R^{i_1}\Big(R^{i_2}(\ldots(R^{i_{m-1}}(\sum\limits^{i_m-1}_{t=1}R^{t}\Big(xR^{i_m-t}(x_1)\Big)x_2)x_3)\ldots)x_m\Big)\\
+R^{i_1}\Big(R^{i_2}(\ldots(R^{i_{m-1}}(R^{i_m}\Big(R(x)x_1+xR(x_1)\Big)x_2)x_3)\ldots)x_m\Big)
+\lambda R^{i_1}\Big(R^{i_2}(\ldots(R^{i_{m-1}}(R^{i_m}(xx_1)x_2)x_3)\ldots)x_m\Big)\\
+\lambda R^{i_1}\Big(xR^{i_2}(\ldots(R^{i_{m-1}}(R^{i_m}(x_1)x_2)x_3)\ldots)x_m\Big).
\end{multline}
We may apply the similar expression for rewriting
$R(x)R^{i_1}\Big(R^{i_2}(\ldots(R^{i_{m-2}}(R^{i_{m-1}}(x_1)x_2)x_3)\ldots)x_{m-1}\Big)x_m.$
\end{lemma}
\begin{proof}

The relations~\eqref{eq1},~\eqref{eq2}, and~\eqref{eq5} follow from the definition of~$L_{B_{m-1}^n}[\vec{x}]$ and~$L_{A_{m-1}^n}[\vec{x}]$.

Let us prove~\eqref{eq3} by induction on $m$. For $m = 3$, we have
\begin{align*}
R\bigg( L_{B_{2}^{n}}[\vec{x}]
 +  L_{A_{2}^{n}} [\vec{x}]\bigg)x_3 
 = R\bigg(R^n(x_1)x_2+\sum\limits^{n-1}_{i=2}R^{n-i}(R^i(x_1)x_2)\bigg)x_3  
 = \sum\limits^{n}_{i=2}R^{n+1-i}(R^i(x_1)x_2)x_3\\
 = L_{B_{3}^{n+1}}[\vec{x},x_3]; \quad |\vec{x}| = 2.
\end{align*}

For $m>3$, we have
\begin{multline*}
 L_{B_{m}^{n+1}} [\vec{x},x_m] 
 =   \\
 \!{=} R\bigg(\!\sum\limits_{B^{n}_{m-1}{\in}\mathcal{B}^n_{m-1}}\!\!\! R^{i_2}(R^{i_3}(\ldots(R^{i_{m-1}}(x_1)x_2)x_3)\ldots)x_{m-1}\bigg)x_m
 {+} R\bigg(\!\sum\limits_{A^{n}_{m-1} {\in} \mathcal{A}^n_{m-1}, i_1{\geq}2}\!\!\!\!\!\! R^{i_1-1}(R^{i_2}(R^{i_3}(\ldots(R^{i_{m-1}}(x_1)x_2)x_3)\ldots)x_{m-1})\bigg)x_m\\
 = R\bigg(L_{B_{m-1}^{n}} [\vec{x}]+L_{A_{m-1}^{n}} [\vec{x}]\bigg)x_m;\quad |\vec{x}|=m-1.
\end{multline*}

We deduce~\eqref{eq4} as follows:
\begin{multline*}
R^{i_1}(R^{i_2}(R^{i_3}(\ldots(R^{i_{m-1}}(x_1)x_2)x_3)\ldots)x_{m-1}).R(x_m) \\
 = \sum\limits^{i_1}_{i=1}R^i(R^{i_1-i+1}(R^{i_2}(R^{i_3}(\ldots(R^{i_{m-1}}(x_1)x_2)x_3)\ldots)x_{m-1})x_m)
 + R^{i_1}\Big(R^{i_2}(R^{i_3}(\ldots(R^{i_{m-1}}(x_1)x_2)x_3)\ldots)(x_{m-1}R(x_m))\Big) \\
 + \lambda \sum\limits^{i_1-1}_{i=1}R^i(R^{i_1-i}(R^{i_2}(R^{i_3}(\ldots(R^{i_{m-1}}(x_1)x_2)x_3)\ldots)x_{m-1})x_m) 
 + \lambda R^{i_1}\Big(R^{i_2}(R^{i_3}(\ldots(R^{i_{m-1}}(x_1)x_2)x_3)\ldots)x_m\Big).
\end{multline*}
Another's relations have been proven similar and
the relations ~\eqref{eq1}--\eqref{3.9} hold true for $C^n_m \in \mathcal{C}^n_m$ and $D^n_m \in \mathcal{D}^n_m$.
\end{proof}

Consider the following relations in the free associative RB-algebra:
\begin{gather*}
(R1)\ R(x)R(y) = R( xR(y) + R(x)y  + \lambda xy ); \quad
(R2)\ R^k( R + \lambda \id)^l(x) = 0; \\
(R3) \ \sum_{t=0}^{l} \lambda^{(t)} \Bigg( 
    R^2 \Big( L_{\mathcal{B}^{k+t-1}_2}[\vec{x}] + L_{\mathcal{A}^{k+t-1}_2}[\vec{x}] \Big)
    + \lambda R \Big( L_{\mathcal{D}^{k+t-1}_2}[\vec{x}] + L_{\mathcal{C}^{k+t-1}_2}[\vec{x}] \Big)
\Bigg) = 0; \quad |\vec{x}| = 2,
\end{gather*}

\vspace{-0.7cm}

\begin{multline*}
(R4) 
\sum^l_{t=0}\lambda^{(t)}
\Bigg(\sum_{i=0}^{m-2}\lambda^i \tbinom{m-2}{i} R^2\Big(L_{\mathcal{B}^{k+t+m-3-i}_m}[\vec{x}]
 + L_{\mathcal{A}^{k+t+m-3-i}_m}[\vec{x}]\Big)   
 + \lambda\sum_{i=0}^{m-3}\lambda^i \tbinom{m-3}{i}R\Big(L_{\mathcal{B}^{k+t+m-3-i}_m}[\vec{x}]\Big)\\
 + \lambda \sum_{i=0}^{m-2}\lambda^i\tbinom{m-2}{i} R\Big(L_{\mathcal{C}^{k+t+m-3-i}_m}[\vec{x}]\Big) 
  + \lambda^2 \sum_{i=0}^{m-3}\lambda^i\tbinom{m-3}{i} R\Big(L_{\mathcal{D}^{k+t+m-4-i}_m}[\vec{x}]\Big)\Bigg) = 0; \ m = |\vec{x}|\geq 3,
\end{multline*}

\vspace{-0.5cm}

$$
(R5) \sum^l_{t=0}\lambda ^{(t)}\Bigg(R^2\Big( L^\op_{\mathcal{B}^{k+t-1}_2}[\vec{x}]+L^{\op}_{\mathcal{A}^{k+t-1}_2}[\vec{x}]\big)
+\lambda R\big( L^{\op}_{\mathcal{B}^{k+t-1}_2}[\vec{x}]+L^\op_{\mathcal{C}^{k+t-1}_2}[\vec{x}]\big)\\
+\lambda^2  R\big(L^\op_{\mathcal{D}^{k+t-2}_m}[\vec{x}]\big)\Bigg)=0;\ |\vec{x}|{=}2, 
$$

\vspace{-0.5cm}

\begin{multline*}
(R6) 
\sum^l_{t=0}\lambda^{(t)}\Bigg(
\sum_{i=0}^{m-2}\lambda^i\tbinom{m-2}{i}R^2\Big(L^\op_{\mathcal{B}^{k+t+m-3-i}_m}[\vec{x}]
 + L^\op_{\mathcal{A}^{k+t+m-3-i}_m}[\vec{x}]\Big)  
 + \lambda \sum_{i=0}^{m-3}\lambda^i\tbinom{m-3}{i} R\Big(L^\op_{\mathcal{B}^{k+t+m-3-i}_m}[\vec{x}]\Big)\\ 
 + \lambda \sum_{i=0}^{m-2}\lambda^i\tbinom{m-2}{i}R\Big(L^\op_{\mathcal{C}^{k+t+m-3-i}_m}[\vec{x}]\Big) 
  +\lambda^2 \sum_{i=0}^{m-3}\lambda^i\tbinom{m-3}{i} R\Big(L^\op_{\mathcal{D}^{k+t+m-4-i}_m}[\vec{x}]\Big)\Bigg)=0, \ m = |\vec{x}|\geq3;
\end{multline*}

\vspace{-0.7cm}

\begin{multline*}
(R7) 
 \sum^l_{t=0}\lambda^{(t)}\Bigg( 
 R^2 \Big(xR\Big( L_{\mathcal{B}^{k+t-1}_2}[\vec{x}]
 + L_{\mathcal{A}^{k+t-1}_2}[\vec{x}]\Big)\Big)
 + \lambda R^2\Big( x\Big( L_{\mathcal{B}^{k+t-1}_2}[\vec{x}]
 + L_{\mathcal{A}^{k+t-1}_2}[\vec{x}]\Big)\Big) \\
 + \lambda R\Big( xR\Big( L_{\mathcal{B}^{k+t-1}_2}[\vec{x}]
 + L_{\mathcal{A}^{k+t-1}_2}[\vec{x}]\Big)\Big)
 + R^2\Big( L^*_{\mathcal{B}^{k+t-1}_3}[x,\vec{x}]+L^*_{\mathcal{A}^{k+t-1}_3}[x,\vec{x}]\Big) \\
 + \lambda R\Big( L^*_{\mathcal{D}^{k+t-1}_3}[x,\vec{x}]
 + L^*_{\mathcal{C}^{k+t-1}_3}[x,\vec{x}]\Big)
 + \lambda^2 R\Big(x\Big( L_{\mathcal{D}^{k+t-1}_2}[\vec{x}]
 + L_{\mathcal{C}^{k+t-1}_2}[\vec{x}]\Big)\Big)\Bigg) {=} 0,\ |\vec{x}| {=}2; 
\end{multline*}

\vspace{-0.5cm}

\begin{multline*}
(R8) 
\sum^l_{t=0}\lambda^{(t)} \Bigg( 
 \sum_{i=0}^{m-2}\lambda^i\tbinom{m-2}{i}R^2\Big(L^*_{\mathcal{B}^{k+t+m-3-i}_m}[x,\vec{x}]
 + L^*_{\mathcal{A}^{k+t+m-3-i}_m}[x,\vec{x}]\Big) 
 + \lambda\sum_{i=0}^{m-3}\lambda^i \tbinom{m-3}{i} R\Big(L^*_{\mathcal{B}^{k+t+m-3-i}_m}[x,\vec{x}]\Big)\\ 
 + \lambda \sum_{i=0}^{m-2}\lambda^i\tbinom{m-2}{i} R\Big(L^*_{\mathcal{C}^{k+t+m-3-i}_m}[x,\vec{x}]\Big) 
 + \lambda^2 \sum_{i=0}^{m-3}\lambda^i\tbinom{m-3}{i} R\Big(L^*_{\mathcal{D}^{k+t+m-4-i}_m}[x,\vec{x}]\Big)
 + \sum_{i=0}^{m-2}\lambda^i\tbinom{m-2}{i}R^2\Big(xR\Big( L_{\mathcal{B}^{k+t+m-3-i}_{m}}[\vec{x}]
 + L_{\mathcal{A}^{k+t+m-3-i}_{m}}[\vec{x}]\Big)\Big) \\
 + \lambda\sum_{i=0}^{m-2}\lambda^i \tbinom{m-2}{i} R^2\Big(x\Big( L_{\mathcal{B}^{k+t+m-3-i}_{m}}[\vec{x}]
 +  L_{\mathcal{A}^{k+t+m-3-i}_{m}}[\vec{x}]\Big)\Big) \allowdisplaybreaks \\
 + \lambda \sum_{i=0}^{m-2}\lambda^i\tbinom{m-2}{i}R\Big(x R\Big(L_{\mathcal{B}^{k+t+m-3-i}_{m}}[\vec{x}]
 +  L_{\mathcal{A}^{k+t+m-3-i}_{m}}[\vec{x}]\Big)\Big)
 + \lambda^2\sum_{i=0}^{m-3}\lambda^i \tbinom{m-3}{i} R\Big(x \Big( L_{\mathcal{B}^{k+t+m-3-i}_{m}}[\vec{x}]\Big)\Big) \\
 + \lambda^2 \sum_{i=0}^{m-2}\lambda^i\tbinom{m-2}{i} R\Big(x \Big( L_{\mathcal{C}^{k+t+m-3-i}_{m}}[\vec{x}]\Big)\Big)
 +\lambda^3\sum_{i=0}^{m-3}\lambda^i \tbinom{m-3}{i}R\Big(x \Big( L_{\mathcal{D}^{k+t+m-4-i}_{m}}[\vec{x}]\Big)\Big) \Bigg) {=} 0,\ m {=} |\vec{x}|{\geq}3; 
\end{multline*}

\vspace{-0.3cm}

\begin{multline*}
(R9) 
\sum^l_{t=0}\lambda^{(t)}\Bigg( R^2 \Big(R\Big( L^\op_{\mathcal{B}^{k+t-1}_2}[\vec{x}]+L^\op_{\mathcal{A}^{n-1}_2}[\vec{x}]\Big)x\Big)
 +\lambda R^2\Big( \Big( L^\op_{\mathcal{B}^{k+t-1}_2}[\vec{x}]+L^\op_{\mathcal{A}^{k+t-1}_2}[\vec{x}]\Big)x\Big) \\
 +\lambda R\Big( R\Big( L^\op_{\mathcal{B}^{k+t-1}_2}[\vec{x}]+L^\op_{\mathcal{A}^{k+t-1}_2}[\vec{x}]\Big)x\Big)
 + R^2\Big( {L^*}^\op_{\mathcal{B}^{k+t-1}_3}[x,\vec{x}]+{L^*}^\op_{\mathcal{A}^{k+t-1}_3}[x,\vec{x}]\Big) \\
 + \lambda R\Big( {L^*}^\op_{\mathcal{B}^{k+t-1}_3}[x,\vec{x}]+{L^*}^\op_{\mathcal{C}^{k+t-1}_3}[x,\vec{x}]\Big)
 + \lambda^2 R\Big(\Big( L^\op_{\mathcal{D}^{k+t-1}_2}[\vec{x}]+L^\op_{\mathcal{C}^{k+t-1}_2}[\vec{x}]\Big)x\Big)\Bigg) {=} 0;\ |\vec{x}|{=}2, 
\end{multline*}

\vspace{-0.3cm}

\begin{multline*}
(R10) 
\sum^l_{t=0}\lambda^{(t)} \Bigg(
\sum_{i=0}^{m-2}\lambda^i\tbinom{m-2}{i}R^2\Big({L^*}^\op_{\mathcal{B}^{k+t+m-3-i}_m}[x,\vec{x}]
+{L^*}^\op_{\mathcal{A}^{k+t+m-3-i}_m}[x,\vec{x}]\Big) 
+\lambda\sum_{i=0}^{m-3}\lambda^i \tbinom{m-3}{i}R\Big({L^*}^\op_{\mathcal{B}^{k+t+m-3-i}_m}[x,\vec{x}]\Big)\\
+\lambda \sum_{i=0}^{m-2}\lambda^i\tbinom{m-2}{i} R\Big({L^*}^\op_{\mathcal{C}^{k+t+m-3-i}_m}[x,\vec{x}]\Big)
+\lambda^2 \sum_{i=0}^{m-3}\lambda^i\tbinom{m-3}{i}R\Big({L^*}^\op_{\mathcal{D}^{k+t+m-4-i}_m}[x,\vec{x}]\Big)\\
 + \sum_{i=0}^{m-2}\lambda^i\tbinom{m-2}{i}R^2\Big(R\Big( L^\op_{\mathcal{B}^{k+t+m-3-i}_{m}}[\vec{x}]
 +  L^\op_{\mathcal{A}^{k+t+m-3-i}_{m}}[\vec{x}]\Big)x\Big)
 + \lambda \sum_{i=0}^{m-2}\lambda^i\tbinom{m-2}{i}R^2\Big(\Big( L^\op_{\mathcal{B}^{k+t+m-3-i}_{m}}[\vec{x}]
 +  L^\op_{\mathcal{A}^{k+t+m-3-i}_{m}}[\vec{x}]\Big)x\Big)\\
 + \lambda\sum_{i=0}^{m-2}\lambda^i \tbinom{m-2}{i}R\Big( R\Big( L^\op_{\mathcal{B}^{k+t+m-3-i}_{m}}[\vec{x}]
 +  L^\op_{\mathcal{A}^{k+t+m-3-i}_{m}}[\vec{x}]\Big)x\Big)
 + \lambda^2 \sum_{i=0}^{m-3}\lambda^i\tbinom{m-3}{i} R\Big( \Big( L^\op_{\mathcal{B}^{k+t+m-3-i}_{m}}[\vec{x}]\Big)x\Big) \\
 + \lambda^2 \sum_{i=0}^{m-2}\lambda^i\tbinom{m-2}{i}R\Big( \Big( L^\op_{\mathcal{C}^{k+t+m-3-i}_{m}}[\vec{x}]\Big)x\Big)
  +\lambda^3 \sum_{i=0}^{m-3}\lambda^i\tbinom{m-3}{i}R\Big( \Big( L^\op_{\mathcal{D}^{k+t+m-4-i}_{m}}[\vec{x}]\Big)x\Big) \Bigg) {=} 0;\ m = |\vec{x}|{\geq}3.
\end{multline*}

Let us denote the relations $(R1)$--$(R10)$ as $f_1,\ldots,f_{10}$, respectively.
The leading words in the relations $(R1)$--$(R10)$ are 
\begin{gather*}
  \bar{f_1}=R(x)R(y), \quad
  \bar{f_2}=R^{k+l}(x), \quad
\bar{f_3}=R^2\Bigg(R^{k+l-1}(x_1)x_2)\Bigg), \quad
\bar{f_4}=R^2\Bigg(R\Big(\ldots(R(R^{k+l-m+1}(x_1)x_2)x_3)\ldots\Big)x_m\Bigg), \\
\bar{f_5}=R^2\Big(x_1R^{k+l-1}(x_2)\Big), \quad
\bar{f_6}=R^2\Bigg(x_1R\Big(\ldots(x_{m-2}R(x_{m-1}R^{k+l-m+1}(x_m))\ldots\Big)\Bigg), \\
\bar{f_7}=R^2\Bigg(xR\Big(R^{k+l-1}(x_1)x_2\Big)\Bigg), \quad
\bar{f_8}=R^2\Bigg(xR\Big(\ldots(R(R^{k+l-m+1}(x_1)x_2)x_3)\ldots\Big)x_m\Bigg), \\
\bar{f_9}=R^2\Bigg(xR\Big(x_1R^{k+l-1}(x_2)\Big)\Bigg), \quad
\bar{f_{10}}=R^2\Bigg(xR\Big(x_1R\Big(\ldots(x_{m-2}R(x_{m-1}R^{k+l-m+1}(x_m))\ldots\Big)\Big)\Bigg).
\end{gather*}
\begin{lemma}\label{lemm3.2}
The relations $(R3)$--$(R10)$ are contained in the ideal generated by $(R1)$ and $(R2)$.
\end{lemma}
\begin{proof}

 We may rewrite $(R1)$ and $(R2)$ as follows,
$$
R^{k+l}(x_1)R(x_2)=R\Big(R^{k+l}(x_1)x_2+R^{k+l-1}(x_1)R(x_2)+\lambda R^{k+l-1}(x_1)x_2\Big), \quad 
R^{k+l}(x_1) =-\sum_{t=0}^{l-1}\lambda^{(t)}R^{k+t}(x_1).
$$

Therefore,
\[
f_1-f_2R(x_2)
=\sum_{t=0}^{l-1}\lambda^{(t)}R^{k+t}(x_1)R(x_2)+R\Big(R^{k+l}(x_1)x_2+R^{k+l-1}(x_1)R(x_2)+\lambda R^{k+l-1}(x_1)x_2\Big).
\]
For $n\geq 1$, we have
\begin{multline*}\label{rel.2}
R^n(x_1)R(x_2)
 = \sum\limits^{n}_{j=1}R^j(R^{n-j+1}(x_1)x_2)+R^n(x_1R(x_2))+\lambda \bigg(\sum\limits^{n-1}_{j=1}R^j(R^{n-j}(x_1)x_2) +R^n(x_1x_2)\bigg)\\
 = R(R^{n}(x_1)x_2)+\sum\limits^{n-1}_{j=2}R^j(R^{n-j+1}(x_1)x_2)+R^n(R(x_1)x_2)+R^n(x_1R(x_2))
 +\lambda \bigg(\sum\limits^{n-1}_{j=1}R^j(R^{n-j}(x_1)x_2) +R^n(x_1x_2)\bigg)\\
 \end{multline*}
 \begin{multline*}
 = R^2\bigg(R^{n-1}(x_1)x_2+R(R^{n-2}(x_1)x_2)+\ldots+R^{n-3}(R^2(x_1)x_2)\bigg)
 +\lambda R\bigg(R^{n-1}(x_1)x_2+R(R^{n-2}(x_1)x_2)+\ldots+R^{n-2}(R(x_1)x_2) \bigg) \\
 +R(R^{n}(x_1)x_2)+R^n(R(x_1)x_2)+R^n(x_1R(x_2))+\lambda R^n(x_1x_2)\\
 = R^2\bigg( L_{B^{n-1}_2}[x_1,x_2]+L_{A^{n-1}_2}[x_1,x_2]\bigg) 
 +\lambda R\bigg( L_{D^{n-1}_2}[x_1,x_2]+L_{C^{n-1}_2}[x_1,x_2]\bigg) \\
 +R(R^{n}(x_1)x_2)+R^n(R(x_1)x_2)+R^n(x_1R(x_2))+\lambda R^n(x_1x_2).
\end{multline*}
Therefore,
\begin{multline*}
f_1-f_2R(x_2)
  = \sum\limits^l_{t=0}\lambda ^{(t)}\Bigg(R^2\Big( L_{\mathcal{B}^{k+t-1}_2}[x_1,x_2]+L_{\mathcal{A}^{k+t-1}_2}[x_1,x_2]\Big)\\
 +\lambda R\Big( L_{\mathcal{D}^{k+t-1}_2}[x_1,x_2]+L_{\mathcal{C}^{k+t-1}_2}[x_1,x_2]\Big)
 +R(R^{k+t}(x_1)x_2)+R^{k+t}(R(x_1)x_2)+R^{k+t}(x_1R(x_2))+\lambda R^{k+t}(x_1x_2)\Bigg) \allowdisplaybreaks \\
 = \sum\limits^l_{t=0}\lambda ^{(t)}\Bigg(R^2\Big( L_{\mathcal{B}^{k+t-1}_2}[x_1,x_2]+L_{\mathcal{A}^{k+t-1}_2}[x_1,x_2]\Big)
+\lambda R\Big( L_{\mathcal{D}^{k+t-1}_2}[x_1,x_2]+L_{C^{k+t-1}_2}[x_1,x_2]\Big)\Bigg)\\
+R(f_1(x_1)x_2)+f_1(R(x_1)x_2+x_1R(x_2))+\lambda f_1(x_1x_2)\\
=f_3+R(f_1(x_1)x_2)+f_1(R(x_1)x_2+x_1R(x_2))+\lambda f_1(x_1x_2).
\end{multline*}
It follows that $\mathrm{f_3} \in \langle \mathrm{f_1}, \mathrm{f_2} \rangle$. The derivations of $(R4)$--$(R10)$ follow similar patterns with increasing combinatorial complexity.
\end{proof}

\begin{theorem}
The relations $(R1)$--$(R10)$ form a Gr\"{o}bner---Shirshov basis.
 \end{theorem} 
\begin{proof}
The relations $(R3)$--$(R6)$ arise from compositions between $(R1)$ and $(R2)$.

Specifically:
\begin{itemize}
    \item Relation $(R3)$ is obtained from the composition of $(R1)$ and $(R2)$ at $w = \bar{f_1} = \bar{f_2}R(x_2)$ (see Lemma~\ref{lemm3.2})
    \item Relation $(R4)$ is derived by

$$
R^2(x)R(y)=R(R^2(x)y)+R^2(R(x)y)+R^2(xR(y))+\lambda R^2(xy)+\lambda R(R(x)y),
$$
we deduce by~\eqref{eq1}--\eqref{eq5}
\begin{multline*}
 R^2\bigg( L_{\mathcal{B}^{k+t-1}_2}[x_1,x_2]+L_{\mathcal{A}^{k+t-1}_2}[x_1,x_2] \bigg)R(x_3)
 = R\Bigg(R^2\Big( L_{\mathcal{B}^{k+t-1}_2}[x_1,x_2]+L_{\mathcal{A}^{k+t-1}_2}[x_1,x_2] \Big)x_3\Bigg)\\
 +R^2\Bigg(R\Big( L_{\mathcal{B}^{k+t-1}_2}[x_1,x_2]+L_{\mathcal{A}^{k+t-1}_2}[x_1,x_2] \Big)x_3\Bigg)
 + R^2\Bigg(\Big( L_{\mathcal{B}^{k+t-1}_2}[x_1,x_2] +  L_{\mathcal{A}^{k+t-1}_2}[x_1,x_2]  \Big)R(x_3)\Bigg)\\
 +\lambda R^2\Bigg(\Big( L_{\mathcal{B}^{k+t-1}_2}[x_1,x_2] + L_{\mathcal{A}^{k+t-1}_2}[x_1,x_2] \Big)x_3\Bigg) 
 +\lambda R\Bigg(R\Big( L_{\mathcal{B}^{k+t-1}_2}[x_1,x_2]+L_{\mathcal{A}^{k+t-1}_2}[x_1,x_2] \Big)x_3\Bigg)
 \\
 = R\Bigg(R^2\Big( L_{\mathcal{B}^{k+t-1}_2}[x_1,x_2]+L_{A^{k+t-1}_2}[x_1,x_2] \Big)x_3\Bigg)
 +R^2\Bigg( L_{\mathcal{B}^{k+t}_3}[x_1,x_2,x_3]\Bigg)\\
 +R^2\Bigg( L_{\mathcal{B}^{k+t-1}_2}[x_1,x_2R(x_3)] + L_{\mathcal{A}^{k+t}_3}[x_1,x_2,x_3]  
 +\lambda L_{\mathcal{A}^{k+t-1}_2}[x_1,x_2x_3]
 + L_{\mathcal{A}^{k+t-1}_2}[x_1,x_2R(x_3)]+\lambda L_{\mathcal{A}^{k+t-1}_3}[x_1,x_2,x_3]  \Bigg) \\
 +\lambda R^2\bigg(L_{\mathcal{B}^{k+t-1}_2}[x_1,x_2x_3] + L_{\mathcal{B}^{k+t-1}_3}[x_1,x_2,x_3]  \bigg)
 +\lambda R\bigg( L_{\mathcal{B}^{k+t}_3}[x_1,x_2,x_3]\bigg).
\end{multline*}

Further,
\begin{multline*}
  \lambda R \Bigg( L_{\mathcal{D}^{k+t-1}_2}[x_1,x_2]+L_{\mathcal{C}^{k+t-1}_2}[x_1,x_2]\Bigg)R(x_3)
  {=}\lambda R\Bigg(R\Big( L_{\mathcal{D}^{k+t-1}_2}[x_1,x_2]+L_{\mathcal{C}^{k+t-1}_2}[x_1,x_2]\Big)x_3\Bigg)\\
  +\lambda R\Bigg(\Big( L_{\mathcal{D}^{k+t-1}_2}[x_1,x_2]+L_{\mathcal{C}^{k+t-1}_2}[x_1,x_2]\Big)R(x_3)\Bigg)
  +\lambda^2 R\Bigg(\Big( L_{\mathcal{D}^{k+t-1}_2}[x_1,x_2]+L_{\mathcal{C}^{k+t-1}_2}[x_1,x_2]\Big)x_3\Big)\\
  = \lambda R\Bigg(R\Big( L_{\mathcal{D}^{k+t-1}_2}[x_1,x_2]+L_{\mathcal{C}^{k+t-1}_2}[x_1,x_2]\Big)x_3\Bigg) \\
  +\lambda R\Bigg( L_{\mathcal{D}^{k+t-1}_2}[x_1,x_2R(x_3)]+L_{\mathcal{C}^{k+t}_3}[x_1,x_2,x_3]
  +L_{C^{k+t-1}_2}[x_1,x_2R(x_3)]+\lambda L_{\mathcal{C}^{k+t-1}_2}[x_1,x_2x_3]
   +\lambda L_{\mathcal{C}^{k+t-1}_3}[x_1,x_2,x_3]\Bigg)\\
  +\lambda^2 R\Bigg( L_{\mathcal{D}^{k+t-1}_2}[x_1,x_2x_3]+L_{\mathcal{D}^{k+t-1}_3}[x_1,x_2,x_3]\Bigg).
\end{multline*}

As a~result, we obtain~$(R4)$ for $m=3$:
\begin{multline*}
    \sum\limits^l_{t=0}\lambda ^{(t)}\Bigg(R^2\Big(L_{\mathcal{B}^{k+t}_3}[\vec{x}]+L_{\mathcal{A}^{k+t}_3}[\vec{x}]
    +\lambda L_{\mathcal{B}^{k+t-1}_3}[\vec{x}]+\lambda L_{\mathcal{A}^{k+t-1}_3}[\vec{x}]\Big)\\
    +\lambda R\Big(L_{\mathcal{B}^{k+t}_3}[\vec{x}]\Big)
    +\lambda R\Big(L_{\mathcal{C}^{k+t}_3}[\vec{x}]+\lambda L_{\mathcal{C}^{k+t-1}_3}[\vec{x}]\Big)
    +\lambda^2 R\Big( L_{\mathcal{D}^{k+t-1}_3}[\vec{x}]\Big)\Bigg)=0;\quad |\vec{x}|=3.
\end{multline*}

\vspace{-0.25cm}

Assume that~$(R4)$ holds for $m$. Below we state the inductive step for~$m+1$ applying~\eqref{eq1}--\eqref{eq5}:
\begin{multline*}
R^2\Bigg(L_{\mathcal{B}^{k+t+m-3-i}_m}[\vec{x}]
 + L_{\mathcal{A}^{k+t+m-3-i}_m}[\vec{x}]\Bigg)R(x_{m+1}) \\
 = R\Bigg(R^2\Big(L_{\mathcal{B}^{k+t+m-3-i}_m}[\vec{x}]
 + L_{\mathcal{A}^{k+t+m-3-i}_m}[\vec{x}]\Big)x_{m+1}\Bigg) 
 + R^2\Bigg(R\Big(L_{\mathcal{B}^{k+t+m-3-i}_m}[\vec{x}] \\
 + L_{\mathcal{A}^{k+t+m-3-i}_m}[\vec{x}]\Big)x_{m+1}\Bigg) 
 + R^2\Bigg(\Big(L_{\mathcal{B}^{k+t+m-3-i}_m}[\vec{x}]
 + L_{\mathcal{A}^{k+t+m-3-i}_m}[\vec{x}]\Big)R(x_{m+1})\Bigg) \\
 + \lambda R^2 \Bigg(\Big(L_{\mathcal{B}^{k+t+m-3-i}_m}[\vec{x}]
 + L_{\mathcal{A}^{k+t+m-3-i}_m}[\vec{x}]\Big)x_{m+1}\Bigg) 
 + \lambda R\Bigg(R\Big(L_{\mathcal{B}^{k+t+m-3-i}_m}[\vec{x}]
 + L_{\mathcal{A}^{k+t+m-3-i}_m}[\vec{x}]\Big)x_{m+1}\Bigg) \\
 = R\Bigg(R^2\Big(L_{\mathcal{B}^{k+t+m-3-i}_m}[\vec{x}]
 + L_{A^{k+t+m-3-i}_m}[\vec{x}]\Big)x_{m+1}\Bigg) 
 + R^2\Bigg(L_{\mathcal{B}^{k+t+m-2-i}_{m+1}}[\vec{x},x_{m+1}]\Bigg) \\   
 + R^2\Bigg(L_{\mathcal{B}^{k+t+m-3-i}_m}[x_1,x_2,\ldots,x_mR(x_{m+1})]
 + L_{\mathcal{A}^{k+t+m-2-i}_{m+1}}[\vec{x},x_{m+1}]\\
 + \lambda L_{\mathcal{A}^{k+t+m-3-i}_m}[x_1,x_2,\ldots,x_mx_{m+1}]
 + \lambda L_{\mathcal{A}^{k+t+m-3-i}_{m+1}}[\vec{x},x_{m+1}] \\
 + L_{\mathcal{A}^{k+t+m-3-i}_m}[x_1,x_2,\ldots,x_mR(x_{m+1})]\Bigg) 
 + \lambda R^2 \Bigg(L_{\mathcal{B}^{k+t+m-3-i}_m}[x_1,x_2,\ldots,x_mx_{m+1}] \\
 + L_{\mathcal{B}^{k+t+m-3-i}_{m+1}}[\vec{x},x_{m+1}]\Bigg) 
 + \lambda R\Bigg(L_{\mathcal{B}^{k+t+m-2-i}_{m+1}}[\vec{x},x_{m+1}]\Bigg).
\end{multline*}

Further,
\begin{multline*}
 R\Bigg(L_{\mathcal{B}^{k+t+m-3-i}_m}[\vec{x}]\Bigg)R(x_{m+1})
  = R\Bigg(R\Big(L_{\mathcal{B}^{k+t+m-3-i}_m}[\vec{x}]\Big)x_{m+1}\Bigg) 
 + R\Bigg(\Big(L_{\mathcal{B}^{k+t+m-3-i}_m}[\vec{x}]\Big)R(x_{m+1})\Bigg) \\
  + \lambda R\Bigg(\Big(L_{\mathcal{B}^{k+t+m-3-i}_m}[\vec{x}]\Big)x_{m+1}\Bigg) 
 = R\Bigg(R\Big(L_{\mathcal{B}^{k+t+m-3-i}_m}[\vec{x}]\Big)x_{m+1}\Bigg) \\
 + R\Bigg(L_{\mathcal{B}^{k+t+m-3-i}_m}[x_1,x_2,\ldots,x_mR(x_{m+1})]\Bigg)
 + \lambda R\Bigg(L_{\mathcal{B}^{k+t+m-3-i}_m}[x_1,x_2,\ldots,x_mx_{m+1}]\Bigg);
\end{multline*}

\vspace{-0.15cm}

\begin{multline*} 
R\Bigg( L_{\mathcal{C}^{k+t+m-3-i}_m}[\vec{x}]\Bigg)R(x_{m+1})
 = R\Bigg(R\Big(L_{\mathcal{C}^{k+t+m-3-i}_m}[\vec{x}]\Big)x_{m+1}\Bigg) 
 + R\Bigg(\Big(L_{\mathcal{C}^{k+t+m-3-i}_m}[\vec{x}]\Big)R(x_{m+1})\Bigg)
 + \lambda R\Bigg(\Big(L_{\mathcal{C}^{k+t+m-3-i}_m}[\vec{x}]\Big)x_{m+1}\Bigg) \allowdisplaybreaks \\
 = R\Bigg(R\Big(L_{\mathcal{C}^{k+t+m-3-i}_m}[\vec{x}]\Big)x_{m+1}\Bigg) 
 + R\Bigg(L_{\mathcal{C}^{k+t+m-2-i}_{m+1}}[\vec{x},x_{m+1}] \\
 + L_{\mathcal{C}^{k+t+m-3-i}_{m}}[x_1,x_2,\ldots,x_{m}R(x_{m+1}] 
 + L_{\mathcal{C}^{k+t+m-3-i}_{m}}[x_1,x_2,\ldots,x_{m}x_{m+1}] \\
 + L_{\mathcal{C}^{k+t+m-3-i}_{m+1}}[\vec{x},x_{m+1}]\Bigg) 
 + \lambda R\Bigg(L_{\mathcal{D}^{k+t+m-3-i}_{m+1}}[\vec{x},x_{m+1}]\Bigg);
\end{multline*}

\vspace{-0.15cm}

\begin{multline*} 
R\Bigg(L_{\mathcal{D}^{k+t+m-4-i}_m}[\vec{x}]\Bigg)R(x_{m+1})
 = R\Bigg(R\Big(L_{\mathcal{D}^{k+t+m-4-i}_m}[\vec{x}]\Big)x_{m+1}\Bigg) 
 + R\Bigg(\Big(L_{\mathcal{D}^{k+t+m-4-i}_m}[\vec{x}]\Big)R(x_{m+1})\Bigg)
 + \lambda R\Bigg(\Big(L_{\mathcal{D}^{k+t+m-4-i}_m}[\vec{x}]\Big)x_{m+1}\Bigg) \\
 = R\Bigg(R\Big(L_{\mathcal{D}^{k+t+m-4-i}_m}[\vec{x}]\Big)x_{m+1}\Bigg) 
 R\Bigg(L_{\mathcal{D}^{k+t+m-4-i}_m}[x_1,x_2,\ldots,x_mR(x_{m+1})]\Bigg)
 + \lambda R\bigg(L_{\mathcal{D}^{k+t+m-4-i}_m}[x_1,x_2,\ldots,x_mx_{m+1}]\Bigg).
\end{multline*}

Denote $\vec{x} = [x_1,\ldots,x_m,x_{m+1}]$,
$\vec{x}\,' = [x_1,\ldots,x_m x_{m+1}]$,
$\vec{x}_0 = [x_1,\ldots,x_m]$,
$\vec{x}_R = [x_1,\ldots,x_m R(x_{m+1})]$.

Therefore
\begin{multline*}
R^2\Bigg(\sum\limits^l_{t=0}\sum\limits^{m-2}_{i=0}\tbinom{m-2}{i} \lambda ^{(t)}\lambda^i\Bigg(\Big( L_{\mathcal{B}^{k+t+m-2-i}_{m+1}}[\vec{x}]+ L_{\mathcal{A}^{k+t+m-2-i}_{m+1}}[\vec{x}]
+\lambda L_{\mathcal{B}^{k+t+m-3-i}_{m+1}}[\vec{x}] \\
+\lambda L_{\mathcal{A}^{k+t+m-3-i}_{m+1}}[\vec{x}]\Big)
 +\lambda R\Big(\sum\limits^l_{t=0}\sum\limits^{m-2}_{i=0}\tbinom{m-2}{i}\lambda ^{(t)}\lambda ^i\bigg(L_{\mathcal{B}^{k+t+m-2-i}_{m+1}}[\vec{x}]\Big)\Bigg)\\
 +\lambda R\Bigg(\sum\limits^l_{t=0}\sum\limits^{m-2}_{i=0}\tbinom{m-2}{i}\lambda ^{(t)}\lambda ^i\Big(L_{\mathcal{C}^{k+t+m-2-i}_{m+1}}[\vec{x}]+\lambda L_{\mathcal{C}^{k+t+m-3-i}_{m+1}}[\vec{x}]\Big)\Bigg) 
 +\lambda^2 R\Bigg(\sum\limits^l_{t=0}\sum\limits^{m-2}_{i=0}\tbinom{m-2}{i}\lambda ^{(t)}\lambda ^i\Big(L_{\mathcal{D}^{k+t+m-3-i}_{m+1}}[\vec{x}]\Big)\Bigg)\\
 +R^2\Bigg(\sum\limits^l_{t=0}\sum\limits^{m-2}_{i=0}\tbinom{m-2}{i}\lambda^{(t)}\lambda^i\Big(L_{\mathcal{B}^{k+t+m-3-i}_m}[\vec{x}_R]+L_{\mathcal{A}^{k+t+m-3-i}_m}[\vec{x}_R]\Big)\Bigg)\\
 +\lambda R\Bigg(\sum\limits^l_{t=0}\sum\limits^{m-3}_{i=0}\tbinom{m-3}{i}\lambda ^{(t)}\lambda ^i\Big(L_{\mathcal{B}^{k+t+m-3-i}_m}[\vec{x}_R]\Big)\Bigg) \\
 +\lambda R\Bigg(\sum\limits^l_{t=0}\sum\limits^{m-2}_{i=0}\tbinom{m-2}{i}\lambda ^{(t)}\lambda^i\Big(L_{\mathcal{C}^{k+t+m-3-i}_m}[\vec{x}_R]\Big)\Bigg)
 +\lambda^2 R\Bigg(\sum\limits^l_{t=0}\sum\limits^{m-3}_{i=0}\tbinom{m-3}{i}\lambda^{(t)}\lambda^i\Big(L_{\mathcal{D}^{k+t+m-4-i}_m}[\vec{x}_R]\Big)\Bigg)\\
 +\lambda R^2\Bigg(\sum\limits^l_{t=0}\sum\limits^{m-2}_{i=0}\tbinom{m-2}{i}\lambda^{(t)}\lambda^i\Big(L_{\mathcal{B}^{k+t+m-3-i}_m}[\vec{x}\,']+L_{\mathcal{A}^{k+t+m-3-i}_m}[\vec{x}\,']\Big)\Bigg)
 +\lambda^2 R\Bigg(\sum\limits^l_{t=0}\sum\limits^{m-3}_{i=0}\tbinom{m-3}{i}\lambda^{(t)}\lambda ^i\Big(L_{\mathcal{B}^{k+t+m-3-i}_m}[\vec{x}\,']\Big)\Bigg) \allowdisplaybreaks \\
 +\lambda^2 R\Bigg(\sum\limits^l_{t=0}\sum\limits^{m-2}_{i=0}\tbinom{m-2}{i}\lambda^{(t)}\lambda^i\Big(L_{\mathcal{C}^{k+t+m-3-i}_m}[\vec{x}\,']\Big)\Bigg)
 +\lambda^3 R\Bigg(\sum\limits^l_{t=0}\sum\limits^{m-3}_{i=0}\tbinom{m-2}{i}\lambda^{(t)}\lambda^i\Big(L_{\mathcal{D}^{k+t+m-4-i}_m}[\vec{x}\,']\Big)\Bigg) \\
 + R\Bigg(R^2\Big(\sum\limits^l_{t=0}\sum\limits^{m-2}_{i=0}   
 \tbinom{m-2}{i}\lambda^{(t)}\lambda^i\Big(L_{\mathcal{B}^{k+t+m-3-i}_m}[\vec{x}_0]+L_{\mathcal{A}^{k+t+m-3-i}_m}[\vec{x}_0]\Big)\Bigg)\\
 +\lambda R\Bigg(\sum\limits^l_{t=0}\sum\limits^{m-3}_{i=0}\tbinom{m-3}{i}\lambda^{(t)}\lambda^i\Big(L_{\mathcal{B}^{k+t+m-3-i}_m}[\vec{x}_0]\Big)\Bigg)
 +\lambda R\Bigg(\sum\limits^l_{t=0}\sum\limits^{m-2}_{i=0}\tbinom{m-2}{i}\lambda^{(t)}\lambda^i\Big(L_{\mathcal{C}^{k+t+m-3-i}_m}[\vec{x}_0]\Big)\Bigg)\\
 +\lambda^2 R\Bigg(\sum\limits^l_{t=0}\sum\limits^{m-3}_{i=0}\tbinom{m-3}{i}\lambda^{(t)}\lambda^i\Big(L_{\mathcal{D}^{k+t+m-4-i}_m}[\vec{x}_0]\Big)\Bigg)\Bigg)x_{m+1} 
 = 0.
\end{multline*}

As a result, we get the relation~$(R4)$:
\begin{multline*}
R^2\Bigg(\sum\limits^l_{t=0}\sum\limits^{m-1}_{i=0}\lambda ^{(t)}\lambda ^i\tbinom{m-1}{i}\Big(L_{\mathcal{B}^{k+t+m-2-i}_{m+1}}[\vec{x}]+L_{\mathcal{A}^{k+t+m-2-i}_{m+1}}[\vec{x}]\Big)\Bigg)\\
 +\lambda R\Bigg(\sum\limits^l_{t=0}\sum\limits^{m-2}_{i=0}\tbinom{m-2}{i}\lambda ^{(t)}\lambda ^i\Big(L_{\mathcal{B}^{k+t+m-2-i}_{m+1}}[\vec{x}]\Big)\Bigg)
 +\lambda R\Bigg(\sum\limits^l_{t=0}\sum\limits^{m-1}_{i=0}\tbinom{m-1}{i}\lambda ^{(t)}\lambda ^i\Big(L_{\mathcal{C}^{k+t+m-2-i}_{m+1}}[\vec{x}]\Big)\Bigg)\\
 +\lambda^2 R\Bigg(\sum\limits^l_{t=0}\sum\limits^{m-2}_{i=0}\tbinom{m-2}{i}\lambda ^{(t)}\lambda ^i\Big(L_{\mathcal{D}^{k+t+m-3-i}_{m+1}}[\vec{x}]\Big)\Bigg)=0;\quad |\vec{x}| = m+1.
\end{multline*}

\item Analogously, we may deduce the relations~$(R5)$ and~$(R6)$. 
\end{itemize}

The relations $(R7)$--$(R8)$ arise from compositions between $(R2)$ and $(R3)$.

Specifically:
\begin{itemize}
    \item It is easy to deduce the relation $(R7)$.
    
\item To derive the relation $(R8)$ by using \eqref{3.10} we have:
\end{itemize}


\begin{multline*}
R(x)\bigg( L_{\mathcal{B}^n_{m}}[\vec{x}]+ L_{\mathcal{A}^n_{m}}[\vec{x}]\bigg) 
 =  L^*_{\mathcal{B}^n_{m}}[x,\vec{x}]+ L^*_{\mathcal{A}^n_{m}}[x,\vec{x}] \\
 +  L_{\mathcal{B}^n_{m}}[R(x)x_1+xR(x_1),x_2,\ldots,x_{m}]+ L_{\mathcal{A}^n_{m}}[R(x)x_1+xR(x_1),x_2,\ldots,x_{m}] \\
 + \lambda  L_{\mathcal{B}^n_{m}}[xx_1,x_2,\ldots,x_{m}]+\lambda  L_{\mathcal{A}^n_{m}}[xx_1,x_2,\ldots,x_{m}].
\end{multline*}

Hence, we get
\begin{multline*}
R(x)R^2\Bigg( L_{\mathcal{B}^n_{m}}[\vec{x}]+ L_{\mathcal{A}^n_{m}}[\vec{x}]\Bigg)
 = R\Bigg(xR^2\Big( L_{\mathcal{B}^n_{m}}[\vec{x}]+ L_{\mathcal{A}^n_{m}}[\vec{x}]\Big)\Bigg) \\ 
 + R^2\Bigg(xR\Big( L_{\mathcal{B}^n_{m}}[\vec{x}]+ L_{\mathcal{A}^n_{m}}[\vec{x}]\Big)\Bigg) 
 + R^2\Bigg(R(x)\Big( L_{\mathcal{B}^n_{m}}[\vec{x}]+L_{\mathcal{A}^n_{m}}[\vec{x}]\Big)\Bigg) \\
 + \lambda R^2\Bigg(x\Big( L_{\mathcal{B}^n_{m}}[\vec{x}]+ L_{\mathcal{A}^n_{m}}[\vec{x}]\Big)\Bigg) 
 + \lambda R\Bigg(x R\Big( L_{\mathcal{B}^n_{m}}[\vec{x}]+ L_{\mathcal{A}^n_{m}}[\vec{x}]\Big)\Bigg) \allowdisplaybreaks \\
 = R\bigg(xR^2\bigg( L_{\mathcal{B}^n_{m}}[\vec{x}]+ L_{\mathcal{A}^n_{m}}[\vec{x}]\bigg)\bigg) 
 + R^2\Bigg(xR\Big( L_{\mathcal{B}^n_{m}}[\vec{x}]+ L_{\mathcal{A}^n_{m}}[\vec{x}]\Big)\Bigg) 
 + R^2\Bigg( L^*_{\mathcal{B}^n_{m}}[x,\vec{x}]+ L^*_{\mathcal{A}^n_{m}}[x,\vec{x}]\Bigg) \\
 + R^2\Bigg( L_{\mathcal{B}^n_{m}}[R(x)x_1+xR(x_1),x_2,\ldots,x_{m}]+ L_{\mathcal{A}^n_{m}}[R(x)x_1+xR(x_1),x_2,\ldots,x_{m}]\Bigg) \\
 + \lambda R^2\Bigg( L_{\mathcal{B}^n_{m}}[xx_1,x_2,\ldots,x_{m}]+ L_{\mathcal{A}^n_{m}}[xx_1,x_2,\ldots,x_{m}]\Bigg) \\
 + \lambda R^2\Bigg(x\Big(L_{\mathcal{B}^n_{m}}[\vec{x}]+ L_{A^n_{m}}[\vec{x}]\Big)\Bigg) 
 + \lambda R\Bigg(x R\Big( L_{\mathcal{B}^n_{m}}[\vec{x}]+ L_{\mathcal{A}^n_{m}}[\vec{x}]\Big)\Bigg).
\end{multline*}

The next equality holds, if one changes all appearances of $ \mathcal{B}^n_m$
on $\mathcal{D}^n_m$ or on $  \mathcal{C}^n_m$:
\begin{multline*}
R(x)R\Bigg( L_{\mathcal{B}^n_{m}}[\vec{x}]\Bigg) 
 = R\Bigg (R(x)\Big(L_{\mathcal{B}^n_{m}}[\vec{x}]\Big)\Bigg)
  +R\Bigg (x R\Big( L_{\mathcal{B}^n_{m}}[\vec{x}]\Big)\Bigg)
 +\lambda R\Bigg (x \Big(L_{\mathcal{B}^n_{m}}[\vec{x}]\Big)\Bigg)\\
  = R\Bigg( L^*_{\mathcal{B}^n_{m}}[x,\vec{x}]
    + L_{\mathcal{B}^n_{m}}[R(x)x_1+xR(x_1),x_2,\ldots,x_{m}]
    +\lambda  L_{\mathcal{B}^n_{m}}[xx_1,x_2,\ldots,x_{m}]\Bigg)
     +R\Bigg (x R\Big( L_{\mathcal{B}^n_{m}}[\vec{x}]\Big)\Bigg)
 +\lambda R\Bigg (x \Big( L_{\mathcal{B}^n_{m}}[\vec{x}]\Big)\Bigg).
\end{multline*}

Taking into account the relation~$(R4)$, we get
$\sum\limits^l_{t=0}\lambda^{(t)} M_t$, where
\begin{multline*}
M_t
 = R^2\Bigg(\sum\limits^{m-2}_{i=0} \tbinom{m-2}{i}\lambda^i\Big(L^*_{\mathcal{B}^{k+t+m-3-i}_m}[x,\vec{x}]+L^*_{\mathcal{A}^{k+t+m-3-i}_m}[x,\vec{x}]\Big)\Bigg) \\
 + \lambda R\Bigg(\sum\limits^{m-3}_{i=0}\tbinom{m-3}{i}\lambda ^i\Big(L^*_{\mathcal{B}^{k+t+m-3-i}_m}[x,\vec{x}]\Big)\Bigg) 
 + \lambda R\Bigg(\sum\limits^{m-2}_{i=0}\tbinom{m-2}{i}\lambda ^i\Big(L^*_{\mathcal{C}^{k+t+m-3-i}_m}[x,\vec{x}]\Big)\Bigg) \\
 + \lambda^2 R\Bigg(\sum\limits^{m-3}_{i=0}\tbinom{m-3}{i}\lambda ^i\Big(L^*_{\mathcal{D}^{k+t+m-4-i}_m}[x,\vec{x}]\Big)\Bigg) 
 + R^2\Bigg(\sum\limits^{m-2}_{i=0}\tbinom{m-2}{i}\lambda^{(t)}\lambda^i\Bigg(xR\Big( L_{\mathcal{B}^{k+t+m-3-i}_{m}}[\vec{x}]+ L_{\mathcal{A}^{k+t+m-3-i}_{m}}[\vec{x}]\Big)\Bigg) \\
 + \lambda R^2\Bigg(\sum\limits^{m-2}_{i=0}\tbinom{m-2}{i}\lambda^i\Bigg(x\Big( L_{\mathcal{B}^{k+t+m-3-i}_{m}}[\vec{x}]+ L_{\mathcal{A}^{k+t+m-3-i}_{m}}[\vec{x}]\Big)\Bigg) 
 + \lambda R\Bigg(\sum\limits^{m-2}_{i=0}\tbinom{m-2}{i}\lambda ^i\Bigg(x R\Big( L_{\mathcal{B}^{k+t+m-3-i}_{m}}[\vec{x}]+ L_{\mathcal{A}^{k+t+m-3-i}_{m}}[\vec{x}]\Big)\Bigg) \allowdisplaybreaks \\
 + \lambda^2 R\Bigg(\sum\limits^{m-3}_{i=0}\tbinom{m-3}{i}\lambda^i\Big(x  L_{\mathcal{B}^{k+t+m-3-i}_{m}}[\vec{x}]\Big)\Bigg) 
 + \lambda^2 R\Bigg(\sum\limits^{m-2}_{i=0}\tbinom{m-2}{i}\lambda^i\Big(xL_{\mathcal{C}^{k+t+m-3-i}_{m}}[\vec{x}]\Big)\Bigg) \\
 + \lambda^3 R\Bigg(\sum\limits^{m-3}_{i=0}\tbinom{m-3}{i}\lambda^i\Big(x  L_{\mathcal{D}^{k+t+m-3-i}_{m}}[\vec{x}]\Big)\Bigg)=0;\quad |\vec{x}|=m.
\end{multline*}
The last expression is exactly the relation $(R8)$. We can prove that the composition between the relations $(R2)$ and $(R6)$ gives us the relations $(R9)$
and $(R10)$ by using \eqref{3.7}-\eqref{3.9} and the same way as above.
With  similar arguments, we do not get new compositions (see Appendix p.15) 

\end{proof}
\begin{example}
Let $\lambda=0$. Then the following relations form a Gr\"{o}bner---Shirshov basis in $R\As\langle X\rangle$.
\[
   \begin{aligned}
&(R1)\quad R(x)R(y) = R( xR(y) + R(x)y),\\
&(R2)\quad R^{k+l}(x) = 0,\\
&(R3)\quad
 R^2\bigg(L_{\mathcal{B}^{k+l+m-3}_m}[\vec{x}]
 + L_{\mathcal{A}^{k+l+m-3}_m}[\vec{x}]\bigg)= 0; \ m = |\vec{x}|\geq 2, \\
&(R4)\quad
R^2\bigg(L^\op_{\mathcal{B}^{k+l+m-3}_m}[\vec{x}]
 + L^\op_{\mathcal{A}^{k+l+m-3}_m}[\vec{x}]\bigg)=0; \ m = |\vec{x}|\geq2,\\
&(R5)\quad
 R^2 \Bigg(xR\Big( L_{\mathcal{B}^{k+l-3}_m}[\vec{x}]
 + L_{\mathcal{A}^{k+l-3}_m}[\vec{x}]\Big)\Bigg)
 + R^2\bigg( L^*_{\mathcal{B}^{k+l-3}_m}[x,\vec{x}]+L^*_{\mathcal{A}^{k+l-1}_m}[x,\vec{x}]\bigg)
  {=} 0;\ m=|\vec{x}| {=}2,\\
&(R6)\quad
R^2 \Bigg(R\Big( L^\op_{\mathcal{B}^{k+l-3}_m}[\vec{x}]+L^\op_{\mathcal{A}^{k+l-3-1}_m}[\vec{x}]\Big)x\Bigg)
 + R^2\Bigg( {L^*}^\op_{\mathcal{B}^{k+l-3}_m}[x,\vec{x}]+{L^*}^\op_{\mathcal{A}^{k+l-3}_m}[x,\vec{x}]\Bigg) {=} 0;\ m=|\vec{x}|{\geq}2. 
\end{aligned}
\]
 In case $k+l=2$, we have $I_2=\{1\}$, so $i_m,j_{m-1}\neq 1$ and $A_m=B_m=\varnothing$. Therefore, GSB is:
    \[
    R(x)R(y)=R(xR(y)+R(x)y),\ R^2(x)=0\\
    \]
  In case $k+l=3$, we have $I_3=\{1,2\}$, and $A_m=\varnothing$,$B_m=\{\underbrace{1,1,\ldots,1}_{m-2},2\}$. Therefore, GSB is: 
\[  
\begin{aligned}
&(R1)\quad R(x)R(y)=R(xR(y)+R(x)y),\ R^3(x)=0,\\
&(R2)\quad    R^2\Bigg(\underbrace{R(\ldots(R}_{m-2}(R^2(x_1)x_2)x_3)\ldots)x_m\Bigg)=0,\\
& (R3)\quad   R^2\Bigg(x_1R(\ldots(x_{m-2}R(x_{m-1}R^2(x_m)))\ldots)\Bigg)=0, \\
&(R4)\quad    R^2\Bigg(x\underbrace{R(\ldots R(R}_{m-2}(R^2(x_1)x_2)x_3)\ldots)+R(x\underbrace{R(\ldots R(R}_{m-3}(R^2(x_1)x_2)x_3)\ldots x_{m-1}))x_m+\\
&    +\ldots+\underbrace{R(\ldots R(R}_{m-2}(xR^2(x_1))x_2)\ldots )x_m\Bigg)=0,\\
 &(R5)\quad R^2\Bigg(R\Big(x_1R(\ldots(x_{m-2}R(x_{m-1}R^2(x_m)))\ldots\Big)x+x_1R\Big(\ldots(x_{m-2}R(x_{m-1}R^2(x_m)))\ldots )x\Big)+\\
 &+\ldots+x_1R\Big(\ldots(x_{m-2}R(x_{m-1}(R^2(x_m)x))\ldots)\Bigg)=0.
\end{aligned}
\]
\end{example}
\begin{example}
Let $\lambda\neq 0$ and $k+l=3$ (for simplicity, we choose $k=2$, $l=1$). 
The following relations form a Gr\"{o}bner---Shirshov basis in $R\As\langle X\rangle$.
\[
\begin{aligned}
&(R1)\quad R(x)R(y)=R(xR(y)+R(x)y+\lambda xy),\\
&(R2)\quad R^3(x)+\lambda R^2(x)=0,\\
&(R3)\quad R^2\Bigg(R^2(x_1)x_2\Bigg)+\lambda R \Bigg( R^2(x_1)x_2+R(R(x_1)x_2)\Bigg)+\lambda^2 R\Bigg(R(x_1)x_2\Bigg)=0,\\
&(R4)\quad   R^2\Bigg(\underbrace{R(\ldots(R}_{m-2}(R^2(x_1)x_2)x_3)\ldots)x_m\Bigg)+\lambda R\Bigg(\underbrace{R(\ldots(R}_{m-2}(R^2(x_1)x_2)x_3)\ldots)x_m\Bigg)\\
&+\lambda R\Bigg(\underbrace{R(\ldots (R(R}_{m}(x_1)x_2)x_3)\ldots)x_m)\Bigg)
+\lambda^2 R\Bigg(\underbrace{R(\ldots (R(R}_{m-1}(x_1)x_2)x_3)\ldots)x_m\Bigg)=0,\\
&(R5)\quad R^2\Bigg(x_1R^2(x_2)\Bigg)+\lambda R \Bigg( x_1R^2(x_2)+R(x_1R(x_2))\Bigg)+\lambda^2 R\Bigg(x_1R(x_1)\Bigg)=0, \\
& (R6)\quad   R^2\Bigg(x_1R(\ldots(x_{m-2}R(x_{m-1}R^2(x_m)))\ldots)\Bigg)+\lambda R\Bigg(x_1R(\ldots(x_{m-2}R(x_{m-1}R(x_m)))\ldots)\Bigg)\\
&+\lambda R\Bigg(R(x_1R(\ldots(x_{m-2}R(x_{m-1}R^2(x_m)))\ldots)\Bigg)+\lambda^2 R\Bigg(x_1R(\ldots(x_{m-2}R(x_{m-1}R(x_m)))\ldots)\Bigg)=0,\\
&
(R7) \quad 
 R^2 \bigg(xR\bigg( R^2(x_1)x_2\bigg)\bigg)
 + \lambda R^2\bigg( x\bigg( R^2(x_1)x_2\bigg)\bigg) 
+ \lambda R\bigg( xR\bigg( R^2(x_1)x_2\bigg)\bigg)
 + R^2\bigg( R\big(x R^2(x_1)\big)x_2\bigg) \\
& + \lambda R\bigg( R\big(xR(x_1))x_2\big)\bigg)
 + \lambda^2 R\bigg(xR(x_1)x_2\big)
 \bigg)\bigg) = 0, \\
&(R8)  \quad
 R^2 \bigg(xR\bigg(\underbrace{R(\ldots(R}_{m-2}(R^2(x_1)x_2)x_3)\ldots)x_m\bigg)\bigg)
 + \lambda R^2\bigg( x\bigg( \underbrace{R(\ldots(R}_{m-2}(R^2(x_1)x_2)x_3)\ldots)x_m\bigg)\bigg) \\
& + \lambda^2 R\bigg( xR\bigg( \underbrace{R(\ldots(R}_{m-2}(R^2(x_1)x_2)x_3)\ldots)x_m\bigg)\bigg)+R^2\bigg( R\bigg(x \underbrace{R(\ldots(R}_{m-2}(R^2(x_1)x_2)x_3)\ldots\bigg)x_m\bigg)\bigg)\\
&+\ldots+ R^2\bigg( R\bigg (\underbrace{R(\ldots(R}_{m-2}(xR^2(x_1))x_2)x_3)\ldots)x_m\bigg)\bigg)
 +\lambda R\bigg( R\bigg(x \underbrace{R(\ldots(R}_{m-2}(R^2(x_1)x_2)x_3)\ldots)x_m\bigg)\bigg)\\
 &+\ldots+\lambda R\bigg( R\bigg (\underbrace{R(\ldots(R}_{m-2}(xR^2(x_1))x_2)x_3)\ldots)x_m\bigg)\bigg) 
+\lambda^2 R\bigg( x\bigg( \underbrace{R(\ldots(R}_{m-2}(R^2(x_1)x_2)x_3)\ldots)x_m\bigg)\bigg) \\
& +\lambda R\bigg( R\bigg(x \underbrace{R(\ldots(R}_{m-1}(R(x_1)x_2)x_3)\ldots)x_m\bigg)\bigg)+\ldots+\lambda R\bigg( R\bigg (\underbrace{R(\ldots(R}_{m-1}(xR(x_1))x_2)x_3)\ldots)x_m\bigg)\bigg)\\
&+\lambda^2 R\bigg( x \underbrace{R\bigg(\ldots(R}_{m-1}(R(x_1)x_2)x_3)\ldots)x_m\bigg)\bigg)  +\lambda^2 R\bigg( R\bigg(x \underbrace{R(\ldots(R}_{m-2}(R(x_1)x_2)x_3)\ldots)x_m\bigg)\bigg)\\
&+\ldots+\lambda^2 R\bigg( R\bigg (\underbrace{R(\ldots(R}_{m-2}(xR(x_1))x_2)x_3)\ldots)x_m\bigg)\bigg)+\lambda^3 R\bigg( x\bigg( \underbrace{R(\ldots(R}_{m-2}(R^2(x_1)x_2)x_3)\ldots)x_m\bigg)\bigg)=0, \\
&(R9) \quad 
 R^2 \bigg(R\bigg( x_1R^2(x_2)\bigg)x\bigg)
 + \lambda R^2\bigg( \bigg( x_1R^2(x_2)\bigg)x\bigg) 
+ \lambda R\bigg( R\bigg( x_1R^2(x_2)\bigg)x\bigg)
 + R^2\bigg(x_1 R\big( R^2(x_2)x\big)\bigg) \\
& + \lambda R\bigg(x_1 R\big(R(x_2)x\big)\bigg)
 + \lambda^2 R\bigg(x_1R(x_2)x\big)
 \bigg)\bigg) = 0, 
 \end{aligned}
\]
\[
\begin{aligned}
 & (R10)\quad   R^2\Bigg(R\Big(x_1R(\ldots(x_{m-2}R(x_{m-1}R^2(x_m)))\Big)x\Bigg)+\lambda   R^2\Bigg(\Big(x_1R(\ldots(x_{m-2}R(x_{m-1}R^2(x_m)))\Big)x\Bigg)\\
 &+\lambda^2   R\Bigg(R\Big(x_1R(\ldots(x_{m-2}R(x_{m-1}R^2(x_m)))\Big)x\Bigg)+R^2\Bigg(R\Big(x_1R(\ldots(x_{m-2}R(x_{m-1}R(R^2(x_m)x))\Big)\Bigg)\\
&+\ldots+R^2\Bigg(R\Big(xR\big(x_1R(\ldots(x_{m-2}R(x_{m-1}R^2(x_m)))\Big)\Bigg)+\lambda R\Bigg(x_1R(\ldots(x_{m-2}R(x_{m-1}R(R^2(x_m)x))\Big)\Bigg) \\
&+\ldots+\lambda R\Bigg(xR\big(x_1R(\ldots(x_{m-2}R(x_{m-1}R^2(x_m)))\Big)\Bigg)+\lambda^2   R\Bigg(\Big(x_1R(\ldots(x_{m-2}R(x_{m-1}R^2(x_m)))\Big)x\Bigg)\\
&+\lambda R\Bigg(R\Big(x_1R(\ldots(x_{m-2}R(x_{m-1}R(R(x_m)x))\Big)\Bigg)+\ldots+\lambda R\Bigg(R\Big(xR\big(x_1R(\ldots(x_{m-2}R(x_{m-1}R(x_m)))\Big)\Bigg)\\
&+\lambda^2   R\Bigg(R\Big(x_1R(\ldots(x_{m-2}R(x_{m-1}R(x_m)))\Big)x\Bigg)+\lambda^2 R\Bigg(x_1R(\ldots(x_{m-2}R(x_{m-1}R(R(x_m)x))\Big)\Bigg)\\
&+\ldots+\lambda^2 R\Bigg(xR\big(x_1R(\ldots(x_{m-2}R(x_{m-1}R(x_m)))\Big)\Bigg)+\lambda^3   R\Bigg(x_1R\Big(\ldots(x_{m-2}R(x_{m-1}R(x_m)))\Big)x\Bigg)=0.
     \end{aligned}
     \]
\end{example}

\newpage

\section{Appendix}

when we found the composition between the relations $(R1)$ and $(R2)$, we noticed that the idea of composition is to rewrite the relation $f.R(x_2)$ via $g$ for $\vec{x}=[x_1,x_2]$ and repeat this process for $\vec{x}=[x_1,x_2,\ldots,x_m]$ and the results are relations $(R4)$-$(R10)$. 

Since each term in the relations $(R4)$-$(R10)$ has the form $R^n(T)$, where $T$ is a sum of words of the form $R^{i_1}\Big(R^{i_2}(\ldots(R^{i_{m-1}}(R^{i_m}(x_1)x_2)x_3)\ldots)x_m\Big)$
or $R^{i_1}\Big(x_1R^{i_2}(\ldots(x_{m-2}R^{i_{m-1}}(x_{m-1}R^{i_m}(x_m)))\ldots)\Big)$ and the leader part in $(R4)$ is $R^{k+l-m+1}\Big(R(\ldots(R(R(x_1)x_2)x_3)\ldots)x_{m-1}\Big)x_m$
and in $(R5)$ is \\$R^{k+l-m+1}\Big(x_1R(\ldots(x_{m-2}R(x_{m-1}R(x_m)))\ldots)\Big)$, so the composition has form $R^n(T)R(T)$ and we get relation $(R4)$ by the same way, therefore the composition between the relations $(R4)$ and $(R5)$ gives a sum of the expressions, which are consequences of $(R4)$ itself.
Below, we show this example for
$\sum\limits_{t=0}^l \sum\limits^{m-1}_{i=0}\lambda^{(t)}\lambda^i
M_{t,i}$, where
\[
\begin{aligned}
M_{t,i}&=\tbinom{m-1}{i} R^2\bigg(L_{\mathcal{B}^{k+t+m-2-i}_{m+1}}[\vec{x},\sum\limits^l_{t=0}\lambda ^{(t)}\bigg(R\bigg( L^\op_{\mathcal{B}^{k+t-1}_2}[\vec{y}]+L^{\op}_{\mathcal{A}^{k+t-1}_2}[\vec{y}]\bigg)
+\lambda \bigg( L^{\op}_{\mathcal{B}^{k+t-1}_2}[\vec{y}]+L^\op_{\mathcal{C}^{k+t-1}_2}[\vec{y}]\bigg)\bigg)]
\\
& + L_{\mathcal{A}^{k+t+m-2-i}_{m+1}}[\vec{x},\sum\limits^l_{t=0}\lambda ^{(t)}\bigg(R\bigg( L^\op_{\mathcal{B}^{k+t-1}_2}[\vec{y}]+L^{\op}_{\mathcal{A}^{k+t-1}_2}[\vec{y}]\bigg)
+\lambda \bigg( L^{\op}_{\mathcal{B}^{k+t-1}_2}[\vec{y}]+L^\op_{\mathcal{C}^{k+t-1}_2}[\vec{y}]\bigg)\bigg)]\bigg) 
 \\
& + \lambda \tbinom{m-2}{i}R\bigg(L_{\mathcal{B}^{k+t+m-2-i}_{m+1}}[\vec{x},\sum\limits^l_{t=0}\lambda ^{(t)}\bigg(R\bigg( L^\op_{\mathcal{B}^{k+t-1}_2}[\vec{y}]+L^{\op}_{\mathcal{A}^{k+t-1}_2}[\vec{y}]\bigg)\\
&+\lambda \bigg( L^{\op}_{\mathcal{B}^{k+t-1}_2}[\vec{y}]+L^\op_{\mathcal{C}^{k+t-1}_2}[\vec{y}]\bigg)\bigg)]\bigg)
+ \lambda \tbinom{m-1}{i} R\bigg(L_{\mathcal{C}^{k+t+m-2-i}_{m+1}}[\vec{x},\sum\limits^l_{t=0}\lambda ^{(t)}\bigg(R\bigg( L^\op_{\mathcal{B}^{k+t-1}_2}[\vec{y}]+L^{\op}_{\mathcal{A}^{k+t-1}_2}[\vec{y}]\bigg)\\
&+\lambda \bigg( L^{\op}_{\mathcal{B}^{k+t-1}_2}[\vec{y}]+L^\op_{\mathcal{C}^{k+t-1}_2}[\vec{y}]\bigg)\bigg)]\bigg) 
  + \lambda^2 \tbinom{m-2}{i} R\bigg(L_{\mathcal{D}^{k+t+m-3-i}_{m+1}}[\vec{x},\sum\limits^l_{t=0}\lambda ^{(t)}\bigg(R\bigg( L^\op_{\mathcal{B}^{k+t-1}_2}[\vec{y}]+L^{\op}_{\mathcal{A}^{k+t-1}_2}[\vec{y}]\bigg)\\
&+\lambda \bigg( L^{\op}_{\mathcal{B}^{k+t-1}_2}[\vec{y}]+L^\op_{\mathcal{C}^{k+t-1}_2}[\vec{y}]\bigg)\bigg)]\bigg) = 0, \ m = |\vec{x}|\geq 3,\ |\vec{y}|=2.
\end{aligned}
\]

Composition between relations $(R4)$ and $(R6)$, we get a sum of the expressions, which are consequences of $(R4)$ itself.
Below, we show this example for  $\sum\limits_{t=0}^l \sum\limits^{m-1}_{i=0}\lambda^{(t)}\lambda^i
M_{t,i}$, where

\[
\begin{aligned}
M_{t,i}&=\tbinom{m-1}{i} R^2\bigg( L_{\mathcal{B}^{k+t+m-2-i}_{m+1}}[\vec{x},\sum\limits^l_{t=0}\sum\limits^{m-2}_{i=0}\lambda^{(t)}\lambda^i\left(
\tbinom{m-2}{i}R\bigg(  L^\op_{\mathcal{B}^{k+t+m-3-i}_m}[\vec{x}]
 +    L^\op_{\mathcal{A}^{k+t+m-3-i}_m}[\vec{x}]\bigg) \right. \\
& + \lambda \tbinom{m-3}{i} \bigg(  L^\op_{\mathcal{B}^{k+t+m-3-i}_m}[\vec{x}]\bigg) 
 + \lambda \tbinom{m-2}{i}\bigg(    L^\op_{\mathcal{C}^{k+t+m-3-i}_m}[\vec{x}]\bigg)
 \left. +\lambda^2 \tbinom{m-3}{i} \bigg(   L^\op_{\mathcal{D}^{k+t+m-4-i}_m}[\vec{x}]\bigg)\right)]
 \\
 &+       L_{\mathcal{A}^{k+t+m-2-i}_{m+1}}[\vec{x},\sum\limits^l_{t=0}\sum\limits^{m-2}_{i=0}\lambda^{(t)}\lambda^i\left(
\tbinom{m-2}{i}R\bigg(  L^\op_{\mathcal{B}^{k+t+m-3-i}_m}[\vec{x}]
 +    L^\op_{\mathcal{A}^{k+t+m-3-i}_m}[\vec{x}]\bigg) \right. \\
& + \lambda \tbinom{m-3}{i} \bigg(  L^\op_{\mathcal{B}^{k+t+m-3-i}_m}[\vec{x}]\bigg) 
 + \lambda \tbinom{m-2}{i}\bigg(    L^\op_{\mathcal{C}^{k+t+m-3-i}_m}[\vec{x}]\bigg)
 \left. +\lambda^2 \tbinom{m-3}{i} \bigg(   L^\op_{\mathcal{D}^{k+t+m-4-i}_m}[\vec{x}]\bigg)\right)]\bigg)  \\
 \end{aligned}
\]
\[
\begin{aligned}
& + \lambda \tbinom{m-2}{i}R\bigg( L_{\mathcal{B}^{k+t+m-2-i}_{m+1}}[\vec{x},\sum\limits^l_{t=0}\sum\limits^{m-2}_{i=0}\lambda^{(t)}\lambda^i\left(
\tbinom{m-2}{i}R\bigg(  L^\op_{\mathcal{B}^{k+t+m-3-i}_m}[\vec{x}]
 +    L^\op_{\mathcal{A}^{k+t+m-3-i}_m}[\vec{x}]\bigg) \right. \\
& + \lambda \tbinom{m-3}{i} \bigg(  L^\op_{\mathcal{B}^{k+t+m-3-i}_m}[\vec{x}]\bigg) 
 + \lambda \tbinom{m-2}{i}\bigg(    L^\op_{\mathcal{C}^{k+t+m-3-i}_m}[\vec{x}]\bigg)
 \left. +\lambda^2 \tbinom{m-3}{i} \bigg(   L^\op_{\mathcal{D}^{k+t+m-4-i}_m}[\vec{x}]\bigg)\right)]\bigg)
 \\
& + \lambda \tbinom{m-1}{i} R\bigg( L_{\mathcal{C}^{k+t+m-2-i}_{m+1}}[\vec{x},\sum\limits^l_{t=0}\sum\limits^{m-2}_{i=0}\lambda^{(t)}\lambda^i\left(
\tbinom{m-2}{i}R\bigg(  L^\op_{\mathcal{B}^{k+t+m-3-i}_m}[\vec{x}]
 +    L^\op_{\mathcal{A}^{k+t+m-3-i}_m}[\vec{x}]\bigg) \right. \\
& + \lambda \tbinom{m-3}{i} \bigg(  L^\op_{\mathcal{B}^{k+t+m-3-i}_m}[\vec{x}]\bigg) 
 + \lambda \tbinom{m-2}{i}\bigg(    L^\op_{\mathcal{C}^{k+t+m-3-i}_m}[\vec{x}]\bigg) 
 \left. +\lambda^2 \tbinom{m-3}{i} \bigg(   L^\op_{\mathcal{D}^{k+t+m-4-i}_m}[\vec{x}]\bigg)\right)]\bigg) 
 \\
 & + \lambda^2 \tbinom{m-2}{i} R\bigg(  L_{\mathcal{D}^{k+t+m-3-i}_{m+1}}[\vec{x},\sum\limits^l_{t=0}\sum\limits^{m-2}_{i=0}\lambda^{(t)}\lambda^i\left(
\tbinom{m-2}{i}R\bigg(  L^\op_{\mathcal{B}^{k+t+m-3-i}_m}[\vec{x}]
 +    L^\op_{\mathcal{A}^{k+t+m-3-i}_m}[\vec{x}]\bigg) \right. \\
& + \lambda \tbinom{m-3}{i} \bigg(  L^\op_{\mathcal{B}^{k+t+m-3-i}_m}[\vec{x}]\bigg) 
 + \lambda \tbinom{m-2}{i}\bigg(    L^\op_{\mathcal{C}^{k+t+m-3-i}_m}[\vec{x}]\bigg) 
 \left. +\lambda^2 \tbinom{m-3}{i} \bigg(   L^\op_{\mathcal{D}^{k+t+m-4-i}_m}[\vec{x}]\bigg)\right)]\bigg)\\
 &= 0, \ m = |\vec{x}|\geq 3.
\end{aligned}
\]

\vspace{-0.3cm}
Composition between relations $(R4)$ and $(R7)$, we get a sum of the expressions, which are consequences of $(R4)$ itself.
Below, we show this example for
  $\sum\limits_{t=0}^l \sum\limits^{m-1}_{i=0}\lambda^{(t)}\lambda^i
M_{t,i}$, where
\[
\begin{aligned}
M_{t,i}&=\tbinom{m-1}{i} R^2\bigg( L_{\mathcal{B}^{k+t+m-2-i}_{m+1}}[\vec{x}, \sum\limits^l_{t=0}\lambda^{(t)}\bigg(
 R \bigg(xR\bigg(  L_{\mathcal{B}^{k+t-1}_2}[\vec{y}]
 +L_{\mathcal{A}^{k+t-1}_2}[\vec{y}]\bigg)\bigg)
 + \lambda R\bigg( x\bigg(  L_{\mathcal{B}^{k+t-1}_2}[\vec{y}]
 +L_{\mathcal{A}^{k+t-1}_2}[\vec{y}]\bigg)\bigg) \\
& + \lambda \bigg( xR\bigg(  L_{\mathcal{B}^{k+t-1}_2}[\vec{y}]
 + L_{\mathcal{A}^{k+t-1}_2}[\vec{y}]\bigg)\bigg)
 + R\bigg(   L^*_{\mathcal{B}^{k+t-1}_3}[x,\vec{y}]+    L^*_{\mathcal{A}^{k+t-1}_3}[x,\vec{y}]\bigg) 
 + \lambda \bigg(   L^*_{\mathcal{B}^{k+t-1}_3}[x,\vec{y}]
 +    L^*_{\mathcal{C}^{k+t-1}_3}[x,\vec{y}]\bigg)\\
 &+ \lambda^2 \bigg(x\bigg(  L_{\mathcal{B}^{k+t-1}_2}[\vec{y}]
 + L_{\mathcal{C}^{k+t-1}_2}[\vec{y}]\bigg)\bigg)\bigg) ]
 +       L_{\mathcal{A}^{k+t+m-2-i}_{m+1}}[\vec{x}, \sum\limits^l_{t=0}\lambda^{(t)}\bigg(
 R \bigg(xR\bigg(  L_{\mathcal{B}^{k+t-1}_2}[\vec{y}]
 + L_{\mathcal{A}^{k+t-1}_2}[\vec{y}]\bigg)\bigg)\\
& + \lambda R\bigg( x\bigg(  L_{\mathcal{B}^{k+t-1}_2}[\vec{y}]
 + L_{\mathcal{A}^{k+t-1}_2}[\vec{y}]\bigg)\bigg) 
 + \lambda \bigg( xR\bigg(  L_{\mathcal{B}^{k+t-1}_2}[\vec{y}]
 + L_{\mathcal{A}^{k+t-1}_2}[\vec{y}]\bigg)\bigg)+ R\bigg(   L^*_{\mathcal{B}^{k+t-1}_3}[x,\vec{y}]+    L^*_{\mathcal{A}^{k+t-1}_3}[x,\vec{y}]\bigg)\\
&  
 + \lambda \bigg(   L^*_{\mathcal{B}^{k+t-1}_3}[x,\vec{y}]
 +    L^*_{\mathcal{C}^{k+t-1}_3}[x,\vec{y}]\bigg)
 + \lambda^2 \bigg(x\bigg(  L_{\mathcal{B}^{k+t-1}_2}[\vec{y}]
 +       L_{\mathcal{C}^{k+t-1}_2}[\vec{y}]\bigg)\bigg)\bigg) ]\bigg)  \\
& + \lambda \tbinom{m-2}{i}R\bigg( L_{\mathcal{B}^{k+t+m-2-i}_{m+1}}[\vec{x},  \sum\limits^l_{t=0}\lambda^{(t)}\bigg(
 R \bigg(xR\bigg(  L_{\mathcal{B}^{k+t-1}_2}[\vec{y}]
 + L_{\mathcal{A}^{k+t-1}_2}[\vec{y}]\bigg)\bigg)
 + \lambda R\bigg( x\bigg(  L_{\mathcal{B}^{k+t-1}_2}[\vec{y}]
 + L_{\mathcal{A}^{k+t-1}_2}[\vec{y}]\bigg)\bigg) \\
& + \lambda \bigg( xR\bigg(  L_{\mathcal{B}^{k+t-1}_2}[\vec{y}]
 + L_{\mathcal{A}^{k+t-1}_2}[\vec{y}]\bigg)\bigg)
 + R\bigg(   L^*_{\mathcal{B}^{k+t-1}_3}[x,\vec{y}]+    L^*_{\mathcal{A}^{k+t-1}_3}[x,\vec{y}]\bigg) \\
& + \lambda \bigg(   L^*_{\mathcal{B}^{k+t-1}_3}[x,\vec{y}]
 +    L^*_{\mathcal{C}^{k+t-1}_3}[x,\vec{y}]\bigg)
 + \lambda^2 \bigg(x\bigg(  L_{\mathcal{B}^{k+t-1}_2}[\vec{y}]
 +       L_{\mathcal{C}^{k+t-1}_2}[\vec{y}]\bigg)\bigg)\bigg) ]\bigg)
\\
& + \lambda \tbinom{m-1}{i} R\bigg( L_{\mathcal{C}^{k+t+m-2-i}_{m+1}}[\vec{x}, \sum\limits^l_{t=0}\lambda^{(t)}\bigg(
 R \bigg(xR\bigg(  L_{\mathcal{B}^{k+t-1}_2}[\vec{y}]
 + L_{\mathcal{A}^{k+t-1}_2}[\vec{y}]\bigg)\bigg)
 + \lambda R\bigg( x\bigg(  L_{\mathcal{B}^{k+t-1}_2}[\vec{y}]
 + L_{\mathcal{A}^{k+t-1}_2}[\vec{y}]\bigg)\bigg) 
 \end{aligned}
\]
\[
\begin{aligned}
& + \lambda \bigg( xR\bigg(  L_{\mathcal{B}^{k+t-1}_2}[\vec{y}]
 + L_{\mathcal{A}^{k+t-1}_2}[\vec{y}]\bigg)\bigg)
 + R\bigg(   L^*_{\mathcal{B}^{k+t-1}_3}[x,\vec{y}]+    L^*_{\mathcal{A}^{k+t-1}_3}[x,\vec{y}]\bigg) 
 + \lambda \bigg(   L^*_{\mathcal{B}^{k+t-1}_3}[x,\vec{y}]
 +    L^*_{\mathcal{C}^{k+t-1}_3}[x,\vec{y}]\bigg)\\
 &+ \lambda^2 \bigg(x\bigg(  L_{\mathcal{B}^{k+t-1}_2}[\vec{y}]
 +       L_{\mathcal{C}^{k+t-1}_2}[\vec{y}]\bigg)\bigg)\bigg) ]\bigg) 
  + \lambda^2 \tbinom{m-2}{i} R\bigg(  L_{\mathcal{D}^{k+t+m-3-i}_{m+1}}[\vec{x}, \sum\limits^l_{t=0}\lambda^{(t)}\bigg(
 R \bigg(xR\bigg(  L_{\mathcal{B}^{k+t-1}_2}[\vec{y}]
 + L_{\mathcal{A}^{k+t-1}_2}[\vec{y}]\bigg)\bigg)\\
 &
 + \lambda R\bigg( x\bigg(  L_{\mathcal{B}^{k+t-1}_2}[\vec{y}]
 + L_{\mathcal{A}^{k+t-1}_2}[\vec{y}]\bigg)\bigg) 
 + \lambda \bigg( xR\bigg(  L_{\mathcal{B}^{k+t-1}_2}[\vec{y}]
 + L_{\mathcal{A}^{k+t-1}_2}[\vec{y}]\bigg)\bigg)
 + R\bigg(   L^*_{\mathcal{B}^{k+t-1}_3}[x,\vec{y}]+    L^*_{\mathcal{A}^{k+t-1}_3}[x,\vec{y}]\bigg) \\
& + \lambda \bigg(   L^*_{\mathcal{B}^{k+t-1}_3}[x,\vec{y}]
 +    L^*_{\mathcal{C}^{k+t-1}_3}[x,\vec{y}]\bigg)
 + \lambda^2 \bigg(x\bigg(  L_{\mathcal{B}^{k+t-1}_2}[\vec{y}]
 +       L_{\mathcal{C}^{k+t-1}_2}[\vec{y}]\bigg)\bigg)\bigg) ]\bigg) = 0; \ m = |\vec{x}|\geq 3, |\vec{y}|=2 .
\end{aligned}
\]

Composition between relations $(R4)$ and $(R8)$, we get a sum of the expressions, which are consequences of $(R4)$ itself.
Below, we show this example for  $\sum\limits_{t=0}^l \sum\limits^{m-1}_{i=0}\lambda^{(t)}\lambda^i
M_{t,i}$, where
\[
\begin{aligned}
M_{t,i}&=\tbinom{m-1}{i} R^2\bigg( L_{\mathcal{B}^{k+t+m-2-i}_{m+1}}[\vec{x},\sum\limits^l_{t=0}\sum\limits^{m-2}_{i=0}\lambda^{(t)}\lambda^i \left( 
 \tbinom{m-2}{i}R\bigg(  L^*_{\mathcal{B}^{k+t+m-3-i}_m}[x,\vec{x}]
 +    L^*_{\mathcal{A}^{k+t+m-3-i}_m}[x,\vec{x}]\bigg)\right. \\
& + \lambda \tbinom{m-3}{i} \bigg(  L^*_{\mathcal{B}^{k+t+m-3-i}_m}[x,\vec{x}]\bigg) 
 + \lambda \tbinom{m-2}{i} \bigg(    L^*_{\mathcal{C}^{k+t+m-3-i}_m}[x,\vec{x}]\bigg) 
 + \lambda^2 \tbinom{m-3}{i} \bigg(   L^*_{\mathcal{D}^{k+t+m-4-i}_m}[x,\vec{x}]\bigg)\\
 &+ \tbinom{m-2}{i}R\bigg(xR\bigg(          L_{\mathcal{B}^{k+t+m-3-i}_{m}}[\vec{x}]
 +        L_{\mathcal{A}^{k+t+m-3-i}_{m}}[\vec{x}]\bigg)\bigg) 
 + \lambda \tbinom{m-2}{i} R\bigg(x\bigg(          L_{\mathcal{B}^{k+t+m-3-i}_{m}}[\vec{x}]
 +            L_{\mathcal{A}^{k+t+m-3-i}_{m}}[\vec{x}]\bigg)\bigg) \allowdisplaybreaks \\
& + \lambda \tbinom{m-2}{i}\bigg(x R\bigg(          L_{\mathcal{B}^{k+t+m-3-i}_{m}}[\vec{x}]
 +        L_{\mathcal{A}^{k+t+m-3-i}_{m}}[\vec{x}]\bigg)\bigg)
 + \lambda^2 \tbinom{m-3}{i} \bigg(x \bigg(          L_{\mathcal{B}^{k+t+m-3-i}_{m}}[\vec{x}]\bigg)\bigg)\bigg) \\
& + \lambda^2 \tbinom{m-2}{i} \bigg(x \bigg(            L_{\mathcal{C}^{k+t+m-3-i}_{m}}[\vec{x}]\bigg)\bigg)\bigg)
\left. +\lambda^3 \tbinom{m-3}{i}\bigg(x \bigg(             L_{\mathcal{D}^{k+t+m-3-i}_{m}}[\vec{x}]\bigg)\bigg)\bigg) \right)]\\
& +       L_{\mathcal{A}^{k+t+m-2-i}_{m+1}}[\vec{x},\sum\limits^l_{t=0}\sum\limits^{m-2}_{i=0}\lambda^{(t)}\lambda^i \left( 
 \tbinom{m-2}{i}R\bigg(  L^*_{\mathcal{B}^{k+t+m-3-i}_m}[x,\vec{x}]
 +    L^*_{\mathcal{A}^{k+t+m-3-i}_m}[x,\vec{x}]\bigg)\right. \\
& + \lambda \tbinom{m-3}{i} \bigg(  L^*_{\mathcal{B}^{k+t+m-3-i}_m}[x,\vec{x}]\bigg) 
 + \lambda \tbinom{m-2}{i} \bigg(    L^*_{\mathcal{C}^{k+t+m-3-i}_m}[x,\vec{x}]\bigg) 
+ \lambda^2 \tbinom{m-3}{i} \bigg(   L^*_{\mathcal{D}^{k+t+m-4-i}_m}[x,\vec{x}]\bigg)\\
 &+ \tbinom{m-2}{i}R\bigg(xR\bigg(          L_{\mathcal{B}^{k+t+m-3-i}_{m}}[\vec{x}]
 +        L_{\mathcal{A}^{k+t+m-3-i}_{m}}[\vec{x}]\bigg)\bigg) 
 + \lambda \tbinom{m-2}{i} R\bigg(x\bigg(          L_{\mathcal{B}^{k+t+m-3-i}_{m}}[\vec{x}]
 +            L_{\mathcal{A}^{k+t+m-3-i}_{m}}[\vec{x}]\bigg)\bigg) \allowdisplaybreaks \\
& + \lambda \tbinom{m-2}{i}\bigg(x R\bigg(          L_{\mathcal{B}^{k+t+m-3-i}_{m}}[\vec{x}]
 +        L_{\mathcal{A}^{k+t+m-3-i}_{m}}[\vec{x}]\bigg)\bigg)
 + \lambda^2 \tbinom{m-3}{i} \bigg(x \bigg(          L_{\mathcal{B}^{k+t+m-3-i}_{m}}[\vec{x}]\bigg)\bigg)\bigg) \\
& + \lambda^2 \tbinom{m-2}{i} \bigg(x \bigg(            L_{\mathcal{C}^{k+t+m-3-i}_{m}}[\vec{x}]\bigg)\bigg)\bigg)
\left. +\lambda^3 \tbinom{m-3}{i}\bigg(x \bigg(             L_{\mathcal{D}^{k+t+m-3-i}_{m}}[\vec{x}]\bigg)\bigg)\bigg) \right)]\bigg)  \\
& + \lambda \tbinom{m-2}{i}R\bigg( L_{\mathcal{B}^{k+t+m-2-i}_{m+1}}[\vec{x},\sum\limits^l_{t=0}\sum\limits^{m-2}_{i=0}\lambda^{(t)}\lambda^i \left( 
 \tbinom{m-2}{i}R\bigg(  L^*_{\mathcal{B}^{k+t+m-3-i}_m}[x,\vec{x}]
 +    L^*_{\mathcal{A}^{k+t+m-3-i}_m}[x,\vec{x}]\bigg)\right. \\
& + \lambda \tbinom{m-3}{i} \bigg(  L^*_{\mathcal{B}^{k+t+m-3-i}_m}[x,\vec{x}]\bigg) 
 + \lambda \tbinom{m-2}{i} \bigg(    L^*_{\mathcal{C}^{k+t+m-3-i}_m}[x,\vec{x}]\bigg) 
 + \lambda^2 \tbinom{m-3}{i} \bigg(   L^*_{\mathcal{D}^{k+t+m-4-i}_m}[x,\vec{x}]\bigg)\\
 &+ \tbinom{m-2}{i}R\bigg(xR\bigg(          L_{\mathcal{B}^{k+t+m-3-i}_{m}}[\vec{x}]
 +        L_{\mathcal{A}^{k+t+m-3-i}_{m}}[\vec{x}]\bigg)\bigg) 
 + \lambda \tbinom{m-2}{i} R\bigg(x\bigg(          L_{\mathcal{B}^{k+t+m-3-i}_{m}}[\vec{x}]
 +            L_{\mathcal{A}^{k+t+m-3-i}_{m}}[\vec{x}]\bigg)\bigg) \allowdisplaybreaks \\
& + \lambda \tbinom{m-2}{i}\bigg(x R\bigg(          L_{\mathcal{B}^{k+t+m-3-i}_{m}}[\vec{x}]
 +        L_{\mathcal{A}^{k+t+m-3-i}_{m}}[\vec{x}]\bigg)\bigg)
 + \lambda^2 \tbinom{m-3}{i} \bigg(x \bigg(          L_{\mathcal{B}^{k+t+m-3-i}_{m}}[\vec{x}]\bigg)\bigg)\bigg) \\
 \end{aligned}
 \]
 \[
 \begin{aligned}
& + \lambda^2 \tbinom{m-2}{i} \bigg(x \bigg(            L_{\mathcal{C}^{k+t+m-3-i}_{m}}[\vec{x}]\bigg)\bigg)\bigg)
\left. +\lambda^3 \tbinom{m-3}{i}\bigg(x \bigg(             L_{\mathcal{D}^{k+t+m-3-i}_{m}}[\vec{x}]\bigg)\bigg)\bigg) \right)]\bigg)
 \\
& + \lambda \tbinom{m-1}{i} R\bigg( L_{\mathcal{C}^{k+t+m-2-i}_{m+1}}[\vec{x},\sum\limits^l_{t=0}\sum\limits^{m-2}_{i=0}\lambda^{(t)}\lambda^i \left( 
 \tbinom{m-2}{i}R\bigg(  L^*_{\mathcal{B}^{k+t+m-3-i}_m}[x,\vec{x}]
 +    L^*_{\mathcal{A}^{k+t+m-3-i}_m}[x,\vec{x}]\bigg)\right. \\
& + \lambda \tbinom{m-3}{i} \bigg(  L^*_{\mathcal{B}^{k+t+m-3-i}_m}[x,\vec{x}]\bigg) 
 + \lambda \tbinom{m-2}{i} \bigg(    L^*_{\mathcal{C}^{k+t+m-3-i}_m}[x,\vec{x}]\bigg) 
 + \lambda^2 \tbinom{m-3}{i} \bigg(   L^*_{\mathcal{D}^{k+t+m-4-i}_m}[x,\vec{x}]\bigg)\\
 &+ \tbinom{m-2}{i}R\bigg(xR\bigg(          L_{\mathcal{B}^{k+t+m-3-i}_{m}}[\vec{x}]
 +        L_{\mathcal{A}^{k+t+m-3-i}_{m}}[\vec{x}]\bigg)\bigg) 
 + \lambda \tbinom{m-2}{i} R\bigg(x\bigg(          L_{\mathcal{B}^{k+t+m-3-i}_{m}}[\vec{x}]
 +            L_{\mathcal{A}^{k+t+m-3-i}_{m}}[\vec{x}]\bigg)\bigg) \allowdisplaybreaks \\
& + \lambda \tbinom{m-2}{i}\bigg(x R\bigg(          L_{\mathcal{B}^{k+t+m-3-i}_{m}}[\vec{x}]
 +        L_{\mathcal{A}^{k+t+m-3-i}_{m}}[\vec{x}]\bigg)\bigg)
 + \lambda^2 \tbinom{m-3}{i} \bigg(x \bigg(          L_{\mathcal{B}^{k+t+m-3-i}_{m}}[\vec{x}]\bigg)\bigg)\bigg) \\
& + \lambda^2 \tbinom{m-2}{i} \bigg(x \bigg(            L_{\mathcal{C}^{k+t+m-3-i}_{m}}[\vec{x}]\bigg)\bigg)\bigg)
\left. +\lambda^3 \tbinom{m-3}{i}\bigg(x \bigg(             L_{\mathcal{D}^{k+t+m-3-i}_{m}}[\vec{x}]\bigg)\bigg)\bigg) \right)]\bigg) 
 \\
&  + \lambda^2 \tbinom{m-2}{i} R\bigg(  L_{\mathcal{D}^{k+t+m-3-i}_{m+1}}[\vec{x},\sum\limits^l_{t=0}\sum\limits^{m-2}_{i=0}\lambda^{(t)}\lambda^i \left( 
 \tbinom{m-2}{i}R\bigg(  L^*_{\mathcal{B}^{k+t+m-3-i}_m}[x,\vec{x}]
 +    L^*_{\mathcal{A}^{k+t+m-3-i}_m}[x,\vec{x}]\bigg)\right. \\
& + \lambda \tbinom{m-3}{i} \bigg(  L^*_{\mathcal{B}^{k+t+m-3-i}_m}[x,\vec{x}]\bigg) 
 + \lambda \tbinom{m-2}{i} \bigg(    L^*_{\mathcal{C}^{k+t+m-3-i}_m}[x,\vec{x}]\bigg) 
 + \lambda^2 \tbinom{m-3}{i} \bigg(   L^*_{\mathcal{D}^{k+t+m-4-i}_m}[x,\vec{x}]\bigg)\\
& + \tbinom{m-2}{i}R\bigg(xR\bigg(          L_{\mathcal{B}^{k+t+m-3-i}_{m}}[\vec{x}]
 +        L_{\mathcal{A}^{k+t+m-3-i}_{m}}[\vec{x}]\bigg)\bigg) 
 + \lambda \tbinom{m-2}{i} R\bigg(x\bigg(          L_{\mathcal{B}^{k+t+m-3-i}_{m}}[\vec{x}]
 +            L_{\mathcal{A}^{k+t+m-3-i}_{m}}[\vec{x}]\bigg)\bigg) \allowdisplaybreaks \\
& + \lambda \tbinom{m-2}{i}\bigg(x R\bigg(          L_{\mathcal{B}^{k+t+m-3-i}_{m}}[\vec{x}]
 +        L_{\mathcal{A}^{k+t+m-3-i}_{m}}[\vec{x}]\bigg)\bigg)
 + \lambda^2 \tbinom{m-3}{i} \bigg(x \bigg(          L_{\mathcal{B}^{k+t+m-3-i}_{m}}[\vec{x}]\bigg)\bigg)\bigg) \\
& + \lambda^2 \tbinom{m-2}{i} \bigg(x \bigg(            L_{\mathcal{C}^{k+t+m-3-i}_{m}}[\vec{x}]\bigg)\bigg)\bigg)
\left. +\lambda^3 \tbinom{m-3}{i}\bigg(x \bigg(             L_{\mathcal{D}^{k+t+m-3-i}_{m}}[\vec{x}]\bigg)\bigg)\bigg) \right)]\bigg) = 0; \ m = |\vec{x}|\geq 3.
\end{aligned}
\]

Composition between relations $(R4)$ and $(R9)$, we get a sum of the expressions, which are consequences of $(R4)$ itself.
Below, we show this example for  $\sum\limits_{t=0}^l \sum\limits^{m-1}_{i=0}\lambda^{(t)}\lambda^i
M_{t,i}$, where
\[
\begin{aligned}
M_{t,i}&=\tbinom{m-1}{i} R^2\bigg( L_{\mathcal{B}^{k+t+m-2-i}_{m+1}}[\vec{x},\sum\limits^l_{t=0}\lambda^{(t)}\bigg( R \bigg(R\bigg(  L^\op_{\mathcal{B}^{k+t-1}_2}[\vec{y}]+              L^\op_{\mathcal{A}^{n-1}_2}[\vec{y}]\bigg)x\bigg)\\
& +\lambda R\bigg( \bigg(  L^\op_{\mathcal{B}^{k+t-1}_2}[\vec{y}]+              L^\op_{\mathcal{A}^{n-1}_2}[\vec{y}]\bigg)x\bigg)
 +\lambda \bigg( R\bigg(  L^\op_{\mathcal{B}^{k+t-1}_2}[\vec{y}]
 +              L^\op_{\mathcal{A}^{n-1}_2}[\vec{y}]\bigg)x\bigg)\\
 &+ R\bigg(   {L^*}^\op_{\mathcal{B}^{k+t-1}_3}[x,\vec{y}]+    {L^*}^\op_{\mathcal{A}^{n-1}_3}[x,\vec{y}]\bigg)
 + \lambda \bigg(   {L^*}^\op_{\mathcal{B}^{k+t-1}_3}[x,\vec{y}]+   {L^*}^\op_{\mathcal{C}^{k+t-1}_3}[x,\vec{y}]\bigg)\\
& + \lambda^2 \bigg(\bigg(  L^\op_{\mathcal{B}^{k+t-1}_2}[\vec{y}]+      L^\op_{\mathcal{C}^{k+t-1}_2}[\vec{y}]\bigg)x\bigg)\bigg)]
+       L_{\mathcal{A}^{k+t+m-2-i}_{m+1}}[\vec{x},\sum\limits^l_{t=0}\lambda^{(t)}\bigg( R \bigg(R\bigg(  L^\op_{\mathcal{B}^{k+t-1}_2}[\vec{y}]+              L^\op_{\mathcal{A}^{n-1}_2}[\vec{y}]\bigg)x\bigg)
\\
 &+\lambda R\bigg( \bigg(  L^\op_{\mathcal{B}^{k+t-1}_2}[\vec{y}]+              L^\op_{\mathcal{A}^{n-1}_2}[\vec{y}]\bigg)x\bigg)
 +\lambda \bigg( R\bigg(  L^\op_{\mathcal{B}^{k+t-1}_2}[\vec{y}]
 +              L^\op_{\mathcal{A}^{n-1}_2}[\vec{y}]\bigg)x\bigg)\\
& + R\bigg(   {L^*}^\op_{\mathcal{B}^{k+t-1}_3}[x,\vec{y}]+    {L^*}^\op_{\mathcal{A}^{n-1}_3}[x,\vec{y}]\bigg)
 + \lambda \bigg(   {L^*}^\op_{\mathcal{B}^{k+t-1}_3}[x,\vec{y}]+   {L^*}^\op_{\mathcal{C}^{k+t-1}_3}[x,\vec{y}]\bigg)\\
& + \lambda^2 \bigg(\bigg(  L^\op_{\mathcal{B}^{k+t-1}_2}[\vec{y}]+      L^\op_{\mathcal{C}^{k+t-1}_2}[\vec{y}]\bigg)x\bigg)\bigg)]\bigg)  \\
 \end{aligned}
 \]
 \[
 \begin{aligned}
 &+ \lambda \tbinom{m-2}{i}R\bigg( L_{\mathcal{B}^{k+t+m-2-i}_{m+1}}[\vec{x}, \sum\limits^l_{t=0}\lambda^{(t)}\bigg( R \bigg(R\bigg(  L^\op_{\mathcal{B}^{k+t-1}_2}[\vec{y}]+              L^\op_{\mathcal{A}^{n-1}_2}[\vec{y}]\bigg)x\bigg)
 +\lambda R\bigg( \bigg(  L^\op_{\mathcal{B}^{k+t-1}_2}[\vec{y}]+              L^\op_{\mathcal{A}^{n-1}_2}[\vec{y}]\bigg)x\bigg)\\&
 +\lambda \bigg( R\bigg(  L^\op_{\mathcal{B}^{k+t-1}_2}[\vec{y}]
 +              L^\op_{\mathcal{A}^{n-1}_2}[\vec{y}]\bigg)x\bigg)
 + R\bigg(   {L^*}^\op_{\mathcal{B}^{k+t-1}_3}[x,\vec{y}]+    {L^*}^\op_{\mathcal{A}^{n-1}_3}[x,\vec{y}]\bigg)
 + \lambda \bigg(   {L^*}^\op_{\mathcal{B}^{k+t-1}_3}[x,\vec{y}]+   {L^*}^\op_{\mathcal{C}^{k+t-1}_3}[x,\vec{y}]\bigg)\\
& + \lambda^2 \bigg(\bigg(  L^\op_{\mathcal{B}^{k+t-1}_2}[\vec{y}]+      L^\op_{\mathcal{C}^{k+t-1}_2}[\vec{y}]\bigg)x\bigg)\bigg)]\bigg)
 + \lambda \tbinom{m-1}{i} R\bigg( L_{\mathcal{C}^{k+t+m-2-i}_{m+1}}[\vec{x},\sum\limits^l_{t=0}\lambda^{(t)}\bigg( R \bigg(R\bigg(  L^\op_{\mathcal{B}^{k+t-1}_2}[\vec{y}]+              L^\op_{\mathcal{A}^{n-1}_2}[\vec{y}]\bigg)x\bigg)\\
& +\lambda R\bigg( \bigg(  L^\op_{\mathcal{B}^{k+t-1}_2}[\vec{y}]+              L^\op_{\mathcal{A}^{n-1}_2}[\vec{y}]\bigg)x\bigg)
 +\lambda \bigg( R\bigg(  L^\op_{\mathcal{B}^{k+t-1}_2}[\vec{y}]
 +              L^\op_{\mathcal{A}^{n-1}_2}[\vec{y}]\bigg)x\bigg)
  + R\bigg(   {L^*}^\op_{\mathcal{B}^{k+t-1}_3}[x,\vec{y}]+    {L^*}^\op_{\mathcal{A}^{n-1}_3}[x,\vec{y}]\bigg)\\&
 + \lambda \bigg(   {L^*}^\op_{\mathcal{B}^{k+t-1}_3}[x,\vec{y}]+   {L^*}^\op_{\mathcal{C}^{k+t-1}_3}[x,\vec{y}]\bigg)
 + \lambda^2 \bigg(\bigg(  L^\op_{\mathcal{B}^{k+t-1}_2}[\vec{y}]+      L^\op_{\mathcal{C}^{k+t-1}_2}[\vec{y}]\bigg)x\bigg)\bigg)]\bigg) 
 \\
 & + \lambda^2 \tbinom{m-2}{i} R\bigg(  L_{\mathcal{D}^{k+t+m-3-i}_{m+1}}[\vec{x},\sum\limits^l_{t=0}\lambda^{(t)}\bigg( R \bigg(R\bigg(  L^\op_{\mathcal{B}^{k+t-1}_2}[\vec{y}]+              L^\op_{\mathcal{A}^{n-1}_2}[\vec{y}]\bigg)x\bigg)
 +\lambda R\bigg( \bigg(  L^\op_{\mathcal{B}^{k+t-1}_2}[\vec{y}]+              L^\op_{\mathcal{A}^{n-1}_2}[\vec{y}]\bigg)x\bigg)\\&
 +\lambda \bigg( R\bigg(  L^\op_{\mathcal{B}^{k+t-1}_2}[\vec{y}]
 +              L^\op_{\mathcal{A}^{n-1}_2}[\vec{y}]\bigg)x\bigg) + R\bigg(   {L^*}^\op_{\mathcal{B}^{k+t-1}_3}[x,\vec{y}]+    {L^*}^\op_{\mathcal{A}^{n-1}_3}[x,\vec{y}]\bigg)
 + \lambda \bigg(   {L^*}^\op_{\mathcal{B}^{k+t-1}_3}[x,\vec{y}]+   {L^*}^\op_{\mathcal{C}^{k+t-1}_3}[x,\vec{y}]\bigg)\\
& + \lambda^2 \bigg(\bigg(  L^\op_{\mathcal{B}^{k+t-1}_2}[\vec{y}]+      L^\op_{\mathcal{C}^{k+t-1}_2}[\vec{y}]\bigg)x\bigg)\bigg)]\bigg) = 0, \ m = |\vec{x}|\geq 3,\ |\vec{y}|=2.
\end{aligned}
\]

Composition between relations $(R4)$ and $(R10)$, we get a sum of the expressions, which are consequences of $(R4)$ itself.
Below, we show this example for  $\sum\limits_{t=0}^l \sum\limits^{m-1}_{i=0}\lambda^{(t)}\lambda^i
M_{t,i}$, where
\[
\begin{aligned}
M_{t,i}&=\tbinom{m-1}{i} R^2\bigg( L_{\mathcal{B}^{k+t+m-2-i}_{m+1}}[\vec{x}, \sum\limits^l_{t=0}\sum\limits^{m-2}_{i=0}\lambda^{(t)}\lambda^i \left(
\tbinom{m-2}{i}R\bigg(  {L^*}^\op_{\mathcal{B}^{k+t+m-3-i}_m}[x,\vec{x}]
+   {L^*}^\op_{\mathcal{A}^{k+t+m-3-i}_m}[x,\vec{x}]\bigg) \right. \\
&+\lambda \tbinom{m-3}{i}\bigg(  {L^*}^\op_{\mathcal{B}^{k+t+m-3-i}_m}[x,\vec{x}]\bigg)
+\lambda \tbinom{m-2}{i} \bigg(    {L^*}^\op_{\mathcal{C}^{k+t+m-3-i}_m}[x,\vec{x}]\bigg)\bigg)\\
&+\lambda^2 \tbinom{m-3}{i}\bigg(   {L^*}^\op_{\mathcal{D}^{k+t+m-4-i}_m}[x,\vec{x}]\bigg)
 + \tbinom{m-2}{i}R\bigg(R\bigg(          L^\op_{\mathcal{B}^{k+t+m-3-i}_{m}}[\vec{x}]
 +            L^\op_{\mathcal{A}^{k+t+m-3-i}_{m}}[\vec{x}]\bigg)x\bigg)\\
& + \lambda \tbinom{m-2}{i}R\bigg(\bigg(          L^\op_{\mathcal{B}^{k+t+m-3-i}_{m}}[\vec{x}]
 +            L^\op_{\mathcal{A}^{k+t+m-3-i}_{m}}[\vec{x}]\bigg)x\bigg)\\
& + \lambda \tbinom{m-2}{i}\bigg( R\bigg(          L^\op_{\mathcal{B}^{k+t+m-3-i}_{m}}[\vec{x}]
 +            L^\op_{\mathcal{A}^{k+t+m-3-i}_{m}}[\vec{x}]\bigg)x\bigg)
 + \lambda^2 \tbinom{m-3}{i} \bigg( \bigg(          L^\op_{\mathcal{B}^{k+t+m-3-i}_{m}}[\vec{x}]\bigg)x\bigg) \\
& + \lambda^2 \tbinom{m-2}{i}\bigg( \bigg(            L^\op_{\mathcal{C}^{k+t+m-3-i}_{m}}[\vec{x}]\bigg)x\bigg)
 \left. +\lambda^3 \tbinom{m-3}{i}\bigg( \bigg(              L^\op_{\mathcal{D}^{k+t+m-3-i}_{m}}[\vec{x}]\bigg)x\bigg) \right)]
\\
& +       L_{\mathcal{A}^{k+t+m-2-i}_{m+1}}[\vec{x},\sum\limits^l_{t=0}\sum\limits^{m-2}_{i=0}\lambda^{(t)}\lambda^i \left(
\tbinom{m-2}{i}R\bigg(  {L^*}^\op_{\mathcal{B}^{k+t+m-3-i}_m}[x,\vec{x}]
+   {L^*}^\op_{\mathcal{A}^{k+t+m-3-i}_m}[x,\vec{x}]\bigg) \right. \\
&+\lambda \tbinom{m-3}{i}\bigg(  {L^*}^\op_{\mathcal{B}^{k+t+m-3-i}_m}[x,\vec{x}]\bigg)
+\lambda \tbinom{m-2}{i} \bigg(    {L^*}^\op_{\mathcal{C}^{k+t+m-3-i}_m}[x,\vec{x}]\bigg)\bigg)\\
&+\lambda^2 \tbinom{m-3}{i}\bigg(   {L^*}^\op_{\mathcal{D}^{k+t+m-4-i}_m}[x,\vec{x}]\bigg)
 + \tbinom{m-2}{i}R\bigg(R\bigg(          L^\op_{\mathcal{B}^{k+t+m-3-i}_{m}}[\vec{x}]
 +            L^\op_{\mathcal{A}^{k+t+m-3-i}_{m}}[\vec{x}]\bigg)x\bigg)\\
& + \lambda \tbinom{m-2}{i}R\bigg(\bigg(          L^\op_{\mathcal{B}^{k+t+m-3-i}_{m}}[\vec{x}]
 +            L^\op_{\mathcal{A}^{k+t+m-3-i}_{m}}[\vec{x}]\bigg)x\bigg)\\
& + \lambda \tbinom{m-2}{i}\bigg( R\bigg(          L^\op_{\mathcal{B}^{k+t+m-3-i}_{m}}[\vec{x}]
 +            L^\op_{\mathcal{A}^{k+t+m-3-i}_{m}}[\vec{x}]\bigg)x\bigg)
 + \lambda^2 \tbinom{m-3}{i} \bigg( \bigg(          L^\op_{\mathcal{B}^{k+t+m-3-i}_{m}}[\vec{x}]\bigg)x\bigg) \\
& + \lambda^2 \tbinom{m-2}{i}\bigg( \bigg(            L^\op_{\mathcal{C}^{k+t+m-3-i}_{m}}[\vec{x}]\bigg)x\bigg)
 \left. +\lambda^3 \tbinom{m-3}{i}\bigg( \bigg(              L^\op_{\mathcal{D}^{k+t+m-3-i}_{m}}[\vec{x}]\bigg)x\bigg) \right)]\bigg) 
 \end{aligned}
 \]
 \[
 \begin{aligned}
 &+ \lambda \tbinom{m-2}{i}R\bigg( L_{\mathcal{B}^{k+t+m-2-i}_{m+1}}[\vec{x}, \sum\limits^l_{t=0}\sum\limits^{m-2}_{i=0}\lambda^{(t)}\lambda^i \left(
\tbinom{m-2}{i}R\bigg(  {L^*}^\op_{\mathcal{B}^{k+t+m-3-i}_m}[x,\vec{x}]
+   {L^*}^\op_{\mathcal{A}^{k+t+m-3-i}_m}[x,\vec{x}]\bigg) \right. \\
&+\lambda \tbinom{m-3}{i}\bigg(  {L^*}^\op_{\mathcal{B}^{k+t+m-3-i}_m}[x,\vec{x}]\bigg)
+\lambda \tbinom{m-2}{i} \bigg(    {L^*}^\op_{\mathcal{C}^{k+t+m-3-i}_m}[x,\vec{x}]\bigg)\bigg)\\
&+\lambda^2 \tbinom{m-3}{i}\bigg(   {L^*}^\op_{\mathcal{D}^{k+t+m-4-i}_m}[x,\vec{x}]\bigg)
 + \tbinom{m-2}{i}R\bigg(R\bigg(          L^\op_{\mathcal{B}^{k+t+m-3-i}_{m}}[\vec{x}]
 +            L^\op_{\mathcal{A}^{k+t+m-3-i}_{m}}[\vec{x}]\bigg)x\bigg)\\
& + \lambda \tbinom{m-2}{i}R\bigg(\bigg(          L^\op_{\mathcal{B}^{k+t+m-3-i}_{m}}[\vec{x}]
 +            L^\op_{\mathcal{A}^{k+t+m-3-i}_{m}}[\vec{x}]\bigg)x\bigg)\\
& + \lambda \tbinom{m-2}{i}\bigg( R\bigg(          L^\op_{\mathcal{B}^{k+t+m-3-i}_{m}}[\vec{x}]
 +            L^\op_{\mathcal{A}^{k+t+m-3-i}_{m}}[\vec{x}]\bigg)x\bigg)
 + \lambda^2 \tbinom{m-3}{i} \bigg( \bigg(          L^\op_{\mathcal{B}^{k+t+m-3-i}_{m}}[\vec{x}]\bigg)x\bigg) \\
& + \lambda^2 \tbinom{m-2}{i}\bigg( \bigg(            L^\op_{\mathcal{C}^{k+t+m-3-i}_{m}}[\vec{x}]\bigg)x\bigg)
 \left. +\lambda^3 \tbinom{m-3}{i}\bigg( \bigg(              L^\op_{\mathcal{D}^{k+t+m-3-i}_{m}}[\vec{x}]\bigg)x\bigg) \right)]\bigg)
 \\
& + \lambda \tbinom{m-1}{i} R\bigg( L_{\mathcal{C}^{k+t+m-2-i}_{m+1}}[\vec{x}, \sum\limits^l_{t=0}\sum\limits^{m-2}_{i=0}\lambda^{(t)}\lambda^i \left(
\tbinom{m-2}{i}R\bigg(  {L^*}^\op_{\mathcal{B}^{k+t+m-3-i}_m}[x,\vec{x}]
+   {L^*}^\op_{\mathcal{A}^{k+t+m-3-i}_m}[x,\vec{x}]\bigg) \right. \\
&+\lambda \tbinom{m-3}{i}\bigg(  {L^*}^\op_{\mathcal{B}^{k+t+m-3-i}_m}[x,\vec{x}]\bigg)
+\lambda \tbinom{m-2}{i} \bigg(    {L^*}^\op_{\mathcal{C}^{k+t+m-3-i}_m}[x,\vec{x}]\bigg)\bigg)+\lambda^2 \tbinom{m-3}{i}\bigg(   {L^*}^\op_{\mathcal{D}^{k+t+m-4-i}_m}[x,\vec{x}]\bigg)\\
&
 + \tbinom{m-2}{i}R\bigg(R\bigg(          L^\op_{\mathcal{B}^{k+t+m-3-i}_{m}}[\vec{x}]
 +            L^\op_{\mathcal{A}^{k+t+m-3-i}_{m}}[\vec{x}]\bigg)x\bigg) + \lambda \tbinom{m-2}{i}R\bigg(\bigg(          L^\op_{\mathcal{B}^{k+t+m-3-i}_{m}}[\vec{x}]
 +            L^\op_{\mathcal{A}^{k+t+m-3-i}_{m}}[\vec{x}]\bigg)x\bigg)\\
& + \lambda \tbinom{m-2}{i}\bigg( R\bigg(          L^\op_{\mathcal{B}^{k+t+m-3-i}_{m}}[\vec{x}]
 +            L^\op_{\mathcal{A}^{k+t+m-3-i}_{m}}[\vec{x}]\bigg)x\bigg)
 + \lambda^2 \tbinom{m-3}{i} \bigg( \bigg(          L^\op_{\mathcal{B}^{k+t+m-3-i}_{m}}[\vec{x}]\bigg)x\bigg) \\
& + \lambda^2 \tbinom{m-2}{i}\bigg( \bigg(            L^\op_{\mathcal{C}^{k+t+m-3-i}_{m}}[\vec{x}]\bigg)x\bigg)
 \left. +\lambda^3 \tbinom{m-3}{i}\bigg( \bigg(              L^\op_{\mathcal{D}^{k+t+m-3-i}_{m}}[\vec{x}]\bigg)x\bigg) \right)]\bigg) 
 \\
& +\lambda^2 \tbinom{m-2}{i} R\bigg(  L_{\mathcal{D}^{k+t+m-3-i}_{m+1}}[\vec{x},\sum\limits^l_{t=0}\sum\limits^{m-2}_{i=0}\lambda^{(t)}\lambda^i \left(
\tbinom{m-2}{i}R\bigg(  {L^*}^\op_{\mathcal{B}^{k+t+m-3-i}_m}[x,\vec{x}]
+   {L^*}^\op_{\mathcal{A}^{k+t+m-3-i}_m}[x,\vec{x}]\bigg) \right. \\
&+\lambda \tbinom{m-3}{i}\bigg(  {L^*}^\op_{\mathcal{B}^{k+t+m-3-i}_m}[x,\vec{x}]\bigg)
+\lambda \tbinom{m-2}{i} \bigg(    {L^*}^\op_{\mathcal{C}^{k+t+m-3-i}_m}[x,\vec{x}]\bigg)\bigg)+\lambda^2 \tbinom{m-3}{i}\bigg(   {L^*}^\op_{\mathcal{D}^{k+t+m-4-i}_m}[x,\vec{x}]\bigg)\\
&
 + \tbinom{m-2}{i}R\bigg(R\bigg(          L^\op_{\mathcal{B}^{k+t+m-3-i}_{m}}[\vec{x}]
 +            L^\op_{\mathcal{A}^{k+t+m-3-i}_{m}}[\vec{x}]\bigg)x\bigg) + \lambda \tbinom{m-2}{i}R\bigg(\bigg(          L^\op_{\mathcal{B}^{k+t+m-3-i}_{m}}[\vec{x}]
 +            L^\op_{\mathcal{A}^{k+t+m-3-i}_{m}}[\vec{x}]\bigg)x\bigg)\\
& + \lambda \tbinom{m-2}{i}\bigg( R\bigg(          L^\op_{\mathcal{B}^{k+t+m-3-i}_{m}}[\vec{x}]
 +            L^\op_{\mathcal{A}^{k+t+m-3-i}_{m}}[\vec{x}]\bigg)x\bigg)
 + \lambda^2 \tbinom{m-3}{i} \bigg( \bigg(          L^\op_{\mathcal{B}^{k+t+m-3-i}_{m}}[\vec{x}]\bigg)x\bigg) \\
& + \lambda^2 \tbinom{m-2}{i}\bigg( \bigg(            L^\op_{\mathcal{C}^{k+t+m-3-i}_{m}}[\vec{x}]\bigg)x\bigg)
 \left. +\lambda^3 \tbinom{m-3}{i}\bigg( \bigg(              L^\op_{\mathcal{D}^{k+t+m-3-i}_{m}}[\vec{x}]\bigg)x\bigg) \right)]\bigg) = 0; \ m = |\vec{x}|\geq 3.
\end{aligned}
\]

Composition between relations $(R5)$ and $(R4)$-$(R10)$, we get a sum of the expressions, which are consequences of $(R6)$ itself.

Composition between relations $(R6)$ and $(R4)$, we get a sum of the expressions, which are consequences of $(R6)$ itself.
Below, we show this example for  $\sum\limits_{t=0}^l \sum\limits^{m-1}_{i=0}\lambda^{(t)}\lambda^i
M_{t,i}$, where
\[
\begin{aligned} 
M_{t,i}&=\tbinom{m-1}{i}R^2\bigg( L^\op_{\mathcal{B}^{k+t+m-2-i}_{m+1}}[\vec{x},R\bigg(\sum\limits^{m-2}_{i=0}\tbinom{m-2}{i}\lambda^i\bigg(  L_{\mathcal{B}^{k+t+m-3-i}_m}[\vec{x}]+   L_{\mathcal{A}^{k+t+m-3-i}_m}[\vec{x}]\bigg)\bigg) \\
\end{aligned}
\]
\[
\begin{aligned}
& +\lambda (\sum\limits^{m-3}_{i=0}\tbinom{m-3}{i}\lambda^i\bigg(  L_{\mathcal{B}^{k+t+m-3-i}_m}[\vec{x}]\bigg)\bigg) 
 +\lambda \bigg(\sum\limits^{m-2}_{i=0}\tbinom{m-2}{i}\lambda^i\bigg(    L_{\mathcal{C}^{k+t+m-3-i}_m}[\vec{x}]\bigg)\bigg) 
  +\lambda^2 \bigg(\sum\limits^{m-3}_{i=0}\tbinom{m-3}{i}\lambda^i\bigg(   L_{\mathcal{D}^{k+t+m-4-i}_m}[\vec{x}]\bigg)\bigg)]
\\
& +       L^\op_{\mathcal{A}^{k+t+m-2-i}_{m+1}}[\vec{x},R\bigg(\sum\limits^{m-2}_{i=0}\tbinom{m-2}{i}\lambda^i\bigg(  L_{\mathcal{B}^{k+t+m-3-i}_m}[\vec{x}]+   L_{\mathcal{A}^{k+t+m-3-i}_m}[\vec{x}]\bigg)\bigg)  +\lambda (\sum\limits^{m-3}_{i=0}\tbinom{m-3}{i}\lambda^i\bigg(  L_{\mathcal{B}^{k+t+m-3-i}_m}[\vec{x}]\bigg)\bigg) \\
 &+\lambda \bigg(\sum\limits^{m-2}_{i=0}\tbinom{m-2}{i}\lambda^i\bigg(    L_{\mathcal{C}^{k+t+m-3-i}_m}[\vec{x}]\bigg)\bigg) 
  +\lambda^2 \bigg(\sum\limits^{m-3}_{i=0}\tbinom{m-3}{i}\lambda^i\bigg(   L_{\mathcal{D}^{k+t+m-4-i}_m}[\vec{x}]\bigg)\bigg)]\bigg)\\
 & + \lambda \tbinom{m-2}{i} R\bigg( L^\op_{\mathcal{B}^{k+t+m-2-i}_m}[\vec{x},R\bigg(\sum\limits^{m-2}_{i=0}\tbinom{m-2}{i}\lambda^i\bigg(  L_{\mathcal{B}^{k+t+m-3-i}_m}[\vec{x}]+   L_{\mathcal{A}^{k+t+m-3-i}_m}[\vec{x}]\bigg)\bigg) \\
& +\lambda (\sum\limits^{m-3}_{i=0}\tbinom{m-3}{i}\lambda^i\bigg(  L_{\mathcal{B}^{k+t+m-3-i}_m}[\vec{x}]\bigg)\bigg) 
 +\lambda \bigg(\sum\limits^{m-2}_{i=0}\tbinom{m-2}{i}\lambda^i\bigg(    L_{\mathcal{C}^{k+t+m-3-i}_m}[\vec{x}]\bigg)\bigg)  +\lambda^2 \bigg(\sum\limits^{m-3}_{i=0}\tbinom{m-3}{i}\lambda^i\bigg(   L_{\mathcal{D}^{k+t+m-4-i}_m}[\vec{x}]\bigg)\bigg)]\bigg)
 \\
& + \lambda \tbinom{m-1}{i}R\bigg( L^\op_{\mathcal{C}^{k+t+m-2-i}_{m+1}}[\vec{x},R\bigg(\sum\limits^{m-2}_{i=0}\tbinom{m-2}{i}\lambda^i\bigg(  L_{\mathcal{B}^{k+t+m-3-i}_m}[\vec{x}]+   L_{\mathcal{A}^{k+t+m-3-i}_m}[\vec{x}]\bigg)\bigg) \\
& +\lambda (\sum\limits^{m-3}_{i=0}\tbinom{m-3}{i}\lambda^i\bigg(  L_{\mathcal{B}^{k+t+m-3-i}_m}[\vec{x}]\bigg)\bigg) 
 +\lambda \bigg(\sum\limits^{m-2}_{i=0}\tbinom{m-2}{i}\lambda^i\bigg(    L_{\mathcal{C}^{k+t+m-3-i}_m}[\vec{x}]\bigg)\bigg) 
 +\lambda^2 \bigg(\sum\limits^{m-3}_{i=0}\tbinom{m-3}{i}\lambda^i\bigg(   L_{\mathcal{D}^{k+t+m-4-i}_m}[\vec{x}]\bigg)\bigg)]\bigg)
\\
& +\lambda^2 \tbinom{m-2}{i} R\bigg(  L^\op_{\mathcal{D}^{k+t+m-3-i}_{m+1}}[\vec{x},R\bigg(\sum\limits^{m-2}_{i=0}\tbinom{m-2}{i}\lambda^i\bigg(  L_{\mathcal{B}^{k+t+m-3-i}_m}[\vec{x}]+   L_{\mathcal{A}^{k+t+m-3-i}_m}[\vec{x}]\bigg)\bigg) \\
& +\lambda (\sum\limits^{m-3}_{i=0}\tbinom{m-3}{i}\lambda^i\bigg(  L_{\mathcal{B}^{k+t+m-3-i}_m}[\vec{x}]\bigg)\bigg) 
 +\lambda \bigg(\sum\limits^{m-2}_{i=0}\tbinom{m-2}{i}\lambda^i\bigg(    L_{\mathcal{C}^{k+t+m-3-i}_m}[\vec{x}]\bigg)\bigg) \\
 & +\lambda^2 \bigg(\sum\limits^{m-3}_{i=0}\tbinom{m-3}{i}\lambda^i\bigg(   L_{\mathcal{D}^{k+t+m-4-i}_m}[\vec{x}]\bigg)\bigg)]\bigg)=0; \ m = |\vec{x}|\geq3.
\end{aligned}
\]

Composition between relations $(R6)$ and $(R5)$, we get a sum of the expressions, which are consequences of $(R6)$ itself.
Below, we show this example for  $\sum\limits_{t=0}^l \sum\limits^{m-1}_{i=0}\lambda^{(t)}\lambda^i
M_{t,i}$, where
\[
\begin{aligned} 
M_{t,i}&=\tbinom{m-1}{i}R^2\bigg( L^\op_{\mathcal{B}^{k+t+m-2-i}_{m+1}}[\vec{x},\sum\limits^l_{t=0}\lambda ^{(t)}\bigg(R\bigg(  L^\op_{\mathcal{B}^{k+t-1}_2}[\vec{y}]+              L^{\op}_{\mathcal{A}^{k+t-1}_2}[\vec{y}]\bigg)
+\lambda \bigg(  L^{\op}_{\mathcal{B}^{k+t-1}_2}[\vec{y}]+ L^\op_{\mathcal{C}^{k+t-1}_2}[\vec{y}]\bigg)\bigg)]\\
 &+       L^\op_{\mathcal{A}^{k+t+m-2-i}_{m+1}}[\vec{x},\sum\limits^l_{t=0}\lambda ^{(t)}\bigg(R\bigg(  L^\op_{\mathcal{B}^{k+t-1}_2}[\vec{y}]+              L^{\op}_{\mathcal{A}^{k+t-1}_2}[\vec{y}]\bigg)
+\lambda \bigg(  L^{\op}_{\mathcal{B}^{k+t-1}_2}[\vec{y}]+      L^\op_{\mathcal{C}^{k+t-1}_2}[\vec{y}]\bigg)\bigg)]\bigg)
 \\
& + \lambda \tbinom{m-2}{i} R\bigg( L^\op_{\mathcal{B}^{k+t+m-2-i}_m}[\vec{x},\sum\limits^l_{t=0}\lambda ^{(t)}\bigg(R\bigg(  L^\op_{\mathcal{B}^{k+t-1}_2}[\vec{y}]+              L^{\op}_{\mathcal{A}^{k+t-1}_2}[\vec{y}]\bigg)
+\lambda \bigg(  L^{\op}_{\mathcal{B}^{k+t-1}_2}[\vec{y}]+      L^\op_{\mathcal{C}^{k+t-1}_2}[\vec{y}]\bigg)\bigg)]\bigg)
 \\
& + \lambda \tbinom{m-1}{i}R\bigg( L^\op_{\mathcal{C}^{k+t+m-2-i}_{m+1}}[\vec{x},\sum\limits^l_{t=0}\lambda ^{(t)}\bigg(R\bigg(  L^\op_{\mathcal{B}^{k+t-1}_2}[\vec{y}]+              L^{\op}_{\mathcal{A}^{k+t-1}_2}[\vec{y}]\bigg)
+\lambda \bigg(  L^{\op}_{\mathcal{B}^{k+t-1}_2}[\vec{y}]+      L^\op_{\mathcal{C}^{k+t-1}_2}[\vec{y}]\bigg)\bigg)]\bigg)
 \\
& +\lambda^2 \tbinom{m-2}{i} R\bigg(  L^\op_{\mathcal{D}^{k+t+m-3-i}_{m+1}}[\vec{x},\sum\limits^l_{t=0}\lambda ^{(t)}\bigg(R\bigg(  L^\op_{\mathcal{B}^{k+t-1}_2}[\vec{y}]+              L^{\op}_{\mathcal{A}^{k+t-1}_2}[\vec{y}]\bigg)
+\lambda \bigg(  L^{\op}_{\mathcal{B}^{k+t-1}_2}[\vec{y}]+      L^\op_{\mathcal{C}^{k+t-1}_2}[\vec{y}]\bigg)\bigg)]\bigg)=0;\\& \ m = |\vec{x}|\geq3;
\end{aligned}
\]

Composition between relations $(R7)$ and $(R4)$-$(R10)$, we get a sum of the expressions, which are consequences of $(R8)$ itself.

Composition between relations $(R8)$ and $(R4)$, we get a sum of the expressions, which are consequences of $(R8)$ itself.
Below, we show this example for  $\sum\limits_{t=0}^l \sum\limits^{m-1}_{i=0}\lambda^{(t)}\lambda^i
M_{t,i}$, where
\[
\begin{aligned}
 M_{t,i}&=\tbinom{m-1}{i}R^2\bigg( L^*_{\mathcal{B}^{k+t+m-2-i}_{m+1}}[x,\vec{x},R\bigg(\sum\limits^{m-2}_{i=0}\tbinom{m-2}{i}\lambda^i\bigg(  L_{\mathcal{B}^{k+t+m-3-i}_m}[\vec{x}]+   L_{\mathcal{A}^{k+t+m-3-i}_m}[\vec{x}]\bigg)\bigg) \\
& +\lambda (\sum\limits^{m-3}_{i=0}\tbinom{m-3}{i}\lambda^i\bigg(  L_{\mathcal{B}^{k+t+m-3-i}_m}[\vec{x}]\bigg)\bigg) 
 +\lambda \bigg(\sum\limits^{m-2}_{i=0}\tbinom{m-2}{i}\lambda^i\bigg(    L_{\mathcal{C}^{k+t+m-3-i}_m}[\vec{x}]\bigg)\bigg) 
  +\lambda^2 \bigg(\sum\limits^{m-3}_{i=0}\tbinom{m-3}{i}\lambda^i\bigg(   L_{\mathcal{D}^{k+t+m-4-i}_m}[\vec{x}]\bigg)\bigg)]
 \\
 &+       L^*_{\mathcal{A}^{k+t+m-2-i}_{m+1}}[x,\vec{x},R\bigg(\sum\limits^{m-2}_{i=0}\tbinom{m-2}{i}\lambda^i\bigg(  L_{\mathcal{B}^{k+t+m-3-i}_m}[\vec{x}]+   L_{\mathcal{A}^{k+t+m-3-i}_m}[\vec{x}]\bigg)\bigg) \\
& +\lambda (\sum\limits^{m-3}_{i=0}\tbinom{m-3}{i}\lambda^i\bigg(  L_{\mathcal{B}^{k+t+m-3-i}_m}[\vec{x}]\bigg)\bigg) 
 +\lambda \bigg(\sum\limits^{m-2}_{i=0}\tbinom{m-2}{i}\lambda^i\bigg(    L_{\mathcal{C}^{k+t+m-3-i}_m}[\vec{x}]\bigg)\bigg)   +\lambda^2 \bigg(\sum\limits^{m-3}_{i=0}\tbinom{m-3}{i}\lambda^i\bigg(   L_{\mathcal{D}^{k+t+m-4-i}_m}[\vec{x}]\bigg)\bigg)]\bigg) 
 \\
& + \lambda \tbinom{m-2}{i} R\bigg( L^*_{\mathcal{B}^{k+t+m-2-i}_{m+1}}[x,\vec{x},R\bigg(\sum\limits^{m-2}_{i=0}\tbinom{m-2}{i}\lambda^i\bigg(  L_{\mathcal{B}^{k+t+m-3-i}_m}[\vec{x}]+   L_{\mathcal{A}^{k+t+m-3-i}_m}[\vec{x}]\bigg)\bigg) \\
& +\lambda (\sum\limits^{m-3}_{i=0}\tbinom{m-3}{i}\lambda^i\bigg(  L_{\mathcal{B}^{k+t+m-3-i}_m}[\vec{x}]\bigg)\bigg) 
 +\lambda \bigg(\sum\limits^{m-2}_{i=0}\tbinom{m-2}{i}\lambda^i\bigg(    L_{\mathcal{C}^{k+t+m-3-i}_m}[\vec{x}]\bigg)\bigg)  +\lambda^2 \bigg(\sum\limits^{m-3}_{i=0}\tbinom{m-3}{i}\lambda^i\bigg(   L_{\mathcal{D}^{k+t+m-4-i}_m}[\vec{x}]\bigg)\bigg)]\bigg) 
\\
& + \lambda \tbinom{m-1}{i} R\bigg( L^*_{\mathcal{C}^{k+t+m-2-i}_{m+1}}[x,\vec{x},R\bigg(\sum\limits^{m-2}_{i=0}\tbinom{m-2}{i}\lambda^i\bigg(  L_{\mathcal{B}^{k+t+m-3-i}_m}[\vec{x}]+   L_{\mathcal{A}^{k+t+m-3-i}_m}[\vec{x}]\bigg)\bigg) \\
& +\lambda (\sum\limits^{m-3}_{i=0}\tbinom{m-3}{i}\lambda^i\bigg(  L_{\mathcal{B}^{k+t+m-3-i}_m}[\vec{x}]\bigg)\bigg) 
 +\lambda \bigg(\sum\limits^{m-2}_{i=0}\tbinom{m-2}{i}\lambda^i\bigg(    L_{\mathcal{C}^{k+t+m-3-i}_m}[\vec{x}]\bigg)\bigg)  +\lambda^2 \bigg(\sum\limits^{m-3}_{i=0}\tbinom{m-3}{i}\lambda^i\bigg(   L_{\mathcal{D}^{k+t+m-4-i}_m}[\vec{x}]\bigg)\bigg)]\bigg)
\\
& + \lambda^2 \tbinom{m-2}{i} R\bigg(  L^*_{\mathcal{D}^{k+t+m-3-i}_{m+1}}[x,\vec{x},R\bigg(\sum\limits^{m-2}_{i=0}\tbinom{m-2}{i}\lambda^i\bigg(  L_{\mathcal{B}^{k+t+m-3-i}_m}[\vec{x}]+   L_{\mathcal{A}^{k+t+m-3-i}_m}[\vec{x}]\bigg)\bigg) \\
& +\lambda (\sum\limits^{m-3}_{i=0}\tbinom{m-3}{i}\lambda^i\bigg(  L_{\mathcal{B}^{k+t+m-3-i}_m}[\vec{x}]\bigg)\bigg) 
 +\lambda \bigg(\sum\limits^{m-2}_{i=0}\tbinom{m-2}{i}\lambda^i\bigg(    L_{\mathcal{C}^{k+t+m-3-i}_m}[\vec{x}]\bigg)\bigg)   +\lambda^2 \bigg(\sum\limits^{m-3}_{i=0}\tbinom{m-3}{i}\lambda^i\bigg(   L_{\mathcal{D}^{k+t+m-4-i}_m}[\vec{x}]\bigg)\bigg)]\bigg)
\\
& + \tbinom{m-1}{i}R^2\bigg(xR\bigg(  L_{\mathcal{B}^{k+t+m-2-i}_{m+1}}[\vec{x},R\bigg(\sum\limits^{m-2}_{i=0}\tbinom{m-2}{i}\lambda^i\bigg(  L_{\mathcal{B}^{k+t+m-3-i}_m}[\vec{x}]+   L_{\mathcal{A}^{k+t+m-3-i}_m}[\vec{x}]\bigg)\bigg) \\
& +\lambda (\sum\limits^{m-3}_{i=0}\tbinom{m-3}{i}\lambda^i\bigg(  L_{\mathcal{B}^{k+t+m-3-i}_m}[\vec{x}]\bigg)\bigg) 
 +\lambda \bigg(\sum\limits^{m-2}_{i=0}\tbinom{m-2}{i}\lambda^i\bigg(    L_{\mathcal{C}^{k+t+m-3-i}_m}[\vec{x}]\bigg)\bigg)   +\lambda^2 \bigg(\sum\limits^{m-3}_{i=0}\tbinom{m-3}{i}\lambda^i\bigg(   L_{\mathcal{D}^{k+t+m-4-i}_m}[\vec{x}]\bigg)\bigg)]
 \\
& +   L_{\mathcal{A}^{k+t+m-2-i}_{m+1}}[\vec{x},R\bigg(\sum\limits^{m-2}_{i=0}\tbinom{m-2}{i}\lambda^i\bigg(  L_{\mathcal{B}^{k+t+m-3-i}_m}[\vec{x}]+   L_{\mathcal{A}^{k+t+m-3-i}_m}[\vec{x}]\bigg)\bigg) \\
& +\lambda (\sum\limits^{m-3}_{i=0}\tbinom{m-3}{i}\lambda^i\bigg(  L_{\mathcal{B}^{k+t+m-3-i}_m}[\vec{x}]\bigg)\bigg) 
 +\lambda \bigg(\sum\limits^{m-2}_{i=0}\tbinom{m-2}{i}\lambda^i\bigg(    L_{\mathcal{C}^{k+t+m-3-i}_m}[\vec{x}]\bigg)\bigg)   +\lambda^2 \bigg(\sum\limits^{m-3}_{i=0}\tbinom{m-3}{i}\lambda^i\bigg(   L_{\mathcal{D}^{k+t+m-4-i}_m}[\vec{x}]\bigg)\bigg)]\bigg)\bigg) \\
 \end{aligned}
\]
\[
\begin{aligned}
 &+ \lambda \tbinom{m-1}{i} R^2\bigg(x\bigg(  L_{\mathcal{B}^{k+t+m-2-i}_{m+1}}[\vec{x},R\bigg(\sum\limits^{m-2}_{i=0}\tbinom{m-2}{i}\lambda^i\bigg(  L_{\mathcal{B}^{k+t+m-3-i}_m}[\vec{x}]+   L_{\mathcal{A}^{k+t+m-3-i}_m}[\vec{x}]\bigg)\bigg) \\
& +\lambda (\sum\limits^{m-3}_{i=0}\tbinom{m-3}{i}\lambda^i\bigg(  L_{\mathcal{B}^{k+t+m-3-i}_m}[\vec{x}]\bigg)\bigg) 
 +\lambda \bigg(\sum\limits^{m-2}_{i=0}\tbinom{m-2}{i}\lambda^i\bigg(    L_{\mathcal{C}^{k+t+m-3-i}_m}[\vec{x}]\bigg)\bigg) 
  +\lambda^2 \bigg(\sum\limits^{m-3}_{i=0}\tbinom{m-3}{i}\lambda^i\bigg(   L_{\mathcal{D}^{k+t+m-4-i}_m}[\vec{x}]\bigg)\bigg)]
\\
 &+        L_{\mathcal{A}^{k+t+m-2-i}_{m+1}}[\vec{x},R\bigg(\sum\limits^{m-2}_{i=0}\tbinom{m-2}{i}\lambda^i\bigg(  L_{\mathcal{B}^{k+t+m-3-i}_m}[\vec{x}]+   L_{\mathcal{A}^{k+t+m-3-i}_m}[\vec{x}]\bigg)\bigg) \\
& +\lambda (\sum\limits^{m-3}_{i=0}\tbinom{m-3}{i}\lambda^i\bigg(  L_{\mathcal{B}^{k+t+m-3-i}_m}[\vec{x}]\bigg)\bigg) 
 +\lambda \bigg(\sum\limits^{m-2}_{i=0}\tbinom{m-2}{i}\lambda^i\bigg(    L_{\mathcal{C}^{k+t+m-3-i}_m}[\vec{x}]\bigg)\bigg)   +\lambda^2 \bigg(\sum\limits^{m-3}_{i=0}\tbinom{m-3}{i}\lambda^i\bigg(   L_{\mathcal{D}^{k+t+m-4-i}_m}[\vec{x}]\bigg)\bigg)]\bigg)\bigg) \allowdisplaybreaks
 \\
 &+ \lambda \tbinom{m-1}{i}R\bigg(x R\bigg(\sum\limits_{\mathcal{B}^{k+t+m-2-i}_{m}\in \mathcal{Q}} L_{\mathcal{B}^{k+t+m-2-i}_{m+1}}[\vec{x},R\bigg(\sum\limits^{m-2}_{i=0}\tbinom{m-2}{i}\lambda^i\bigg(  L_{\mathcal{B}^{k+t+m-3-i}_m}[\vec{x}]+   L_{\mathcal{A}^{k+t+m-3-i}_m}[\vec{x}]\bigg)\bigg) \\
& +\lambda (\sum\limits^{m-3}_{i=0}\tbinom{m-3}{i}\lambda^i\bigg(  L_{\mathcal{B}^{k+t+m-3-i}_m}[\vec{x}]\bigg)\bigg) 
 +\lambda \bigg(\sum\limits^{m-2}_{i=0}\tbinom{m-2}{i}\lambda^i\bigg(    L_{\mathcal{C}^{k+t+m-3-i}_m}[\vec{x}]\bigg)\bigg)   +\lambda^2 \bigg(\sum\limits^{m-3}_{i=0}\tbinom{m-3}{i}\lambda^i\bigg(   L_{\mathcal{D}^{k+t+m-4-i}_m}[\vec{x}]\bigg)\bigg)]
 \\
& +   L_{\mathcal{A}^{k+t+m-2-i}_{m+1}}[\vec{x},R\bigg(\sum\limits^{m-2}_{i=0}\tbinom{m-2}{i}\lambda^i\bigg(  L_{\mathcal{B}^{k+t+m-3-i}_m}[\vec{x}]+   L_{\mathcal{A}^{k+t+m-3-i}_m}[\vec{x}]\bigg)\bigg) \\
& +\lambda (\sum\limits^{m-3}_{i=0}\tbinom{m-3}{i}\lambda^i\bigg(  L_{\mathcal{B}^{k+t+m-3-i}_m}[\vec{x}]\bigg)\bigg) 
 +\lambda \bigg(\sum\limits^{m-2}_{i=0}\tbinom{m-2}{i}\lambda^i\bigg(    L_{\mathcal{C}^{k+t+m-3-i}_m}[\vec{x}]\bigg)\bigg)  +\lambda^2 \bigg(\sum\limits^{m-3}_{i=0}\tbinom{m-3}{i}\lambda^i\bigg(   L_{\mathcal{D}^{k+t+m-4-i}_m}[\vec{x}]\bigg)\bigg)]\bigg)\bigg)
 \\
 &+ \lambda^2 \tbinom{m-2}{i} R\bigg(x \bigg(  L_{\mathcal{B}^{k+t+m-2-i}_{m+1}}[\vec{x},R\bigg(\sum\limits^{m-2}_{i=0}\tbinom{m-2}{i}\lambda^i\bigg(  L_{\mathcal{B}^{k+t+m-3-i}_m}[\vec{x}]+   L_{\mathcal{A}^{k+t+m-3-i}_m}[\vec{x}]\bigg)\bigg) \\
& +\lambda (\sum\limits^{m-3}_{i=0}\tbinom{m-3}{i}\lambda^i\bigg(  L_{\mathcal{B}^{k+t+m-3-i}_m}[\vec{x}]\bigg)\bigg) 
 +\lambda \bigg(\sum\limits^{m-2}_{i=0}\tbinom{m-2}{i}\lambda^i\bigg(    L_{\mathcal{C}^{k+t+m-3-i}_m}[\vec{x}]\bigg)\bigg) 
  +\lambda^2 \bigg(\sum\limits^{m-3}_{i=0}\tbinom{m-3}{i}\lambda^i\bigg(   L_{\mathcal{D}^{k+t+m-4-i}_m}[\vec{x}]\bigg)\bigg)]\bigg)\bigg)\bigg)
\\
 &+ \lambda^2 \tbinom{m-1}{i} R\bigg(x \bigg(  L_{\mathcal{C}^{k+t+m-2-i}_{m+1}}[\vec{x},R\bigg(\sum\limits^{m-2}_{i=0}\tbinom{m-2}{i}\lambda^i\bigg(  L_{\mathcal{B}^{k+t+m-3-i}_m}[\vec{x}]+   L_{\mathcal{A}^{k+t+m-3-i}_m}[\vec{x}]\bigg)\bigg) \\
& +\lambda (\sum\limits^{m-3}_{i=0}\tbinom{m-3}{i}\lambda^i\bigg(  L_{\mathcal{B}^{k+t+m-3-i}_m}[\vec{x}]\bigg)\bigg) 
 +\lambda \bigg(\sum\limits^{m-2}_{i=0}\tbinom{m-2}{i}\lambda^i\bigg(    L_{\mathcal{C}^{k+t+m-3-i}_m}[\vec{x}]\bigg)\bigg)  +\lambda^2 \bigg(\sum\limits^{m-3}_{i=0}\tbinom{m-3}{i}\lambda^i\bigg(   L_{\mathcal{D}^{k+t+m-4-i}_m}[\vec{x}]\bigg)\bigg)]\bigg)\bigg)\bigg)
 \\
&+\lambda^3 \tbinom{m-2}{i}R\bigg(x \bigg(  L_{\mathcal{D}^{k+t+m-2-i}_{m+1}}[\vec{x},R\bigg(\sum\limits^{m-2}_{i=0}\tbinom{m-2}{i}\lambda^i\bigg(  L_{\mathcal{B}^{k+t+m-3-i}_m}[\vec{x}]+   L_{\mathcal{A}^{k+t+m-3-i}_m}[\vec{x}]\bigg)\bigg) \\
& +\lambda (\sum\limits^{m-3}_{i=0}\tbinom{m-3}{i}\lambda^i\bigg(  L_{\mathcal{B}^{k+t+m-3-i}_m}[\vec{x}]\bigg)\bigg) 
 +\lambda \bigg(\sum\limits^{m-2}_{i=0}\tbinom{m-2}{i}\lambda^i\bigg(    L_{\mathcal{C}^{k+t+m-3-i}_m}[\vec{x}]\bigg)\bigg) \\
 &  +\lambda^2 \bigg(\sum\limits^{m-3}_{i=0}\tbinom{m-3}{i}\lambda^i\bigg(   L_{\mathcal{D}^{k+t+m-4-i}_m}[\vec{x}]\bigg)\bigg)]\bigg)\bigg)\bigg) {=} 0;\ m {=} |\vec{x}|{\geq}3. 
\end{aligned}
\]
 Composition between relations $(R8)$ and $(R5)$, we get a sum of the expressions, which are consequences of~$(R8)$ itself.
Below, we show this example for  $\sum\limits_{t=0}^l \sum\limits^{m-1}_{i=0}\lambda^{(t)}\lambda^i
M_{t,i}$, where
\[
\begin{aligned}
 M_{t,i}&=\tbinom{m-1}{i}R^2\bigg( L^*_{\mathcal{B}^{k+t+m-2-i}_{m+1}}[x,\vec{x},\sum\limits^l_{t=0}\lambda ^{(t)}\bigg(R\bigg(  L^\op_{\mathcal{B}^{k+t-1}_2}[\vec{y}]+              L^{\op}_{\mathcal{A}^{k+t-1}_2}[\vec{y}]\bigg)+\lambda \bigg(  L^{\op}_{\mathcal{B}^{k+t-1}_2}[\vec{y}]+      L^\op_{\mathcal{C}^{k+t-1}_2}[\vec{y}]\bigg)\bigg)]
 \\
 &+       L^*_{\mathcal{A}^{k+t+m-2-i}_{m+1}}[x,\vec{x},\sum\limits^l_{t=0}\lambda ^{(t)}\bigg(R\bigg(  L^\op_{\mathcal{B}^{k+t-1}_2}[\vec{y}]+              L^{\op}_{\mathcal{A}^{k+t-1}_2}[\vec{y}]\bigg)+\lambda \bigg(  L^{\op}_{\mathcal{B}^{k+t-1}_2}[\vec{y}]+      L^\op_{\mathcal{C}^{k+t-1}_2}[\vec{y}]\bigg)\bigg)]\bigg) 
 \\
& + \lambda \tbinom{m-2}{i} R\bigg( L^*_{\mathcal{B}^{k+t+m-2-i}_{m+1}}[x,\vec{x},\sum\limits^l_{t=0}\lambda ^{(t)}\bigg(R\bigg(  L^\op_{\mathcal{B}^{k+t-1}_2}[\vec{y}]+              L^{\op}_{\mathcal{A}^{k+t-1}_2}[\vec{y}]\bigg)+\lambda \bigg(  L^{\op}_{\mathcal{B}^{k+t-1}_2}[\vec{y}]+      L^\op_{\mathcal{C}^{k+t-1}_2}[\vec{y}]\bigg)\bigg)]\bigg) 
\\
& + \lambda \tbinom{m-1}{i} R\bigg( L^*_{\mathcal{C}^{k+t+m-2-i}_{m+1}}[x,\vec{x},\sum\limits^l_{t=0}\lambda ^{(t)}\bigg(R\bigg(  L^\op_{\mathcal{B}^{k+t-1}_2}[\vec{y}]+              L^{\op}_{\mathcal{A}^{k+t-1}_2}[\vec{y}]\bigg)+\lambda \bigg(  L^{\op}_{\mathcal{B}^{k+t-1}_2}[\vec{y}]+      L^\op_{\mathcal{C}^{k+t-1}_2}[\vec{y}]\bigg)\bigg)]\bigg)\\
& + \lambda^2 \tbinom{m-2}{i} R\bigg(  L^*_{\mathcal{D}^{k+t+m-3-i}_{m+1}}[x,\vec{x},\sum\limits^l_{t=0}\lambda ^{(t)}\bigg(R\bigg(  L^\op_{\mathcal{B}^{k+t-1}_2}[\vec{y}]+              L^{\op}_{\mathcal{A}^{k+t-1}_2}[\vec{y}]\bigg)+\lambda \bigg(  L^{\op}_{\mathcal{B}^{k+t-1}_2}[\vec{y}]+      L^\op_{\mathcal{C}^{k+t-1}_2}[\vec{y}]\bigg)\bigg)]\bigg)
\\
& + \tbinom{m-1}{i}R^2\bigg(xR\bigg(  L_{\mathcal{B}^{k+t+m-2-i}_{m+1}}[\vec{x},\sum\limits^l_{t=0}\lambda ^{(t)}\bigg(R\bigg(  L^\op_{\mathcal{B}^{k+t-1}_2}[\vec{y}]+              L^{\op}_{\mathcal{A}^{k+t-1}_2}[\vec{y}]\bigg)+\lambda \bigg(  L^{\op}_{\mathcal{B}^{k+t-1}_2}[\vec{y}]+      L^\op_{\mathcal{C}^{k+t-1}_2}[\vec{y}]\bigg)\bigg)]
\\
& +   L_{\mathcal{A}^{k+t+m-2-i}_{m+1}}[\vec{x},\sum\limits^l_{t=0}\lambda ^{(t)}\bigg(R\bigg(  L^\op_{\mathcal{B}^{k+t-1}_2}[\vec{y}]+              L^{\op}_{\mathcal{A}^{k+t-1}_2}[\vec{y}]\bigg)+\lambda \bigg(  L^{\op}_{\mathcal{B}^{k+t-1}_2}[\vec{y}]+      L^\op_{\mathcal{C}^{k+t-1}_2}[\vec{y}]\bigg)\bigg)]\bigg)\bigg) 
 \\
 &+ \lambda \tbinom{m-1}{i} R^2\bigg(x\bigg(  L_{\mathcal{B}^{k+t+m-2-i}_{m+1}}[\vec{x},\sum\limits^l_{t=0}\lambda ^{(t)}\bigg(R\bigg(  L^\op_{\mathcal{B}^{k+t-1}_2}[\vec{y}]+              L^{\op}_{\mathcal{A}^{k+t-1}_2}[\vec{y}]\bigg)+\lambda \bigg(  L^{\op}_{\mathcal{B}^{k+t-1}_2}[\vec{y}]+      L^\op_{\mathcal{C}^{k+t-1}_2}[\vec{y}]\bigg)\bigg)]
 \\
 &+        L_{\mathcal{A}^{k+t+m-2-i}_{m+1}}[\vec{x},\sum\limits^l_{t=0}\lambda ^{(t)}\bigg(R\bigg(  L^\op_{\mathcal{B}^{k+t-1}_2}[\vec{y}]+              L^{\op}_{\mathcal{A}^{k+t-1}_2}[\vec{y}]\bigg)+\lambda \bigg(  L^{\op}_{\mathcal{B}^{k+t-1}_2}[\vec{y}]+      L^\op_{\mathcal{C}^{k+t-1}_2}[\vec{y}]\bigg)\bigg)]\bigg)\bigg) \allowdisplaybreaks
 \\
 &+ \lambda \tbinom{m-1}{i}R\bigg(x R\bigg(\sum\limits_{\mathcal{B}^{k+t+m-2-i}_{m}\in \mathcal{Q}} L_{\mathcal{B}^{k+t+m-2-i}_{m+1}}[\vec{x},\sum\limits^l_{t=0}\lambda ^{(t)}\bigg(R\bigg(  L^\op_{\mathcal{B}^{k+t-1}_2}[\vec{y}]+              L^{\op}_{\mathcal{A}^{k+t-1}_2}[\vec{y}]\bigg)+\lambda \bigg(  L^{\op}_{\mathcal{B}^{k+t-1}_2}[\vec{y}]+      L^\op_{\mathcal{C}^{k+t-1}_2}[\vec{y}]\bigg)\bigg)]
\\
& +   L_{\mathcal{A}^{k+t+m-2-i}_{m+1}}[\vec{x},\sum\limits^l_{t=0}\lambda ^{(t)}\bigg(R\bigg(  L^\op_{\mathcal{B}^{k+t-1}_2}[\vec{y}]+              L^{\op}_{\mathcal{A}^{k+t-1}_2}[\vec{y}]\bigg)+\lambda \bigg(  L^{\op}_{\mathcal{B}^{k+t-1}_2}[\vec{y}]+      L^\op_{\mathcal{C}^{k+t-1}_2}[\vec{y}]\bigg)\bigg)]\bigg)\bigg)
 \\
 &+ \lambda^2 \tbinom{m-2}{i} R\bigg(x \bigg(  L_{\mathcal{B}^{k+t+m-2-i}_{m+1}}[\vec{x},\sum\limits^l_{t=0}\lambda ^{(t)}\bigg(R\bigg(  L^\op_{\mathcal{B}^{k+t-1}_2}[\vec{y}]+              L^{\op}_{\mathcal{A}^{k+t-1}_2}[\vec{y}]\bigg)+\lambda \bigg(  L^{\op}_{\mathcal{B}^{k+t-1}_2}[\vec{y}]+      L^\op_{\mathcal{C}^{k+t-1}_2}[\vec{y}]\bigg)\bigg)]\bigg)\bigg)\bigg)
 \\
 &+ \lambda^2 \tbinom{m-1}{i} R\bigg(x \bigg(  L_{\mathcal{C}^{k+t+m-2-i}_{m+1}}[\vec{x},\sum\limits^l_{t=0}\lambda ^{(t)}\bigg(R\bigg(  L^\op_{\mathcal{B}^{k+t-1}_2}[\vec{y}]+              L^{\op}_{\mathcal{A}^{k+t-1}_2}[\vec{y}]\bigg)+\lambda \bigg(  L^{\op}_{\mathcal{B}^{k+t-1}_2}[\vec{y}]+      L^\op_{\mathcal{C}^{k+t-1}_2}[\vec{y}]\bigg)\bigg)]\bigg)\bigg)\bigg)
\\
&+\lambda^3 \tbinom{m-2}{i}R\bigg(x \bigg(  L_{\mathcal{D}^{k+t+m-2-i}_{m+1}}[\vec{x},\sum\limits^l_{t=0}\lambda ^{(t)}\bigg(R\bigg(  L^\op_{\mathcal{B}^{k+t-1}_2}[\vec{y}]+              L^{\op}_{\mathcal{A}^{k+t-1}_2}[\vec{y}]\bigg)+\lambda \bigg(  L^{\op}_{\mathcal{B}^{k+t-1}_2}[\vec{y}]+      L^\op_{\mathcal{C}^{k+t-1}_2}[\vec{y}]\bigg)\bigg)]\bigg)\bigg)\bigg) {=} 0;\\&\ m {=} |\vec{x}|{\geq}3,\ |\vec{y}|=2. 
\end{aligned}
\]
Composition between relations $(R8)$ and $(R6)$, we get a sum of the expressions, which are consequences of~$(R8)$ itself.
Below, we show this example for  $\sum\limits_{t=0}^l \sum\limits^{m-1}_{i=0}\lambda^{(t)}\lambda^i
M_{t,i}$, where
\[
\begin{aligned}
 M_{t,i}&=\tbinom{m-1}{i}R^2\bigg( L^*_{\mathcal{B}^{k+t+m-2-i}_{m+1}}[x,\vec{x},\sum\limits^l_{t=0}\sum\limits^{m-2}_{i=0}\lambda^{(t)}\lambda^i\left(
\tbinom{m-2}{i}R\bigg(  L^\op_{\mathcal{B}^{k+t+m-3-i}_m}[\vec{x}]
 +    L^\op_{\mathcal{A}^{k+t+m-3-i}_m}[\vec{x}]\bigg) \right. \\
 &+ \lambda \tbinom{m-3}{i} \bigg(  L^\op_{\mathcal{B}^{k+t+m-3-i}_m}[\vec{x}]\bigg) 
 + \lambda \tbinom{m-2}{i}\bigg(    L^\op_{\mathcal{C}^{k+t+m-3-i}_m}[\vec{x}]\bigg)
 \left. +\lambda^2 \tbinom{m-3}{i} \bigg(   L^\op_{\mathcal{D}^{k+t+m-4-i}_m}[\vec{x}]\bigg)\right)]
 \\
 &+       L^*_{\mathcal{A}^{k+t+m-2-i}_{m+1}}[x,\vec{x},\sum\limits^l_{t=0}\sum\limits^{m-2}_{i=0}\lambda^{(t)}\lambda^i\left(
\tbinom{m-2}{i}R\bigg(  L^\op_{\mathcal{B}^{k+t+m-3-i}_m}[\vec{x}]
 +    L^\op_{\mathcal{A}^{k+t+m-3-i}_m}[\vec{x}]\bigg) \right. \\
 &+ \lambda \tbinom{m-3}{i} \bigg(  L^\op_{\mathcal{B}^{k+t+m-3-i}_m}[\vec{x}]\bigg) 
 + \lambda \tbinom{m-2}{i}\bigg(    L^\op_{\mathcal{C}^{k+t+m-3-i}_m}[\vec{x}]\bigg)
 \left. +\lambda^2 \tbinom{m-3}{i} \bigg(   L^\op_{\mathcal{D}^{k+t+m-4-i}_m}[\vec{x}]\bigg)\right)]\bigg) 
\\
& + \lambda \tbinom{m-2}{i} R\bigg( L^*_{\mathcal{B}^{k+t+m-2-i}_{m+1}}[x,\vec{x},\sum\limits^l_{t=0}\sum\limits^{m-2}_{i=0}\lambda^{(t)}\lambda^i\left(
\tbinom{m-2}{i}R\bigg(  L^\op_{\mathcal{B}^{k+t+m-3-i}_m}[\vec{x}]
 +    L^\op_{\mathcal{A}^{k+t+m-3-i}_m}[\vec{x}]\bigg) \right. \\
 &+ \lambda \tbinom{m-3}{i} \bigg(  L^\op_{\mathcal{B}^{k+t+m-3-i}_m}[\vec{x}]\bigg) 
 + \lambda \tbinom{m-2}{i}\bigg(    L^\op_{\mathcal{C}^{k+t+m-3-i}_m}[\vec{x}]\bigg)
 \left. +\lambda^2 \tbinom{m-3}{i} \bigg(   L^\op_{\mathcal{D}^{k+t+m-4-i}_m}[\vec{x}]\bigg)\right)]\bigg) 
\\
& + \lambda \tbinom{m-1}{i} R\bigg( L^*_{\mathcal{C}^{k+t+m-2-i}_{m+1}}[x,\vec{x},\sum\limits^l_{t=0}\sum\limits^{m-2}_{i=0}\lambda^{(t)}\lambda^i\left(
\tbinom{m-2}{i}R\bigg(  L^\op_{\mathcal{B}^{k+t+m-3-i}_m}[\vec{x}]
 +    L^\op_{\mathcal{A}^{k+t+m-3-i}_m}[\vec{x}]\bigg) \right. \\
 &+ \lambda \tbinom{m-3}{i} \bigg(  L^\op_{\mathcal{B}^{k+t+m-3-i}_m}[\vec{x}]\bigg) 
 + \lambda \tbinom{m-2}{i}\bigg(    L^\op_{\mathcal{C}^{k+t+m-3-i}_m}[\vec{x}]\bigg)
 \left. +\lambda^2 \tbinom{m-3}{i} \bigg(   L^\op_{\mathcal{D}^{k+t+m-4-i}_m}[\vec{x}]\bigg)\right)]\bigg)
\\
& + \lambda^2 \tbinom{m-2}{i} R\bigg(  L^*_{\mathcal{D}^{k+t+m-3-i}_{m+1}}[x,\vec{x},\sum\limits^l_{t=0}\sum\limits^{m-2}_{i=0}\lambda^{(t)}\lambda^i\left(
\tbinom{m-2}{i}R\bigg(  L^\op_{\mathcal{B}^{k+t+m-3-i}_m}[\vec{x}]
 +    L^\op_{\mathcal{A}^{k+t+m-3-i}_m}[\vec{x}]\bigg) \right. \\
 &+ \lambda \tbinom{m-3}{i} \bigg(  L^\op_{\mathcal{B}^{k+t+m-3-i}_m}[\vec{x}]\bigg) 
 + \lambda \tbinom{m-2}{i}\bigg(    L^\op_{\mathcal{C}^{k+t+m-3-i}_m}[\vec{x}]\bigg)
 \left. +\lambda^2 \tbinom{m-3}{i} \bigg(   L^\op_{\mathcal{D}^{k+t+m-4-i}_m}[\vec{x}]\bigg)\right)]\bigg)
\\
& + \tbinom{m-1}{i}R^2\bigg(xR\bigg(  L_{\mathcal{B}^{k+t+m-2-i}_{m+1}}[\vec{x},\sum\limits^l_{t=0}\sum\limits^{m-2}_{i=0}\lambda^{(t)}\lambda^i\left(
\tbinom{m-2}{i}R\bigg(  L^\op_{\mathcal{B}^{k+t+m-3-i}_m}[\vec{x}]
 +    L^\op_{\mathcal{A}^{k+t+m-3-i}_m}[\vec{x}]\bigg) \right. \\
 &+ \lambda \tbinom{m-3}{i} \bigg(  L^\op_{\mathcal{B}^{k+t+m-3-i}_m}[\vec{x}]\bigg) 
 + \lambda \tbinom{m-2}{i}\bigg(    L^\op_{\mathcal{C}^{k+t+m-3-i}_m}[\vec{x}]\bigg)
 \left. +\lambda^2 \tbinom{m-3}{i} \bigg(   L^\op_{\mathcal{D}^{k+t+m-4-i}_m}[\vec{x}]\bigg)\right)]
 \\
& +   L_{\mathcal{A}^{k+t+m-2-i}_{m+1}}[\vec{x},\sum\limits^l_{t=0}\sum\limits^{m-2}_{i=0}\lambda^{(t)}\lambda^i\left(
\tbinom{m-2}{i}R\bigg(  L^\op_{\mathcal{B}^{k+t+m-3-i}_m}[\vec{x}]
 +    L^\op_{\mathcal{A}^{k+t+m-3-i}_m}[\vec{x}]\bigg) \right. \\
 &+ \lambda \tbinom{m-3}{i} \bigg(  L^\op_{\mathcal{B}^{k+t+m-3-i}_m}[\vec{x}]\bigg) 
 + \lambda \tbinom{m-2}{i}\bigg(    L^\op_{\mathcal{C}^{k+t+m-3-i}_m}[\vec{x}]\bigg)
 \left. +\lambda^2 \tbinom{m-3}{i} \bigg(   L^\op_{\mathcal{D}^{k+t+m-4-i}_m}[\vec{x}]\bigg)\right)]\bigg)\bigg) 
 \\
 &+ \lambda \tbinom{m-1}{i} R^2\bigg(x\bigg(  L_{\mathcal{B}^{k+t+m-2-i}_{m+1}}[\vec{x},\sum\limits^l_{t=0}\sum\limits^{m-2}_{i=0}\lambda^{(t)}\lambda^i\left(
\tbinom{m-2}{i}R\bigg(  L^\op_{\mathcal{B}^{k+t+m-3-i}_m}[\vec{x}]
 +    L^\op_{\mathcal{A}^{k+t+m-3-i}_m}[\vec{x}]\bigg) \right. \\
 &+ \lambda \tbinom{m-3}{i} \bigg(  L^\op_{\mathcal{B}^{k+t+m-3-i}_m}[\vec{x}]\bigg) 
 + \lambda \tbinom{m-2}{i}\bigg(    L^\op_{\mathcal{C}^{k+t+m-3-i}_m}[\vec{x}]\bigg)
 \left. +\lambda^2 \tbinom{m-3}{i} \bigg(   L^\op_{\mathcal{D}^{k+t+m-4-i}_m}[\vec{x}]\bigg)\right)]\\
 &+        L_{\mathcal{A}^{k+t+m-2-i}_{m+1}}[\vec{x},\sum\limits^l_{t=0}\sum\limits^{m-2}_{i=0}\lambda^{(t)}\lambda^i\left(
\tbinom{m-2}{i}R\bigg(  L^\op_{\mathcal{B}^{k+t+m-3-i}_m}[\vec{x}]
 +    L^\op_{\mathcal{A}^{k+t+m-3-i}_m}[\vec{x}]\bigg) \right. \\
 &+ \lambda \tbinom{m-3}{i} \bigg(  L^\op_{\mathcal{B}^{k+t+m-3-i}_m}[\vec{x}]\bigg) 
 + \lambda \tbinom{m-2}{i}\bigg(    L^\op_{\mathcal{C}^{k+t+m-3-i}_m}[\vec{x}]\bigg)
 \left. +\lambda^2 \tbinom{m-3}{i} \bigg(   L^\op_{\mathcal{D}^{k+t+m-4-i}_m}[\vec{x}]\bigg)\right)]\bigg)\bigg) \allowdisplaybreaks
 \end{aligned}
 \]
 \[
 \begin{aligned}
 &+ \lambda \tbinom{m-1}{i}R\bigg(x R\bigg(\sum\limits_{\mathcal{B}^{k+t+m-2-i}_{m}\in \mathcal{Q}} L_{\mathcal{B}^{k+t+m-2-i}_{m+1}}[\vec{x},\sum\limits^l_{t=0}\sum\limits^{m-2}_{i=0}\lambda^{(t)}\lambda^i\left(
\tbinom{m-2}{i}R\bigg(  L^\op_{\mathcal{B}^{k+t+m-3-i}_m}[\vec{x}]
 +    L^\op_{\mathcal{A}^{k+t+m-3-i}_m}[\vec{x}]\bigg) \right. \\
 &+ \lambda \tbinom{m-3}{i} \bigg(  L^\op_{\mathcal{B}^{k+t+m-3-i}_m}[\vec{x}]\bigg) 
 + \lambda \tbinom{m-2}{i}\bigg(    L^\op_{\mathcal{C}^{k+t+m-3-i}_m}[\vec{x}]\bigg)
 \left. +\lambda^2 \tbinom{m-3}{i} \bigg(   L^\op_{\mathcal{D}^{k+t+m-4-i}_m}[\vec{x}]\bigg)\right)]
 \\
& +   L_{\mathcal{A}^{k+t+m-2-i}_{m+1}}[\vec{x},\sum\limits^l_{t=0}\sum\limits^{m-2}_{i=0}\lambda^{(t)}\lambda^i\left(
\tbinom{m-2}{i}R\bigg(  L^\op_{\mathcal{B}^{k+t+m-3-i}_m}[\vec{x}]
 +    L^\op_{\mathcal{A}^{k+t+m-3-i}_m}[\vec{x}]\bigg) \right. \\
 &+ \lambda \tbinom{m-3}{i} \bigg(  L^\op_{\mathcal{B}^{k+t+m-3-i}_m}[\vec{x}]\bigg) 
 + \lambda \tbinom{m-2}{i}\bigg(    L^\op_{\mathcal{C}^{k+t+m-3-i}_m}[\vec{x}]\bigg)
 \left. +\lambda^2 \tbinom{m-3}{i} \bigg(   L^\op_{\mathcal{D}^{k+t+m-4-i}_m}[\vec{x}]\bigg)\right)]\bigg)\bigg)
 \\
 &+ \lambda^2 \tbinom{m-2}{i} R\bigg(x \bigg(  L_{\mathcal{B}^{k+t+m-2-i}_{m+1}}[\vec{x},\sum\limits^l_{t=0}\sum\limits^{m-2}_{i=0}\lambda^{(t)}\lambda^i\left(
\tbinom{m-2}{i}R\bigg(  L^\op_{\mathcal{B}^{k+t+m-3-i}_m}[\vec{x}]
 +    L^\op_{\mathcal{A}^{k+t+m-3-i}_m}[\vec{x}]\bigg) \right. \\
 &+ \lambda \tbinom{m-3}{i} \bigg(  L^\op_{\mathcal{B}^{k+t+m-3-i}_m}[\vec{x}]\bigg) 
 + \lambda \tbinom{m-2}{i}\bigg(    L^\op_{\mathcal{C}^{k+t+m-3-i}_m}[\vec{x}]\bigg)
 \left. +\lambda^2 \tbinom{m-3}{i} \bigg(   L^\op_{\mathcal{D}^{k+t+m-4-i}_m}[\vec{x}]\bigg)\right)]\bigg)\bigg)\bigg)
 \\
 &+ \lambda^2 \tbinom{m-1}{i} R\bigg(x \bigg(  L_{\mathcal{C}^{k+t+m-2-i}_{m+1}}[\vec{x},\sum\limits^l_{t=0}\sum\limits^{m-2}_{i=0}\lambda^{(t)}\lambda^i\left(
\tbinom{m-2}{i}R\bigg(  L^\op_{\mathcal{B}^{k+t+m-3-i}_m}[\vec{x}]
 +    L^\op_{\mathcal{A}^{k+t+m-3-i}_m}[\vec{x}]\bigg) \right. \\
 &+ \lambda \tbinom{m-3}{i} \bigg(  L^\op_{\mathcal{B}^{k+t+m-3-i}_m}[\vec{x}]\bigg) 
 + \lambda \tbinom{m-2}{i}\bigg(    L^\op_{\mathcal{C}^{k+t+m-3-i}_m}[\vec{x}]\bigg)
 \left. +\lambda^2 \tbinom{m-3}{i} \bigg(   L^\op_{\mathcal{D}^{k+t+m-4-i}_m}[\vec{x}]\bigg)\right)]\bigg)\bigg)\bigg)
 \\
&+\lambda^3 \tbinom{m-2}{i}R\bigg(x \bigg(  L_{\mathcal{D}^{k+t+m-2-i}_{m+1}}[\vec{x},\sum\limits^l_{t=0}\sum\limits^{m-2}_{i=0}\lambda^{(t)}\lambda^i\left(
\tbinom{m-2}{i}R\bigg(  L^\op_{\mathcal{B}^{k+t+m-3-i}_m}[\vec{x}]
 +    L^\op_{\mathcal{A}^{k+t+m-3-i}_m}[\vec{x}]\bigg) \right. \\
 &+ \lambda \tbinom{m-3}{i} \bigg(  L^\op_{\mathcal{B}^{k+t+m-3-i}_m}[\vec{x}]\bigg) 
 + \lambda \tbinom{m-2}{i}\bigg(    L^\op_{\mathcal{C}^{k+t+m-3-i}_m}[\vec{x}]\bigg)
 \left. +\lambda^2 \tbinom{m-3}{i} \bigg(   L^\op_{\mathcal{D}^{k+t+m-4-i}_m}[\vec{x}]\bigg)\right)]\bigg)\bigg)\bigg) {=} 0; m {=} |\vec{x}|{\geq}3. 
\end{aligned}
\]
Composition between relations $(R8)$ and $(R7)$, we get a sum of the expressions, which are consequences of~$(R8)$ itself.
Below, we show this example for  $\sum\limits_{t=0}^l \sum\limits^{m-1}_{i=0}\lambda^{(t)}\lambda^i
M_{t,i}$, where
\[
\begin{aligned}
 M_{t,i}&=\tbinom{m-1}{i}R^2\bigg( L^*_{\mathcal{B}^{k+t+m-2-i}_{m+1}}[x,\vec{x}, \sum\limits^l_{t=0}\lambda^{(t)}\bigg(
 R \bigg(xR\bigg(  L_{\mathcal{B}^{k+t-1}_2}[\vec{y}] + L_{\mathcal{A}^{k+t-1}_2}[\vec{y}]\bigg)\bigg)
 + \lambda R\bigg( x\bigg(  L_{\mathcal{B}^{k+t-1}_2}[\vec{y}]
 + L_{\mathcal{A}^{k+t-1}_2}[\vec{y}]\bigg)\bigg) \\
& + \lambda \bigg( xR\bigg(  L_{\mathcal{B}^{k+t-1}_2}[\vec{y}]
 + L_{\mathcal{A}^{k+t-1}_2}[\vec{y}]\bigg)\bigg)
 + R\bigg(   L^*_{\mathcal{B}^{k+t-1}_3}[x,\vec{y}]+    L^*_{\mathcal{A}^{k+t-1}_3}[x,\vec{y}]\bigg)  + \lambda \bigg(   L^*_{\mathcal{B}^{k+t-1}_3}[x,\vec{y}]
 +    L^*_{\mathcal{C}^{k+t-1}_3}[x,\vec{y}]\bigg)\\
 &+ \lambda^2 \bigg(x\bigg(  L_{\mathcal{B}^{k+t-1}_2}[\vec{y}]
 +       L_{\mathcal{C}^{k+t-1}_2}[\vec{y}]\bigg)\bigg)\bigg)]
 +       L^*_{\mathcal{A}^{k+t+m-2-i}_{m+1}}[x,\vec{x}, \sum\limits^l_{t=0}\lambda^{(t)}\bigg(
 R \bigg(xR\bigg(  L_{\mathcal{B}^{k+t-1}_2}[\vec{y}]
 + L_{\mathcal{A}^{k+t-1}_2}[\vec{y}]\bigg)\bigg)\\
 &+ \lambda R\bigg( x\bigg(  L_{\mathcal{B}^{k+t-1}_2}[\vec{y}]
 + L_{\mathcal{A}^{k+t-1}_2}[\vec{y}]\bigg)\bigg) 
 + \lambda \bigg( xR\bigg(  L_{\mathcal{B}^{k+t-1}_2}[\vec{y}]
 + L_{\mathcal{A}^{k+t-1}_2}[\vec{y}]\bigg)\bigg)
 + R\bigg(   L^*_{\mathcal{B}^{k+t-1}_3}[x,\vec{y}]+    L^*_{\mathcal{A}^{k+t-1}_3}[x,\vec{y}]\bigg) \\
& + \lambda \bigg(   L^*_{\mathcal{B}^{k+t-1}_3}[x,\vec{y}]
 +    L^*_{\mathcal{C}^{k+t-1}_3}[x,\vec{y}]\bigg)
 + \lambda^2 \bigg(x\bigg(  L_{\mathcal{B}^{k+t-1}_2}[\vec{y}]
 +       L_{\mathcal{C}^{k+t-1}_2}[\vec{y}]\bigg)\bigg)\bigg)]\bigg) 
 \\
& + \lambda \tbinom{m-2}{i} R\bigg( L^*_{\mathcal{B}^{k+t+m-2-i}_{m+1}}[x,\vec{x}, \sum\limits^l_{t=0}\lambda^{(t)}\bigg(
 R \bigg(xR\bigg(  L_{\mathcal{B}^{k+t-1}_2}[\vec{y}] + L_{\mathcal{A}^{k+t-1}_2}[\vec{y}]\bigg)\bigg)
 + \lambda R\bigg( x\bigg(  L_{\mathcal{B}^{k+t-1}_2}[\vec{y}]
 + L_{\mathcal{A}^{k+t-1}_2}[\vec{y}]\bigg)\bigg) \\
 \end{aligned}
\]
\[
\begin{aligned}
& + \lambda \bigg( xR\bigg(  L_{\mathcal{B}^{k+t-1}_2}[\vec{y}]
 + L_{\mathcal{A}^{k+t-1}_2}[\vec{y}]\bigg)\bigg)
 + R\bigg(   L^*_{\mathcal{B}^{k+t-1}_3}[x,\vec{y}]+    L^*_{\mathcal{A}^{k+t-1}_3}[x,\vec{y}]\bigg) \\
& + \lambda \bigg(   L^*_{\mathcal{B}^{k+t-1}_3}[x,\vec{y}]
 +    L^*_{\mathcal{C}^{k+t-1}_3}[x,\vec{y}]\bigg)
 + \lambda^2 \bigg(x\bigg(  L_{\mathcal{B}^{k+t-1}_2}[\vec{y}]
 +       L_{\mathcal{C}^{k+t-1}_2}[\vec{y}]\bigg)\bigg)\bigg)]\bigg) 
\\
& + \lambda \tbinom{m-1}{i} R\bigg( L^*_{\mathcal{C}^{k+t+m-2-i}_{m+1}}[x,\vec{x}, \sum\limits^l_{t=0}\lambda^{(t)}\bigg(
 R \bigg(xR\bigg(  L_{\mathcal{B}^{k+t-1}_2}[\vec{y}] + L_{\mathcal{A}^{k+t-1}_2}[\vec{y}]\bigg)\bigg)
 + \lambda R\bigg( x\bigg(  L_{\mathcal{B}^{k+t-1}_2}[\vec{y}]
 + L_{\mathcal{A}^{k+t-1}_2}[\vec{y}]\bigg)\bigg) \\
& + \lambda \bigg( xR\bigg(  L_{\mathcal{B}^{k+t-1}_2}[\vec{y}]
 + L_{\mathcal{A}^{k+t-1}_2}[\vec{y}]\bigg)\bigg)
 + R\bigg(   L^*_{\mathcal{B}^{k+t-1}_3}[x,\vec{y}]+    L^*_{\mathcal{A}^{k+t-1}_3}[x,\vec{y}]\bigg)  + \lambda \bigg(   L^*_{\mathcal{B}^{k+t-1}_3}[x,\vec{y}]
 +    L^*_{\mathcal{C}^{k+t-1}_3}[x,\vec{y}]\bigg)\\
 &+ \lambda^2 \bigg(x\bigg(  L_{\mathcal{B}^{k+t-1}_2}[\vec{y}]
 +       L_{\mathcal{C}^{k+t-1}_2}[\vec{y}]\bigg)\bigg)\bigg)]\bigg)
+ \lambda^2 \tbinom{m-2}{i} R\bigg(  L^*_{\mathcal{D}^{k+t+m-3-i}_{m+1}}[x,\vec{x}, \sum\limits^l_{t=0}\lambda^{(t)}\bigg(
 R \bigg(xR\bigg(  L_{\mathcal{B}^{k+t-1}_2}[\vec{y}] + L_{\mathcal{A}^{k+t-1}_2}[\vec{y}]\bigg)\bigg)\\
 &+ \lambda R\bigg( x\bigg(  L_{\mathcal{B}^{k+t-1}_2}[\vec{y}]
 + L_{\mathcal{A}^{k+t-1}_2}[\vec{y}]\bigg)\bigg) + \lambda \bigg( xR\bigg(  L_{\mathcal{B}^{k+t-1}_2}[\vec{y}]
 + L_{\mathcal{A}^{k+t-1}_2}[\vec{y}]\bigg)\bigg)
 + R\bigg(   L^*_{\mathcal{B}^{k+t-1}_3}[x,\vec{y}]+    L^*_{\mathcal{A}^{k+t-1}_3}[x,\vec{y}]\bigg) \\
& + \lambda \bigg(   L^*_{\mathcal{B}^{k+t-1}_3}[x,\vec{y}]
 +    L^*_{\mathcal{C}^{k+t-1}_3}[x,\vec{y}]\bigg)
 + \lambda^2 \bigg(x\bigg(  L_{\mathcal{B}^{k+t-1}_2}[\vec{y}]
 +       L_{\mathcal{C}^{k+t-1}_2}[\vec{y}]\bigg)\bigg)\bigg)]\bigg)
\\
& + \tbinom{m-1}{i}R^2\bigg(xR\bigg(  L_{\mathcal{B}^{k+t+m-2-i}_{m+1}}[\vec{x}, \sum\limits^l_{t=0}\lambda^{(t)}\bigg(
 R \bigg(xR\bigg(  L_{\mathcal{B}^{k+t-1}_2}[\vec{y}] + L_{\mathcal{A}^{k+t-1}_2}[\vec{y}]\bigg)\bigg)
 + \lambda R\bigg( x\bigg(  L_{\mathcal{B}^{k+t-1}_2}[\vec{y}]
 + L_{\mathcal{A}^{k+t-1}_2}[\vec{y}]\bigg)\bigg) \\
& + \lambda \bigg( xR\bigg(  L_{\mathcal{B}^{k+t-1}_2}[\vec{y}]
 + L_{\mathcal{A}^{k+t-1}_2}[\vec{y}]\bigg)\bigg)
 + R\bigg(   L^*_{\mathcal{B}^{k+t-1}_3}[x,\vec{y}]+    L^*_{\mathcal{A}^{k+t-1}_3}[x,\vec{y}]\bigg)  + \lambda \bigg(   L^*_{\mathcal{B}^{k+t-1}_3}[x,\vec{y}]
 +    L^*_{\mathcal{C}^{k+t-1}_3}[x,\vec{y}]\bigg)\\
& + \lambda^2 \bigg(x\bigg(  L_{\mathcal{B}^{k+t-1}_2}[\vec{y}]
 +       L_{\mathcal{C}^{k+t-1}_2}[\vec{y}]\bigg)\bigg)\bigg)]
 +   L_{\mathcal{A}^{k+t+m-2-i}_{m+1}}[\vec{x}, \sum\limits^l_{t=0}\lambda^{(t)}\bigg(
 R \bigg(xR\bigg(  L_{\mathcal{B}^{k+t-1}_2}[\vec{y}]
 + L_{\mathcal{A}^{k+t-1}_2}[\vec{y}]\bigg)\bigg)\\
 &+ \lambda R\bigg( x\bigg(  L_{\mathcal{B}^{k+t-1}_2}[\vec{y}]
 + L_{\mathcal{A}^{k+t-1}_2}[\vec{y}]\bigg)\bigg)  + \lambda \bigg( xR\bigg(  L_{\mathcal{B}^{k+t-1}_2}[\vec{y}]
 + L_{\mathcal{A}^{k+t-1}_2}[\vec{y}]\bigg)\bigg)
 + R\bigg(   L^*_{\mathcal{B}^{k+t-1}_3}[x,\vec{y}]+    L^*_{\mathcal{A}^{k+t-1}_3}[x,\vec{y}]\bigg) \\
& + \lambda \bigg(   L^*_{\mathcal{B}^{k+t-1}_3}[x,\vec{y}]
 +    L^*_{\mathcal{C}^{k+t-1}_3}[x,\vec{y}]\bigg)
 + \lambda^2 \bigg(x\bigg(  L_{\mathcal{B}^{k+t-1}_2}[\vec{y}]
 +       L_{\mathcal{C}^{k+t-1}_2}[\vec{y}]\bigg)\bigg)\bigg)]\bigg)\bigg) 
 \\
 &+ \lambda \tbinom{m-1}{i} R^2\bigg(x\bigg(  L_{\mathcal{B}^{k+t+m-2-i}_{m+1}}[\vec{x}, \sum\limits^l_{t=0}\lambda^{(t)}\bigg(
 R \bigg(xR\bigg(  L_{\mathcal{B}^{k+t-1}_2}[\vec{y}] + L_{\mathcal{A}^{k+t-1}_2}[\vec{y}]\bigg)\bigg)
 + \lambda R\bigg( x\bigg(  L_{\mathcal{B}^{k+t-1}_2}[\vec{y}]
 + L_{\mathcal{A}^{k+t-1}_2}[\vec{y}]\bigg)\bigg) 
 \\
& + \lambda \bigg( xR\bigg(  L_{\mathcal{B}^{k+t-1}_2}[\vec{y}]
 + L_{\mathcal{A}^{k+t-1}_2}[\vec{y}]\bigg)\bigg)
 + R\bigg(   L^*_{\mathcal{B}^{k+t-1}_3}[x,\vec{y}]+    L^*_{\mathcal{A}^{k+t-1}_3}[x,\vec{y}]\bigg)  + \lambda \bigg(   L^*_{\mathcal{B}^{k+t-1}_3}[x,\vec{y}]
 +    L^*_{\mathcal{C}^{k+t-1}_3}[x,\vec{y}]\bigg)\\
& + \lambda^2 \bigg(x\bigg(  L_{\mathcal{B}^{k+t-1}_2}[\vec{y}]
 +       L_{\mathcal{C}^{k+t-1}_2}[\vec{y}]\bigg)\bigg)\bigg)]
 +        L_{\mathcal{A}^{k+t+m-2-i}_{m+1}}[\vec{x}, \sum\limits^l_{t=0}\lambda^{(t)}\bigg(
 R \bigg(xR\bigg(  L_{\mathcal{B}^{k+t-1}_2}[\vec{y}] + L_{\mathcal{A}^{k+t-1}_2}[\vec{y}]\bigg)\bigg)\\
 &+ \lambda R\bigg( x\bigg(  L_{\mathcal{B}^{k+t-1}_2}[\vec{y}]
 + L_{\mathcal{A}^{k+t-1}_2}[\vec{y}]\bigg)\bigg)  + \lambda \bigg( xR\bigg(  L_{\mathcal{B}^{k+t-1}_2}[\vec{y}]
 + L_{\mathcal{A}^{k+t-1}_2}[\vec{y}]\bigg)\bigg)
 + R\bigg(   L^*_{\mathcal{B}^{k+t-1}_3}[x,\vec{y}]+    L^*_{\mathcal{A}^{k+t-1}_3}[x,\vec{y}]\bigg) \\
& + \lambda \bigg(   L^*_{\mathcal{B}^{k+t-1}_3}[x,\vec{y}]
 +    L^*_{\mathcal{C}^{k+t-1}_3}[x,\vec{y}]\bigg)
 + \lambda^2 \bigg(x\bigg(  L_{\mathcal{B}^{k+t-1}_2}[\vec{y}]
 +       L_{\mathcal{C}^{k+t-1}_2}[\vec{y}]\bigg)\bigg)\bigg)]\bigg)\bigg) \allowdisplaybreaks
 \\
 &+ \lambda \tbinom{m-1}{i}R\bigg(x R\bigg(\sum\limits_{\mathcal{B}^{k+t+m-2-i}_{m}\in \mathcal{Q}} L_{\mathcal{B}^{k+t+m-2-i}_{m+1}}[\vec{x}, \sum\limits^l_{t=0}\lambda^{(t)}\bigg(
 R \bigg(xR\bigg(  L_{\mathcal{B}^{k+t-1}_2}[\vec{y}] + L_{\mathcal{A}^{k+t-1}_2}[\vec{y}]\bigg)\bigg)\\
 &+ \lambda R\bigg( x\bigg(  L_{\mathcal{B}^{k+t-1}_2}[\vec{y}]
 + L_{\mathcal{A}^{k+t-1}_2}[\vec{y}]\bigg)\bigg)  + \lambda \bigg( xR\bigg(  L_{\mathcal{B}^{k+t-1}_2}[\vec{y}]
 + L_{\mathcal{A}^{k+t-1}_2}[\vec{y}]\bigg)\bigg)
 + R\bigg(   L^*_{\mathcal{B}^{k+t-1}_3}[x,\vec{y}]+    L^*_{\mathcal{A}^{k+t-1}_3}[x,\vec{y}]\bigg) \\
 \end{aligned}
\]
\[
\begin{aligned}
& + \lambda \bigg(   L^*_{\mathcal{B}^{k+t-1}_3}[x,\vec{y}]
 +    L^*_{\mathcal{C}^{k+t-1}_3}[x,\vec{y}]\bigg)
 + \lambda^2 \bigg(x\bigg(  L_{\mathcal{B}^{k+t-1}_2}[\vec{y}]
 +       L_{\mathcal{C}^{k+t-1}_2}[\vec{y}]\bigg)\bigg)\bigg)]
  \\
& +   L_{\mathcal{A}^{k+t+m-2-i}_{m+1}}[\vec{x}, \sum\limits^l_{t=0}\lambda^{(t)}\bigg(
 R \bigg(xR\bigg(  L_{\mathcal{B}^{k+t-1}_2}[\vec{y}] + L_{\mathcal{A}^{k+t-1}_2}[\vec{y}]\bigg)\bigg)
 + \lambda R\bigg( x\bigg(  L_{\mathcal{B}^{k+t-1}_2}[\vec{y}]
 + L_{\mathcal{A}^{k+t-1}_2}[\vec{y}]\bigg)\bigg) \\
& + \lambda \bigg( xR\bigg(  L_{\mathcal{B}^{k+t-1}_2}[\vec{y}]
 + L_{\mathcal{A}^{k+t-1}_2}[\vec{y}]\bigg)\bigg)
 + R\bigg(   L^*_{\mathcal{B}^{k+t-1}_3}[x,\vec{y}]+    L^*_{\mathcal{A}^{k+t-1}_3}[x,\vec{y}]\bigg) \\
& + \lambda \bigg(   L^*_{\mathcal{B}^{k+t-1}_3}[x,\vec{y}]
 +    L^*_{\mathcal{C}^{k+t-1}_3}[x,\vec{y}]\bigg)
 + \lambda^2 \bigg(x\bigg(  L_{\mathcal{B}^{k+t-1}_2}[\vec{y}]
 +       L_{\mathcal{C}^{k+t-1}_2}[\vec{y}]\bigg)\bigg)\bigg)]\bigg)\bigg)
 \\
 &+ \lambda^2 \tbinom{m-2}{i} R\bigg(x \bigg(  L_{\mathcal{B}^{k+t+m-2-i}_{m+1}}[\vec{x}, \sum\limits^l_{t=0}\lambda^{(t)}\bigg(
 R \bigg(xR\bigg(  L_{\mathcal{B}^{k+t-1}_2}[\vec{y}]
+ L_{\mathcal{A}^{k+t-1}_2}[\vec{y}]\bigg)\bigg)
 + \lambda R\bigg( x\bigg(  L_{\mathcal{B}^{k+t-1}_2}[\vec{y}]
 + L_{\mathcal{A}^{k+t-1}_2}[\vec{y}]\bigg)\bigg) \\
& + \lambda \bigg( xR\bigg(  L_{\mathcal{B}^{k+t-1}_2}[\vec{y}]
 + L_{\mathcal{A}^{k+t-1}_2}[\vec{y}]\bigg)\bigg)
 + R\bigg(   L^*_{\mathcal{B}^{k+t-1}_3}[x,\vec{y}]+    L^*_{\mathcal{A}^{k+t-1}_3}[x,\vec{y}]\bigg) \\
& + \lambda \bigg(   L^*_{\mathcal{B}^{k+t-1}_3}[x,\vec{y}]
 +    L^*_{\mathcal{C}^{k+t-1}_3}[x,\vec{y}]\bigg)
 + \lambda^2 \bigg(x\bigg(  L_{\mathcal{B}^{k+t-1}_2}[\vec{y}]
 +       L_{\mathcal{C}^{k+t-1}_2}[\vec{y}]\bigg)\bigg)\bigg)]\bigg)\bigg)\bigg)
 \\
 &+ \lambda^2 \tbinom{m-1}{i} R\bigg(x \bigg(  L_{\mathcal{C}^{k+t+m-2-i}_{m+1}}[\vec{x}, \sum\limits^l_{t=0}\lambda^{(t)}\bigg(
 R \bigg(xR\bigg(  L_{\mathcal{B}^{k+t-1}_2}[\vec{y}] + L_{\mathcal{A}^{k+t-1}_2}[\vec{y}]\bigg)\bigg)
 + \lambda R\bigg( x\bigg(  L_{\mathcal{B}^{k+t-1}_2}[\vec{y}]
 + L_{\mathcal{A}^{k+t-1}_2}[\vec{y}]\bigg)\bigg) \\
& + \lambda \bigg( xR\bigg(  L_{\mathcal{B}^{k+t-1}_2}[\vec{y}]
 + L_{\mathcal{A}^{k+t-1}_2}[\vec{y}]\bigg)\bigg)
 + R\bigg(   L^*_{\mathcal{B}^{k+t-1}_3}[x,\vec{y}]+    L^*_{\mathcal{A}^{k+t-1}_3}[x,\vec{y}]\bigg) \\
& + \lambda \bigg(   L^*_{\mathcal{B}^{k+t-1}_3}[x,\vec{y}]
 +    L^*_{\mathcal{C}^{k+t-1}_3}[x,\vec{y}]\bigg)
 + \lambda^2 \bigg(x\bigg(  L_{\mathcal{B}^{k+t-1}_2}[\vec{y}]
 +       L_{\mathcal{C}^{k+t-1}_2}[\vec{y}]\bigg)\bigg)\bigg)]\bigg)\bigg)\bigg)
 \\
&+\lambda^3 \tbinom{m-2}{i}R\bigg(x \bigg(  L_{\mathcal{D}^{k+t+m-2-i}_{m+1}}[\vec{x}, \sum\limits^l_{t=0}\lambda^{(t)}\bigg(
 R \bigg(xR\bigg(  L_{\mathcal{B}^{k+t-1}_2}[\vec{y}] + L_{\mathcal{A}^{k+t-1}_2}[\vec{y}]\bigg)\bigg)
 + \lambda R\bigg( x\bigg(  L_{\mathcal{B}^{k+t-1}_2}[\vec{y}]
 + L_{\mathcal{A}^{k+t-1}_2}[\vec{y}]\bigg)\bigg) \\
& + \lambda \bigg( xR\bigg(  L_{\mathcal{B}^{k+t-1}_2}[\vec{y}]
 + L_{\mathcal{A}^{k+t-1}_2}[\vec{y}]\bigg)\bigg)
 + R\bigg(   L^*_{\mathcal{B}^{k+t-1}_3}[x,\vec{y}]+    L^*_{\mathcal{A}^{k+t-1}_3}[x,\vec{y}]\bigg) \\
& + \lambda \bigg(   L^*_{\mathcal{B}^{k+t-1}_3}[x,\vec{y}]
 +    L^*_{\mathcal{C}^{k+t-1}_3}[x,\vec{y}]\bigg)
 + \lambda^2 \bigg(x\bigg(  L_{\mathcal{B}^{k+t-1}_2}[\vec{y}]
 +       L_{\mathcal{C}^{k+t-1}_2}[\vec{y}]\bigg)\bigg)\bigg)]\bigg)\bigg)\bigg) {=} 0;\ m {=} |\vec{x}|{\geq}3, |\vec{y}|=2. 
\end{aligned}
\]
Composition between relations $(R8)$ and $(R8)$, we get a sum of the expressions, which are consequences of~$(R8)$ itself.
Below, we show this example for  $\sum\limits_{t=0}^l \sum\limits^{m-1}_{i=0}\lambda^{(t)}\lambda^i
M_{t,i}$, where
\[
\begin{aligned}
 M_{t,i}&=\tbinom{m-1}{i}R^2\bigg( L^*_{\mathcal{B}^{k+t+m-2-i}_{m+1}}[x,\vec{x},\sum\limits^l_{t=0}\sum\limits^{m-2}_{i=0}\lambda^{(t)}\lambda^i \left( 
 \tbinom{m-2}{i}R\bigg(  L^*_{\mathcal{B}^{k+t+m-3-i}_m}[x,\vec{x}]
 +    L^*_{\mathcal{A}^{k+t+m-3-i}_m}[x,\vec{x}]\bigg)\right. \\
& + \lambda \tbinom{m-3}{i} \bigg(  L^*_{\mathcal{B}^{k+t+m-3-i}_m}[x,\vec{x}]\bigg) 
 + \lambda \tbinom{m-2}{i} \bigg(    L^*_{\mathcal{C}^{k+t+m-3-i}_m}[x,\vec{x}]\bigg) + \lambda^2 \tbinom{m-3}{i} \bigg(   L^*_{\mathcal{D}^{k+t+m-4-i}_m}[x,\vec{x}]\bigg)\\
& + \tbinom{m-2}{i}R\bigg(xR\bigg(          L_{\mathcal{B}^{k+t+m-3-i}_{m}}[\vec{x}]
 +        L_{\mathcal{A}^{k+t+m-3-i}_{m}}[\vec{x}]\bigg)\bigg)  + \lambda \tbinom{m-2}{i} R\bigg(x\bigg(          L_{\mathcal{B}^{k+t+m-3-i}_{m}}[\vec{x}]
 +            L_{\mathcal{A}^{k+t+m-3-i}_{m}}[\vec{x}]\bigg)\bigg) \allowdisplaybreaks \\
& + \lambda \tbinom{m-2}{i}\bigg(x R\bigg(          L_{\mathcal{B}^{k+t+m-3-i}_{m}}[\vec{x}]
 +        L_{\mathcal{A}^{k+t+m-3-i}_{m}}[\vec{x}]\bigg)\bigg)
 + \lambda^2 \tbinom{m-3}{i} \bigg(x \bigg(          L_{\mathcal{B}^{k+t+m-3-i}_{m}}[\vec{x}]\bigg)\bigg)\bigg) \\
& + \lambda^2 \tbinom{m-2}{i} \bigg(x \bigg(            L_{\mathcal{C}^{k+t+m-3-i}_{m}}[\vec{x}]\bigg)\bigg)\bigg)
\left. +\lambda^3 \tbinom{m-3}{i}\bigg(x \bigg(             L_{\mathcal{D}^{k+t+m-3-i}_{m}}[\vec{x}]\bigg)\bigg)\bigg) \right)]
 \end{aligned}
\]
\[
\begin{aligned}
 &+       L^*_{\mathcal{A}^{k+t+m-2-i}_{m+1}}[x,\vec{x},\sum\limits^l_{t=0}\sum\limits^{m-2}_{i=0}\lambda^{(t)}\lambda^i \left( 
 \tbinom{m-2}{i}R\bigg(  L^*_{\mathcal{B}^{k+t+m-3-i}_m}[x,\vec{x}]
 +    L^*_{\mathcal{A}^{k+t+m-3-i}_m}[x,\vec{x}]\bigg)\right. \\
& + \lambda \tbinom{m-3}{i} \bigg(  L^*_{\mathcal{B}^{k+t+m-3-i}_m}[x,\vec{x}]\bigg) 
 + \lambda \tbinom{m-2}{i} \bigg(    L^*_{\mathcal{C}^{k+t+m-3-i}_m}[x,\vec{x}]\bigg)  + \lambda^2 \tbinom{m-3}{i} \bigg(   L^*_{\mathcal{D}^{k+t+m-4-i}_m}[x,\vec{x}]\bigg)\\
 &+ \tbinom{m-2}{i}R\bigg(xR\bigg(          L_{\mathcal{B}^{k+t+m-3-i}_{m}}[\vec{x}]
 +        L_{\mathcal{A}^{k+t+m-3-i}_{m}}[\vec{x}]\bigg)\bigg)  + \lambda \tbinom{m-2}{i} R\bigg(x\bigg(          L_{\mathcal{B}^{k+t+m-3-i}_{m}}[\vec{x}]
 +            L_{\mathcal{A}^{k+t+m-3-i}_{m}}[\vec{x}]\bigg)\bigg) \allowdisplaybreaks \\
& + \lambda \tbinom{m-2}{i}\bigg(x R\bigg(          L_{\mathcal{B}^{k+t+m-3-i}_{m}}[\vec{x}]
 +        L_{\mathcal{A}^{k+t+m-3-i}_{m}}[\vec{x}]\bigg)\bigg)
 + \lambda^2 \tbinom{m-3}{i} \bigg(x \bigg(          L_{\mathcal{B}^{k+t+m-3-i}_{m}}[\vec{x}]\bigg)\bigg)\bigg) \\
& + \lambda^2 \tbinom{m-2}{i} \bigg(x \bigg(            L_{\mathcal{C}^{k+t+m-3-i}_{m}}[\vec{x}]\bigg)\bigg)\bigg)
\left. +\lambda^3 \tbinom{m-3}{i}\bigg(x \bigg(             L_{\mathcal{D}^{k+t+m-3-i}_{m}}[\vec{x}]\bigg)\bigg)\bigg) \right)]\bigg) 
 \\
& + \lambda \tbinom{m-2}{i} R\bigg( L^*_{\mathcal{B}^{k+t+m-2-i}_{m+1}}[x,\vec{x},\sum\limits^l_{t=0}\sum\limits^{m-2}_{i=0}\lambda^{(t)}\lambda^i \left( 
 \tbinom{m-2}{i}R\bigg(  L^*_{\mathcal{B}^{k+t+m-3-i}_m}[x,\vec{x}]
 +    L^*_{\mathcal{A}^{k+t+m-3-i}_m}[x,\vec{x}]\bigg)\right. \\
& + \lambda \tbinom{m-3}{i} \bigg(  L^*_{\mathcal{B}^{k+t+m-3-i}_m}[x,\vec{x}]\bigg)+ \lambda \tbinom{m-2}{i} \bigg(    L^*_{\mathcal{C}^{k+t+m-3-i}_m}[x,\vec{x}]\bigg)  + \lambda^2 \tbinom{m-3}{i} \bigg(   L^*_{\mathcal{D}^{k+t+m-4-i}_m}[x,\vec{x}]\bigg)\\
 &+ \tbinom{m-2}{i}R\bigg(xR\bigg(          L_{\mathcal{B}^{k+t+m-3-i}_{m}}[\vec{x}]
 +        L_{\mathcal{A}^{k+t+m-3-i}_{m}}[\vec{x}]\bigg)\bigg) + \lambda \tbinom{m-2}{i} R\bigg(x\bigg(          L_{\mathcal{B}^{k+t+m-3-i}_{m}}[\vec{x}]
 +            L_{\mathcal{A}^{k+t+m-3-i}_{m}}[\vec{x}]\bigg)\bigg) \allowdisplaybreaks \\
& + \lambda \tbinom{m-2}{i}\bigg(x R\bigg(          L_{\mathcal{B}^{k+t+m-3-i}_{m}}[\vec{x}]
 +        L_{\mathcal{A}^{k+t+m-3-i}_{m}}[\vec{x}]\bigg)\bigg)
 + \lambda^2 \tbinom{m-3}{i} \bigg(x \bigg(          L_{\mathcal{B}^{k+t+m-3-i}_{m}}[\vec{x}]\bigg)\bigg)\bigg) \\
& + \lambda^2 \tbinom{m-2}{i} \bigg(x \bigg(            L_{\mathcal{C}^{k+t+m-3-i}_{m}}[\vec{x}]\bigg)\bigg)\bigg)
\left. +\lambda^3 \tbinom{m-3}{i}\bigg(x \bigg(             L_{\mathcal{D}^{k+t+m-3-i}_{m}}[\vec{x}]\bigg)\bigg)\bigg) \right)]\bigg) 
 \\
& + \lambda \tbinom{m-1}{i} R\bigg( L^*_{\mathcal{C}^{k+t+m-2-i}_{m+1}}[x,\vec{x},\sum\limits^l_{t=0}\sum\limits^{m-2}_{i=0}\lambda^{(t)}\lambda^i \left( 
 \tbinom{m-2}{i}R\bigg(  L^*_{\mathcal{B}^{k+t+m-3-i}_m}[x,\vec{x}]
 +    L^*_{\mathcal{A}^{k+t+m-3-i}_m}[x,\vec{x}]\bigg)\right. \\
& + \lambda \tbinom{m-3}{i} \bigg(  L^*_{\mathcal{B}^{k+t+m-3-i}_m}[x,\vec{x}]\bigg) 
 + \lambda \tbinom{m-2}{i} \bigg(    L^*_{\mathcal{C}^{k+t+m-3-i}_m}[x,\vec{x}]\bigg)  + \lambda^2 \tbinom{m-3}{i} \bigg(   L^*_{\mathcal{D}^{k+t+m-4-i}_m}[x,\vec{x}]\bigg)\\
& + \tbinom{m-2}{i}R\bigg(xR\bigg(          L_{\mathcal{B}^{k+t+m-3-i}_{m}}[\vec{x}]
 +        L_{\mathcal{A}^{k+t+m-3-i}_{m}}[\vec{x}]\bigg)\bigg)  + \lambda \tbinom{m-2}{i} R\bigg(x\bigg(          L_{\mathcal{B}^{k+t+m-3-i}_{m}}[\vec{x}]
 +            L_{\mathcal{A}^{k+t+m-3-i}_{m}}[\vec{x}]\bigg)\bigg) \allowdisplaybreaks \\
& + \lambda \tbinom{m-2}{i}\bigg(x R\bigg(          L_{\mathcal{B}^{k+t+m-3-i}_{m}}[\vec{x}]
 +        L_{\mathcal{A}^{k+t+m-3-i}_{m}}[\vec{x}]\bigg)\bigg)
 + \lambda^2 \tbinom{m-3}{i} \bigg(x \bigg(          L_{\mathcal{B}^{k+t+m-3-i}_{m}}[\vec{x}]\bigg)\bigg)\bigg) \\
& + \lambda^2 \tbinom{m-2}{i} \bigg(x \bigg(            L_{\mathcal{C}^{k+t+m-3-i}_{m}}[\vec{x}]\bigg)\bigg)\bigg)
\left. +\lambda^3 \tbinom{m-3}{i}\bigg(x \bigg(             L_{\mathcal{D}^{k+t+m-3-i}_{m}}[\vec{x}]\bigg)\bigg)\bigg) \right)]\bigg)
\\
& + \lambda^2 \tbinom{m-2}{i} R\bigg(  L^*_{\mathcal{D}^{k+t+m-3-i}_{m+1}}[x,\vec{x},\sum\limits^l_{t=0}\sum\limits^{m-2}_{i=0}\lambda^{(t)}\lambda^i \left( 
 \tbinom{m-2}{i}R\bigg(  L^*_{\mathcal{B}^{k+t+m-3-i}_m}[x,\vec{x}]
 +    L^*_{\mathcal{A}^{k+t+m-3-i}_m}[x,\vec{x}]\bigg)\right. \\
& + \lambda \tbinom{m-3}{i} \bigg(  L^*_{\mathcal{B}^{k+t+m-3-i}_m}[x,\vec{x}]\bigg) 
 + \lambda \tbinom{m-2}{i} \bigg(    L^*_{\mathcal{C}^{k+t+m-3-i}_m}[x,\vec{x}]\bigg) 
 + \lambda^2 \tbinom{m-3}{i} \bigg(   L^*_{\mathcal{D}^{k+t+m-4-i}_m}[x,\vec{x}]\bigg)\\
& + \tbinom{m-2}{i}R\bigg(xR\bigg(          L_{\mathcal{B}^{k+t+m-3-i}_{m}}[\vec{x}]
 +        L_{\mathcal{A}^{k+t+m-3-i}_{m}}[\vec{x}]\bigg)\bigg)  + \lambda \tbinom{m-2}{i} R\bigg(x\bigg(          L_{\mathcal{B}^{k+t+m-3-i}_{m}}[\vec{x}]
 +            L_{\mathcal{A}^{k+t+m-3-i}_{m}}[\vec{x}]\bigg)\bigg) \allowdisplaybreaks  \\
 &+ \lambda \tbinom{m-2}{i}\bigg(x R\bigg(          L_{\mathcal{B}^{k+t+m-3-i}_{m}}[\vec{x}]
 +        L_{\mathcal{A}^{k+t+m-3-i}_{m}}[\vec{x}]\bigg)\bigg) + \lambda^2 \tbinom{m-3}{i} \bigg(x \bigg(          L_{\mathcal{B}^{k+t+m-3-i}_{m}}[\vec{x}]\bigg)\bigg)\bigg) \\
& + \lambda^2 \tbinom{m-2}{i} \bigg(x \bigg(            L_{\mathcal{C}^{k+t+m-3-i}_{m}}[\vec{x}]\bigg)\bigg)\bigg)
\left. +\lambda^3 \tbinom{m-3}{i}\bigg(x \bigg(             L_{\mathcal{D}^{k+t+m-3-i}_{m}}[\vec{x}]\bigg)\bigg)\bigg) \right)]\bigg)
\\
& + \tbinom{m-1}{i}R^2\bigg(xR\bigg(  L_{\mathcal{B}^{k+t+m-2-i}_{m+1}}[\vec{x},\sum\limits^l_{t=0}\sum\limits^{m-2}_{i=0}\lambda^{(t)}\lambda^i \left( 
 \tbinom{m-2}{i}R\bigg(  L^*_{\mathcal{B}^{k+t+m-3-i}_m}[x,\vec{x}]
 +    L^*_{\mathcal{A}^{k+t+m-3-i}_m}[x,\vec{x}]\bigg)\right. \\
 \end{aligned}
\]
\[
\begin{aligned}
& + \lambda \tbinom{m-3}{i} \bigg(  L^*_{\mathcal{B}^{k+t+m-3-i}_m}[x,\vec{x}]\bigg) + \lambda \tbinom{m-2}{i} \bigg(    L^*_{\mathcal{C}^{k+t+m-3-i}_m}[x,\vec{x}]\bigg)  + \lambda^2 \tbinom{m-3}{i} \bigg(   L^*_{\mathcal{D}^{k+t+m-4-i}_m}[x,\vec{x}]\bigg)\\
 &+ \tbinom{m-2}{i}R\bigg(xR\bigg(          L_{\mathcal{B}^{k+t+m-3-i}_{m}}[\vec{x}]
 +        L_{\mathcal{A}^{k+t+m-3-i}_{m}}[\vec{x}]\bigg)\bigg) + \lambda \tbinom{m-2}{i} R\bigg(x\bigg(          L_{\mathcal{B}^{k+t+m-3-i}_{m}}[\vec{x}]
 +            L_{\mathcal{A}^{k+t+m-3-i}_{m}}[\vec{x}]\bigg)\bigg) \allowdisplaybreaks \\
& + \lambda \tbinom{m-2}{i}\bigg(x R\bigg(          L_{\mathcal{B}^{k+t+m-3-i}_{m}}[\vec{x}]
 +        L_{\mathcal{A}^{k+t+m-3-i}_{m}}[\vec{x}]\bigg)\bigg)
 + \lambda^2 \tbinom{m-3}{i} \bigg(x \bigg(          L_{\mathcal{B}^{k+t+m-3-i}_{m}}[\vec{x}]\bigg)\bigg)\bigg) \\
& + \lambda^2 \tbinom{m-2}{i} \bigg(x \bigg(            L_{\mathcal{C}^{k+t+m-3-i}_{m}}[\vec{x}]\bigg)\bigg)\bigg)
\left. +\lambda^3 \tbinom{m-3}{i}\bigg(x \bigg(             L_{\mathcal{D}^{k+t+m-3-i}_{m}}[\vec{x}]\bigg)\bigg)\bigg) \right)]
 \\
& +   L_{\mathcal{A}^{k+t+m-2-i}_{m+1}}[\vec{x},,\sum\limits^l_{t=0}\sum\limits^{m-2}_{i=0}\lambda^{(t)}\lambda^i \left( 
 \tbinom{m-2}{i}R\bigg(  L^*_{\mathcal{B}^{k+t+m-3-i}_m}[x,\vec{x}]
 +    L^*_{\mathcal{A}^{k+t+m-3-i}_m}[x,\vec{x}]\bigg)\right. \\
& + \lambda \tbinom{m-3}{i} \bigg(  L^*_{\mathcal{B}^{k+t+m-3-i}_m}[x,\vec{x}]\bigg) 
 + \lambda \tbinom{m-2}{i} \bigg(    L^*_{\mathcal{C}^{k+t+m-3-i}_m}[x,\vec{x}]\bigg)  + \lambda^2 \tbinom{m-3}{i} \bigg(   L^*_{\mathcal{D}^{k+t+m-4-i}_m}[x,\vec{x}]\bigg)\\
 &+ \tbinom{m-2}{i}R\bigg(xR\bigg(          L_{\mathcal{B}^{k+t+m-3-i}_{m}}[\vec{x}]
 +        L_{\mathcal{A}^{k+t+m-3-i}_{m}}[\vec{x}]\bigg)\bigg)  + \lambda \tbinom{m-2}{i} R\bigg(x\bigg(          L_{\mathcal{B}^{k+t+m-3-i}_{m}}[\vec{x}]
 +            L_{\mathcal{A}^{k+t+m-3-i}_{m}}[\vec{x}]\bigg)\bigg) \allowdisplaybreaks \\
& + \lambda \tbinom{m-2}{i}\bigg(x R\bigg(          L_{\mathcal{B}^{k+t+m-3-i}_{m}}[\vec{x}]
 +        L_{\mathcal{A}^{k+t+m-3-i}_{m}}[\vec{x}]\bigg)\bigg)
 + \lambda^2 \tbinom{m-3}{i} \bigg(x \bigg(          L_{\mathcal{B}^{k+t+m-3-i}_{m}}[\vec{x}]\bigg)\bigg)\bigg) \\
& + \lambda^2 \tbinom{m-2}{i} \bigg(x \bigg(            L_{\mathcal{C}^{k+t+m-3-i}_{m}}[\vec{x}]\bigg)\bigg)\bigg)
\left. +\lambda^3 \tbinom{m-3}{i}\bigg(x \bigg(             L_{\mathcal{D}^{k+t+m-3-i}_{m}}[\vec{x}]\bigg)\bigg)\bigg) \right)]\bigg)\bigg) 
 \\
 &+ \lambda \tbinom{m-1}{i} R^2\bigg(x\bigg(  L_{\mathcal{B}^{k+t+m-2-i}_{m+1}}[\vec{x},\sum\limits^l_{t=0}\sum\limits^{m-2}_{i=0}\lambda^{(t)}\lambda^i \left( 
 \tbinom{m-2}{i}R\bigg(  L^*_{\mathcal{B}^{k+t+m-3-i}_m}[x,\vec{x}]
 +    L^*_{\mathcal{A}^{k+t+m-3-i}_m}[x,\vec{x}]\bigg)\right. \\
& + \lambda \tbinom{m-3}{i} \bigg(  L^*_{\mathcal{B}^{k+t+m-3-i}_m}[x,\vec{x}]\bigg) 
 + \lambda \tbinom{m-2}{i} \bigg(    L^*_{\mathcal{C}^{k+t+m-3-i}_m}[x,\vec{x}]\bigg)  + \lambda^2 \tbinom{m-3}{i} \bigg(   L^*_{\mathcal{D}^{k+t+m-4-i}_m}[x,\vec{x}]\bigg)\\
&
 + \tbinom{m-2}{i}R\bigg(xR\bigg(          L_{\mathcal{B}^{k+t+m-3-i}_{m}}[\vec{x}]
 +        L_{\mathcal{A}^{k+t+m-3-i}_{m}}[\vec{x}]\bigg)\bigg) 
 + \lambda \tbinom{m-2}{i} R\bigg(x\bigg(          L_{\mathcal{B}^{k+t+m-3-i}_{m}}[\vec{x}]
 +            L_{\mathcal{A}^{k+t+m-3-i}_{m}}[\vec{x}]\bigg)\bigg) \allowdisplaybreaks  \\&+ \lambda \tbinom{m-2}{i}\bigg(x R\bigg(          L_{\mathcal{B}^{k+t+m-3-i}_{m}}[\vec{x}]
 +        L_{\mathcal{A}^{k+t+m-3-i}_{m}}[\vec{x}]\bigg)\bigg)
 + \lambda^2 \tbinom{m-3}{i} \bigg(x \bigg(          L_{\mathcal{B}^{k+t+m-3-i}_{m}}[\vec{x}]\bigg)\bigg)\bigg) 
 + \lambda^2 \tbinom{m-2}{i} \bigg(x \bigg(            L_{\mathcal{C}^{k+t+m-3-i}_{m}}[\vec{x}]\bigg)\bigg)\bigg)\\
 &
\left. +\lambda^3 \tbinom{m-3}{i}\bigg(x \bigg(             L_{\mathcal{D}^{k+t+m-3-i}_{m}}[\vec{x}]\bigg)\bigg)\bigg) \right)]
 +        L_{\mathcal{A}^{k+t+m-2-i}_{m+1}}[\vec{x},\sum\limits^l_{t=0}\sum\limits^{m-2}_{i=0}\lambda^{(t)}\lambda^i \left( 
 \tbinom{m-2}{i}R\bigg(  L^*_{\mathcal{B}^{k+t+m-3-i}_m}[x,\vec{x}]
 +    L^*_{\mathcal{A}^{k+t+m-3-i}_m}[x,\vec{x}]\bigg)\right. \\
& + \lambda \tbinom{m-3}{i} \bigg(  L^*_{\mathcal{B}^{k+t+m-3-i}_m}[x,\vec{x}]\bigg) 
 + \lambda \tbinom{m-2}{i} \bigg(    L^*_{\mathcal{C}^{k+t+m-3-i}_m}[x,\vec{x}]\bigg)  + \lambda^2 \tbinom{m-3}{i} \bigg(   L^*_{\mathcal{D}^{k+t+m-4-i}_m}[x,\vec{x}]\bigg)\\
 &+ \tbinom{m-2}{i}R\bigg(xR\bigg(          L_{\mathcal{B}^{k+t+m-3-i}_{m}}[\vec{x}]
 +        L_{\mathcal{A}^{k+t+m-3-i}_{m}}[\vec{x}]\bigg)\bigg)  + \lambda \tbinom{m-2}{i} R\bigg(x\bigg(          L_{\mathcal{B}^{k+t+m-3-i}_{m}}[\vec{x}]
 +            L_{\mathcal{A}^{k+t+m-3-i}_{m}}[\vec{x}]\bigg)\bigg) \allowdisplaybreaks \\
& + \lambda \tbinom{m-2}{i}\bigg(x R\bigg(          L_{\mathcal{B}^{k+t+m-3-i}_{m}}[\vec{x}]
 + L_{\mathcal{A}^{k+t+m-3-i}_{m}}[\vec{x}]\bigg)\bigg)
 + \lambda^2 \tbinom{m-3}{i} \bigg(x \bigg(          L_{\mathcal{B}^{k+t+m-3-i}_{m}}[\vec{x}]\bigg)\bigg)\bigg) \\
& + \lambda^2 \tbinom{m-2}{i} \bigg(x \bigg(            L_{\mathcal{C}^{k+t+m-3-i}_{m}}[\vec{x}]\bigg)\bigg)\bigg)
\left. +\lambda^3 \tbinom{m-3}{i}\bigg(x \bigg(             L_{\mathcal{D}^{k+t+m-3-i}_{m}}[\vec{x}]\bigg)\bigg)\bigg) \right)]\bigg)\bigg) \allowdisplaybreaks
 \\
 &+ \lambda \tbinom{m-1}{i}R\bigg(x R\bigg( L_{\mathcal{B}^{k+t+m-2-i}_{m+1}}[\vec{x},\sum\limits^l_{t=0}\sum\limits^{m-2}_{i=0}\lambda^{(t)}\lambda^i \left( 
 \tbinom{m-2}{i}R\bigg(  L^*_{\mathcal{B}^{k+t+m-3-i}_m}[x,\vec{x}]
 +    L^*_{\mathcal{A}^{k+t+m-3-i}_m}[x,\vec{x}]\bigg)\right. \\
& + \lambda \tbinom{m-3}{i} \bigg(  L^*_{\mathcal{B}^{k+t+m-3-i}_m}[x,\vec{x}]\bigg) 
 + \lambda \tbinom{m-2}{i} \bigg(    L^*_{\mathcal{C}^{k+t+m-3-i}_m}[x,\vec{x}]\bigg)  + \lambda^2 \tbinom{m-3}{i} \bigg(   L^*_{\mathcal{D}^{k+t+m-4-i}_m}[x,\vec{x}]\bigg)\\
&
 + \tbinom{m-2}{i}R\bigg(xR\bigg(          L_{\mathcal{B}^{k+t+m-3-i}_{m}}[\vec{x}]
 +        L_{\mathcal{A}^{k+t+m-3-i}_{m}}[\vec{x}]\bigg)\bigg)  + \lambda \tbinom{m-2}{i} R\bigg(x\bigg(          L_{\mathcal{B}^{k+t+m-3-i}_{m}}[\vec{x}]
 +            L_{\mathcal{A}^{k+t+m-3-i}_{m}}[\vec{x}]\bigg)\bigg) \allowdisplaybreaks 
 \end{aligned}
\]
\[
\begin{aligned}
& + \lambda \tbinom{m-2}{i}\bigg(x R\bigg(          L_{\mathcal{B}^{k+t+m-3-i}_{m}}[\vec{x}]
 +        L_{\mathcal{A}^{k+t+m-3-i}_{m}}[\vec{x}]\bigg)\bigg)
 + \lambda^2 \tbinom{m-3}{i} \bigg(x \bigg(          L_{\mathcal{B}^{k+t+m-3-i}_{m}}[\vec{x}]\bigg)\bigg)\bigg)  + \lambda^2 \tbinom{m-2}{i} \bigg(x \bigg(            L_{\mathcal{C}^{k+t+m-3-i}_{m}}[\vec{x}]\bigg)\bigg)\bigg)
 \\
 &
\left. +\lambda^3 \tbinom{m-3}{i}\bigg(x \bigg(             L_{\mathcal{D}^{k+t+m-3-i}_{m}}[\vec{x}]\bigg)\bigg)\bigg) \right)]
 +   L_{\mathcal{A}^{k+t+m-2-i}_{m+1}}[\vec{x},\sum\limits^l_{t=0}\sum\limits^{m-2}_{i=0}\lambda^{(t)}\lambda^i \left( 
 \tbinom{m-2}{i}R\bigg(  L^*_{\mathcal{B}^{k+t+m-3-i}_m}[x,\vec{x}]
 +    L^*_{\mathcal{A}^{k+t+m-3-i}_m}[x,\vec{x}]\bigg)\right. \\
& + \lambda \tbinom{m-3}{i} \bigg(  L^*_{\mathcal{B}^{k+t+m-3-i}_m}[x,\vec{x}]\bigg) 
 + \lambda \tbinom{m-2}{i} \bigg(    L^*_{\mathcal{C}^{k+t+m-3-i}_m}[x,\vec{x}]\bigg)  + \lambda^2 \tbinom{m-3}{i} \bigg(   L^*_{\mathcal{D}^{k+t+m-4-i}_m}[x,\vec{x}]\bigg)\\
&
 + \tbinom{m-2}{i}R\bigg(xR\bigg(          L_{\mathcal{B}^{k+t+m-3-i}_{m}}[\vec{x}]
 +        L_{\mathcal{A}^{k+t+m-3-i}_{m}}[\vec{x}]\bigg)\bigg) + \lambda \tbinom{m-2}{i} R\bigg(x\bigg(          L_{\mathcal{B}^{k+t+m-3-i}_{m}}[\vec{x}]
 +L_{\mathcal{A}^{k+t+m-3-i}_{m}}[\vec{x}]\bigg)\bigg) \allowdisplaybreaks \\
& + \lambda \tbinom{m-2}{i}\bigg(x R\bigg(          L_{\mathcal{B}^{k+t+m-3-i}_{m}}[\vec{x}]
 +        L_{\mathcal{A}^{k+t+m-3-i}_{m}}[\vec{x}]\bigg)\bigg)
 + \lambda^2 \tbinom{m-3}{i} \bigg(x \bigg(          L_{\mathcal{B}^{k+t+m-3-i}_{m}}[\vec{x}]\bigg)\bigg)\bigg) \\
& + \lambda^2 \tbinom{m-2}{i} \bigg(x \bigg(            L_{\mathcal{C}^{k+t+m-3-i}_{m}}[\vec{x}]\bigg)\bigg)\bigg)
\left. +\lambda^3 \tbinom{m-3}{i}\bigg(x \bigg(             L_{\mathcal{D}^{k+t+m-3-i}_{m}}[\vec{x}]\bigg)\bigg)\bigg) \right)]\bigg)\bigg)
 + \lambda^2 \tbinom{m-2}{i} R\bigg(x \bigg(  L_{\mathcal{B}^{k+t+m-2-i}_{m+1}}[\vec{x}]\bigg)\bigg)\bigg)
 \\
 &+ \lambda^2 \tbinom{m-1}{i} R\bigg(x \bigg(  L_{\mathcal{C}^{k+t+m-2-i}_{m+1}}[\vec{x},\sum\limits^l_{t=0}\sum\limits^{m-2}_{i=0}\lambda^{(t)}\lambda^i \left( 
 \tbinom{m-2}{i}R\bigg(  L^*_{\mathcal{B}^{k+t+m-3-i}_m}[x,\vec{x}]
 +    L^*_{\mathcal{A}^{k+t+m-3-i}_m}[x,\vec{x}]\bigg)\right. \\
& + \lambda \tbinom{m-3}{i} \bigg(  L^*_{\mathcal{B}^{k+t+m-3-i}_m}[x,\vec{x}]\bigg) 
 + \lambda \tbinom{m-2}{i} \bigg(    L^*_{\mathcal{C}^{k+t+m-3-i}_m}[x,\vec{x}]\bigg)  + \lambda^2 \tbinom{m-3}{i} \bigg(   L^*_{\mathcal{D}^{k+t+m-4-i}_m}[x,\vec{x}]\bigg)\\
 &
 + \tbinom{m-2}{i}R\bigg(xR\bigg(          L_{\mathcal{B}^{k+t+m-3-i}_{m}}[\vec{x}]
 +        L_{\mathcal{A}^{k+t+m-3-i}_{m}}[\vec{x}]\bigg)\bigg)  + \lambda \tbinom{m-2}{i} R\bigg(x\bigg(          L_{\mathcal{B}^{k+t+m-3-i}_{m}}[\vec{x}]
 +            L_{\mathcal{A}^{k+t+m-3-i}_{m}}[\vec{x}]\bigg)\bigg) \allowdisplaybreaks \\
 &+ \lambda \tbinom{m-2}{i}\bigg(x R\bigg(          L_{\mathcal{B}^{k+t+m-3-i}_{m}}[\vec{x}]
 +        L_{\mathcal{A}^{k+t+m-3-i}_{m}}[\vec{x}]\bigg)\bigg)
 + \lambda^2 \tbinom{m-3}{i} \bigg(x \bigg(          L_{\mathcal{B}^{k+t+m-3-i}_{m}}[\vec{x}]\bigg)\bigg)\bigg) \\
& + \lambda^2 \tbinom{m-2}{i} \bigg(x \bigg(            L_{\mathcal{C}^{k+t+m-3-i}_{m}}[\vec{x}]\bigg)\bigg)\bigg)
\left. +\lambda^3 \tbinom{m-3}{i}\bigg(x \bigg(             L_{\mathcal{D}^{k+t+m-3-i}_{m}}[\vec{x}]\bigg)\bigg)\bigg) \right)]\bigg)\bigg)\bigg)
 \\
&+\lambda^3 \tbinom{m-2}{i}R\bigg(x \bigg(  L_{\mathcal{D}^{k+t+m-2-i}_{m+1}}[\vec{x},\sum\limits^l_{t=0}\sum\limits^{m-2}_{i=0}\lambda^{(t)}\lambda^i \left( 
 \tbinom{m-2}{i}R\bigg(  L^*_{\mathcal{B}^{k+t+m-3-i}_m}[x,\vec{x}]
 +    L^*_{\mathcal{A}^{k+t+m-3-i}_m}[x,\vec{x}]\bigg)\right. \\
& + \lambda \tbinom{m-3}{i} \bigg(  L^*_{\mathcal{B}^{k+t+m-3-i}_m}[x,\vec{x}]\bigg) 
 + \lambda \tbinom{m-2}{i} \bigg(    L^*_{\mathcal{C}^{k+t+m-3-i}_m}[x,\vec{x}]\bigg)  + \lambda^2 \tbinom{m-3}{i} \bigg(   L^*_{\mathcal{D}^{k+t+m-4-i}_m}[x,\vec{x}]\bigg)\\
 &+ \tbinom{m-2}{i}R\bigg(xR\bigg(          L_{\mathcal{B}^{k+t+m-3-i}_{m}}[\vec{x}]
 +        L_{\mathcal{A}^{k+t+m-3-i}_{m}}[\vec{x}]\bigg)\bigg)  + \lambda \tbinom{m-2}{i} R\bigg(x\bigg(          L_{\mathcal{B}^{k+t+m-3-i}_{m}}[\vec{x}]
 +            L_{\mathcal{A}^{k+t+m-3-i}_{m}}[\vec{x}]\bigg)\bigg) \allowdisplaybreaks 
 \\
 &+ \lambda \tbinom{m-2}{i}\bigg(x R\bigg(          L_{\mathcal{B}^{k+t+m-3-i}_{m}}[\vec{x}]
 +        L_{\mathcal{A}^{k+t+m-3-i}_{m}}[\vec{x}]\bigg)\bigg)+ \lambda^2 \tbinom{m-3}{i} \bigg(x \bigg(          L_{\mathcal{B}^{k+t+m-3-i}_{m}}[\vec{x}]\bigg)\bigg)\bigg) \\
 &+ \lambda^2 \tbinom{m-2}{i} \bigg(x \bigg(            L_{\mathcal{C}^{k+t+m-3-i}_{m}}[\vec{x}]\bigg)\bigg)\bigg)
\left. +\lambda^3 \tbinom{m-3}{i}\bigg(x \bigg(             L_{\mathcal{D}^{k+t+m-3-i}_{m}}[\vec{x}]\bigg)\bigg)\bigg) \right)]\bigg)\bigg)\bigg) {=} 0;\ m {=} |\vec{x}|{\geq}3. 
\end{aligned}
\]
Composition between relations $(R8)$ and $(R9)$, we get a sum of the expressions, which are consequences of~$(R8)$ itself.
Below, we show this example for  $\sum\limits_{t=0}^l \sum\limits^{m-1}_{i=0}\lambda^{(t)}\lambda^i
M_{t,i}$, where
\[
\begin{aligned}
 M_{t,i}&=\tbinom{m-1}{i}R^2\bigg( L^*_{\mathcal{B}^{k+t+m-2-i}_{m+1}}[x,\vec{x},\sum\limits^l_{t=0}\lambda^{(t)}\bigg( R \bigg(R\bigg(  L^\op_{\mathcal{B}^{k+t-1}_2}[\vec{y}]+              L^\op_{\mathcal{A}^{n-1}_2}[\vec{y}]\bigg)x\bigg) +\lambda R\bigg( \bigg(  L^\op_{\mathcal{B}^{k+t-1}_2}[\vec{y}]+              L^\op_{\mathcal{A}^{n-1}_2}[\vec{y}]\bigg)x\bigg)\\
& +\lambda \bigg( R\bigg(  L^\op_{\mathcal{B}^{k+t-1}_2}[\vec{y}]
 +              L^\op_{\mathcal{A}^{n-1}_2}[\vec{y}]\bigg)x\bigg)+ R\bigg(   {L^*}^\op_{\mathcal{B}^{k+t-1}_3}[x,\vec{y}]+    {L^*}^\op_{\mathcal{A}^{n-1}_3}[x,\vec{y}]\bigg)
 + \lambda \bigg(   {L^*}^\op_{\mathcal{B}^{k+t-1}_3}[x,\vec{y}]+   {L^*}^\op_{\mathcal{C}^{k+t-1}_3}[x,\vec{y}]\bigg)\\
& + \lambda^2 \bigg(\bigg(  L^\op_{\mathcal{B}^{k+t-1}_2}[\vec{y}]+      L^\op_{\mathcal{C}^{k+t-1}_2}[\vec{y}]\bigg)x\bigg)\bigg)]
+       L^*_{\mathcal{A}^{k+t+m-2-i}_{m+1}}[x,\vec{x},\sum\limits^l_{t=0}\lambda^{(t)}\bigg( R \bigg(R\bigg(  L^\op_{\mathcal{B}^{k+t-1}_2}[\vec{y}]+              L^\op_{\mathcal{A}^{n-1}_2}[\vec{y}]\bigg)x\bigg)\\
& +\lambda R\bigg( \bigg(  L^\op_{\mathcal{B}^{k+t-1}_2}[\vec{y}]+              L^\op_{\mathcal{A}^{n-1}_2}[\vec{y}]\bigg)x\bigg)
 +\lambda \bigg( R\bigg(  L^\op_{\mathcal{B}^{k+t-1}_2}[\vec{y}]
 +              L^\op_{\mathcal{A}^{n-1}_2}[\vec{y}]\bigg)x\bigg)+ R\bigg(   {L^*}^\op_{\mathcal{B}^{k+t-1}_3}[x,\vec{y}]+    {L^*}^\op_{\mathcal{A}^{n-1}_3}[x,\vec{y}]\bigg)\\
 \end{aligned}
\]
\[
\begin{aligned}
 &+ \lambda \bigg(   {L^*}^\op_{\mathcal{B}^{k+t-1}_3}[x,\vec{y}]+   {L^*}^\op_{\mathcal{C}^{k+t-1}_3}[x,\vec{y}]\bigg) + \lambda^2 \bigg(\bigg(  L^\op_{\mathcal{B}^{k+t-1}_2}[\vec{y}]+      L^\op_{\mathcal{C}^{k+t-1}_2}[\vec{y}]\bigg)x\bigg)\bigg)]\bigg) 
 \\
& + \lambda \tbinom{m-2}{i} R\bigg( L^*_{\mathcal{B}^{k+t+m-2-i}_{m+1}}[x,\vec{x},\sum\limits^l_{t=0}\lambda^{(t)}\bigg( R \bigg(R\bigg(  L^\op_{\mathcal{B}^{k+t-1}_2}[\vec{y}]+              L^\op_{\mathcal{A}^{n-1}_2}[\vec{y}]\bigg)x\bigg) +\lambda R\bigg( \bigg(  L^\op_{\mathcal{B}^{k+t-1}_2}[\vec{y}]+              L^\op_{\mathcal{A}^{n-1}_2}[\vec{y}]\bigg)x\bigg)\\
 &+\lambda \bigg( R\bigg(  L^\op_{\mathcal{B}^{k+t-1}_2}[\vec{y}]
 +              L^\op_{\mathcal{A}^{n-1}_2}[\vec{y}]\bigg)x\bigg)+ R\bigg(   {L^*}^\op_{\mathcal{B}^{k+t-1}_3}[x,\vec{y}]+    {L^*}^\op_{\mathcal{A}^{n-1}_3}[x,\vec{y}]\bigg)
 + \lambda \bigg(   {L^*}^\op_{\mathcal{B}^{k+t-1}_3}[x,\vec{y}]+   {L^*}^\op_{\mathcal{C}^{k+t-1}_3}[x,\vec{y}]\bigg)\\
& + \lambda^2 \bigg(\bigg(  L^\op_{\mathcal{B}^{k+t-1}_2}[\vec{y}]+      L^\op_{\mathcal{C}^{k+t-1}_2}[\vec{y}]\bigg)x\bigg)\bigg)]\bigg) 
 + \lambda \tbinom{m-1}{i} R\bigg( L^*_{\mathcal{C}^{k+t+m-2-i}_{m+1}}[x,\vec{x},\sum\limits^l_{t=0}\lambda^{(t)}\bigg( R \bigg(R\bigg(  L^\op_{\mathcal{B}^{k+t-1}_2}[\vec{y}]+              L^\op_{\mathcal{A}^{n-1}_2}[\vec{y}]\bigg)x\bigg)\\
& +\lambda R\bigg( \bigg(  L^\op_{\mathcal{B}^{k+t-1}_2}[\vec{y}]+              L^\op_{\mathcal{A}^{n-1}_2}[\vec{y}]\bigg)x\bigg)
 +\lambda \bigg( R\bigg(  L^\op_{\mathcal{B}^{k+t-1}_2}[\vec{y}]
 +              L^\op_{\mathcal{A}^{n-1}_2}[\vec{y}]\bigg)x\bigg)+ R\bigg(   {L^*}^\op_{\mathcal{B}^{k+t-1}_3}[x,\vec{y}]+    {L^*}^\op_{\mathcal{A}^{n-1}_3}[x,\vec{y}]\bigg)\\
 &+ \lambda \bigg(   {L^*}^\op_{\mathcal{B}^{k+t-1}_3}[x,\vec{y}]+   {L^*}^\op_{\mathcal{C}^{k+t-1}_3}[x,\vec{y}]\bigg) +\lambda R\bigg( \bigg(  L^\op_{\mathcal{B}^{k+t-1}_2}[\vec{y}]+              L^\op_{\mathcal{A}^{n-1}_2}[\vec{y}]\bigg)x\bigg)
 +\lambda \bigg( R\bigg(  L^\op_{\mathcal{B}^{k+t-1}_2}[\vec{y}]
 +              L^\op_{\mathcal{A}^{n-1}_2}[\vec{y}]\bigg)x\bigg)\\
& + \lambda^2 \bigg(\bigg(  L^\op_{\mathcal{B}^{k+t-1}_2}[\vec{y}]+      L^\op_{\mathcal{C}^{k+t-1}_2}[\vec{y}]\bigg)x\bigg)\bigg)]\bigg)
 + \lambda^2 \tbinom{m-2}{i} R\bigg(  L^*_{\mathcal{D}^{k+t+m-3-i}_{m+1}}[x,\vec{x},\sum\limits^l_{t=0}\lambda^{(t)}\bigg( R \bigg(R\bigg(  L^\op_{\mathcal{B}^{k+t-1}_2}[\vec{y}]+              L^\op_{\mathcal{A}^{n-1}_2}[\vec{y}]\bigg)x\bigg)\\
 &+ R\bigg(   {L^*}^\op_{\mathcal{B}^{k+t-1}_3}[x,\vec{y}]+    {L^*}^\op_{\mathcal{A}^{n-1}_3}[x,\vec{y}]\bigg)
 + \lambda \bigg(   {L^*}^\op_{\mathcal{B}^{k+t-1}_3}[x,\vec{y}]+   {L^*}^\op_{\mathcal{C}^{k+t-1}_3}[x,\vec{y}]\bigg) + \lambda^2 \bigg(\bigg(  L^\op_{\mathcal{B}^{k+t-1}_2}[\vec{y}]+      L^\op_{\mathcal{C}^{k+t-1}_2}[\vec{y}]\bigg)x\bigg)\bigg)]\bigg)
 \\
& + \tbinom{m-1}{i}R^2\bigg(xR\bigg(  L_{\mathcal{B}^{k+t+m-2-i}_{m+1}}[\vec{x},\sum\limits^l_{t=0}\lambda^{(t)}\bigg( R \bigg(R\bigg(  L^\op_{\mathcal{B}^{k+t-1}_2}[\vec{y}]+              L^\op_{\mathcal{A}^{n-1}_2}[\vec{y}]\bigg)x\bigg) +\lambda R\bigg( \bigg(  L^\op_{\mathcal{B}^{k+t-1}_2}[\vec{y}]+              L^\op_{\mathcal{A}^{n-1}_2}[\vec{y}]\bigg)x\bigg)\\
 &+\lambda \bigg( R\bigg(  L^\op_{\mathcal{B}^{k+t-1}_2}[\vec{y}]
 +              L^\op_{\mathcal{A}^{n-1}_2}[\vec{y}]\bigg)x\bigg)+ R\bigg(   {L^*}^\op_{\mathcal{B}^{k+t-1}_3}[x,\vec{y}]+    {L^*}^\op_{\mathcal{A}^{n-1}_3}[x,\vec{y}]\bigg)
 + \lambda \bigg(   {L^*}^\op_{\mathcal{B}^{k+t-1}_3}[x,\vec{y}]+   {L^*}^\op_{\mathcal{C}^{k+t-1}_3}[x,\vec{y}]\bigg)\\
& + \lambda^2 \bigg(\bigg(  L^\op_{\mathcal{B}^{k+t-1}_2}[\vec{y}]+      L^\op_{\mathcal{C}^{k+t-1}_2}[\vec{y}]\bigg)x\bigg)\bigg)]
  +   L_{\mathcal{A}^{k+t+m-2-i}_{m+1}}[\vec{x},\sum\limits^l_{t=0}\lambda^{(t)}\bigg( R \bigg(R\bigg(  L^\op_{\mathcal{B}^{k+t-1}_2}[\vec{y}]+              L^\op_{\mathcal{A}^{n-1}_2}[\vec{y}]\bigg)x\bigg)\\
& +\lambda R\bigg( \bigg(  L^\op_{\mathcal{B}^{k+t-1}_2}[\vec{y}]+              L^\op_{\mathcal{A}^{n-1}_2}[\vec{y}]\bigg)x\bigg)
 +\lambda \bigg( R\bigg(  L^\op_{\mathcal{B}^{k+t-1}_2}[\vec{y}]
 +              L^\op_{\mathcal{A}^{n-1}_2}[\vec{y}]\bigg)x\bigg)+ R\bigg(   {L^*}^\op_{\mathcal{B}^{k+t-1}_3}[x,\vec{y}]+    {L^*}^\op_{\mathcal{A}^{n-1}_3}[x,\vec{y}]\bigg)\\
 &+ \lambda \bigg(   {L^*}^\op_{\mathcal{B}^{k+t-1}_3}[x,\vec{y}]+   {L^*}^\op_{\mathcal{C}^{k+t-1}_3}[x,\vec{y}]\bigg) + \lambda^2 \bigg(\bigg(  L^\op_{\mathcal{B}^{k+t-1}_2}[\vec{y}]+      L^\op_{\mathcal{C}^{k+t-1}_2}[\vec{y}]\bigg)x\bigg)\bigg)]\bigg)\bigg) 
  \\
 &+ \lambda \tbinom{m-1}{i} R^2\bigg(x\bigg(  L_{\mathcal{B}^{k+t+m-2-i}_{m+1}}[\vec{x},\sum\limits^l_{t=0}\lambda^{(t)}\bigg( R \bigg(R\bigg(  L^\op_{\mathcal{B}^{k+t-1}_2}[\vec{y}]+              L^\op_{\mathcal{A}^{n-1}_2}[\vec{y}]\bigg)x\bigg) +\lambda R\bigg( \bigg(  L^\op_{\mathcal{B}^{k+t-1}_2}[\vec{y}]+              L^\op_{\mathcal{A}^{n-1}_2}[\vec{y}]\bigg)x\bigg)\\
 &+\lambda \bigg( R\bigg(  L^\op_{\mathcal{B}^{k+t-1}_2}[\vec{y}]
 +              L^\op_{\mathcal{A}^{n-1}_2}[\vec{y}]\bigg)x\bigg)+ R\bigg(   {L^*}^\op_{\mathcal{B}^{k+t-1}_3}[x,\vec{y}]+    {L^*}^\op_{\mathcal{A}^{n-1}_3}[x,\vec{y}]\bigg)
 + \lambda \bigg(   {L^*}^\op_{\mathcal{B}^{k+t-1}_3}[x,\vec{y}]+   {L^*}^\op_{\mathcal{C}^{k+t-1}_3}[x,\vec{y}]\bigg)\\
& + \lambda^2 \bigg(\bigg(  L^\op_{\mathcal{B}^{k+t-1}_2}[\vec{y}]+      L^\op_{\mathcal{C}^{k+t-1}_2}[\vec{y}]\bigg)x\bigg)\bigg)]
 +        L_{\mathcal{A}^{k+t+m-2-i}_{m+1}}[\vec{x},\sum\limits^l_{t=0}\lambda^{(t)}\bigg( R \bigg(R\bigg(  L^\op_{\mathcal{B}^{k+t-1}_2}[\vec{y}]+              L^\op_{\mathcal{A}^{n-1}_2}[\vec{y}]\bigg)x\bigg)\\
& +\lambda R\bigg( \bigg(  L^\op_{\mathcal{B}^{k+t-1}_2}[\vec{y}]+              L^\op_{\mathcal{A}^{n-1}_2}[\vec{y}]\bigg)x\bigg)
 +\lambda \bigg( R\bigg(  L^\op_{\mathcal{B}^{k+t-1}_2}[\vec{y}]
 +              L^\op_{\mathcal{A}^{n-1}_2}[\vec{y}]\bigg)x\bigg)
 + R\bigg(   {L^*}^\op_{\mathcal{B}^{k+t-1}_3}[x,\vec{y}]+    {L^*}^\op_{\mathcal{A}^{n-1}_3}[x,\vec{y}]\bigg)\\
 &+ \lambda \bigg(   {L^*}^\op_{\mathcal{B}^{k+t-1}_3}[x,\vec{y}]+   {L^*}^\op_{\mathcal{C}^{k+t-1}_3}[x,\vec{y}]\bigg)+ \lambda^2 \bigg(\bigg(  L^\op_{\mathcal{B}^{k+t-1}_2}[\vec{y}]+      L^\op_{\mathcal{C}^{k+t-1}_2}[\vec{y}]\bigg)x\bigg)\bigg)]\bigg)\bigg) \allowdisplaybreaks
\\
 &+ \lambda \tbinom{m-1}{i}R\bigg(x R\bigg(\sum\limits_{\mathcal{B}^{k+t+m-2-i}_{m}\in \mathcal{Q}} L_{\mathcal{B}^{k+t+m-2-i}_{m+1}}[\vec{x},\sum\limits^l_{t=0}\lambda^{(t)}\bigg( R \bigg(R\bigg(  L^\op_{\mathcal{B}^{k+t-1}_2}[\vec{y}]+              L^\op_{\mathcal{A}^{n-1}_2}[\vec{y}]\bigg)x\bigg)\\
& +\lambda R\bigg( \bigg(  L^\op_{\mathcal{B}^{k+t-1}_2}[\vec{y}]+              L^\op_{\mathcal{A}^{n-1}_2}[\vec{y}]\bigg)x\bigg)
 +\lambda \bigg( R\bigg(  L^\op_{\mathcal{B}^{k+t-1}_2}[\vec{y}]
 +              L^\op_{\mathcal{A}^{n-1}_2}[\vec{y}]\bigg)x\bigg)+ R\bigg(   {L^*}^\op_{\mathcal{B}^{k+t-1}_3}[x,\vec{y}]+    {L^*}^\op_{\mathcal{A}^{n-1}_3}[x,\vec{y}]\bigg)
 \end{aligned}
 \]
 \[
 \begin{aligned}
 &+ \lambda \bigg(   {L^*}^\op_{\mathcal{B}^{k+t-1}_3}[x,\vec{y}]+   {L^*}^\op_{\mathcal{C}^{k+t-1}_3}[x,\vec{y}]\bigg) + \lambda^2 \bigg(\bigg(  L^\op_{\mathcal{B}^{k+t-1}_2}[\vec{y}]+      L^\op_{\mathcal{C}^{k+t-1}_2}[\vec{y}]\bigg)x\bigg)\bigg)]
 \\
& +   L_{\mathcal{A}^{k+t+m-2-i}_{m+1}}[\vec{x},\sum\limits^l_{t=0}\lambda^{(t)}\bigg( R \bigg(R\bigg(  L^\op_{\mathcal{B}^{k+t-1}_2}[\vec{y}]+              L^\op_{\mathcal{A}^{n-1}_2}[\vec{y}]\bigg)x\bigg)+\lambda R\bigg( \bigg(  L^\op_{\mathcal{B}^{k+t-1}_2}[\vec{y}]+              L^\op_{\mathcal{A}^{n-1}_2}[\vec{y}]\bigg)x\bigg)\\
& +\lambda \bigg( R\bigg(  L^\op_{\mathcal{B}^{k+t-1}_2}[\vec{y}]
 +              L^\op_{\mathcal{A}^{n-1}_2}[\vec{y}]\bigg)x\bigg)
 + R\bigg(   {L^*}^\op_{\mathcal{B}^{k+t-1}_3}[x,\vec{y}]+    {L^*}^\op_{\mathcal{A}^{n-1}_3}[x,\vec{y}]\bigg)
 + \lambda \bigg(   {L^*}^\op_{\mathcal{B}^{k+t-1}_3}[x,\vec{y}]+   {L^*}^\op_{\mathcal{C}^{k+t-1}_3}[x,\vec{y}]\bigg)\\
& + \lambda^2 \bigg(\bigg(  L^\op_{\mathcal{B}^{k+t-1}_2}[\vec{y}]+      L^\op_{\mathcal{C}^{k+t-1}_2}[\vec{y}]\bigg)x\bigg)\bigg)]\bigg)\bigg)
 + \lambda^2 \tbinom{m-2}{i} R\bigg(x \bigg(  L_{\mathcal{B}^{k+t+m-2-i}_{m+1}}[\vec{x},\sum\limits^l_{t=0}\lambda^{(t)}\bigg( R \bigg(R\bigg(  L^\op_{\mathcal{B}^{k+t-1}_2}[\vec{y}]+              L^\op_{\mathcal{A}^{n-1}_2}[\vec{y}]\bigg)x\bigg)\\
& +\lambda R\bigg( \bigg(  L^\op_{\mathcal{B}^{k+t-1}_2}[\vec{y}]+              L^\op_{\mathcal{A}^{n-1}_2}[\vec{y}]\bigg)x\bigg)
 +\lambda \bigg( R\bigg(  L^\op_{\mathcal{B}^{k+t-1}_2}[\vec{y}]
 +              L^\op_{\mathcal{A}^{n-1}_2}[\vec{y}]\bigg)x\bigg)
 + R\bigg(   {L^*}^\op_{\mathcal{B}^{k+t-1}_3}[x,\vec{y}]+    {L^*}^\op_{\mathcal{A}^{n-1}_3}[x,\vec{y}]\bigg)\\
 &+ \lambda \bigg(   {L^*}^\op_{\mathcal{B}^{k+t-1}_3}[x,\vec{y}]+   {L^*}^\op_{\mathcal{C}^{k+t-1}_3}[x,\vec{y}]\bigg) + \lambda^2 \bigg(\bigg(  L^\op_{\mathcal{B}^{k+t-1}_2}[\vec{y}]+      L^\op_{\mathcal{C}^{k+t-1}_2}[\vec{y}]\bigg)x\bigg)\bigg)]\bigg)\bigg)\bigg)
   \\
 &+ \lambda^2 \tbinom{m-1}{i} R\bigg(x \bigg(  L_{\mathcal{C}^{k+t+m-2-i}_{m+1}}[\vec{x},\sum\limits^l_{t=0}\lambda^{(t)}\bigg( R \bigg(R\bigg(  L^\op_{\mathcal{B}^{k+t-1}_2}[\vec{y}]+              L^\op_{\mathcal{A}^{n-1}_2}[\vec{y}]\bigg)x\bigg)
 +\lambda R\bigg( \bigg(  L^\op_{\mathcal{B}^{k+t-1}_2}[\vec{y}]+              L^\op_{\mathcal{A}^{n-1}_2}[\vec{y}]\bigg)x\bigg)\\
& +\lambda \bigg( R\bigg(  L^\op_{\mathcal{B}^{k+t-1}_2}[\vec{y}]
 +              L^\op_{\mathcal{A}^{n-1}_2}[\vec{y}]\bigg)x\bigg)+ R\bigg(   {L^*}^\op_{\mathcal{B}^{k+t-1}_3}[x,\vec{y}]+    {L^*}^\op_{\mathcal{A}^{n-1}_3}[x,\vec{y}]\bigg)
 + \lambda \bigg(   {L^*}^\op_{\mathcal{B}^{k+t-1}_3}[x,\vec{y}]+   {L^*}^\op_{\mathcal{C}^{k+t-1}_3}[x,\vec{y}]\bigg)\\
& + \lambda^2 \bigg(\bigg(  L^\op_{\mathcal{B}^{k+t-1}_2}[\vec{y}]+      L^\op_{\mathcal{C}^{k+t-1}_2}[\vec{y}]\bigg)x\bigg)\bigg)]\bigg)\bigg)\bigg)
+\lambda^3 \tbinom{m-2}{i}R\bigg(x \bigg(  L_{\mathcal{D}^{k+t+m-2-i}_{m+1}}[\vec{x},\sum\limits^l_{t=0}\lambda^{(t)}\bigg( R \bigg(R\bigg(  L^\op_{\mathcal{B}^{k+t-1}_2}[\vec{y}]+              L^\op_{\mathcal{A}^{n-1}_2}[\vec{y}]\bigg)x\bigg)\\
& +\lambda R\bigg( \bigg(  L^\op_{\mathcal{B}^{k+t-1}_2}[\vec{y}]+              L^\op_{\mathcal{A}^{n-1}_2}[\vec{y}]\bigg)x\bigg)
 +\lambda \bigg( R\bigg(  L^\op_{\mathcal{B}^{k+t-1}_2}[\vec{y}]
 +              L^\op_{\mathcal{A}^{n-1}_2}[\vec{y}]\bigg)x\bigg)+ R\bigg(   {L^*}^\op_{\mathcal{B}^{k+t-1}_3}[x,\vec{y}]+    {L^*}^\op_{\mathcal{A}^{n-1}_3}[x,\vec{y}]\bigg)\\
 &+ \lambda \bigg(   {L^*}^\op_{\mathcal{B}^{k+t-1}_3}[x,\vec{y}]+   {L^*}^\op_{\mathcal{C}^{k+t-1}_3}[x,\vec{y}]\bigg) + \lambda^2 \bigg(\bigg(  L^\op_{\mathcal{B}^{k+t-1}_2}[\vec{y}]+      L^\op_{\mathcal{C}^{k+t-1}_2}[\vec{y}]\bigg)x\bigg)\bigg)]\bigg)\bigg)\bigg) {=} 0;\ m {=} |\vec{x}|{\geq}3,\ |\vec{y}|=2. 
\end{aligned}
\]
Composition between relations $(R8)$ and $(R10)$, we get a sum of the expressions, which are consequences of~$(R8)$ itself.
Below, we show this example for  $\sum\limits_{t=0}^l \sum\limits^{m-1}_{i=0}\lambda^{(t)}\lambda^i
M_{t,i}$, where 
\[
\begin{aligned}
 M_{t,i}&=\tbinom{m-1}{i}R^2\bigg( L^*_{\mathcal{B}^{k+t+m-2-i}_{m+1}}[x,\vec{x},\sum\limits^l_{t=0}\sum\limits^{m-2}_{i=0}\lambda^{(t)}\lambda^i \left(
\tbinom{m-2}{i}R\bigg(  {L^*}^\op_{\mathcal{B}^{k+t+m-3-i}_m}[x,\vec{x}]
+   {L^*}^\op_{\mathcal{A}^{k+t+m-3-i}_m}[x,\vec{x}]\bigg) \right. \\
&+\lambda \tbinom{m-3}{i}\bigg(  {L^*}^\op_{\mathcal{B}^{k+t+m-3-i}_m}[x,\vec{x}]\bigg)
+\lambda \tbinom{m-2}{i} \bigg(    {L^*}^\op_{\mathcal{C}^{k+t+m-3-i}_m}[x,\vec{x}]\bigg)\bigg)
+\lambda^2 \tbinom{m-3}{i}\bigg(   {L^*}^\op_{\mathcal{D}^{k+t+m-4-i}_m}[x,\vec{x}]\bigg)\\
& + \tbinom{m-2}{i}R\bigg(R\bigg(          L^\op_{\mathcal{B}^{k+t+m-3-i}_{m}}[\vec{x}]
 +            L^\op_{\mathcal{A}^{k+t+m-3-i}_{m}}[\vec{x}]\bigg)x\bigg) + \lambda \tbinom{m-2}{i}R\bigg(\bigg(          L^\op_{\mathcal{B}^{k+t+m-3-i}_{m}}[\vec{x}]
 +            L^\op_{\mathcal{A}^{k+t+m-3-i}_{m}}[\vec{x}]\bigg)x\bigg)\\
& + \lambda \tbinom{m-2}{i}\bigg( R\bigg(          L^\op_{\mathcal{B}^{k+t+m-3-i}_{m}}[\vec{x}]
 +            L^\op_{\mathcal{A}^{k+t+m-3-i}_{m}}[\vec{x}]\bigg)x\bigg)
 + \lambda^2 \tbinom{m-3}{i} \bigg( \bigg(          L^\op_{\mathcal{B}^{k+t+m-3-i}_{m}}[\vec{x}]\bigg)x\bigg) \\
& + \lambda^2 \tbinom{m-2}{i}\bigg( \bigg(            L^\op_{\mathcal{C}^{k+t+m-3-i}_{m}}[\vec{x}]\bigg)x\bigg)
 \left. +\lambda^3 \tbinom{m-3}{i}\bigg( \bigg(              L^\op_{\mathcal{D}^{k+t+m-3-i}_{m}}[\vec{x}]\bigg)x\bigg) \right)]
 \\
 &+       L^*_{\mathcal{A}^{k+t+m-2-i}_{m+1}}[x,\vec{x},\sum\limits^l_{t=0}\sum\limits^{m-2}_{i=0}\lambda^{(t)}\lambda^i \left(
\tbinom{m-2}{i}R\bigg(  {L^*}^\op_{\mathcal{B}^{k+t+m-3-i}_m}[x,\vec{x}]
+   {L^*}^\op_{\mathcal{A}^{k+t+m-3-i}_m}[x,\vec{x}]\bigg) \right. \\
&+\lambda \tbinom{m-3}{i}\bigg(  {L^*}^\op_{\mathcal{B}^{k+t+m-3-i}_m}[x,\vec{x}]\bigg)
+\lambda \tbinom{m-2}{i} \bigg(    {L^*}^\op_{\mathcal{C}^{k+t+m-3-i}_m}[x,\vec{x}]\bigg)\bigg)
 +\lambda^2 \tbinom{m-3}{i}\bigg(   {L^*}^\op_{\mathcal{D}^{k+t+m-4-i}_m}[x,\vec{x}]\bigg)\\&
 + \tbinom{m-2}{i}R\bigg(R\bigg(          L^\op_{\mathcal{B}^{k+t+m-3-i}_{m}}[\vec{x}]
 +            L^\op_{\mathcal{A}^{k+t+m-3-i}_{m}}[\vec{x}]\bigg)x\bigg)
 + \lambda \tbinom{m-2}{i}R\bigg(\bigg(          L^\op_{\mathcal{B}^{k+t+m-3-i}_{m}}[\vec{x}]
 +            L^\op_{\mathcal{A}^{k+t+m-3-i}_{m}}[\vec{x}]\bigg)x\bigg)\\
 \end{aligned}
\]
\[
\begin{aligned}
& + \lambda \tbinom{m-2}{i}\bigg( R\bigg(          L^\op_{\mathcal{B}^{k+t+m-3-i}_{m}}[\vec{x}]
 +            L^\op_{\mathcal{A}^{k+t+m-3-i}_{m}}[\vec{x}]\bigg)x\bigg)
 + \lambda^2 \tbinom{m-3}{i} \bigg( \bigg(          L^\op_{\mathcal{B}^{k+t+m-3-i}_{m}}[\vec{x}]\bigg)x\bigg) \\
& + \lambda^2 \tbinom{m-2}{i}\bigg( \bigg(            L^\op_{\mathcal{C}^{k+t+m-3-i}_{m}}[\vec{x}]\bigg)x\bigg)
 \left. +\lambda^3 \tbinom{m-3}{i}\bigg( \bigg(              L^\op_{\mathcal{D}^{k+t+m-3-i}_{m}}[\vec{x}]\bigg)x\bigg) \right)]\bigg)\\ 
& + \lambda \tbinom{m-2}{i} R\bigg( L^*_{\mathcal{B}^{k+t+m-2-i}_{m+1}}[x,\vec{x},\sum\limits^l_{t=0}\sum\limits^{m-2}_{i=0}\lambda^{(t)}\lambda^i \left(
\tbinom{m-2}{i}R\bigg(  {L^*}^\op_{\mathcal{B}^{k+t+m-3-i}_m}[x,\vec{x}]
+   {L^*}^\op_{\mathcal{A}^{k+t+m-3-i}_m}[x,\vec{x}]\bigg) \right. \\
&+\lambda \tbinom{m-3}{i}\bigg(  {L^*}^\op_{\mathcal{B}^{k+t+m-3-i}_m}[x,\vec{x}]\bigg)
+\lambda \tbinom{m-2}{i} \bigg(    {L^*}^\op_{\mathcal{C}^{k+t+m-3-i}_m}[x,\vec{x}]\bigg)\bigg)+\lambda^2 \tbinom{m-3}{i}\bigg(   {L^*}^\op_{\mathcal{D}^{k+t+m-4-i}_m}[x,\vec{x}]\bigg)\\
&+ \tbinom{m-2}{i}R\bigg(R\bigg(          L^\op_{\mathcal{B}^{k+t+m-3-i}_{m}}[\vec{x}]
 +            L^\op_{\mathcal{A}^{k+t+m-3-i}_{m}}[\vec{x}]\bigg)x\bigg) + \lambda \tbinom{m-2}{i}R\bigg(\bigg(          L^\op_{\mathcal{B}^{k+t+m-3-i}_{m}}[\vec{x}]
 +            L^\op_{\mathcal{A}^{k+t+m-3-i}_{m}}[\vec{x}]\bigg)x\bigg)\\
& + \lambda \tbinom{m-2}{i}\bigg( R\bigg(          L^\op_{\mathcal{B}^{k+t+m-3-i}_{m}}[\vec{x}]
 +            L^\op_{\mathcal{A}^{k+t+m-3-i}_{m}}[\vec{x}]\bigg)x\bigg)
 + \lambda^2 \tbinom{m-3}{i} \bigg( \bigg(          L^\op_{\mathcal{B}^{k+t+m-3-i}_{m}}[\vec{x}]\bigg)x\bigg) \\
& + \lambda^2 \tbinom{m-2}{i}\bigg( \bigg(            L^\op_{\mathcal{C}^{k+t+m-3-i}_{m}}[\vec{x}]\bigg)x\bigg)
 \left. +\lambda^3 \tbinom{m-3}{i}\bigg( \bigg(              L^\op_{\mathcal{D}^{k+t+m-3-i}_{m}}[\vec{x}]\bigg)x\bigg) \right)]\bigg) 
\\
& + \lambda \tbinom{m-1}{i} R\bigg( L^*_{\mathcal{C}^{k+t+m-2-i}_{m+1}}[x,\vec{x},\sum\limits^l_{t=0}\sum\limits^{m-2}_{i=0}\lambda^{(t)}\lambda^i \left(
\tbinom{m-2}{i}R\bigg(  {L^*}^\op_{\mathcal{B}^{k+t+m-3-i}_m}[x,\vec{x}]
+   {L^*}^\op_{\mathcal{A}^{k+t+m-3-i}_m}[x,\vec{x}]\bigg) \right. \\
&+\lambda \tbinom{m-3}{i}\bigg(  {L^*}^\op_{\mathcal{B}^{k+t+m-3-i}_m}[x,\vec{x}]\bigg)
+\lambda \tbinom{m-2}{i} \bigg(    {L^*}^\op_{\mathcal{C}^{k+t+m-3-i}_m}[x,\vec{x}]\bigg)\bigg)
+\lambda^2 \tbinom{m-3}{i}\bigg(   {L^*}^\op_{\mathcal{D}^{k+t+m-4-i}_m}[x,\vec{x}]\bigg)\\
 &+ \tbinom{m-2}{i}R\bigg(R\bigg(          L^\op_{\mathcal{B}^{k+t+m-3-i}_{m}}[\vec{x}]
 +            L^\op_{\mathcal{A}^{k+t+m-3-i}_{m}}[\vec{x}]\bigg)x\bigg)
 + \lambda \tbinom{m-2}{i}R\bigg(\bigg(          L^\op_{\mathcal{B}^{k+t+m-3-i}_{m}}[\vec{x}]
 +            L^\op_{\mathcal{A}^{k+t+m-3-i}_{m}}[\vec{x}]\bigg)x\bigg)\\
& + \lambda \tbinom{m-2}{i}\bigg( R\bigg(          L^\op_{\mathcal{B}^{k+t+m-3-i}_{m}}[\vec{x}]
 +            L^\op_{\mathcal{A}^{k+t+m-3-i}_{m}}[\vec{x}]\bigg)x\bigg)
 + \lambda^2 \tbinom{m-3}{i} \bigg( \bigg(          L^\op_{\mathcal{B}^{k+t+m-3-i}_{m}}[\vec{x}]\bigg)x\bigg) \\
& + \lambda^2 \tbinom{m-2}{i}\bigg( \bigg(            L^\op_{\mathcal{C}^{k+t+m-3-i}_{m}}[\vec{x}]\bigg)x\bigg)
 \left. +\lambda^3 \tbinom{m-3}{i}\bigg( \bigg(              L^\op_{\mathcal{D}^{k+t+m-3-i}_{m}}[\vec{x}]\bigg)x\bigg) \right)]\bigg)
\\
& + \lambda^2 \tbinom{m-2}{i} R\bigg(  L^*_{\mathcal{D}^{k+t+m-3-i}_{m+1}}[x,\vec{x},\sum\limits^l_{t=0}\sum\limits^{m-2}_{i=0}\lambda^{(t)}\lambda^i \left(
\tbinom{m-2}{i}R\bigg(  {L^*}^\op_{\mathcal{B}^{k+t+m-3-i}_m}[x,\vec{x}]
+   {L^*}^\op_{\mathcal{A}^{k+t+m-3-i}_m}[x,\vec{x}]\bigg) \right. \\
&+\lambda \tbinom{m-3}{i}\bigg(  {L^*}^\op_{\mathcal{B}^{k+t+m-3-i}_m}[x,\vec{x}]\bigg)
+\lambda \tbinom{m-2}{i} \bigg(    {L^*}^\op_{\mathcal{C}^{k+t+m-3-i}_m}[x,\vec{x}]\bigg)\bigg)
+\lambda^2 \tbinom{m-3}{i}\bigg(   {L^*}^\op_{\mathcal{D}^{k+t+m-4-i}_m}[x,\vec{x}]\bigg)\\
 &+ \tbinom{m-2}{i}R\bigg(R\bigg(          L^\op_{\mathcal{B}^{k+t+m-3-i}_{m}}[\vec{x}]
 +            L^\op_{\mathcal{A}^{k+t+m-3-i}_{m}}[\vec{x}]\bigg)x\bigg)
 + \lambda \tbinom{m-2}{i}R\bigg(\bigg(          L^\op_{\mathcal{B}^{k+t+m-3-i}_{m}}[\vec{x}]
 +            L^\op_{\mathcal{A}^{k+t+m-3-i}_{m}}[\vec{x}]\bigg)x\bigg)\\
& + \lambda \tbinom{m-2}{i}\bigg( R\bigg(          L^\op_{\mathcal{B}^{k+t+m-3-i}_{m}}[\vec{x}]
 +            L^\op_{\mathcal{A}^{k+t+m-3-i}_{m}}[\vec{x}]\bigg)x\bigg)
 + \lambda^2 \tbinom{m-3}{i} \bigg( \bigg(          L^\op_{\mathcal{B}^{k+t+m-3-i}_{m}}[\vec{x}]\bigg)x\bigg) \\
& + \lambda^2 \tbinom{m-2}{i}\bigg( \bigg(            L^\op_{\mathcal{C}^{k+t+m-3-i}_{m}}[\vec{x}]\bigg)x\bigg)
 \left. +\lambda^3 \tbinom{m-3}{i}\bigg( \bigg(              L^\op_{\mathcal{D}^{k+t+m-3-i}_{m}}[\vec{x}]\bigg)x\bigg) \right)]\bigg)+ \lambda^2 \tbinom{m-2}{i}\bigg( \bigg(            L^\op_{\mathcal{C}^{k+t+m-3-i}_{m}}[\vec{x}]\bigg)x\bigg)
\\
& + \tbinom{m-1}{i}R^2\bigg(xR\bigg(  L_{\mathcal{B}^{k+t+m-2-i}_{m+1}}[\vec{x},\sum\limits^l_{t=0}\sum\limits^{m-2}_{i=0}\lambda^{(t)}\lambda^i \left(
\tbinom{m-2}{i}R\bigg(  {L^*}^\op_{\mathcal{B}^{k+t+m-3-i}_m}[x,\vec{x}]
+   {L^*}^\op_{\mathcal{A}^{k+t+m-3-i}_m}[x,\vec{x}]\bigg) \right. \\
&+\lambda \tbinom{m-3}{i}\bigg(  {L^*}^\op_{\mathcal{B}^{k+t+m-3-i}_m}[x,\vec{x}]\bigg)
+\lambda \tbinom{m-2}{i} \bigg(    {L^*}^\op_{\mathcal{C}^{k+t+m-3-i}_m}[x,\vec{x}]\bigg)\bigg)
+\lambda^2 \tbinom{m-3}{i}\bigg(   {L^*}^\op_{\mathcal{D}^{k+t+m-4-i}_m}[x,\vec{x}]\bigg)\\
& + \tbinom{m-2}{i}R\bigg(R\bigg(  L^\op_{\mathcal{B}^{k+t+m-3-i}_{m}}[\vec{x}]
 +            L^\op_{\mathcal{A}^{k+t+m-3-i}_{m}}[\vec{x}]\bigg)x\bigg)
  + \lambda \tbinom{m-2}{i}R\bigg(\bigg(          L^\op_{\mathcal{B}^{k+t+m-3-i}_{m}}[\vec{x}]
 +            L^\op_{\mathcal{A}^{k+t+m-3-i}_{m}}[\vec{x}]\bigg)x\bigg)\\
 &+ \lambda \tbinom{m-2}{i}\bigg( R\bigg(          L^\op_{\mathcal{B}^{k+t+m-3-i}_{m}}[\vec{x}]
 +            L^\op_{\mathcal{A}^{k+t+m-3-i}_{m}}[\vec{x}]\bigg)x\bigg)+ \lambda^2 \tbinom{m-3}{i} \bigg( \bigg(          L^\op_{\mathcal{B}^{k+t+m-3-i}_{m}}[\vec{x}]\bigg)x\bigg)
 \left. +\lambda^3 \tbinom{m-3}{i}\bigg( \bigg(              L^\op_{\mathcal{D}^{k+t+m-3-i}_{m}}[\vec{x}]\bigg)x\bigg) \right)]
 \end{aligned}
\]
\[
\begin{aligned}
& +   L_{\mathcal{A}^{k+t+m-2-i}_{m+1}}[\vec{x},\sum\limits^l_{t=0}\sum\limits^{m-2}_{i=0}\lambda^{(t)}\lambda^i \left(
\tbinom{m-2}{i}R\bigg(  {L^*}^\op_{\mathcal{B}^{k+t+m-3-i}_m}[x,\vec{x}]
+   {L^*}^\op_{\mathcal{A}^{k+t+m-3-i}_m}[x,\vec{x}]\bigg) \right. 
 \\
&+\lambda \tbinom{m-3}{i}\bigg(  {L^*}^\op_{\mathcal{B}^{k+t+m-3-i}_m}[x,\vec{x}]\bigg)
+\lambda \tbinom{m-2}{i} \bigg(    {L^*}^\op_{\mathcal{C}^{k+t+m-3-i}_m}[x,\vec{x}]\bigg)\bigg)+\lambda^2 \tbinom{m-3}{i}\bigg(   {L^*}^\op_{\mathcal{D}^{k+t+m-4-i}_m}[x,\vec{x}]\bigg)\\
&
 + \tbinom{m-2}{i}R\bigg(R\bigg(          L^\op_{\mathcal{B}^{k+t+m-3-i}_{m}}[\vec{x}]
 +            L^\op_{\mathcal{A}^{k+t+m-3-i}_{m}}[\vec{x}]\bigg)x\bigg) + \lambda \tbinom{m-2}{i}R\bigg(\bigg(          L^\op_{\mathcal{B}^{k+t+m-3-i}_{m}}[\vec{x}]
 +            L^\op_{\mathcal{A}^{k+t+m-3-i}_{m}}[\vec{x}]\bigg)x\bigg)\\
& + \lambda \tbinom{m-2}{i}\bigg( R\bigg(          L^\op_{\mathcal{B}^{k+t+m-3-i}_{m}}[\vec{x}]
 +            L^\op_{\mathcal{A}^{k+t+m-3-i}_{m}}[\vec{x}]\bigg)x\bigg)
 + \lambda^2 \tbinom{m-3}{i} \bigg( \bigg(          L^\op_{\mathcal{B}^{k+t+m-3-i}_{m}}[\vec{x}]\bigg)x\bigg) \\
& + \lambda^2 \tbinom{m-2}{i}\bigg( \bigg(            L^\op_{\mathcal{C}^{k+t+m-3-i}_{m}}[\vec{x}]\bigg)x\bigg)
 \left. +\lambda^3 \tbinom{m-3}{i}\bigg( \bigg(              L^\op_{\mathcal{D}^{k+t+m-3-i}_{m}}[\vec{x}]\bigg)x\bigg) \right)]\bigg)\bigg) 
  \\
 &+ \lambda \tbinom{m-1}{i} R^2\bigg(x\bigg(  L_{\mathcal{B}^{k+t+m-2-i}_{m+1}}[\vec{x},\sum\limits^l_{t=0}\sum\limits^{m-2}_{i=0}\lambda^{(t)}\lambda^i \left(
\tbinom{m-2}{i}R\bigg(  {L^*}^\op_{\mathcal{B}^{k+t+m-3-i}_m}[x,\vec{x}]
+   {L^*}^\op_{\mathcal{A}^{k+t+m-3-i}_m}[x,\vec{x}]\bigg) \right. \\
&+\lambda \tbinom{m-3}{i}\bigg(  {L^*}^\op_{\mathcal{B}^{k+t+m-3-i}_m}[x,\vec{x}]\bigg)
+\lambda \tbinom{m-2}{i} \bigg(    {L^*}^\op_{\mathcal{C}^{k+t+m-3-i}_m}[x,\vec{x}]\bigg)\bigg)+\lambda^2 \tbinom{m-3}{i}\bigg(   {L^*}^\op_{\mathcal{D}^{k+t+m-4-i}_m}[x,\vec{x}]\bigg)\\
&
 + \tbinom{m-2}{i}R\bigg(R\bigg(          L^\op_{\mathcal{B}^{k+t+m-3-i}_{m}}[\vec{x}]
 +            L^\op_{\mathcal{A}^{k+t+m-3-i}_{m}}[\vec{x}]\bigg)x\bigg)
 + \lambda \tbinom{m-2}{i}R\bigg(\bigg(          L^\op_{\mathcal{B}^{k+t+m-3-i}_{m}}[\vec{x}]
 +            L^\op_{\mathcal{A}^{k+t+m-3-i}_{m}}[\vec{x}]\bigg)x\bigg)\\
& + \lambda \tbinom{m-2}{i}\bigg( R\bigg(          L^\op_{\mathcal{B}^{k+t+m-3-i}_{m}}[\vec{x}]
 +            L^\op_{\mathcal{A}^{k+t+m-3-i}_{m}}[\vec{x}]\bigg)x\bigg)
 + \lambda^2 \tbinom{m-3}{i} \bigg( \bigg(          L^\op_{\mathcal{B}^{k+t+m-3-i}_{m}}[\vec{x}]\bigg)x\bigg) \\
& + \lambda^2 \tbinom{m-2}{i}\bigg( \bigg(            L^\op_{\mathcal{C}^{k+t+m-3-i}_{m}}[\vec{x}]\bigg)x\bigg)
 \left. +\lambda^3 \tbinom{m-3}{i}\bigg( \bigg(              L^\op_{\mathcal{D}^{k+t+m-3-i}_{m}}[\vec{x}]\bigg)x\bigg) \right)]
\\&+        L_{\mathcal{A}^{k+t+m-2-i}_{m+1}}[\vec{x},\sum\limits^l_{t=0}\sum\limits^{m-2}_{i=0}\lambda^{(t)}\lambda^i \left(
\tbinom{m-2}{i}R\bigg(  {L^*}^\op_{\mathcal{B}^{k+t+m-3-i}_m}[x,\vec{x}]
+   {L^*}^\op_{\mathcal{A}^{k+t+m-3-i}_m}[x,\vec{x}]\bigg) \right. \\
&+\lambda \tbinom{m-3}{i}\bigg(  {L^*}^\op_{\mathcal{B}^{k+t+m-3-i}_m}[x,\vec{x}]\bigg)
+\lambda \tbinom{m-2}{i} \bigg(    {L^*}^\op_{\mathcal{C}^{k+t+m-3-i}_m}[x,\vec{x}]\bigg)\bigg)+\lambda^2 \tbinom{m-3}{i}\bigg(   {L^*}^\op_{\mathcal{D}^{k+t+m-4-i}_m}[x,\vec{x}]\bigg)\\
&
 + \tbinom{m-2}{i}R\bigg(R\bigg(          L^\op_{\mathcal{B}^{k+t+m-3-i}_{m}}[\vec{x}]
 +            L^\op_{\mathcal{A}^{k+t+m-3-i}_{m}}[\vec{x}]\bigg)x\bigg)+ \lambda \tbinom{m-2}{i}R\bigg(\bigg(          L^\op_{\mathcal{B}^{k+t+m-3-i}_{m}}[\vec{x}]
 +            L^\op_{\mathcal{A}^{k+t+m-3-i}_{m}}[\vec{x}]\bigg)x\bigg)\\
& + \lambda \tbinom{m-2}{i}R\bigg(\bigg(          L^\op_{\mathcal{B}^{k+t+m-3-i}_{m}}[\vec{x}]
 +            L^\op_{\mathcal{A}^{k+t+m-3-i}_{m}}[\vec{x}]\bigg)x\bigg) + \lambda \tbinom{m-2}{i}\bigg( R\bigg(          L^\op_{\mathcal{B}^{k+t+m-3-i}_{m}}[\vec{x}]
 +            L^\op_{\mathcal{A}^{k+t+m-3-i}_{m}}[\vec{x}]\bigg)x\bigg)\\
 &+ \lambda^2 \tbinom{m-3}{i} \bigg( \bigg(          L^\op_{\mathcal{B}^{k+t+m-3-i}_{m}}[\vec{x}]\bigg)x\bigg)  + \lambda^2 \tbinom{m-2}{i}\bigg( \bigg(            L^\op_{\mathcal{C}^{k+t+m-3-i}_{m}}[\vec{x}]\bigg)x\bigg)
 \left. +\lambda^3 \tbinom{m-3}{i}\bigg( \bigg(              L^\op_{\mathcal{D}^{k+t+m-3-i}_{m}}[\vec{x}]\bigg)x\bigg) \right)]\bigg)\bigg) \allowdisplaybreaks
  \\
 &+ \lambda \tbinom{m-1}{i}R\bigg(x R\bigg(\sum\limits_{\mathcal{B}^{k+t+m-2-i}_{m}\in \mathcal{Q}} L_{\mathcal{B}^{k+t+m-2-i}_{m+1}}[\vec{x},\sum\limits^l_{t=0}\sum\limits^{m-2}_{i=0}\lambda^{(t)}\lambda^i \left(
\tbinom{m-2}{i}R\bigg(  {L^*}^\op_{\mathcal{B}^{k+t+m-3-i}_m}[x,\vec{x}]
+   {L^*}^\op_{\mathcal{A}^{k+t+m-3-i}_m}[x,\vec{x}]\bigg) \right. \\
&+\lambda \tbinom{m-3}{i}\bigg(  {L^*}^\op_{\mathcal{B}^{k+t+m-3-i}_m}[x,\vec{x}]\bigg)
+\lambda \tbinom{m-2}{i} \bigg(    {L^*}^\op_{\mathcal{C}^{k+t+m-3-i}_m}[x,\vec{x}]\bigg)\bigg)
+\lambda^2 \tbinom{m-3}{i}\bigg(   {L^*}^\op_{\mathcal{D}^{k+t+m-4-i}_m}[x,\vec{x}]\bigg)\\
 &+ \tbinom{m-2}{i}R\bigg(R\bigg(          L^\op_{\mathcal{B}^{k+t+m-3-i}_{m}}[\vec{x}]
 +            L^\op_{\mathcal{A}^{k+t+m-3-i}_{m}}[\vec{x}]\bigg)x\bigg) + \lambda \tbinom{m-2}{i}\bigg( R\bigg(          L^\op_{\mathcal{B}^{k+t+m-3-i}_{m}}[\vec{x}]
 +            L^\op_{\mathcal{A}^{k+t+m-3-i}_{m}}[\vec{x}]\bigg)x\bigg)\\
 & 
 + \lambda^2 \tbinom{m-3}{i} \bigg( \bigg(          L^\op_{\mathcal{B}^{k+t+m-3-i}_{m}}[\vec{x}]\bigg)x\bigg)  + \lambda^2 \tbinom{m-2}{i}\bigg( \bigg(            L^\op_{\mathcal{C}^{k+t+m-3-i}_{m}}[\vec{x}]\bigg)x\bigg)
 \left. +\lambda^3 \tbinom{m-3}{i}\bigg( \bigg(              L^\op_{\mathcal{D}^{k+t+m-3-i}_{m}}[\vec{x}]\bigg)x\bigg) \right)]
 \\
& +   L_{\mathcal{A}^{k+t+m-2-i}_{m+1}}[\vec{x},\sum\limits^l_{t=0}\sum\limits^{m-2}_{i=0}\lambda^{(t)}\lambda^i \left(
\tbinom{m-2}{i}R\bigg(  {L^*}^\op_{\mathcal{B}^{k+t+m-3-i}_m}[x,\vec{x}]
+   {L^*}^\op_{\mathcal{A}^{k+t+m-3-i}_m}[x,\vec{x}]\bigg) \right. +\lambda \tbinom{m-3}{i}\bigg(  {L^*}^\op_{\mathcal{B}^{k+t+m-3-i}_m}[x,\vec{x}]\bigg)\\
&+\lambda \tbinom{m-2}{i} \bigg(    {L^*}^\op_{\mathcal{C}^{k+t+m-3-i}_m}[x,\vec{x}]\bigg)\bigg)+\lambda^2 \tbinom{m-3}{i}\bigg(   {L^*}^\op_{\mathcal{D}^{k+t+m-4-i}_m}[x,\vec{x}]\bigg)
 + \tbinom{m-2}{i}R\bigg(R\bigg(          L^\op_{\mathcal{B}^{k+t+m-3-i}_{m}}[\vec{x}]
 +            L^\op_{\mathcal{A}^{k+t+m-3-i}_{m}}[\vec{x}]\bigg)x\bigg)\\
& + \lambda \tbinom{m-2}{i}R\bigg(\bigg(          L^\op_{\mathcal{B}^{k+t+m-3-i}_{m}}[\vec{x}]
 +            L^\op_{\mathcal{A}^{k+t+m-3-i}_{m}}[\vec{x}]\bigg)x\bigg) + \lambda \tbinom{m-2}{i}\bigg( R\bigg(          L^\op_{\mathcal{B}^{k+t+m-3-i}_{m}}[\vec{x}]
 +            L^\op_{\mathcal{A}^{k+t+m-3-i}_{m}}[\vec{x}]\bigg)x\bigg)\\
 \end{aligned}
\]
\[
\begin{aligned}
 &+ \lambda^2 \tbinom{m-3}{i} \bigg( \bigg(          L^\op_{\mathcal{B}^{k+t+m-3-i}_{m}}[\vec{x}]\bigg)x\bigg) 
 + \lambda^2 \tbinom{m-2}{i}\bigg( \bigg(            L^\op_{\mathcal{C}^{k+t+m-3-i}_{m}}[\vec{x}]\bigg)x\bigg)
 \left. +\lambda^3 \tbinom{m-3}{i}\bigg( \bigg(              L^\op_{\mathcal{D}^{k+t+m-3-i}_{m}}[\vec{x}]\bigg)x\bigg) \right)]\bigg)\bigg)
 \\
 &+ \lambda^2 \tbinom{m-2}{i} R\bigg(x \bigg(  L_{\mathcal{B}^{k+t+m-2-i}_{m+1}}[\vec{x},\sum\limits^l_{t=0}\sum\limits^{m-2}_{i=0}\lambda^{(t)}\lambda^i \left(
\tbinom{m-2}{i}R\bigg(  {L^*}^\op_{\mathcal{B}^{k+t+m-3-i}_m}[x,\vec{x}]
+   {L^*}^\op_{\mathcal{A}^{k+t+m-3-i}_m}[x,\vec{x}]\bigg) \right. \\
&+\lambda \tbinom{m-3}{i}\bigg(  {L^*}^\op_{\mathcal{B}^{k+t+m-3-i}_m}[x,\vec{x}]\bigg)
+\lambda \tbinom{m-2}{i} \bigg(    {L^*}^\op_{\mathcal{C}^{k+t+m-3-i}_m}[x,\vec{x}]\bigg)\bigg)+\lambda^2 \tbinom{m-3}{i}\bigg(   {L^*}^\op_{\mathcal{D}^{k+t+m-4-i}_m}[x,\vec{x}]\bigg)\\
 &+ \tbinom{m-2}{i}R\bigg(R\bigg(          L^\op_{\mathcal{B}^{k+t+m-3-i}_{m}}[\vec{x}]
 +            L^\op_{\mathcal{A}^{k+t+m-3-i}_{m}}[\vec{x}]\bigg)x\bigg) + \lambda \tbinom{m-2}{i}R\bigg(\bigg(          L^\op_{\mathcal{B}^{k+t+m-3-i}_{m}}[\vec{x}]
 +            L^\op_{\mathcal{A}^{k+t+m-3-i}_{m}}[\vec{x}]\bigg)x\bigg)\\
& + \lambda \tbinom{m-2}{i}\bigg( R\bigg(          L^\op_{\mathcal{B}^{k+t+m-3-i}_{m}}[\vec{x}]
 +            L^\op_{\mathcal{A}^{k+t+m-3-i}_{m}}[\vec{x}]\bigg)x\bigg)
 + \lambda^2 \tbinom{m-3}{i} \bigg( \bigg(          L^\op_{\mathcal{B}^{k+t+m-3-i}_{m}}[\vec{x}]\bigg)x\bigg) \\
& + \lambda^2 \tbinom{m-2}{i}\bigg( \bigg(            L^\op_{\mathcal{C}^{k+t+m-3-i}_{m}}[\vec{x}]\bigg)x\bigg)
 \left. +\lambda^3 \tbinom{m-3}{i}\bigg( \bigg(              L^\op_{\mathcal{D}^{k+t+m-3-i}_{m}}[\vec{x}]\bigg)x\bigg) \right)]\bigg)\bigg)\bigg)
 \\
 &+ \lambda^2 \tbinom{m-1}{i} R\bigg(x \bigg(  L_{\mathcal{C}^{k+t+m-2-i}_{m+1}}[\vec{x},\sum\limits^l_{t=0}\sum\limits^{m-2}_{i=0}\lambda^{(t)}\lambda^i \left(
\tbinom{m-2}{i}R\bigg(  {L^*}^\op_{\mathcal{B}^{k+t+m-3-i}_m}[x,\vec{x}]
+   {L^*}^\op_{\mathcal{A}^{k+t+m-3-i}_m}[x,\vec{x}]\bigg)
\right.\\
&+\lambda \tbinom{m-3}{i}\bigg(  {L^*}^\op_{\mathcal{B}^{k+t+m-3-i}_m}[x,\vec{x}]\bigg)+\lambda \tbinom{m-2}{i} \bigg(    {L^*}^\op_{\mathcal{C}^{k+t+m-3-i}_m}[x,\vec{x}]\bigg)\bigg)+\lambda^2 \tbinom{m-3}{i}\bigg(   {L^*}^\op_{\mathcal{D}^{k+t+m-4-i}_m}[x,\vec{x}]\bigg)\\
&
 + \tbinom{m-2}{i}R\bigg(R\bigg(          L^\op_{\mathcal{B}^{k+t+m-3-i}_{m}}[\vec{x}]
 +            L^\op_{\mathcal{A}^{k+t+m-3-i}_{m}}[\vec{x}]\bigg)x\bigg)+ \lambda \tbinom{m-2}{i}R\bigg(\bigg(          L^\op_{\mathcal{B}^{k+t+m-3-i}_{m}}[\vec{x}]
 +            L^\op_{\mathcal{A}^{k+t+m-3-i}_{m}}[\vec{x}]\bigg)x\bigg)\\
 &+ \lambda \tbinom{m-2}{i}\bigg( R\bigg(          L^\op_{\mathcal{B}^{k+t+m-3-i}_{m}}[\vec{x}]
 +            L^\op_{\mathcal{A}^{k+t+m-3-i}_{m}}[\vec{x}]\bigg)x\bigg)
 + \lambda^2 \tbinom{m-3}{i} \bigg( \bigg(          L^\op_{\mathcal{B}^{k+t+m-3-i}_{m}}[\vec{x}]\bigg)x\bigg)  + \lambda^2 \tbinom{m-2}{i}\bigg( \bigg(            L^\op_{\mathcal{C}^{k+t+m-3-i}_{m}}[\vec{x}]\bigg)x\bigg)\\
 &
 +\lambda^3 \tbinom{m-2}{i}R\bigg(x \bigg(  L_{\mathcal{D}^{k+t+m-2-i}_{m+1}}[\vec{x},\sum\limits^l_{t=0}\sum\limits^{m-2}_{i=0}\lambda^{(t)}\lambda^i \left(
\tbinom{m-2}{i}R\bigg(  {L^*}^\op_{\mathcal{B}^{k+t+m-3-i}_m}[x,\vec{x}]
+   {L^*}^\op_{\mathcal{A}^{k+t+m-3-i}_m}[x,\vec{x}]\bigg) \right. \\
&\left. +\lambda^3 \tbinom{m-3}{i}\bigg( \bigg(              L^\op_{\mathcal{D}^{k+t+m-3-i}_{m}}[\vec{x}]\bigg)x\bigg) \right)]\bigg)\bigg)\bigg)+\lambda \tbinom{m-3}{i}\bigg(  {L^*}^\op_{\mathcal{B}^{k+t+m-3-i}_m}[x,\vec{x}]\bigg)
+\lambda \tbinom{m-2}{i} \bigg(    {L^*}^\op_{\mathcal{C}^{k+t+m-3-i}_m}[x,\vec{x}]\bigg)\bigg)\\
&+\lambda^2 \tbinom{m-3}{i}\bigg(   {L^*}^\op_{\mathcal{D}^{k+t+m-4-i}_m}[x,\vec{x}]\bigg)
 + \tbinom{m-2}{i}R\bigg(R\bigg(          L^\op_{\mathcal{B}^{k+t+m-3-i}_{m}}[\vec{x}]
 +            L^\op_{\mathcal{A}^{k+t+m-3-i}_{m}}[\vec{x}]\bigg)x\bigg)\\
& + \lambda \tbinom{m-2}{i}R\bigg(\bigg(          L^\op_{\mathcal{B}^{k+t+m-3-i}_{m}}[\vec{x}]
 +            L^\op_{\mathcal{A}^{k+t+m-3-i}_{m}}[\vec{x}]\bigg)x\bigg)+ \lambda \tbinom{m-2}{i}\bigg( R\bigg(          L^\op_{\mathcal{B}^{k+t+m-3-i}_{m}}[\vec{x}]
 +            L^\op_{\mathcal{A}^{k+t+m-3-i}_{m}}[\vec{x}]\bigg)x\bigg)\\
 &+ \lambda^2 \tbinom{m-3}{i} \bigg( \bigg(          L^\op_{\mathcal{B}^{k+t+m-3-i}_{m}}[\vec{x}]\bigg)x\bigg) 
 + \lambda^2 \tbinom{m-2}{i}\bigg( \bigg(            L^\op_{\mathcal{C}^{k+t+m-3-i}_{m}}[\vec{x}]\bigg)x\bigg)
 \left. +\lambda^3 \tbinom{m-3}{i}\bigg( \bigg(              L^\op_{\mathcal{D}^{k+t+m-3-i}_{m}}[\vec{x}]\bigg)x\bigg) \right)]\bigg)\bigg)\bigg) {=} 0;\ m {=} |\vec{x}|{\geq}3. 
\end{aligned}
\]
Composition between relations $(R9)$ and $(R4)$-$(R10)$, we get a sum of the expressions, which are consequences of~$(R10)$.

Composition between relations $(R10)$ and $(R4)$, we get a sum of the expressions, which are consequences of~$(R10)$ itself.
Below, we show this example for  $\sum\limits_{t=0}^l \sum\limits^{m-1}_{i=0}\lambda^{(t)}\lambda^i
M_{t,i}$, where
\[
\begin{aligned}
M_{t,i}&=
\tbinom{m-1}{i}R^2\bigg( {L^*}^\op_{\mathcal{B}^{k+t+m-2-i}_{m+1}}[x,\vec{x},R\bigg(\sum\limits^{m-2}_{i=0}\tbinom{m-2}{i}\lambda^i\bigg(  L_{\mathcal{B}^{k+t+m-3-i}_m}[\vec{x}]+   L_{\mathcal{A}^{k+t+m-3-i}_m}[\vec{x}]\bigg)\bigg)\\ 
&+\lambda (\sum\limits^{m-3}_{i=0}\tbinom{m-3}{i}\lambda^i\bigg(  L_{\mathcal{B}^{k+t+m-3-i}_m}[\vec{x}]\bigg)\bigg) 
 +\lambda \bigg(\sum\limits^{m-2}_{i=0}\tbinom{m-2}{i}\lambda^i\bigg(    L_{\mathcal{C}^{k+t+m-3-i}_m}[\vec{x}]\bigg)\bigg)+\lambda^2 \bigg(\sum\limits^{m-3}_{i=0}\tbinom{m-3}{i}\lambda^i\bigg(   L_{\mathcal{D}^{k+t+m-4-i}_m}[\vec{x}]\bigg)\bigg)]
\\
&+ {L^*}^\op_{\mathcal{A}^{k+t+m-2-i}_{m+1}}[x,\vec{x},R\bigg(\sum\limits^{m-2}_{i=0}\tbinom{m-2}{i}\lambda^i\bigg(  L_{\mathcal{B}^{k+t+m-3-i}_m}[\vec{x}]+   L_{\mathcal{A}^{k+t+m-3-i}_m}[\vec{x}]\bigg)\bigg) 
\end{aligned}
\]
\[
\begin{aligned}
& +\lambda (\sum\limits^{m-3}_{i=0}\tbinom{m-3}{i}\lambda^i\bigg(  L_{\mathcal{B}^{k+t+m-3-i}_m}[\vec{x}]\bigg)\bigg) 
 +\lambda \bigg(\sum\limits^{m-2}_{i=0}\tbinom{m-2}{i}\lambda^i\bigg(    L_{\mathcal{C}^{k+t+m-3-i}_m}[\vec{x}]\bigg)\bigg)  +\lambda^2 \bigg(\sum\limits^{m-3}_{i=0}\tbinom{m-3}{i}\lambda^i\bigg(   L_{\mathcal{D}^{k+t+m-4-i}_m}[\vec{x}]\bigg)\bigg)]\bigg) 
\\
&+\lambda \tbinom{m-2}{i}R\bigg( {L^*}^\op_{\mathcal{B}^{k+t+m-2-i}_{m+1}}[x,\vec{x},R\bigg(\sum\limits^{m-2}_{i=0}\tbinom{m-2}{i}\lambda^i\bigg(  L_{\mathcal{B}^{k+t+m-3-i}_m}[\vec{x}]+   L_{\mathcal{A}^{k+t+m-3-i}_m}[\vec{x}]\bigg)\bigg) \\
& +\lambda (\sum\limits^{m-3}_{i=0}\tbinom{m-3}{i}\lambda^i\bigg(  L_{\mathcal{B}^{k+t+m-3-i}_m}[\vec{x}]\bigg)\bigg) 
 +\lambda \bigg(\sum\limits^{m-2}_{i=0}\tbinom{m-2}{i}\lambda^i\bigg(    L_{\mathcal{C}^{k+t+m-3-i}_m}[\vec{x}]\bigg)\bigg)   +\lambda^2 \bigg(\sum\limits^{m-3}_{i=0}\tbinom{m-3}{i}\lambda^i\bigg(   L_{\mathcal{D}^{k+t+m-4-i}_m}[\vec{x}]\bigg)\bigg)]\bigg)
\\
&+\lambda \tbinom{m-1}{i} R\bigg( {L^*}^\op_{\mathcal{C}^{k+t+m-2-i}_{m+1}}[x,\vec{x},R\bigg(\sum\limits^{m-2}_{i=0}\tbinom{m-2}{i}\lambda^i\bigg(  L_{\mathcal{B}^{k+t+m-3-i}_m}[\vec{x}]+   L_{\mathcal{A}^{k+t+m-3-i}_m}[\vec{x}]\bigg)\bigg) \\
& +\lambda (\sum\limits^{m-3}_{i=0}\tbinom{m-3}{i}\lambda^i\bigg(  L_{\mathcal{B}^{k+t+m-3-i}_m}[\vec{x}]\bigg)\bigg) 
 +\lambda \bigg(\sum\limits^{m-2}_{i=0}\tbinom{m-2}{i}\lambda^i\bigg(    L_{\mathcal{C}^{k+t+m-3-i}_m}[\vec{x}]\bigg)\bigg)  +\lambda^2 \bigg(\sum\limits^{m-3}_{i=0}\tbinom{m-3}{i}\lambda^i\bigg(   L_{\mathcal{D}^{k+t+m-4-i}_m}[\vec{x}]\bigg)\bigg)]\bigg)\bigg)
 \\
&+\lambda^2 \tbinom{m-2}{i}R\bigg(  {L^*}^\op_{\mathcal{D}^{k+t+m-3-i}_{m+1}}[x,\vec{x},R\bigg(\sum\limits^{m-2}_{i=0}\tbinom{m-2}{i}\lambda^i\bigg(  L_{\mathcal{B}^{k+t+m-3-i}_m}[\vec{x}]+   L_{\mathcal{A}^{k+t+m-3-i}_m}[\vec{x}]\bigg)\bigg) \\
& +\lambda (\sum\limits^{m-3}_{i=0}\tbinom{m-3}{i}\lambda^i\bigg(  L_{\mathcal{B}^{k+t+m-3-i}_m}[\vec{x}]\bigg)\bigg) 
 +\lambda \bigg(\sum\limits^{m-2}_{i=0}\tbinom{m-2}{i}\lambda^i\bigg(    L_{\mathcal{C}^{k+t+m-3-i}_m}[\vec{x}]\bigg)\bigg)   +\lambda^2 \bigg(\sum\limits^{m-3}_{i=0}\tbinom{m-3}{i}\lambda^i\bigg(   L_{\mathcal{D}^{k+t+m-4-i}_m}[\vec{x}]\bigg)\bigg)]\bigg)
\\
 &+ \tbinom{m-1}{i}R^2\bigg(R\bigg(  L^\op_{\mathcal{B}^{k+t+m-2-i}_{m+1}}[\vec{x},R\bigg(\sum\limits^{m-2}_{i=0}\tbinom{m-2}{i}\lambda^i\bigg(  L_{\mathcal{B}^{k+t+m-3-i}_m}[\vec{x}]+   L_{\mathcal{A}^{k+t+m-3-i}_m}[\vec{x}]\bigg)\bigg) 
 \\
& +\lambda (\sum\limits^{m-3}_{i=0}\tbinom{m-3}{i}\lambda^i\bigg(  L_{\mathcal{B}^{k+t+m-3-i}_m}[\vec{x}]\bigg)\bigg) 
 +\lambda \bigg(\sum\limits^{m-2}_{i=0}\tbinom{m-2}{i}\lambda^i\bigg(    L_{\mathcal{C}^{k+t+m-3-i}_m}[\vec{x}]\bigg)\bigg) 
 +\lambda^2 \bigg(\sum\limits^{m-3}_{i=0}\tbinom{m-3}{i}\lambda^i\bigg(   L_{\mathcal{D}^{k+t+m-4-i}_m}[\vec{x}]\bigg)\bigg)]
 \\
& +        L^\op_{\mathcal{A}^{k+t+m-2-i}_{m+1}}[\vec{x},R\bigg(\sum\limits^{m-2}_{i=0}\tbinom{m-2}{i}\lambda^i\bigg(  L_{\mathcal{B}^{k+t+m-3-i}_m}[\vec{x}]+   L_{\mathcal{A}^{k+t+m-3-i}_m}[\vec{x}]\bigg)\bigg)  +\lambda (\sum\limits^{m-3}_{i=0}\tbinom{m-3}{i}\lambda^i\bigg(  L_{\mathcal{B}^{k+t+m-3-i}_m}[\vec{x}]\bigg)\bigg) \\
 &+\lambda \bigg(\sum\limits^{m-2}_{i=0}\tbinom{m-2}{i}\lambda^i\bigg(    L_{\mathcal{C}^{k+t+m-3-i}_m}[\vec{x}]\bigg)\bigg)   +\lambda^2 \bigg(\sum\limits^{m-3}_{i=0}\tbinom{m-3}{i}\lambda^i\bigg(   L_{\mathcal{D}^{k+t+m-4-i}_m}[\vec{x}]\bigg)\bigg)]\bigg)x\bigg)
\\
 &+ \lambda \tbinom{m-1}{i}R^2\bigg(\bigg(  L^\op_{\mathcal{B}^{k+t+m-2-i}_{m+1}}[\vec{x},R\bigg(\sum\limits^{m-2}_{i=0}\tbinom{m-2}{i}\lambda^i\bigg(  L_{\mathcal{B}^{k+t+m-3-i}_m}[\vec{x}]+   L_{\mathcal{A}^{k+t+m-3-i}_m}[\vec{x}]\bigg)\bigg) \\
& +\lambda (\sum\limits^{m-3}_{i=0}\tbinom{m-3}{i}\lambda^i\bigg(  L_{\mathcal{B}^{k+t+m-3-i}_m}[\vec{x}]\bigg)\bigg) 
 +\lambda \bigg(\sum\limits^{m-2}_{i=0}\tbinom{m-2}{i}\lambda^i\bigg(    L_{\mathcal{C}^{k+t+m-3-i}_m}[\vec{x}]\bigg)\bigg)  +\lambda^2 \bigg(\sum\limits^{m-3}_{i=0}\tbinom{m-3}{i}\lambda^i\bigg(   L_{\mathcal{D}^{k+t+m-4-i}_m}[\vec{x}]\bigg)\bigg)]
 \\
& +        L^\op_{\mathcal{A}^{k+t+m-2-i}_{m+1}}[\vec{x},R\bigg(\sum\limits^{m-2}_{i=0}\tbinom{m-2}{i}\lambda^i\bigg(  L_{\mathcal{B}^{k+t+m-3-i}_m}[\vec{x}]+   L_{\mathcal{A}^{k+t+m-3-i}_m}[\vec{x}]\bigg)\bigg) +\lambda (\sum\limits^{m-3}_{i=0}\tbinom{m-3}{i}\lambda^i\bigg(  L_{\mathcal{B}^{k+t+m-3-i}_m}[\vec{x}]\bigg)\bigg)\\ 
& +\lambda \bigg(\sum\limits^{m-2}_{i=0}\tbinom{m-2}{i}\lambda^i\bigg(    L_{\mathcal{C}^{k+t+m-3-i}_m}[\vec{x}]\bigg)\bigg)   +\lambda^2 \bigg(\sum\limits^{m-3}_{i=0}\tbinom{m-3}{i}\lambda^i\bigg(   L_{\mathcal{D}^{k+t+m-4-i}_m}[\vec{x}]\bigg)\bigg)]\bigg)x\bigg)
\\
& + \lambda \tbinom{m-1}{i}R\bigg( R\bigg(  L^\op_{\mathcal{B}^{k+t+m-2-i}_{m+1}}[\vec{x},R\bigg(\sum\limits^{m-2}_{i=0}\tbinom{m-2}{i}\lambda^i\bigg(  L_{\mathcal{B}^{k+t+m-3-i}_m}[\vec{x}]+   L_{\mathcal{A}^{k+t+m-3-i}_m}[\vec{x}]\bigg)\bigg) \\
& +\lambda (\sum\limits^{m-3}_{i=0}\tbinom{m-3}{i}\lambda^i\bigg(  L_{\mathcal{B}^{k+t+m-3-i}_m}[\vec{x}]\bigg)\bigg) 
 +\lambda \bigg(\sum\limits^{m-2}_{i=0}\tbinom{m-2}{i}\lambda^i\bigg(    L_{\mathcal{C}^{k+t+m-3-i}_m}[\vec{x}]\bigg)\bigg)  +\lambda^2 \bigg(\sum\limits^{m-3}_{i=0}\tbinom{m-3}{i}\lambda^i\bigg(   L_{\mathcal{D}^{k+t+m-4-i}_m}[\vec{x}]\bigg)\bigg)]
\end{aligned}
\]
\[
\begin{aligned}
& +        L^\op_{\mathcal{A}^{k+t+m-2-i}_{m+1}}[\vec{x},R\bigg(\sum\limits^{m-2}_{i=0}\tbinom{m-2}{i}\lambda^i\bigg(  L_{\mathcal{B}^{k+t+m-3-i}_m}[\vec{x}]+   L_{\mathcal{A}^{k+t+m-3-i}_m}[\vec{x}]\bigg)\bigg)  +\lambda (\sum\limits^{m-3}_{i=0}\tbinom{m-3}{i}\lambda^i\bigg(  L_{\mathcal{B}^{k+t+m-3-i}_m}[\vec{x}]\bigg)\bigg)\\ 
 &+\lambda \bigg(\sum\limits^{m-2}_{i=0}\tbinom{m-2}{i}\lambda^i\bigg(    L_{\mathcal{C}^{k+t+m-3-i}_m}[\vec{x}]\bigg)\bigg)   +\lambda^2 \bigg(\sum\limits^{m-3}_{i=0}\tbinom{m-3}{i}\lambda^i\bigg(   L_{\mathcal{D}^{k+t+m-4-i}_m}[\vec{x}]\bigg)\bigg)]\bigg)x\bigg)
\\
 &+ \lambda^2 \tbinom{m-2}{i} R\bigg( \bigg(  L^\op_{\mathcal{B}^{k+t+m-2-i}_{m+1}}[\vec{x},R\bigg(\sum\limits^{m-2}_{i=0}\tbinom{m-2}{i}\lambda^i\bigg(  L_{\mathcal{B}^{k+t+m-3-i}_m}[\vec{x}]+   L_{\mathcal{A}^{k+t+m-3-i}_m}[\vec{x}]\bigg)\bigg) \\
& +\lambda (\sum\limits^{m-3}_{i=0}\tbinom{m-3}{i}\lambda^i\bigg(  L_{\mathcal{B}^{k+t+m-3-i}_m}[\vec{x}]\bigg)\bigg) 
 +\lambda \bigg(\sum\limits^{m-2}_{i=0}\tbinom{m-2}{i}\lambda^i\bigg(    L_{\mathcal{C}^{k+t+m-3-i}_m}[\vec{x}]\bigg)\bigg)  +\lambda^2 \bigg(\sum\limits^{m-3}_{i=0}\tbinom{m-3}{i}\lambda^i\bigg(   L_{\mathcal{D}^{k+t+m-4-i}_m}[\vec{x}]\bigg)\bigg)]\bigg)x\bigg) 
\\
& + \lambda^2 \tbinom{m-1}{i}R\bigg( \bigg(  L^\op_{\mathcal{C}^{k+t+m-2-i}_{m+1}}[\vec{x},R\bigg(\sum\limits^{m-2}_{i=0}\tbinom{m-2}{i}\lambda^i\bigg(  L_{\mathcal{B}^{k+t+m-3-i}_m}[\vec{x}]+   L_{\mathcal{A}^{k+t+m-3-i}_m}[\vec{x}]\bigg)\bigg) 
 +\lambda (\sum\limits^{m-3}_{i=0}\tbinom{m-3}{i}\lambda^i\bigg(  L_{\mathcal{B}^{k+t+m-3-i}_m}[\vec{x}]\bigg)\bigg)\\ 
 &+\lambda \bigg(\sum\limits^{m-2}_{i=0}\tbinom{m-2}{i}\lambda^i\bigg(    L_{\mathcal{C}^{k+t+m-3-i}_m}[\vec{x}]\bigg)\bigg)  +\lambda^2 \bigg(\sum\limits^{m-3}_{i=0}\tbinom{m-3}{i}\lambda^i\bigg(   L_{\mathcal{D}^{k+t+m-4-i}_m}[\vec{x}]\bigg)\bigg)]\bigg)x\bigg)
\\
  &+\lambda^3 \tbinom{m-2}{i}R\bigg( \bigg(  L^\op_{\mathcal{D}^{k+t+m-2-i}_{m+1}}[\vec{x},R\bigg(\sum\limits^{m-2}_{i=0}\tbinom{m-2}{i}\lambda^i\bigg(  L_{\mathcal{B}^{k+t+m-3-i}_m}[\vec{x}]+   L_{\mathcal{A}^{k+t+m-3-i}_m}[\vec{x}]\bigg)\bigg)  +\lambda (\sum\limits^{m-3}_{i=0}\tbinom{m-3}{i}\lambda^i\bigg(  L_{\mathcal{B}^{k+t+m-3-i}_m}[\vec{x}]\bigg)\bigg)\\
 &+\lambda \bigg(\sum\limits^{m-2}_{i=0}\tbinom{m-2}{i}\lambda^i\bigg(    L_{\mathcal{C}^{k+t+m-3-i}_m}[\vec{x}]\bigg)\bigg)   +\lambda^2 \bigg(\sum\limits^{m-3}_{i=0}\tbinom{m-3}{i}\lambda^i\bigg(   L_{\mathcal{D}^{k+t+m-4-i}_m}[\vec{x}]\bigg)\bigg)]\bigg)x\bigg) ) {=} 0;\ m = |\vec{x}|{\geq}3.
\end{aligned}
\]
Composition between relations $(R10)$ and $(R5)$, we get a sum of the expressions, which are consequences of~$(R10)$ itself.
Below, we show this example for  $\sum\limits_{t=0}^l \sum\limits^{m-1}_{i=0}\lambda^{(t)}\lambda^i
M_{t,i}$, where
\[
\begin{aligned}
M_{t,i}&= 
\tbinom{m-1}{i}R^2\bigg( {L^*}^\op_{\mathcal{B}^{k+t+m-2-i}_{m+1}}[x,\vec{x},\sum\limits^l_{t=0}\lambda ^{(t)}\bigg(R\bigg(  L^\op_{\mathcal{B}^{k+t-1}_2}[\vec{y}]+              L^{\op}_{\mathcal{A}^{k+t-1}_2}[\vec{y}]\bigg)
+\lambda \bigg(  L^{\op}_{\mathcal{B}^{k+t-1}_2}[\vec{y}]+      L^\op_{\mathcal{C}^{k+t-1}_2}[\vec{y}]\bigg)\bigg)]\\
&+ {L^*}^\op_{\mathcal{A}^{k+t+m-2-i}_{m+1}}[x,\vec{x},\sum\limits^l_{t=0}\lambda ^{(t)}\bigg(R\bigg(  L^\op_{\mathcal{B}^{k+t-1}_2}[\vec{y}]+              L^{\op}_{\mathcal{A}^{k+t-1}_2}[\vec{y}]\bigg)+\lambda \bigg(  L^{\op}_{\mathcal{B}^{k+t-1}_2}[\vec{y}]+      L^\op_{\mathcal{C}^{k+t-1}_2}[\vec{y}]\bigg)\bigg)]\bigg) 
\\
&+\lambda \tbinom{m-2}{i}R\bigg( {L^*}^\op_{\mathcal{B}^{k+t+m-2-i}_{m+1}}[x,\vec{x},\sum\limits^l_{t=0}\lambda ^{(t)}\bigg(R\bigg(  L^\op_{\mathcal{B}^{k+t-1}_2}[\vec{y}]+              L^{\op}_{\mathcal{A}^{k+t-1}_2}[\vec{y}]\bigg)+\lambda \bigg(  L^{\op}_{\mathcal{B}^{k+t-1}_2}[\vec{y}]+      L^\op_{\mathcal{C}^{k+t-1}_2}[\vec{y}]\bigg)\bigg)]\bigg)
 \\
&+\lambda \tbinom{m-1}{i} R\bigg( {L^*}^\op_{\mathcal{C}^{k+t+m-2-i}_{m+1}}[x,\vec{x},\sum\limits^l_{t=0}\lambda ^{(t)}\bigg(R\bigg(  L^\op_{\mathcal{B}^{k+t-1}_2}[\vec{y}]+              L^{\op}_{\mathcal{A}^{k+t-1}_2}[\vec{y}]\bigg)
+\lambda \bigg(  L^{\op}_{\mathcal{B}^{k+t-1}_2}[\vec{y}]+      L^\op_{\mathcal{C}^{k+t-1}_2}[\vec{y}]\bigg)\bigg)]\bigg)\bigg)
\\
&+\lambda^2 \tbinom{m-2}{i}R\bigg(  {L^*}^\op_{\mathcal{D}^{k+t+m-3-i}_{m+1}}[x,\vec{x},\sum\limits^l_{t=0}\lambda ^{(t)}\bigg(R\bigg(  L^\op_{\mathcal{B}^{k+t-1}_2}[\vec{y}]+              L^{\op}_{\mathcal{A}^{k+t-1}_2}[\vec{y}]\bigg)
+\lambda \bigg(  L^{\op}_{\mathcal{B}^{k+t-1}_2}[\vec{y}]+      L^\op_{\mathcal{C}^{k+t-1}_2}[\vec{y}]\bigg)\bigg)]\bigg)
 \\
 &+ \tbinom{m-1}{i}R^2\bigg(R\bigg(  L^\op_{\mathcal{B}^{k+t+m-2-i}_{m+1}}[\vec{x},\sum\limits^l_{t=0}\lambda ^{(t)}\bigg(R\bigg(  L^\op_{\mathcal{B}^{k+t-1}_2}[\vec{y}]+              L^{\op}_{\mathcal{A}^{k+t-1}_2}[\vec{y}]\bigg)+\lambda \bigg(  L^{\op}_{\mathcal{B}^{k+t-1}_2}[\vec{y}]+      L^\op_{\mathcal{C}^{k+t-1}_2}[\vec{y}]\bigg)\bigg)]
 \\
& +        L^\op_{\mathcal{A}^{k+t+m-2-i}_{m+1}}[\vec{x},\sum\limits^l_{t=0}\lambda ^{(t)}\bigg(R\bigg(  L^\op_{\mathcal{B}^{k+t-1}_2}[\vec{y}]+              L^{\op}_{\mathcal{A}^{k+t-1}_2}[\vec{y}]\bigg)+\lambda \bigg(  L^{\op}_{\mathcal{B}^{k+t-1}_2}[\vec{y}]+      L^\op_{\mathcal{C}^{k+t-1}_2}[\vec{y}]\bigg)\bigg)]\bigg)x\bigg)
\\
\end{aligned}
 \]
 \[
 \begin{aligned}
 &+ \lambda \tbinom{m-1}{i}R^2\bigg(\bigg(  L^\op_{\mathcal{B}^{k+t+m-2-i}_{m+1}}[\vec{x},\sum\limits^l_{t=0}\lambda ^{(t)}\bigg(R\bigg(  L^\op_{\mathcal{B}^{k+t-1}_2}[\vec{y}]+              L^{\op}_{\mathcal{A}^{k+t-1}_2}[\vec{y}]\bigg)+\lambda \bigg(  L^{\op}_{\mathcal{B}^{k+t-1}_2}[\vec{y}]+      L^\op_{\mathcal{C}^{k+t-1}_2}[\vec{y}]\bigg)\bigg)]
 \\
& +        L^\op_{\mathcal{A}^{k+t+m-2-i}_{m+1}}[\vec{x},\sum\limits^l_{t=0}\lambda ^{(t)}\bigg(R\bigg(  L^\op_{\mathcal{B}^{k+t-1}_2}[\vec{y}]+              L^{\op}_{\mathcal{A}^{k+t-1}_2}[\vec{y}]\bigg)+\lambda \bigg(  L^{\op}_{\mathcal{B}^{k+t-1}_2}[\vec{y}]+      L^\op_{\mathcal{C}^{k+t-1}_2}[\vec{y}]\bigg)\bigg)]\bigg)x\bigg)
\\
& + \lambda \tbinom{m-1}{i}R\bigg( R\bigg(  L^\op_{\mathcal{B}^{k+t+m-2-i}_{m+1}}[\vec{x},\sum\limits^l_{t=0}\lambda ^{(t)}\bigg(R\bigg(  L^\op_{\mathcal{B}^{k+t-1}_2}[\vec{y}]+              L^{\op}_{\mathcal{A}^{k+t-1}_2}[\vec{y}]\bigg)+\lambda \bigg(  L^{\op}_{\mathcal{B}^{k+t-1}_2}[\vec{y}]+      L^\op_{\mathcal{C}^{k+t-1}_2}[\vec{y}]\bigg)\bigg)]
 \\
& +        L^\op_{\mathcal{A}^{k+t+m-2-i}_{m+1}}[\vec{x},\sum\limits^l_{t=0}\lambda ^{(t)}\bigg(R\bigg(  L^\op_{\mathcal{B}^{k+t-1}_2}[\vec{y}]+              L^{\op}_{\mathcal{A}^{k+t-1}_2}[\vec{y}]\bigg)+\lambda \bigg(  L^{\op}_{\mathcal{B}^{k+t-1}_2}[\vec{y}]+      L^\op_{\mathcal{C}^{k+t-1}_2}[\vec{y}]\bigg)\bigg)]\bigg)x\bigg)
 \\
 &+ \lambda^2 \tbinom{m-2}{i} R\bigg( \bigg(  L^\op_{\mathcal{B}^{k+t+m-2-i}_{m+1}}[\vec{x},\sum\limits^l_{t=0}\lambda ^{(t)}\bigg(R\bigg(  L^\op_{\mathcal{B}^{k+t-1}_2}[\vec{y}]+              L^{\op}_{\mathcal{A}^{k+t-1}_2}[\vec{y}]\bigg)+\lambda \bigg(  L^{\op}_{\mathcal{B}^{k+t-1}_2}[\vec{y}]+      L^\op_{\mathcal{C}^{k+t-1}_2}[\vec{y}]\bigg)\bigg)]\bigg)x\bigg) 
 \\
& + \lambda^2 \tbinom{m-1}{i}R\bigg( \bigg(  L^\op_{\mathcal{C}^{k+t+m-2-i}_{m+1}}[\vec{x},\sum\limits^l_{t=0}\lambda ^{(t)}\bigg(R\bigg(  L^\op_{\mathcal{B}^{k+t-1}_2}[\vec{y}]+              L^{\op}_{\mathcal{A}^{k+t-1}_2}[\vec{y}]\bigg)+\lambda \bigg(  L^{\op}_{\mathcal{B}^{k+t-1}_2}[\vec{y}]+      L^\op_{\mathcal{C}^{k+t-1}_2}[\vec{y}]\bigg)\bigg)]\bigg)x\bigg)
\\
  &+\lambda^3 \tbinom{m-2}{i}R\bigg( \bigg(  L^\op_{\mathcal{D}^{k+t+m-2-i}_{m+1}}[\vec{x},\sum\limits^l_{t=0}\lambda ^{(t)}\bigg(R\bigg(  L^\op_{\mathcal{B}^{k+t-1}_2}[\vec{y}]+              L^{\op}_{\mathcal{A}^{k+t-1}_2}[\vec{y}]\bigg)+\lambda \bigg(  L^{\op}_{\mathcal{B}^{k+t-1}_2}[\vec{y}]+      L^\op_{\mathcal{C}^{k+t-1}_2}[\vec{y}]\bigg)\bigg)]\bigg)x\bigg) ) {=} 0;\\&  m = |\vec{x}|{\geq}3,\  |\vec{y}|=2.
\end{aligned}
\]
Composition between relations $(R10)$ and $(R6)$, we get a sum of the expressions, which are consequences of~$(R10)$ itself.
Below, we show this example for  $\sum\limits_{t=0}^l \sum\limits^{m-1}_{i=0}\lambda^{(t)}\lambda^i
M_{t,i}$, where
\[
\begin{aligned}
M_{t,i}=&
\tbinom{m-1}{i}R^2\bigg( {L^*}^\op_{\mathcal{B}^{k+t+m-2-i}_{m+1}}[x,\vec{x},\sum\limits^l_{t=0}\sum\limits^{m-2}_{i=0}\lambda^{(t)}\lambda^i\left(
\tbinom{m-2}{i}R\bigg(  L^\op_{\mathcal{B}^{k+t+m-3-i}_m}[\vec{x}]
 +    L^\op_{\mathcal{A}^{k+t+m-3-i}_m}[\vec{x}]\bigg) \right. \\
 &+ \lambda \tbinom{m-3}{i} \bigg(  L^\op_{\mathcal{B}^{k+t+m-3-i}_m}[\vec{x}]\bigg) 
 + \lambda \tbinom{m-2}{i}\bigg(    L^\op_{\mathcal{C}^{k+t+m-3-i}_m}[\vec{x}]\bigg)
 \left. +\lambda^2 \tbinom{m-3}{i} \bigg(   L^\op_{\mathcal{D}^{k+t+m-4-i}_m}[\vec{x}]\bigg)\right)]
\\
&+ {L^*}^\op_{\mathcal{A}^{k+t+m-2-i}_{m+1}}[x,\vec{x},\sum\limits^l_{t=0}\sum\limits^{m-2}_{i=0}\lambda^{(t)}\lambda^i\left(
\tbinom{m-2}{i}R\bigg(  L^\op_{\mathcal{B}^{k+t+m-3-i}_m}[\vec{x}]
 +    L^\op_{\mathcal{A}^{k+t+m-3-i}_m}[\vec{x}]\bigg) \right. \\
 &+ \lambda \tbinom{m-3}{i} \bigg(  L^\op_{\mathcal{B}^{k+t+m-3-i}_m}[\vec{x}]\bigg) 
 + \lambda \tbinom{m-2}{i}\bigg(    L^\op_{\mathcal{C}^{k+t+m-3-i}_m}[\vec{x}]\bigg)
 \left. +\lambda^2 \tbinom{m-3}{i} \bigg(   L^\op_{\mathcal{D}^{k+t+m-4-i}_m}[\vec{x}]\bigg)\right)]\bigg) 
\\
&+\lambda \tbinom{m-2}{i}R\bigg( {L^*}^\op_{\mathcal{B}^{k+t+m-2-i}_{m+1}}[x,\vec{x},\sum\limits^l_{t=0}\sum\limits^{m-2}_{i=0}\lambda^{(t)}\lambda^i\left(
\tbinom{m-2}{i}R\bigg(  L^\op_{\mathcal{B}^{k+t+m-3-i}_m}[\vec{x}]
 +    L^\op_{\mathcal{A}^{k+t+m-3-i}_m}[\vec{x}]\bigg) \right. \\
 &+ \lambda \tbinom{m-3}{i} \bigg(  L^\op_{\mathcal{B}^{k+t+m-3-i}_m}[\vec{x}]\bigg) 
 + \lambda \tbinom{m-2}{i}\bigg(    L^\op_{\mathcal{C}^{k+t+m-3-i}_m}[\vec{x}]\bigg)
 \left. +\lambda^2 \tbinom{m-3}{i} \bigg(   L^\op_{\mathcal{D}^{k+t+m-4-i}_m}[\vec{x}]\bigg)\right)]\bigg)
\\
&+\lambda \tbinom{m-1}{i} R\bigg( {L^*}^\op_{\mathcal{C}^{k+t+m-2-i}_{m+1}}[x,\vec{x},\sum\limits^l_{t=0}\sum\limits^{m-2}_{i=0}\lambda^{(t)}\lambda^i\left(
\tbinom{m-2}{i}R\bigg(  L^\op_{\mathcal{B}^{k+t+m-3-i}_m}[\vec{x}]
 +    L^\op_{\mathcal{A}^{k+t+m-3-i}_m}[\vec{x}]\bigg) \right. \\
 &+ \lambda \tbinom{m-3}{i} \bigg(  L^\op_{\mathcal{B}^{k+t+m-3-i}_m}[\vec{x}]\bigg) 
 + \lambda \tbinom{m-2}{i}\bigg(    L^\op_{\mathcal{C}^{k+t+m-3-i}_m}[\vec{x}]\bigg)
 \left. +\lambda^2 \tbinom{m-3}{i} \bigg(   L^\op_{\mathcal{D}^{k+t+m-4-i}_m}[\vec{x}]\bigg)\right)]\bigg)\bigg)
\\
&+\lambda^2 \tbinom{m-2}{i}R\bigg(  {L^*}^\op_{\mathcal{D}^{k+t+m-3-i}_{m+1}}[x,\vec{x},\sum\limits^l_{t=0}\sum\limits^{m-2}_{i=0}\lambda^{(t)}\lambda^i\left(
\tbinom{m-2}{i}R\bigg(  L^\op_{\mathcal{B}^{k+t+m-3-i}_m}[\vec{x}]
 +    L^\op_{\mathcal{A}^{k+t+m-3-i}_m}[\vec{x}]\bigg) \right. \\
 \end{aligned}
 \]
 \[
 \begin{aligned}
 &+ \lambda \tbinom{m-3}{i} \bigg(  L^\op_{\mathcal{B}^{k+t+m-3-i}_m}[\vec{x}]\bigg) 
 + \lambda \tbinom{m-2}{i}\bigg(    L^\op_{\mathcal{C}^{k+t+m-3-i}_m}[\vec{x}]\bigg)
 \left. +\lambda^2 \tbinom{m-3}{i} \bigg(   L^\op_{\mathcal{D}^{k+t+m-4-i}_m}[\vec{x}]\bigg)\right)]\bigg)
\\
 &+ \tbinom{m-1}{i}R^2\bigg(R\bigg(  L^\op_{\mathcal{B}^{k+t+m-2-i}_{m+1}}[\vec{x},\sum\limits^l_{t=0}\sum\limits^{m-2}_{i=0}\lambda^{(t)}\lambda^i\left(
\tbinom{m-2}{i}R\bigg(  L^\op_{\mathcal{B}^{k+t+m-3-i}_m}[\vec{x}]
 +    L^\op_{\mathcal{A}^{k+t+m-3-i}_m}[\vec{x}]\bigg) \right. \\
 &+ \lambda \tbinom{m-3}{i} \bigg(  L^\op_{\mathcal{B}^{k+t+m-3-i}_m}[\vec{x}]\bigg) 
 + \lambda \tbinom{m-2}{i}\bigg(    L^\op_{\mathcal{C}^{k+t+m-3-i}_m}[\vec{x}]\bigg)
 \left. +\lambda^2 \tbinom{m-3}{i} \bigg(   L^\op_{\mathcal{D}^{k+t+m-4-i}_m}[\vec{x}]\bigg)\right)]
 \\
& +        L^\op_{\mathcal{A}^{k+t+m-2-i}_{m+1}}[\vec{x},,\sum\limits^l_{t=0}\sum\limits^{m-2}_{i=0}\lambda^{(t)}\lambda^i\left(
\tbinom{m-2}{i}R\bigg(  L^\op_{\mathcal{B}^{k+t+m-3-i}_m}[\vec{x}]
 +    L^\op_{\mathcal{A}^{k+t+m-3-i}_m}[\vec{x}]\bigg) \right. \\
 &+ \lambda \tbinom{m-3}{i} \bigg(  L^\op_{\mathcal{B}^{k+t+m-3-i}_m}[\vec{x}]\bigg) 
 + \lambda \tbinom{m-2}{i}\bigg(    L^\op_{\mathcal{C}^{k+t+m-3-i}_m}[\vec{x}]\bigg)
 \left. +\lambda^2 \tbinom{m-3}{i} \bigg(   L^\op_{\mathcal{D}^{k+t+m-4-i}_m}[\vec{x}]\bigg)\right)]\bigg)x\bigg)
 \\
 &+ \lambda \tbinom{m-1}{i}R^2\bigg(\bigg(  L^\op_{\mathcal{B}^{k+t+m-2-i}_{m+1}}[\vec{x},\sum\limits^l_{t=0}\sum\limits^{m-2}_{i=0}\lambda^{(t)}\lambda^i\left(
\tbinom{m-2}{i}R\bigg(  L^\op_{\mathcal{B}^{k+t+m-3-i}_m}[\vec{x}]
 +    L^\op_{\mathcal{A}^{k+t+m-3-i}_m}[\vec{x}]\bigg) \right. \\
 &+ \lambda \tbinom{m-3}{i} \bigg(  L^\op_{\mathcal{B}^{k+t+m-3-i}_m}[\vec{x}]\bigg) 
 + \lambda \tbinom{m-2}{i}\bigg(    L^\op_{\mathcal{C}^{k+t+m-3-i}_m}[\vec{x}]\bigg)
 \left. +\lambda^2 \tbinom{m-3}{i} \bigg(   L^\op_{\mathcal{D}^{k+t+m-4-i}_m}[\vec{x}]\bigg)\right)]
 \\
& +        L^\op_{\mathcal{A}^{k+t+m-2-i}_{m+1}}[\vec{x},\sum\limits^l_{t=0}\sum\limits^{m-2}_{i=0}\lambda^{(t)}\lambda^i\left(
\tbinom{m-2}{i}R\bigg(  L^\op_{\mathcal{B}^{k+t+m-3-i}_m}[\vec{x}]
 +    L^\op_{\mathcal{A}^{k+t+m-3-i}_m}[\vec{x}]\bigg) \right. \\
 &+ \lambda \tbinom{m-3}{i} \bigg(  L^\op_{\mathcal{B}^{k+t+m-3-i}_m}[\vec{x}]\bigg) 
 + \lambda \tbinom{m-2}{i}\bigg(    L^\op_{\mathcal{C}^{k+t+m-3-i}_m}[\vec{x}]\bigg)
 \left. +\lambda^2 \tbinom{m-3}{i} \bigg(   L^\op_{\mathcal{D}^{k+t+m-4-i}_m}[\vec{x}]\bigg)\right)]\bigg)x\bigg)
\\
& + \lambda \tbinom{m-1}{i}R\bigg( R\bigg(  L^\op_{\mathcal{B}^{k+t+m-2-i}_{m+1}}[\vec{x},\sum\limits^l_{t=0}\sum\limits^{m-2}_{i=0}\lambda^{(t)}\lambda^i\left(
\tbinom{m-2}{i}R\bigg(  L^\op_{\mathcal{B}^{k+t+m-3-i}_m}[\vec{x}]
 +    L^\op_{\mathcal{A}^{k+t+m-3-i}_m}[\vec{x}]\bigg) \right. \\
 &+ \lambda \tbinom{m-3}{i} \bigg(  L^\op_{\mathcal{B}^{k+t+m-3-i}_m}[\vec{x}]\bigg) 
 + \lambda \tbinom{m-2}{i}\bigg(    L^\op_{\mathcal{C}^{k+t+m-3-i}_m}[\vec{x}]\bigg)
 \left. +\lambda^2 \tbinom{m-3}{i} \bigg(   L^\op_{\mathcal{D}^{k+t+m-4-i}_m}[\vec{x}]\bigg)\right)]
 \\
& +        L^\op_{\mathcal{A}^{k+t+m-2-i}_{m+1}}[\vec{x},\sum\limits^l_{t=0}\sum\limits^{m-2}_{i=0}\lambda^{(t)}\lambda^i\left(
\tbinom{m-2}{i}R\bigg(  L^\op_{\mathcal{B}^{k+t+m-3-i}_m}[\vec{x}]
 +    L^\op_{\mathcal{A}^{k+t+m-3-i}_m}[\vec{x}]\bigg) \right. \\
 &+ \lambda \tbinom{m-3}{i} \bigg(  L^\op_{\mathcal{B}^{k+t+m-3-i}_m}[\vec{x}]\bigg) 
 + \lambda \tbinom{m-2}{i}\bigg(    L^\op_{\mathcal{C}^{k+t+m-3-i}_m}[\vec{x}]\bigg)
 \left. +\lambda^2 \tbinom{m-3}{i} \bigg(   L^\op_{\mathcal{D}^{k+t+m-4-i}_m}[\vec{x}]\bigg)\right)]\bigg)x\bigg)
\\
 &+ \lambda^2 \tbinom{m-2}{i} R\bigg( \bigg(  L^\op_{\mathcal{B}^{k+t+m-2-i}_{m+1}}[\vec{x},\sum\limits^l_{t=0}\sum\limits^{m-2}_{i=0}\lambda^{(t)}\lambda^i\left(
\tbinom{m-2}{i}R\bigg(  L^\op_{\mathcal{B}^{k+t+m-3-i}_m}[\vec{x}]
 +    L^\op_{\mathcal{A}^{k+t+m-3-i}_m}[\vec{x}]\bigg) \right. \\
 &+ \lambda \tbinom{m-3}{i} \bigg(  L^\op_{\mathcal{B}^{k+t+m-3-i}_m}[\vec{x}]\bigg) 
 + \lambda \tbinom{m-2}{i}\bigg(    L^\op_{\mathcal{C}^{k+t+m-3-i}_m}[\vec{x}]\bigg)
 \left. +\lambda^2 \tbinom{m-3}{i} \bigg(   L^\op_{\mathcal{D}^{k+t+m-4-i}_m}[\vec{x}]\bigg)\right)]\bigg)x\bigg) 
 \\
& + \lambda^2 \tbinom{m-1}{i}R\bigg( \bigg(  L^\op_{\mathcal{C}^{k+t+m-2-i}_{m+1}}[\vec{x},\sum\limits^l_{t=0}\sum\limits^{m-2}_{i=0}\lambda^{(t)}\lambda^i\left(
\tbinom{m-2}{i}R\bigg(  L^\op_{\mathcal{B}^{k+t+m-3-i}_m}[\vec{x}]
 +    L^\op_{\mathcal{A}^{k+t+m-3-i}_m}[\vec{x}]\bigg) \right. \\
 &+ \lambda \tbinom{m-3}{i} \bigg(  L^\op_{\mathcal{B}^{k+t+m-3-i}_m}[\vec{x}]\bigg) 
 + \lambda \tbinom{m-2}{i}\bigg(    L^\op_{\mathcal{C}^{k+t+m-3-i}_m}[\vec{x}]\bigg)
 \left. +\lambda^2 \tbinom{m-3}{i} \bigg(   L^\op_{\mathcal{D}^{k+t+m-4-i}_m}[\vec{x}]\bigg)\right)]\bigg)x\bigg)
\\
  &+\lambda^3 \tbinom{m-2}{i}R\bigg( \bigg(  L^\op_{\mathcal{D}^{k+t+m-2-i}_{m+1}}[\vec{x},\sum\limits^l_{t=0}\sum\limits^{m-2}_{i=0}\lambda^{(t)}\lambda^i\left(
\tbinom{m-2}{i}R\bigg(  L^\op_{\mathcal{B}^{k+t+m-3-i}_m}[\vec{x}]
 +    L^\op_{\mathcal{A}^{k+t+m-3-i}_m}[\vec{x}]\bigg) \right. \\
 &+ \lambda \tbinom{m-3}{i} \bigg(  L^\op_{\mathcal{B}^{k+t+m-3-i}_m}[\vec{x}]\bigg) 
 + \lambda \tbinom{m-2}{i}\bigg(    L^\op_{\mathcal{C}^{k+t+m-3-i}_m}[\vec{x}]\bigg)
 \left. +\lambda^2 \tbinom{m-3}{i} \bigg(   L^\op_{\mathcal{D}^{k+t+m-4-i}_m}[\vec{x}]\bigg)\right)]\bigg)x\bigg) ) {=} 0;
 \ m = |\vec{x}|{\geq}3.
\end{aligned}
\]

Composition between relations $(R10)$ and $(R7)$, we get a sum of the expressions, which are consequences of~$(R10)$ itself.
Below, we show this example for  $\sum\limits_{t=0}^l \sum\limits^{m-1}_{i=0}\lambda^{(t)}\lambda^i
M_{t,i}$, where
\[
\begin{aligned}
M_{t,i}&=\sum\limits^l_{t=0}\sum\limits^{m-1}_{i=0}\lambda^{(t)}\lambda^i (
\tbinom{m-1}{i}R^2\bigg( {L^*}^\op_{\mathcal{B}^{k+t+m-2-i}_{m+1}}[x,\vec{x}, \sum\limits^l_{t=0}\lambda^{(t)}\bigg(
 R \bigg(xR\bigg(  L_{\mathcal{B}^{k+t-1}_2}[\vec{y}]
 + L_{\mathcal{A}^{k+t-1}_2}[\vec{y}]\bigg)\bigg)\\
 &+ \lambda R\bigg( x\bigg(  L_{\mathcal{B}^{k+t-1}_2}[\vec{y}]
 + L_{\mathcal{A}^{k+t-1}_2}[\vec{y}]\bigg)\bigg) 
+ \lambda \bigg( xR\bigg(  L_{\mathcal{B}^{k+t-1}_2}[\vec{y}]
 + L_{\mathcal{A}^{k+t-1}_2}[\vec{y}]\bigg)\bigg)
 + R\bigg(   L^*_{\mathcal{B}^{k+t-1}_3}[x,\vec{y}]+    L^*_{\mathcal{A}^{k+t-1}_3}[x,\vec{y}]\bigg) \\
& + \lambda \bigg(   L^*_{\mathcal{B}^{k+t-1}_3}[x,\vec{y}]
 +    L^*_{\mathcal{C}^{k+t-1}_3}[x,\vec{y}]\bigg)
 + \lambda^2 \bigg(x\bigg(  L_{\mathcal{B}^{k+t-1}_2}[\vec{y}]
 +       L_{\mathcal{C}^{k+t-1}_2}[\vec{y}]\bigg)\bigg)\bigg)]\\
&+ {L^*}^\op_{\mathcal{A}^{k+t+m-2-i}_{m+1}}[x,\vec{x}, \sum\limits^l_{t=0}\lambda^{(t)}\bigg(
 R \bigg(xR\bigg(  L_{\mathcal{B}^{k+t-1}_2}[\vec{y}] + L_{\mathcal{A}^{k+t-1}_2}[\vec{y}]\bigg)\bigg)
 + \lambda R\bigg( x\bigg(  L_{\mathcal{B}^{k+t-1}_2}[\vec{y}]
 + L_{\mathcal{A}^{k+t-1}_2}[\vec{y}]\bigg)\bigg) \\
& + \lambda \bigg( xR\bigg(  L_{\mathcal{B}^{k+t-1}_2}[\vec{y}]
 + L_{\mathcal{A}^{k+t-1}_2}[\vec{y}]\bigg)\bigg)
 + R\bigg(   L^*_{\mathcal{B}^{k+t-1}_3}[x,\vec{y}]+    L^*_{\mathcal{A}^{k+t-1}_3}[x,\vec{y}]\bigg)  + \lambda \bigg(   L^*_{\mathcal{B}^{k+t-1}_3}[x,\vec{y}]
 +    L^*_{\mathcal{C}^{k+t-1}_3}[x,\vec{y}]\bigg)\\
 &+ \lambda^2 \bigg(x\bigg(  L_{\mathcal{B}^{k+t-1}_2}[\vec{y}]
 +       L_{\mathcal{C}^{k+t-1}_2}[\vec{y}]\bigg)\bigg)\bigg)]\bigg) 
 +\lambda \tbinom{m-2}{i}R\bigg( {L^*}^\op_{\mathcal{B}^{k+t+m-2-i}_{m+1}}[x,\vec{x}, \sum\limits^l_{t=0}\lambda^{(t)}\bigg(
 R \bigg(xR\bigg(  L_{\mathcal{B}^{k+t-1}_2}[\vec{y}]
 + L_{\mathcal{A}^{k+t-1}_2}[\vec{y}]\bigg)\bigg)\\
 &+ \lambda R\bigg( x\bigg(  L_{\mathcal{B}^{k+t-1}_2}[\vec{y}]
 + L_{\mathcal{A}^{k+t-1}_2}[\vec{y}]\bigg)\bigg)  + \lambda \bigg( xR\bigg(  L_{\mathcal{B}^{k+t-1}_2}[\vec{y}]
 + L_{\mathcal{A}^{k+t-1}_2}[\vec{y}]\bigg)\bigg)
 + R\bigg(   L^*_{\mathcal{B}^{k+t-1}_3}[x,\vec{y}]+    L^*_{\mathcal{A}^{k+t-1}_3}[x,\vec{y}]\bigg) \\
& + \lambda \bigg(   L^*_{\mathcal{B}^{k+t-1}_3}[x,\vec{y}]
 +    L^*_{\mathcal{C}^{k+t-1}_3}[x,\vec{y}]\bigg)
 + \lambda^2 \bigg(x\bigg(  L_{\mathcal{B}^{k+t-1}_2}[\vec{y}]
 +       L_{\mathcal{C}^{k+t-1}_2}[\vec{y}]\bigg)\bigg)\bigg)]\bigg)
\\
&+\lambda \tbinom{m-1}{i} R\bigg( {L^*}^\op_{\mathcal{C}^{k+t+m-2-i}_{m+1}}[x,\vec{x}, \sum\limits^l_{t=0}\lambda^{(t)}\bigg(
 R \bigg(xR\bigg(  L_{\mathcal{B}^{k+t-1}_2}[\vec{y}] + L_{\mathcal{A}^{k+t-1}_2}[\vec{y}]\bigg)\bigg)
 + \lambda R\bigg( x\bigg(  L_{\mathcal{B}^{k+t-1}_2}[\vec{y}]
 + L_{\mathcal{A}^{k+t-1}_2}[\vec{y}]\bigg)\bigg) \\
& + \lambda \bigg( xR\bigg(  L_{\mathcal{B}^{k+t-1}_2}[\vec{y}]
 + L_{\mathcal{A}^{k+t-1}_2}[\vec{y}]\bigg)\bigg)
 + R\bigg(   L^*_{\mathcal{B}^{k+t-1}_3}[x,\vec{y}]+    L^*_{\mathcal{A}^{k+t-1}_3}[x,\vec{y}]\bigg) \\
& + \lambda \bigg(   L^*_{\mathcal{B}^{k+t-1}_3}[x,\vec{y}]
 +    L^*_{\mathcal{C}^{k+t-1}_3}[x,\vec{y}]\bigg)
 + \lambda^2 \bigg(x\bigg(  L_{\mathcal{B}^{k+t-1}_2}[\vec{y}]
 +       L_{\mathcal{C}^{k+t-1}_2}[\vec{y}]\bigg)\bigg)\bigg)]\bigg)\bigg)
\\
&+\lambda^2 \tbinom{m-2}{i}R\bigg(  {L^*}^\op_{\mathcal{D}^{k+t+m-3-i}_{m+1}}[x,\vec{x}, \sum\limits^l_{t=0}\lambda^{(t)}\bigg(
 R \bigg(xR\bigg(  L_{\mathcal{B}^{k+t-1}_2}[\vec{y}] + L_{\mathcal{A}^{k+t-1}_2}[\vec{y}]\bigg)\bigg)
 + \lambda R\bigg( x\bigg(  L_{\mathcal{B}^{k+t-1}_2}[\vec{y}]
 + L_{\mathcal{A}^{k+t-1}_2}[\vec{y}]\bigg)\bigg) \\
& + \lambda \bigg( xR\bigg(  L_{\mathcal{B}^{k+t-1}_2}[\vec{y}]
 + L_{\mathcal{A}^{k+t-1}_2}[\vec{y}]\bigg)\bigg)
 + R\bigg(   L^*_{\mathcal{B}^{k+t-1}_3}[x,\vec{y}]+    L^*_{\mathcal{A}^{k+t-1}_3}[x,\vec{y}]\bigg)  + \lambda \bigg(   L^*_{\mathcal{B}^{k+t-1}_3}[x,\vec{y}]
 +    L^*_{\mathcal{C}^{k+t-1}_3}[x,\vec{y}]\bigg)\\
& + \lambda^2 \bigg(x\bigg(  L_{\mathcal{B}^{k+t-1}_2}[\vec{y}]
 +       L_{\mathcal{C}^{k+t-1}_2}[\vec{y}]\bigg)\bigg)\bigg)]\bigg)
+ \tbinom{m-1}{i}R^2\bigg(R\bigg(  L^\op_{\mathcal{B}^{k+t+m-2-i}_{m+1}}[\vec{x}, \sum\limits^l_{t=0}\lambda^{(t)}\bigg(
 R \bigg(xR\bigg(  L_{\mathcal{B}^{k+t-1}_2}[\vec{y}] + L_{\mathcal{A}^{k+t-1}_2}[\vec{y}]\bigg)\bigg)\\
&
 + \lambda R\bigg( x\bigg(  L_{\mathcal{B}^{k+t-1}_2}[\vec{y}]
 + L_{\mathcal{A}^{k+t-1}_2}[\vec{y}]\bigg)\bigg)  + \lambda \bigg( xR\bigg(  L_{\mathcal{B}^{k+t-1}_2}[\vec{y}]
 + L_{\mathcal{A}^{k+t-1}_2}[\vec{y}]\bigg)\bigg)
 + R\bigg(   L^*_{\mathcal{B}^{k+t-1}_3}[x,\vec{y}]+    L^*_{\mathcal{A}^{k+t-1}_3}[x,\vec{y}]\bigg) \\
& + \lambda \bigg(   L^*_{\mathcal{B}^{k+t-1}_3}[x,\vec{y}]
 +    L^*_{\mathcal{C}^{k+t-1}_3}[x,\vec{y}]\bigg)
 + \lambda^2 \bigg(x\bigg(  L_{\mathcal{B}^{k+t-1}_2}[\vec{y}]
 +       L_{\mathcal{C}^{k+t-1}_2}[\vec{y}]\bigg)\bigg)\bigg)]
 \\
& +        L^\op_{\mathcal{A}^{k+t+m-2-i}_{m+1}}[\vec{x}, \sum\limits^l_{t=0}\lambda^{(t)}\bigg(
 R \bigg(xR\bigg(  L_{\mathcal{B}^{k+t-1}_2}[\vec{y}] + L_{\mathcal{A}^{k+t-1}_2}[\vec{y}]\bigg)\bigg)
 + \lambda R\bigg( x\bigg(  L_{\mathcal{B}^{k+t-1}_2}[\vec{y}]
 + L_{\mathcal{A}^{k+t-1}_2}[\vec{y}]\bigg)\bigg) \\
& + \lambda \bigg( xR\bigg(  L_{\mathcal{B}^{k+t-1}_2}[\vec{y}]
 + L_{\mathcal{A}^{k+t-1}_2}[\vec{y}]\bigg)\bigg)
 + R\bigg(   L^*_{\mathcal{B}^{k+t-1}_3}[x,\vec{y}]+    L^*_{\mathcal{A}^{k+t-1}_3}[x,\vec{y}]\bigg) \\
& + \lambda \bigg(   L^*_{\mathcal{B}^{k+t-1}_3}[x,\vec{y}]
 +    L^*_{\mathcal{C}^{k+t-1}_3}[x,\vec{y}]\bigg)
 + \lambda^2 \bigg(x\bigg(  L_{\mathcal{B}^{k+t-1}_2}[\vec{y}]
 +       L_{\mathcal{C}^{k+t-1}_2}[\vec{y}]\bigg)\bigg)\bigg)]\bigg)x\bigg)
\end{aligned}
 \]
 \[
 \begin{aligned}
 &+ \lambda \tbinom{m-1}{i}R^2\bigg(\bigg(  L^\op_{\mathcal{B}^{k+t+m-2-i}_{m+1}}[\vec{x}, \sum\limits^l_{t=0}\lambda^{(t)}\bigg(
 R \bigg(xR\bigg(  L_{\mathcal{B}^{k+t-1}_2}[\vec{y}] + L_{\mathcal{A}^{k+t-1}_2}[\vec{y}]\bigg)\bigg)
 + \lambda R\bigg( x\bigg(  L_{\mathcal{B}^{k+t-1}_2}[\vec{y}]
 + L_{\mathcal{A}^{k+t-1}_2}[\vec{y}]\bigg)\bigg) \\
& + \lambda \bigg( xR\bigg(  L_{\mathcal{B}^{k+t-1}_2}[\vec{y}]
 + L_{\mathcal{A}^{k+t-1}_2}[\vec{y}]\bigg)\bigg)
 + R\bigg(   L^*_{\mathcal{B}^{k+t-1}_3}[x,\vec{y}]+    L^*_{\mathcal{A}^{k+t-1}_3}[x,\vec{y}]\bigg) \\
& + \lambda \bigg(   L^*_{\mathcal{B}^{k+t-1}_3}[x,\vec{y}]
 +    L^*_{\mathcal{C}^{k+t-1}_3}[x,\vec{y}]\bigg)
 + \lambda^2 \bigg(x\bigg(  L_{\mathcal{B}^{k+t-1}_2}[\vec{y}]
 +       L_{\mathcal{C}^{k+t-1}_2}[\vec{y}]\bigg)\bigg)\bigg)]
 \\
& +        L^\op_{\mathcal{A}^{k+t+m-2-i}_{m+1}}[\vec{x}, \sum\limits^l_{t=0}\lambda^{(t)}\bigg(
 R \bigg(xR\bigg(  L_{\mathcal{B}^{k+t-1}_2}[\vec{y}] + L_{\mathcal{A}^{k+t-1}_2}[\vec{y}]\bigg)\bigg)
 + \lambda R\bigg( x\bigg(  L_{\mathcal{B}^{k+t-1}_2}[\vec{y}]
 + L_{\mathcal{A}^{k+t-1}_2}[\vec{y}]\bigg)\bigg) \\
& + \lambda \bigg( xR\bigg(  L_{\mathcal{B}^{k+t-1}_2}[\vec{y}]
 + L_{\mathcal{A}^{k+t-1}_2}[\vec{y}]\bigg)\bigg)
 + R\bigg(   L^*_{\mathcal{B}^{k+t-1}_3}[x,\vec{y}]+    L^*_{\mathcal{A}^{k+t-1}_3}[x,\vec{y}]\bigg)  + \lambda \bigg(   L^*_{\mathcal{B}^{k+t-1}_3}[x,\vec{y}]
 +    L^*_{\mathcal{C}^{k+t-1}_3}[x,\vec{y}]\bigg)\\
 &+ \lambda^2 \bigg(x\bigg(  L_{\mathcal{B}^{k+t-1}_2}[\vec{y}]
 +       L_{\mathcal{C}^{k+t-1}_2}[\vec{y}]\bigg)\bigg)\bigg)]\bigg)x\bigg)
 + \lambda \tbinom{m-1}{i}R\bigg( R\bigg(  L^\op_{\mathcal{B}^{k+t+m-2-i}_{m+1}}[\vec{x}, \sum\limits^l_{t=0}\lambda^{(t)}\bigg(
 R \bigg(xR\bigg(  L_{\mathcal{B}^{k+t-1}_2}[\vec{y}] + L_{\mathcal{A}^{k+t-1}_2}[\vec{y}]\bigg)\bigg)\\
 &+ \lambda R\bigg( x\bigg(  L_{\mathcal{B}^{k+t-1}_2}[\vec{y}]
 + L_{\mathcal{A}^{k+t-1}_2}[\vec{y}]\bigg)\bigg)  + \lambda \bigg( xR\bigg(  L_{\mathcal{B}^{k+t-1}_2}[\vec{y}]
 + L_{\mathcal{A}^{k+t-1}_2}[\vec{y}]\bigg)\bigg)
 + R\bigg(   L^*_{\mathcal{B}^{k+t-1}_3}[x,\vec{y}]+    L^*_{\mathcal{A}^{k+t-1}_3}[x,\vec{y}]\bigg) \\
& + \lambda \bigg(   L^*_{\mathcal{B}^{k+t-1}_3}[x,\vec{y}]
 +    L^*_{\mathcal{C}^{k+t-1}_3}[x,\vec{y}]\bigg)
 + \lambda^2 \bigg(x\bigg(  L_{\mathcal{B}^{k+t-1}_2}[\vec{y}]
 +       L_{\mathcal{C}^{k+t-1}_2}[\vec{y}]\bigg)\bigg)\bigg)]
 \\
& +        L^\op_{\mathcal{A}^{k+t+m-2-i}_{m+1}}[\vec{x}, \sum\limits^l_{t=0}\lambda^{(t)}\bigg(
 R \bigg(xR\bigg(  L_{\mathcal{B}^{k+t-1}_2}[\vec{y}]
+ L_{\mathcal{A}^{k+t-1}_2}[\vec{y}]\bigg)\bigg) + \lambda R\bigg( x\bigg(  L_{\mathcal{B}^{k+t-1}_2}[\vec{y}]
 + L_{\mathcal{A}^{k+t-1}_2}[\vec{y}]\bigg)\bigg) \\
& + \lambda \bigg( xR\bigg(  L_{\mathcal{B}^{k+t-1}_2}[\vec{y}]
 + L_{\mathcal{A}^{k+t-1}_2}[\vec{y}]\bigg)\bigg)
 + R\bigg(   L^*_{\mathcal{B}^{k+t-1}_3}[x,\vec{y}]+    L^*_{\mathcal{A}^{k+t-1}_3}[x,\vec{y}]\bigg)  + \lambda \bigg(   L^*_{\mathcal{B}^{k+t-1}_3}[x,\vec{y}]
 +    L^*_{\mathcal{C}^{k+t-1}_3}[x,\vec{y}]\bigg)\\
 &+ \lambda^2 \bigg(x\bigg(  L_{\mathcal{B}^{k+t-1}_2}[\vec{y}]
 +       L_{\mathcal{C}^{k+t-1}_2}[\vec{y}]\bigg)\bigg)\bigg)]\bigg)x\bigg)
 + \lambda^2 \tbinom{m-2}{i} R\bigg( \bigg(  L^\op_{\mathcal{B}^{k+t+m-2-i}_{m+1}}[\vec{x}, \sum\limits^l_{t=0}\lambda^{(t)}\bigg(
 R \bigg(xR\bigg(  L_{\mathcal{B}^{k+t-1}_2}[\vec{y}]
 + L_{\mathcal{A}^{k+t-1}_2}[\vec{y}]\bigg)\bigg)\\
 &+ \lambda R\bigg( x\bigg(  L_{\mathcal{B}^{k+t-1}_2}[\vec{y}]
 + L_{\mathcal{A}^{k+t-1}_2}[\vec{y}]\bigg)\bigg) + \lambda \bigg( xR\bigg(  L_{\mathcal{B}^{k+t-1}_2}[\vec{y}]
 + L_{\mathcal{A}^{k+t-1}_2}[\vec{y}]\bigg)\bigg)
 + R\bigg(   L^*_{\mathcal{B}^{k+t-1}_3}[x,\vec{y}]+    L^*_{\mathcal{A}^{k+t-1}_3}[x,\vec{y}]\bigg) \\
& + \lambda \bigg(   L^*_{\mathcal{B}^{k+t-1}_3}[x,\vec{y}]
 +    L^*_{\mathcal{C}^{k+t-1}_3}[x,\vec{y}]\bigg)
 + \lambda^2 \bigg(x\bigg(  L_{\mathcal{B}^{k+t-1}_2}[\vec{y}]
 +       L_{\mathcal{C}^{k+t-1}_2}[\vec{y}]\bigg)\bigg)\bigg)]\bigg)x\bigg) 
 \\
& + \lambda^2 \tbinom{m-1}{i}R\bigg( \bigg(  L^\op_{\mathcal{C}^{k+t+m-2-i}_{m+1}}[\vec{x}, \sum\limits^l_{t=0}\lambda^{(t)}\bigg(
 R \bigg(xR\bigg(  L_{\mathcal{B}^{k+t-1}_2}[\vec{y}]
 + L_{\mathcal{A}^{k+t-1}_2}[\vec{y}]\bigg)\bigg)
 + \lambda R\bigg( x\bigg(  L_{\mathcal{B}^{k+t-1}_2}[\vec{y}]
 + L_{\mathcal{A}^{k+t-1}_2}[\vec{y}]\bigg)\bigg) \\
 &+ \lambda \bigg( xR\bigg(  L_{\mathcal{B}^{k+t-1}_2}[\vec{y}]
 + L_{\mathcal{A}^{k+t-1}_2}[\vec{y}]\bigg)\bigg)
 + R\bigg(   L^*_{\mathcal{B}^{k+t-1}_3}[x,\vec{y}]+    L^*_{\mathcal{A}^{k+t-1}_3}[x,\vec{y}]\bigg)  + \lambda \bigg(   L^*_{\mathcal{B}^{k+t-1}_3}[x,\vec{y}]
 +    L^*_{\mathcal{C}^{k+t-1}_3}[x,\vec{y}]\bigg)\\
 &+ \lambda^2 \bigg(x\bigg(  L_{\mathcal{B}^{k+t-1}_2}[\vec{y}]
 +       L_{\mathcal{C}^{k+t-1}_2}[\vec{y}]\bigg)\bigg)\bigg)]\bigg)x\bigg)
+\lambda^3 \tbinom{m-2}{i}R\bigg( \bigg(  L^\op_{\mathcal{D}^{k+t+m-2-i}_{m+1}}[\vec{x}, \sum\limits^l_{t=0}\lambda^{(t)}\bigg(
 R \bigg(xR\bigg(  L_{\mathcal{B}^{k+t-1}_2}[\vec{y}]
 + L_{\mathcal{A}^{k+t-1}_2}[\vec{y}]\bigg)\bigg)\\
& + \lambda R\bigg( x\bigg(  L_{\mathcal{B}^{k+t-1}_2}[\vec{y}]
 + L_{\mathcal{A}^{k+t-1}_2}[\vec{y}]\bigg)\bigg) 
 + \lambda \bigg( xR\bigg(  L_{\mathcal{B}^{k+t-1}_2}[\vec{y}]
 + L_{\mathcal{A}^{k+t-1}_2}[\vec{y}]\bigg)\bigg)
 + R\bigg(   L^*_{\mathcal{B}^{k+t-1}_3}[x,\vec{y}]+    L^*_{\mathcal{A}^{k+t-1}_3}[x,\vec{y}]\bigg) \\
&+ \lambda \bigg(   L^*_{\mathcal{B}^{k+t-1}_3}[x,\vec{y}]
 +    L^*_{\mathcal{C}^{k+t-1}_3}[x,\vec{y}]\bigg) + \lambda^2 \bigg(x\bigg(  L_{\mathcal{B}^{k+t-1}_2}[\vec{y}]
 +       L_{\mathcal{C}^{k+t-1}_2}[\vec{y}]\bigg)\bigg)\bigg)]\bigg)x\bigg) ) {=} 0;\ m = |\vec{x}|{\geq}3, |\vec{y}|=2.
\end{aligned}
\]
Composition between relations $(R10)$ and $(R8)$, we get a sum of the expressions, which are consequences of~$(R10)$ itself.
Below, we show this example for  $\sum\limits_{t=0}^l \sum\limits^{m-1}_{i=0}\lambda^{(t)}\lambda^i
M_{t,i}$, where
\[
\begin{aligned}
M_{t,i}&=
\tbinom{m-1}{i}R^2\bigg( {L^*}^\op_{\mathcal{B}^{k+t+m-2-i}_{m+1}}[x,\vec{x},\sum\limits^l_{t=0}\sum\limits^{m-2}_{i=0}\lambda^{(t)}\lambda^i \left( 
 \tbinom{m-2}{i}R\bigg(  L^*_{\mathcal{B}^{k+t+m-3-i}_m}[x,\vec{x}]
 +    L^*_{\mathcal{A}^{k+t+m-3-i}_m}[x,\vec{x}]\bigg)\right. \\
& + \lambda \tbinom{m-3}{i} \bigg(  L^*_{\mathcal{B}^{k+t+m-3-i}_m}[x,\vec{x}]\bigg) 
 + \lambda \tbinom{m-2}{i} \bigg(    L^*_{\mathcal{C}^{k+t+m-3-i}_m}[x,\vec{x}]\bigg)  + \lambda^2 \tbinom{m-3}{i} \bigg(   L^*_{\mathcal{D}^{k+t+m-4-i}_m}[x,\vec{x}]\bigg)\\
 &+ \tbinom{m-2}{i}R\bigg(xR\bigg(          L_{\mathcal{B}^{k+t+m-3-i}_{m}}[\vec{x}]
 +        L_{\mathcal{A}^{k+t+m-3-i}_{m}}[\vec{x}]\bigg)\bigg)  + \lambda \tbinom{m-2}{i} R\bigg(x\bigg(          L_{\mathcal{B}^{k+t+m-3-i}_{m}}[\vec{x}]
 +            L_{\mathcal{A}^{k+t+m-3-i}_{m}}[\vec{x}]\bigg)\bigg) \allowdisplaybreaks \\
& + \lambda \tbinom{m-2}{i}\bigg(x R\bigg(          L_{\mathcal{B}^{k+t+m-3-i}_{m}}[\vec{x}]
 +        L_{\mathcal{A}^{k+t+m-3-i}_{m}}[\vec{x}]\bigg)\bigg)
 + \lambda^2 \tbinom{m-3}{i} \bigg(x \bigg(          L_{\mathcal{B}^{k+t+m-3-i}_{m}}[\vec{x}]\bigg)\bigg)\bigg) \\
& + \lambda^2 \tbinom{m-2}{i} \bigg(x \bigg(            L_{\mathcal{C}^{k+t+m-3-i}_{m}}[\vec{x}]\bigg)\bigg)\bigg)
\left. +\lambda^3 \tbinom{m-3}{i}\bigg(x \bigg(             L_{\mathcal{D}^{k+t+m-3-i}_{m}}[\vec{x}]\bigg)\bigg)\bigg) \right)]
\\
&+ {L^*}^\op_{\mathcal{A}^{k+t+m-2-i}_{m+1}}[x,\vec{x},\sum\limits^l_{t=0}\sum\limits^{m-2}_{i=0}\lambda^{(t)}\lambda^i \left( 
 \tbinom{m-2}{i}R\bigg(  L^*_{\mathcal{B}^{k+t+m-3-i}_m}[x,\vec{x}]
 +    L^*_{\mathcal{A}^{k+t+m-3-i}_m}[x,\vec{x}]\bigg)\right. \\
& + \lambda \tbinom{m-3}{i} \bigg(  L^*_{\mathcal{B}^{k+t+m-3-i}_m}[x,\vec{x}]\bigg) 
 + \lambda \tbinom{m-2}{i} \bigg(    L^*_{\mathcal{C}^{k+t+m-3-i}_m}[x,\vec{x}]\bigg)  + \lambda^2 \tbinom{m-3}{i} \bigg(   L^*_{\mathcal{D}^{k+t+m-4-i}_m}[x,\vec{x}]\bigg)\\
& + \tbinom{m-2}{i}R\bigg(xR\bigg(          L_{\mathcal{B}^{k+t+m-3-i}_{m}}[\vec{x}]
 +        L_{\mathcal{A}^{k+t+m-3-i}_{m}}[\vec{x}]\bigg)\bigg)  + \lambda \tbinom{m-2}{i} R\bigg(x\bigg(          L_{\mathcal{B}^{k+t+m-3-i}_{m}}[\vec{x}]
 +            L_{\mathcal{A}^{k+t+m-3-i}_{m}}[\vec{x}]\bigg)\bigg) \allowdisplaybreaks \\
& + \lambda \tbinom{m-2}{i}\bigg(x R\bigg(          L_{\mathcal{B}^{k+t+m-3-i}_{m}}[\vec{x}]
 +        L_{\mathcal{A}^{k+t+m-3-i}_{m}}[\vec{x}]\bigg)\bigg)
 + \lambda^2 \tbinom{m-3}{i} \bigg(x \bigg(          L_{\mathcal{B}^{k+t+m-3-i}_{m}}[\vec{x}]\bigg)\bigg)\bigg) \\
& + \lambda^2 \tbinom{m-2}{i} \bigg(x \bigg(            L_{\mathcal{C}^{k+t+m-3-i}_{m}}[\vec{x}]\bigg)\bigg)\bigg)
\left. +\lambda^3 \tbinom{m-3}{i}\bigg(x \bigg(             L_{\mathcal{D}^{k+t+m-3-i}_{m}}[\vec{x}]\bigg)\bigg)\bigg) \right)]\bigg) 
\\
&+\lambda \tbinom{m-2}{i}R\bigg( {L^*}^\op_{\mathcal{B}^{k+t+m-2-i}_{m+1}}[x,\vec{x},\sum\limits^l_{t=0}\sum\limits^{m-2}_{i=0}\lambda^{(t)}\lambda^i \left( 
 \tbinom{m-2}{i}R\bigg(  L^*_{\mathcal{B}^{k+t+m-3-i}_m}[x,\vec{x}]
 +    L^*_{\mathcal{A}^{k+t+m-3-i}_m}[x,\vec{x}]\bigg)\right. \\
& + \lambda \tbinom{m-3}{i} \bigg(  L^*_{\mathcal{B}^{k+t+m-3-i}_m}[x,\vec{x}]\bigg) 
 + \lambda \tbinom{m-2}{i} \bigg(    L^*_{\mathcal{C}^{k+t+m-3-i}_m}[x,\vec{x}]\bigg)  + \lambda^2 \tbinom{m-3}{i} \bigg(   L^*_{\mathcal{D}^{k+t+m-4-i}_m}[x,\vec{x}]\bigg)\\
 &+ \tbinom{m-2}{i}R\bigg(xR\bigg(          L_{\mathcal{B}^{k+t+m-3-i}_{m}}[\vec{x}]
 +        L_{\mathcal{A}^{k+t+m-3-i}_{m}}[\vec{x}]\bigg)\bigg) 
 + \lambda \tbinom{m-2}{i} R\bigg(x\bigg(          L_{\mathcal{B}^{k+t+m-3-i}_{m}}[\vec{x}]
 +            L_{\mathcal{A}^{k+t+m-3-i}_{m}}[\vec{x}]\bigg)\bigg) \allowdisplaybreaks \\
& + \lambda \tbinom{m-2}{i}\bigg(x R\bigg(          L_{\mathcal{B}^{k+t+m-3-i}_{m}}[\vec{x}]
 +        L_{\mathcal{A}^{k+t+m-3-i}_{m}}[\vec{x}]\bigg)\bigg)
 + \lambda^2 \tbinom{m-3}{i} \bigg(x \bigg(          L_{\mathcal{B}^{k+t+m-3-i}_{m}}[\vec{x}]\bigg)\bigg)\bigg) \\
& + \lambda^2 \tbinom{m-2}{i} \bigg(x \bigg(            L_{\mathcal{C}^{k+t+m-3-i}_{m}}[\vec{x}]\bigg)\bigg)\bigg)
\left. +\lambda^3 \tbinom{m-3}{i}\bigg(x \bigg(             L_{\mathcal{D}^{k+t+m-3-i}_{m}}[\vec{x}]\bigg)\bigg)\bigg) \right)]\bigg)
\\
&+\lambda \tbinom{m-1}{i} R\bigg( {L^*}^\op_{\mathcal{C}^{k+t+m-2-i}_{m+1}}[x,\vec{x},\sum\limits^l_{t=0}\sum\limits^{m-2}_{i=0}\lambda^{(t)}\lambda^i \left( 
 \tbinom{m-2}{i}R\bigg(  L^*_{\mathcal{B}^{k+t+m-3-i}_m}[x,\vec{x}]
 +    L^*_{\mathcal{A}^{k+t+m-3-i}_m}[x,\vec{x}]\bigg)\right. \\
& + \lambda \tbinom{m-3}{i} \bigg(  L^*_{\mathcal{B}^{k+t+m-3-i}_m}[x,\vec{x}]\bigg) 
 + \lambda \tbinom{m-2}{i} \bigg(    L^*_{\mathcal{C}^{k+t+m-3-i}_m}[x,\vec{x}]\bigg)  + \lambda^2 \tbinom{m-3}{i} \bigg(   L^*_{\mathcal{D}^{k+t+m-4-i}_m}[x,\vec{x}]\bigg)\\
 &+ \tbinom{m-2}{i}R\bigg(xR\bigg(          L_{\mathcal{B}^{k+t+m-3-i}_{m}}[\vec{x}]
 +        L_{\mathcal{A}^{k+t+m-3-i}_{m}}[\vec{x}]\bigg)\bigg)  + \lambda \tbinom{m-2}{i} R\bigg(x\bigg(          L_{\mathcal{B}^{k+t+m-3-i}_{m}}[\vec{x}]
 +            L_{\mathcal{A}^{k+t+m-3-i}_{m}}[\vec{x}]\bigg)\bigg) \allowdisplaybreaks \\
& + \lambda \tbinom{m-2}{i}\bigg(x R\bigg(          L_{\mathcal{B}^{k+t+m-3-i}_{m}}[\vec{x}]
 +        L_{\mathcal{A}^{k+t+m-3-i}_{m}}[\vec{x}]\bigg)\bigg)
 + \lambda^2 \tbinom{m-3}{i} \bigg(x \bigg(          L_{\mathcal{B}^{k+t+m-3-i}_{m}}[\vec{x}]\bigg)\bigg)\bigg) \\
& + \lambda^2 \tbinom{m-2}{i} \bigg(x \bigg(            L_{\mathcal{C}^{k+t+m-3-i}_{m}}[\vec{x}]\bigg)\bigg)\bigg)
\left. +\lambda^3 \tbinom{m-3}{i}\bigg(x \bigg(             L_{\mathcal{D}^{k+t+m-3-i}_{m}}[\vec{x}]\bigg)\bigg)\bigg) \right)]\bigg)\bigg)
\end{aligned}
 \]
 \[
 \begin{aligned}
&+\lambda^2 \tbinom{m-2}{i}R\bigg(  {L^*}^\op_{\mathcal{D}^{k+t+m-3-i}_{m+1}}[x,\vec{x},\sum\limits^l_{t=0}\sum\limits^{m-2}_{i=0}\lambda^{(t)}\lambda^i \left( 
 \tbinom{m-2}{i}R\bigg(  L^*_{\mathcal{B}^{k+t+m-3-i}_m}[x,\vec{x}]
 +    L^*_{\mathcal{A}^{k+t+m-3-i}_m}[x,\vec{x}]\bigg)\right. \\
& + \lambda \tbinom{m-3}{i} \bigg(  L^*_{\mathcal{B}^{k+t+m-3-i}_m}[x,\vec{x}]\bigg) 
 + \lambda \tbinom{m-2}{i} \bigg(    L^*_{\mathcal{C}^{k+t+m-3-i}_m}[x,\vec{x}]\bigg)  + \lambda^2 \tbinom{m-3}{i} \bigg(   L^*_{\mathcal{D}^{k+t+m-4-i}_m}[x,\vec{x}]\bigg)\\
& + \tbinom{m-2}{i}R\bigg(xR\bigg(          L_{\mathcal{B}^{k+t+m-3-i}_{m}}[\vec{x}]
 +        L_{\mathcal{A}^{k+t+m-3-i}_{m}}[\vec{x}]\bigg)\bigg) + \lambda \tbinom{m-2}{i} R\bigg(x\bigg(          L_{\mathcal{B}^{k+t+m-3-i}_{m}}[\vec{x}]
 +            L_{\mathcal{A}^{k+t+m-3-i}_{m}}[\vec{x}]\bigg)\bigg) \allowdisplaybreaks \\
& + \lambda \tbinom{m-2}{i}\bigg(x R\bigg(          L_{\mathcal{B}^{k+t+m-3-i}_{m}}[\vec{x}]
 +        L_{\mathcal{A}^{k+t+m-3-i}_{m}}[\vec{x}]\bigg)\bigg)
 + \lambda^2 \tbinom{m-3}{i} \bigg(x \bigg(          L_{\mathcal{B}^{k+t+m-3-i}_{m}}[\vec{x}]\bigg)\bigg)\bigg) \\
& + \lambda^2 \tbinom{m-2}{i} \bigg(x \bigg(            L_{\mathcal{C}^{k+t+m-3-i}_{m}}[\vec{x}]\bigg)\bigg)\bigg)
\left. +\lambda^3 \tbinom{m-3}{i}\bigg(x \bigg(             L_{\mathcal{D}^{k+t+m-3-i}_{m}}[\vec{x}]\bigg)\bigg)\bigg) \right)]\bigg)
\\
 &+ \tbinom{m-1}{i}R^2\bigg(R\bigg(  L^\op_{\mathcal{B}^{k+t+m-2-i}_{m+1}}[\vec{x},\sum\limits^l_{t=0}\sum\limits^{m-2}_{i=0}\lambda^{(t)}\lambda^i \left( 
 \tbinom{m-2}{i}R\bigg(  L^*_{\mathcal{B}^{k+t+m-3-i}_m}[x,\vec{x}]
 +    L^*_{\mathcal{A}^{k+t+m-3-i}_m}[x,\vec{x}]\bigg)\right. \\
& + \lambda \tbinom{m-3}{i} \bigg(  L^*_{\mathcal{B}^{k+t+m-3-i}_m}[x,\vec{x}]\bigg) 
 + \lambda \tbinom{m-2}{i} \bigg(    L^*_{\mathcal{C}^{k+t+m-3-i}_m}[x,\vec{x}]\bigg)  + \lambda^2 \tbinom{m-3}{i} \bigg(   L^*_{\mathcal{D}^{k+t+m-4-i}_m}[x,\vec{x}]\bigg)\\
 &+ \tbinom{m-2}{i}R\bigg(xR\bigg(          L_{\mathcal{B}^{k+t+m-3-i}_{m}}[\vec{x}]
 +        L_{\mathcal{A}^{k+t+m-3-i}_{m}}[\vec{x}]\bigg)\bigg)  + \lambda \tbinom{m-2}{i} R\bigg(x\bigg(          L_{\mathcal{B}^{k+t+m-3-i}_{m}}[\vec{x}]
 +            L_{\mathcal{A}^{k+t+m-3-i}_{m}}[\vec{x}]\bigg)\bigg) \allowdisplaybreaks \\
& + \lambda \tbinom{m-2}{i}\bigg(x R\bigg(          L_{\mathcal{B}^{k+t+m-3-i}_{m}}[\vec{x}]
 +        L_{\mathcal{A}^{k+t+m-3-i}_{m}}[\vec{x}]\bigg)\bigg)
 + \lambda^2 \tbinom{m-3}{i} \bigg(x \bigg(          L_{\mathcal{B}^{k+t+m-3-i}_{m}}[\vec{x}]\bigg)\bigg)\bigg) \\
& + \lambda^2 \tbinom{m-2}{i} \bigg(x \bigg(            L_{\mathcal{C}^{k+t+m-3-i}_{m}}[\vec{x}]\bigg)\bigg)\bigg)
\left. +\lambda^3 \tbinom{m-3}{i}\bigg(x \bigg(             L_{\mathcal{D}^{k+t+m-3-i}_{m}}[\vec{x}]\bigg)\bigg)\bigg) \right)]
 \\
& +        L^\op_{\mathcal{A}^{k+t+m-2-i}_{m+1}}[\vec{x},\sum\limits^l_{t=0}\sum\limits^{m-2}_{i=0}\lambda^{(t)}\lambda^i \left( 
 \tbinom{m-2}{i}R\bigg(  L^*_{\mathcal{B}^{k+t+m-3-i}_m}[x,\vec{x}]
 +    L^*_{\mathcal{A}^{k+t+m-3-i}_m}[x,\vec{x}]\bigg)\right. \\
& + \lambda \tbinom{m-3}{i} \bigg(  L^*_{\mathcal{B}^{k+t+m-3-i}_m}[x,\vec{x}]\bigg) 
 + \lambda \tbinom{m-2}{i} \bigg(    L^*_{\mathcal{C}^{k+t+m-3-i}_m}[x,\vec{x}]\bigg)  + \lambda^2 \tbinom{m-3}{i} \bigg(   L^*_{\mathcal{D}^{k+t+m-4-i}_m}[x,\vec{x}]\bigg)\\
& + \tbinom{m-2}{i}R\bigg(xR\bigg(          L_{\mathcal{B}^{k+t+m-3-i}_{m}}[\vec{x}]
 +        L_{\mathcal{A}^{k+t+m-3-i}_{m}}[\vec{x}]\bigg)\bigg)  + \lambda \tbinom{m-2}{i} R\bigg(x\bigg(          L_{\mathcal{B}^{k+t+m-3-i}_{m}}[\vec{x}]
 +            L_{\mathcal{A}^{k+t+m-3-i}_{m}}[\vec{x}]\bigg)\bigg) \allowdisplaybreaks \\
& + \lambda \tbinom{m-2}{i}\bigg(x R\bigg(          L_{\mathcal{B}^{k+t+m-3-i}_{m}}[\vec{x}]
 +        L_{\mathcal{A}^{k+t+m-3-i}_{m}}[\vec{x}]\bigg)\bigg)
 + \lambda^2 \tbinom{m-3}{i} \bigg(x \bigg(          L_{\mathcal{B}^{k+t+m-3-i}_{m}}[\vec{x}]\bigg)\bigg)\bigg) \\
& + \lambda^2 \tbinom{m-2}{i} \bigg(x \bigg(            L_{\mathcal{C}^{k+t+m-3-i}_{m}}[\vec{x}]\bigg)\bigg)\bigg)
\left. +\lambda^3 \tbinom{m-3}{i}\bigg(x \bigg(             L_{\mathcal{D}^{k+t+m-3-i}_{m}}[\vec{x}]\bigg)\bigg)\bigg) \right)]\bigg)x\bigg)
\\
 &+ \lambda \tbinom{m-1}{i}R^2\bigg(\bigg(  L^\op_{\mathcal{B}^{k+t+m-2-i}_{m+1}}[\vec{x},\sum\limits^l_{t=0}\sum\limits^{m-2}_{i=0}\lambda^{(t)}\lambda^i \left( 
 \tbinom{m-2}{i}R\bigg(  L^*_{\mathcal{B}^{k+t+m-3-i}_m}[x,\vec{x}]
 +    L^*_{\mathcal{A}^{k+t+m-3-i}_m}[x,\vec{x}]\bigg)\right. \\
& + \lambda \tbinom{m-3}{i} \bigg(  L^*_{\mathcal{B}^{k+t+m-3-i}_m}[x,\vec{x}]\bigg) 
 + \lambda \tbinom{m-2}{i} \bigg(    L^*_{\mathcal{C}^{k+t+m-3-i}_m}[x,\vec{x}]\bigg)  + \lambda^2 \tbinom{m-3}{i} \bigg(   L^*_{\mathcal{D}^{k+t+m-4-i}_m}[x,\vec{x}]\bigg)\\
 &+ \tbinom{m-2}{i}R\bigg(xR\bigg(          L_{\mathcal{B}^{k+t+m-3-i}_{m}}[\vec{x}]
 +        L_{\mathcal{A}^{k+t+m-3-i}_{m}}[\vec{x}]\bigg)\bigg)
 + \lambda \tbinom{m-2}{i} R\bigg(x\bigg(          L_{\mathcal{B}^{k+t+m-3-i}_{m}}[\vec{x}]
 +            L_{\mathcal{A}^{k+t+m-3-i}_{m}}[\vec{x}]\bigg)\bigg) \allowdisplaybreaks \\
& + \lambda \tbinom{m-2}{i}\bigg(x R\bigg(          L_{\mathcal{B}^{k+t+m-3-i}_{m}}[\vec{x}]
 +        L_{\mathcal{A}^{k+t+m-3-i}_{m}}[\vec{x}]\bigg)\bigg)
 + \lambda^2 \tbinom{m-3}{i} \bigg(x \bigg(          L_{\mathcal{B}^{k+t+m-3-i}_{m}}[\vec{x}]\bigg)\bigg)\bigg) \\
& + \lambda^2 \tbinom{m-2}{i} \bigg(x \bigg(            L_{\mathcal{C}^{k+t+m-3-i}_{m}}[\vec{x}]\bigg)\bigg)\bigg)
\left. +\lambda^3 \tbinom{m-3}{i}\bigg(x \bigg(             L_{\mathcal{D}^{k+t+m-3-i}_{m}}[\vec{x}]\bigg)\bigg)\bigg) \right)]
 \\
& +        L^\op_{\mathcal{A}^{k+t+m-2-i}_{m+1}}[\vec{x},\sum\limits^l_{t=0}\sum\limits^{m-2}_{i=0}\lambda^{(t)}\lambda^i \left( 
 \tbinom{m-2}{i}R\bigg(  L^*_{\mathcal{B}^{k+t+m-3-i}_m}[x,\vec{x}]
 +    L^*_{\mathcal{A}^{k+t+m-3-i}_m}[x,\vec{x}]\bigg)\right. \\
 \end{aligned}
 \]
 \[
 \begin{aligned}
& + \lambda \tbinom{m-3}{i} \bigg(  L^*_{\mathcal{B}^{k+t+m-3-i}_m}[x,\vec{x}]\bigg) 
 + \lambda \tbinom{m-2}{i} \bigg(    L^*_{\mathcal{C}^{k+t+m-3-i}_m}[x,\vec{x}]\bigg)  + \lambda^2 \tbinom{m-3}{i} \bigg(   L^*_{\mathcal{D}^{k+t+m-4-i}_m}[x,\vec{x}]\bigg)\\
 &+ \tbinom{m-2}{i}R\bigg(xR\bigg(          L_{\mathcal{B}^{k+t+m-3-i}_{m}}[\vec{x}]
 +        L_{\mathcal{A}^{k+t+m-3-i}_{m}}[\vec{x}]\bigg)\bigg) 
 + \lambda \tbinom{m-2}{i} R\bigg(x\bigg(          L_{\mathcal{B}^{k+t+m-3-i}_{m}}[\vec{x}]
 +            L_{\mathcal{A}^{k+t+m-3-i}_{m}}[\vec{x}]\bigg)\bigg) \allowdisplaybreaks \\
& + \lambda \tbinom{m-2}{i}\bigg(x R\bigg(          L_{\mathcal{B}^{k+t+m-3-i}_{m}}[\vec{x}]
 +        L_{\mathcal{A}^{k+t+m-3-i}_{m}}[\vec{x}]\bigg)\bigg)
 + \lambda^2 \tbinom{m-3}{i} \bigg(x \bigg(          L_{\mathcal{B}^{k+t+m-3-i}_{m}}[\vec{x}]\bigg)\bigg)\bigg) \\
& + \lambda^2 \tbinom{m-2}{i} \bigg(x \bigg(            L_{\mathcal{C}^{k+t+m-3-i}_{m}}[\vec{x}]\bigg)\bigg)\bigg)
\left. +\lambda^3 \tbinom{m-3}{i}\bigg(x \bigg(             L_{\mathcal{D}^{k+t+m-3-i}_{m}}[\vec{x}]\bigg)\bigg)\bigg) \right)]\bigg)x\bigg)
\\
& + \lambda \tbinom{m-1}{i}R\bigg( R\bigg(  L^\op_{\mathcal{B}^{k+t+m-2-i}_{m+1}}[\vec{x},\sum\limits^l_{t=0}\sum\limits^{m-2}_{i=0}\lambda^{(t)}\lambda^i \left( 
 \tbinom{m-2}{i}R\bigg(  L^*_{\mathcal{B}^{k+t+m-3-i}_m}[x,\vec{x}]
 +    L^*_{\mathcal{A}^{k+t+m-3-i}_m}[x,\vec{x}]\bigg)\right. \\
& + \lambda \tbinom{m-3}{i} \bigg(  L^*_{\mathcal{B}^{k+t+m-3-i}_m}[x,\vec{x}]\bigg) 
 + \lambda \tbinom{m-2}{i} \bigg(    L^*_{\mathcal{C}^{k+t+m-3-i}_m}[x,\vec{x}]\bigg)  + \lambda^2 \tbinom{m-3}{i} \bigg(   L^*_{\mathcal{D}^{k+t+m-4-i}_m}[x,\vec{x}]\bigg)\\
& + \tbinom{m-2}{i}R\bigg(xR\bigg(          L_{\mathcal{B}^{k+t+m-3-i}_{m}}[\vec{x}]
 +        L_{\mathcal{A}^{k+t+m-3-i}_{m}}[\vec{x}]\bigg)\bigg)  + \lambda \tbinom{m-2}{i} R\bigg(x\bigg(          L_{\mathcal{B}^{k+t+m-3-i}_{m}}[\vec{x}]
 +            L_{\mathcal{A}^{k+t+m-3-i}_{m}}[\vec{x}]\bigg)\bigg) \allowdisplaybreaks 
 \\
& + \lambda \tbinom{m-2}{i}\bigg(x R\bigg(          L_{\mathcal{B}^{k+t+m-3-i}_{m}}[\vec{x}]
 +        L_{\mathcal{A}^{k+t+m-3-i}_{m}}[\vec{x}]\bigg)\bigg)
 + \lambda^2 \tbinom{m-3}{i} \bigg(x \bigg(          L_{\mathcal{B}^{k+t+m-3-i}_{m}}[\vec{x}]\bigg)\bigg)\bigg) \\
& + \lambda^2 \tbinom{m-2}{i} \bigg(x \bigg(            L_{\mathcal{C}^{k+t+m-3-i}_{m}}[\vec{x}]\bigg)\bigg)\bigg)
\left. +\lambda^3 \tbinom{m-3}{i}\bigg(x \bigg(             L_{\mathcal{D}^{k+t+m-3-i}_{m}}[\vec{x}]\bigg)\bigg)\bigg) \right)]
 \\
& +        L^\op_{\mathcal{A}^{k+t+m-2-i}_{m+1}}[\vec{x},\sum\limits^l_{t=0}\sum\limits^{m-2}_{i=0}\lambda^{(t)}\lambda^i \left( 
 \tbinom{m-2}{i}R\bigg(  L^*_{\mathcal{B}^{k+t+m-3-i}_m}[x,\vec{x}]
 +    L^*_{\mathcal{A}^{k+t+m-3-i}_m}[x,\vec{x}]\bigg)\right. \\
& + \lambda \tbinom{m-3}{i} \bigg(  L^*_{\mathcal{B}^{k+t+m-3-i}_m}[x,\vec{x}]\bigg) 
 + \lambda \tbinom{m-2}{i} \bigg(    L^*_{\mathcal{C}^{k+t+m-3-i}_m}[x,\vec{x}]\bigg) 
 + \lambda^2 \tbinom{m-3}{i} \bigg(   L^*_{\mathcal{D}^{k+t+m-4-i}_m}[x,\vec{x}]\bigg)\\
& + \tbinom{m-2}{i}R\bigg(xR\bigg(          L_{\mathcal{B}^{k+t+m-3-i}_{m}}[\vec{x}]
 +        L_{\mathcal{A}^{k+t+m-3-i}_{m}}[\vec{x}]\bigg)\bigg) 
+ \lambda \tbinom{m-2}{i} R\bigg(x\bigg(          L_{\mathcal{B}^{k+t+m-3-i}_{m}}[\vec{x}]
 +            L_{\mathcal{A}^{k+t+m-3-i}_{m}}[\vec{x}]\bigg)\bigg) \allowdisplaybreaks \\
& + \lambda \tbinom{m-2}{i}\bigg(x R\bigg(          L_{\mathcal{B}^{k+t+m-3-i}_{m}}[\vec{x}]
 +        L_{\mathcal{A}^{k+t+m-3-i}_{m}}[\vec{x}]\bigg)\bigg)
 + \lambda^2 \tbinom{m-3}{i} \bigg(x \bigg(          L_{\mathcal{B}^{k+t+m-3-i}_{m}}[\vec{x}]\bigg)\bigg)\bigg) 
 \\
& + \lambda^2 \tbinom{m-2}{i} \bigg(x \bigg(            L_{\mathcal{C}^{k+t+m-3-i}_{m}}[\vec{x}]\bigg)\bigg)\bigg)
\left. +\lambda^3 \tbinom{m-3}{i}\bigg(x \bigg(             L_{\mathcal{D}^{k+t+m-3-i}_{m}}[\vec{x}]\bigg)\bigg)\bigg) \right)]\bigg)x\bigg)
\\
 &+ \lambda^2 \tbinom{m-2}{i} R\bigg( \bigg(  L^\op_{\mathcal{B}^{k+t+m-2-i}_{m+1}}[\vec{x},\sum\limits^l_{t=0}\sum\limits^{m-2}_{i=0}\lambda^{(t)}\lambda^i \left( 
 \tbinom{m-2}{i}R\bigg(  L^*_{\mathcal{B}^{k+t+m-3-i}_m}[x,\vec{x}]
 +    L^*_{\mathcal{A}^{k+t+m-3-i}_m}[x,\vec{x}]\bigg)\right. \\
& + \lambda \tbinom{m-3}{i} \bigg(  L^*_{\mathcal{B}^{k+t+m-3-i}_m}[x,\vec{x}]\bigg) 
 + \lambda \tbinom{m-2}{i} \bigg(    L^*_{\mathcal{C}^{k+t+m-3-i}_m}[x,\vec{x}]\bigg) 
 + \lambda^2 \tbinom{m-3}{i} \bigg(   L^*_{\mathcal{D}^{k+t+m-4-i}_m}[x,\vec{x}]\bigg)\\
 &+ \tbinom{m-2}{i}R\bigg(xR\bigg(          L_{\mathcal{B}^{k+t+m-3-i}_{m}}[\vec{x}]
 +        L_{\mathcal{A}^{k+t+m-3-i}_{m}}[\vec{x}]\bigg)\bigg) 
+ \lambda \tbinom{m-2}{i} R\bigg(x\bigg(          L_{\mathcal{B}^{k+t+m-3-i}_{m}}[\vec{x}]
 +            L_{\mathcal{A}^{k+t+m-3-i}_{m}}[\vec{x}]\bigg)\bigg) \allowdisplaybreaks \\
& + \lambda \tbinom{m-2}{i}\bigg(x R\bigg(          L_{\mathcal{B}^{k+t+m-3-i}_{m}}[\vec{x}]
 +        L_{\mathcal{A}^{k+t+m-3-i}_{m}}[\vec{x}]\bigg)\bigg)
 + \lambda^2 \tbinom{m-3}{i} \bigg(x \bigg(          L_{\mathcal{B}^{k+t+m-3-i}_{m}}[\vec{x}]\bigg)\bigg)\bigg) \\
& + \lambda^2 \tbinom{m-2}{i} \bigg(x \bigg(            L_{\mathcal{C}^{k+t+m-3-i}_{m}}[\vec{x}]\bigg)\bigg)\bigg)
\left. +\lambda^3 \tbinom{m-3}{i}\bigg(x \bigg(             L_{\mathcal{D}^{k+t+m-3-i}_{m}}[\vec{x}]\bigg)\bigg)\bigg) \right)]\bigg)x\bigg) 
 \\
& + \lambda^2 \tbinom{m-1}{i}R\bigg( \bigg(  L^\op_{\mathcal{C}^{k+t+m-2-i}_{m+1}}[\vec{x},\sum\limits^l_{t=0}\sum\limits^{m-2}_{i=0}\lambda^{(t)}\lambda^i \left( 
 \tbinom{m-2}{i}R\bigg(  L^*_{\mathcal{B}^{k+t+m-3-i}_m}[x,\vec{x}]
 +    L^*_{\mathcal{A}^{k+t+m-3-i}_m}[x,\vec{x}]\bigg)\right. \\
& + \lambda \tbinom{m-3}{i} \bigg(  L^*_{\mathcal{B}^{k+t+m-3-i}_m}[x,\vec{x}]\bigg) 
 + \lambda \tbinom{m-2}{i} \bigg(    L^*_{\mathcal{C}^{k+t+m-3-i}_m}[x,\vec{x}]\bigg) 
 + \lambda^2 \tbinom{m-3}{i} \bigg(   L^*_{\mathcal{D}^{k+t+m-4-i}_m}[x,\vec{x}]\bigg)\\
& + \tbinom{m-2}{i}R\bigg(xR\bigg(          L_{\mathcal{B}^{k+t+m-3-i}_{m}}[\vec{x}]
 +        L_{\mathcal{A}^{k+t+m-3-i}_{m}}[\vec{x}]\bigg)\bigg) + \lambda \tbinom{m-2}{i} R\bigg(x\bigg(          L_{\mathcal{B}^{k+t+m-3-i}_{m}}[\vec{x}]
 +            L_{\mathcal{A}^{k+t+m-3-i}_{m}}[\vec{x}]\bigg)\bigg) \allowdisplaybreaks \\
 \end{aligned}
 \]
 \[
 \begin{aligned}
& + \lambda \tbinom{m-2}{i}\bigg(x R\bigg(          L_{\mathcal{B}^{k+t+m-3-i}_{m}}[\vec{x}]
 +        L_{\mathcal{A}^{k+t+m-3-i}_{m}}[\vec{x}]\bigg)\bigg)
 + \lambda^2 \tbinom{m-3}{i} \bigg(x \bigg(          L_{\mathcal{B}^{k+t+m-3-i}_{m}}[\vec{x}]\bigg)\bigg)\bigg) \\
& + \lambda^2 \tbinom{m-2}{i} \bigg(x \bigg(            L_{\mathcal{C}^{k+t+m-3-i}_{m}}[\vec{x}]\bigg)\bigg)\bigg)
\left. +\lambda^3 \tbinom{m-3}{i}\bigg(x \bigg(             L_{\mathcal{D}^{k+t+m-3-i}_{m}}[\vec{x}]\bigg)\bigg)\bigg) \right)]\bigg)x\bigg)
\\
  &+\lambda^3 \tbinom{m-2}{i}R\bigg( \bigg(  L^\op_{\mathcal{D}^{k+t+m-2-i}_{m+1}}[\vec{x},\sum\limits^l_{t=0}\sum\limits^{m-2}_{i=0}\lambda^{(t)}\lambda^i \left( 
 \tbinom{m-2}{i}R\bigg(  L^*_{\mathcal{B}^{k+t+m-3-i}_m}[x,\vec{x}]
 +    L^*_{\mathcal{A}^{k+t+m-3-i}_m}[x,\vec{x}]\bigg)\right. \\
& + \lambda \tbinom{m-3}{i} \bigg(  L^*_{\mathcal{B}^{k+t+m-3-i}_m}[x,\vec{x}]\bigg) 
 + \lambda \tbinom{m-2}{i} \bigg(    L^*_{\mathcal{C}^{k+t+m-3-i}_m}[x,\vec{x}]\bigg) 
 + \lambda^2 \tbinom{m-3}{i} \bigg(   L^*_{\mathcal{D}^{k+t+m-4-i}_m}[x,\vec{x}]\bigg)\\
& + \tbinom{m-2}{i}R\bigg(xR\bigg(          L_{\mathcal{B}^{k+t+m-3-i}_{m}}[\vec{x}]
 +        L_{\mathcal{A}^{k+t+m-3-i}_{m}}[\vec{x}]\bigg)\bigg) 
 + \lambda \tbinom{m-2}{i} R\bigg(x\bigg(          L_{\mathcal{B}^{k+t+m-3-i}_{m}}[\vec{x}]
 +            L_{\mathcal{A}^{k+t+m-3-i}_{m}}[\vec{x}]\bigg)\bigg) \allowdisplaybreaks \\
& + \lambda \tbinom{m-2}{i}\bigg(x R\bigg(          L_{\mathcal{B}^{k+t+m-3-i}_{m}}[\vec{x}]
 +        L_{\mathcal{A}^{k+t+m-3-i}_{m}}[\vec{x}]\bigg)\bigg)
 + \lambda^2 \tbinom{m-3}{i} \bigg(x \bigg(          L_{\mathcal{B}^{k+t+m-3-i}_{m}}[\vec{x}]\bigg)\bigg)\bigg) \\
& + \lambda^2 \tbinom{m-2}{i} \bigg(x \bigg(            L_{\mathcal{C}^{k+t+m-3-i}_{m}}[\vec{x}]\bigg)\bigg)\bigg)
\left. +\lambda^3 \tbinom{m-3}{i}\bigg(x \bigg(             L_{\mathcal{D}^{k+t+m-3-i}_{m}}[\vec{x}]\bigg)\bigg)\bigg) \right)]\bigg)x\bigg) ) {=} 0;\ m = |\vec{x}|{\geq}3.
\end{aligned}
\]
Composition between relations $(R10)$ and $(R9)$, we get a sum of the expressions, which are consequences of~$(R10)$ itself.
Below, we show this example for  $\sum\limits_{t=0}^l \sum\limits^{m-1}_{i=0}\lambda^{(t)}\lambda^i
M_{t,i}$, where
\[
\begin{aligned}
M_{t,i}&=\sum\limits^l_{t=0}\sum\limits^{m-1}_{i=0}\lambda^{(t)}\lambda^i (
\tbinom{m-1}{i}R^2\bigg( {L^*}^\op_{\mathcal{B}^{k+t+m-2-i}_{m+1}}[x,\vec{x},\sum\limits^l_{t=0}\lambda^{(t)}\bigg( R \bigg(R\bigg(  L^\op_{\mathcal{B}^{k+t-1}_2}[\vec{y}]+              L^\op_{\mathcal{A}^{n-1}_2}[\vec{y}]\bigg)x\bigg)\\
& +\lambda R\bigg( \bigg(  L^\op_{\mathcal{B}^{k+t-1}_2}[\vec{y}]+              L^\op_{\mathcal{A}^{n-1}_2}[\vec{y}]\bigg)x\bigg)
 +\lambda \bigg( R\bigg(  L^\op_{\mathcal{B}^{k+t-1}_2}[\vec{y}]
 +              L^\op_{\mathcal{A}^{n-1}_2}[\vec{y}]\bigg)x\bigg)\\
 &+ R\bigg(   {L^*}^\op_{\mathcal{B}^{k+t-1}_3}[x,\vec{y}]+    {L^*}^\op_{\mathcal{A}^{n-1}_3}[x,\vec{y}]\bigg)
 + \lambda \bigg(   {L^*}^\op_{\mathcal{B}^{k+t-1}_3}[x,\vec{y}]+   {L^*}^\op_{\mathcal{C}^{k+t-1}_3}[x,\vec{y}]\bigg)\\
& + \lambda^2 \bigg(\bigg(  L^\op_{\mathcal{B}^{k+t-1}_2}[\vec{y}]+      L^\op_{\mathcal{C}^{k+t-1}_2}[\vec{y}]\bigg)x\bigg)\bigg)]
+ {L^*}^\op_{\mathcal{A}^{k+t+m-2-i}_{m+1}}[x,\vec{x},\sum\limits^l_{t=0}\lambda^{(t)}\bigg( R \bigg(R\bigg(  L^\op_{\mathcal{B}^{k+t-1}_2}[\vec{y}]+              L^\op_{\mathcal{A}^{n-1}_2}[\vec{y}]\bigg)x\bigg)\\
& +\lambda R\bigg( \bigg(  L^\op_{\mathcal{B}^{k+t-1}_2}[\vec{y}]+              L^\op_{\mathcal{A}^{n-1}_2}[\vec{y}]\bigg)x\bigg)
 +\lambda \bigg( R\bigg(  L^\op_{\mathcal{B}^{k+t-1}_2}[\vec{y}]
 +              L^\op_{\mathcal{A}^{n-1}_2}[\vec{y}]\bigg)x\bigg)\\
 &+ R\bigg(   {L^*}^\op_{\mathcal{B}^{k+t-1}_3}[x,\vec{y}]+    {L^*}^\op_{\mathcal{A}^{n-1}_3}[x,\vec{y}]\bigg)
 + \lambda \bigg(   {L^*}^\op_{\mathcal{B}^{k+t-1}_3}[x,\vec{y}]+   {L^*}^\op_{\mathcal{C}^{k+t-1}_3}[x,\vec{y}]\bigg)\\
& + \lambda^2 \bigg(\bigg(  L^\op_{\mathcal{B}^{k+t-1}_2}[\vec{y}]+      L^\op_{\mathcal{C}^{k+t-1}_2}[\vec{y}]\bigg)x\bigg)\bigg)]\bigg) 
+\lambda \tbinom{m-2}{i}R\bigg( {L^*}^\op_{\mathcal{B}^{k+t+m-2-i}_{m+1}}[x,\vec{x},\sum\limits^l_{t=0}\lambda^{(t)}\bigg( R \bigg(R\bigg(  L^\op_{\mathcal{B}^{k+t-1}_2}[\vec{y}]+              L^\op_{\mathcal{A}^{n-1}_2}[\vec{y}]\bigg)x\bigg)\\
& +\lambda R\bigg( \bigg(  L^\op_{\mathcal{B}^{k+t-1}_2}[\vec{y}]+              L^\op_{\mathcal{A}^{n-1}_2}[\vec{y}]\bigg)x\bigg)
 +\lambda \bigg( R\bigg(  L^\op_{\mathcal{B}^{k+t-1}_2}[\vec{y}]
 +              L^\op_{\mathcal{A}^{n-1}_2}[\vec{y}]\bigg)x\bigg)\\
 &+ R\bigg(   {L^*}^\op_{\mathcal{B}^{k+t-1}_3}[x,\vec{y}]+    {L^*}^\op_{\mathcal{A}^{n-1}_3}[x,\vec{y}]\bigg)
 + \lambda \bigg(   {L^*}^\op_{\mathcal{B}^{k+t-1}_3}[x,\vec{y}]+   {L^*}^\op_{\mathcal{C}^{k+t-1}_3}[x,\vec{y}]\bigg)\\
& + \lambda^2 \bigg(\bigg(  L^\op_{\mathcal{B}^{k+t-1}_2}[\vec{y}]+      L^\op_{\mathcal{C}^{k+t-1}_2}[\vec{y}]\bigg)x\bigg)\bigg)]\bigg)
+\lambda \tbinom{m-1}{i} R\bigg( {L^*}^\op_{\mathcal{C}^{k+t+m-2-i}_{m+1}}[x,\vec{x},\sum\limits^l_{t=0}\lambda^{(t)}\bigg( R \bigg(R\bigg(  L^\op_{\mathcal{B}^{k+t-1}_2}[\vec{y}]+              L^\op_{\mathcal{A}^{n-1}_2}[\vec{y}]\bigg)x\bigg)\\
& +\lambda R\bigg( \bigg(  L^\op_{\mathcal{B}^{k+t-1}_2}[\vec{y}]+              L^\op_{\mathcal{A}^{n-1}_2}[\vec{y}]\bigg)x\bigg)
 +\lambda \bigg( R\bigg(  L^\op_{\mathcal{B}^{k+t-1}_2}[\vec{y}]
 +              L^\op_{\mathcal{A}^{n-1}_2}[\vec{y}]\bigg)x\bigg)\\
 &+ R\bigg(   {L^*}^\op_{\mathcal{B}^{k+t-1}_3}[x,\vec{y}]+    {L^*}^\op_{\mathcal{A}^{n-1}_3}[x,\vec{y}]\bigg)
 + \lambda \bigg(   {L^*}^\op_{\mathcal{B}^{k+t-1}_3}[x,\vec{y}]+   {L^*}^\op_{\mathcal{C}^{k+t-1}_3}[x,\vec{y}]\bigg)\\
 \end{aligned}
 \]
 \[
 \begin{aligned}
& + \lambda^2 \bigg(\bigg(  L^\op_{\mathcal{B}^{k+t-1}_2}[\vec{y}]+      L^\op_{\mathcal{C}^{k+t-1}_2}[\vec{y}]\bigg)x\bigg)\bigg)]\bigg)\bigg)
+\lambda^2 \tbinom{m-2}{i}R\bigg(  {L^*}^\op_{\mathcal{D}^{k+t+m-3-i}_{m+1}}[x,\vec{x},\sum\limits^l_{t=0}\lambda^{(t)}\bigg( R \bigg(R\bigg(  L^\op_{\mathcal{B}^{k+t-1}_2}[\vec{y}]+              L^\op_{\mathcal{A}^{n-1}_2}[\vec{y}]\bigg)x\bigg)\\
& +\lambda R\bigg( \bigg(  L^\op_{\mathcal{B}^{k+t-1}_2}[\vec{y}]+              L^\op_{\mathcal{A}^{n-1}_2}[\vec{y}]\bigg)x\bigg)
 +\lambda \bigg( R\bigg(  L^\op_{\mathcal{B}^{k+t-1}_2}[\vec{y}]
 +              L^\op_{\mathcal{A}^{n-1}_2}[\vec{y}]\bigg)x\bigg)\\
 &+ R\bigg(   {L^*}^\op_{\mathcal{B}^{k+t-1}_3}[x,\vec{y}]+    {L^*}^\op_{\mathcal{A}^{n-1}_3}[x,\vec{y}]\bigg)
 + \lambda \bigg(   {L^*}^\op_{\mathcal{B}^{k+t-1}_3}[x,\vec{y}]+   {L^*}^\op_{\mathcal{C}^{k+t-1}_3}[x,\vec{y}]\bigg) + \lambda^2 \bigg(\bigg(  L^\op_{\mathcal{B}^{k+t-1}_2}[\vec{y}]+      L^\op_{\mathcal{C}^{k+t-1}_2}[\vec{y}]\bigg)x\bigg)\bigg)]\bigg)
\\
 &+ \tbinom{m-1}{i}R^2\bigg(R\bigg(  L^\op_{\mathcal{B}^{k+t+m-2-i}_{m+1}}[\vec{x},\sum\limits^l_{t=0}\lambda^{(t)}\bigg( R \bigg(R\bigg(  L^\op_{\mathcal{B}^{k+t-1}_2}[\vec{y}]+              L^\op_{\mathcal{A}^{n-1}_2}[\vec{y}]\bigg)x\bigg)
 +\lambda R\bigg( \bigg(  L^\op_{\mathcal{B}^{k+t-1}_2}[\vec{y}]+              L^\op_{\mathcal{A}^{n-1}_2}[\vec{y}]\bigg)x\bigg)\\
 &
 +\lambda \bigg( R\bigg(  L^\op_{\mathcal{B}^{k+t-1}_2}[\vec{y}]
 +              L^\op_{\mathcal{A}^{n-1}_2}[\vec{y}]\bigg)x\bigg)
 + R\bigg(   {L^*}^\op_{\mathcal{B}^{k+t-1}_3}[x,\vec{y}]+    {L^*}^\op_{\mathcal{A}^{n-1}_3}[x,\vec{y}]\bigg)
 + \lambda \bigg(   {L^*}^\op_{\mathcal{B}^{k+t-1}_3}[x,\vec{y}]+   {L^*}^\op_{\mathcal{C}^{k+t-1}_3}[x,\vec{y}]\bigg)\\
& + \lambda^2 \bigg(\bigg(  L^\op_{\mathcal{B}^{k+t-1}_2}[\vec{y}]+      L^\op_{\mathcal{C}^{k+t-1}_2}[\vec{y}]\bigg)x\bigg)\bigg)]
 +        L^\op_{\mathcal{A}^{k+t+m-2-i}_{m+1}}[\vec{x},\sum\limits^l_{t=0}\lambda^{(t)}\bigg( R \bigg(R\bigg(  L^\op_{\mathcal{B}^{k+t-1}_2}[\vec{y}]+              L^\op_{\mathcal{A}^{n-1}_2}[\vec{y}]\bigg)x\bigg)\\
& +\lambda R\bigg( \bigg(  L^\op_{\mathcal{B}^{k+t-1}_2}[\vec{y}]+              L^\op_{\mathcal{A}^{n-1}_2}[\vec{y}]\bigg)x\bigg)
 +\lambda \bigg( R\bigg(  L^\op_{\mathcal{B}^{k+t-1}_2}[\vec{y}]
 +              L^\op_{\mathcal{A}^{n-1}_2}[\vec{y}]\bigg)x\bigg)
 + R\bigg(   {L^*}^\op_{\mathcal{B}^{k+t-1}_3}[x,\vec{y}]+    {L^*}^\op_{\mathcal{A}^{n-1}_3}[x,\vec{y}]\bigg)\\
 &+ \lambda \bigg(   {L^*}^\op_{\mathcal{B}^{k+t-1}_3}[x,\vec{y}]+   {L^*}^\op_{\mathcal{C}^{k+t-1}_3}[x,\vec{y}]\bigg)
 + \lambda^2 \bigg(\bigg(  L^\op_{\mathcal{B}^{k+t-1}_2}[\vec{y}]+      L^\op_{\mathcal{C}^{k+t-1}_2}[\vec{y}]\bigg)x\bigg)\bigg)]\bigg)x\bigg)
\\
 &+ \lambda \tbinom{m-1}{i}R^2\bigg(\bigg(  L^\op_{\mathcal{B}^{k+t+m-2-i}_{m+1}}[\vec{x},\sum\limits^l_{t=0}\lambda^{(t)}\bigg( R \bigg(R\bigg(  L^\op_{\mathcal{B}^{k+t-1}_2}[\vec{y}]+              L^\op_{\mathcal{A}^{n-1}_2}[\vec{y}]\bigg)x\bigg) +\lambda R\bigg( \bigg(  L^\op_{\mathcal{B}^{k+t-1}_2}[\vec{y}]+              L^\op_{\mathcal{A}^{n-1}_2}[\vec{y}]\bigg)x\bigg)\\
& +\lambda \bigg( R\bigg(  L^\op_{\mathcal{B}^{k+t-1}_2}[\vec{y}]
 +              L^\op_{\mathcal{A}^{n-1}_2}[\vec{y}]\bigg)x\bigg)
 + R\bigg(   {L^*}^\op_{\mathcal{B}^{k+t-1}_3}[x,\vec{y}]+    {L^*}^\op_{\mathcal{A}^{n-1}_3}[x,\vec{y}]\bigg)
 + \lambda \bigg(   {L^*}^\op_{\mathcal{B}^{k+t-1}_3}[x,\vec{y}]+   {L^*}^\op_{\mathcal{C}^{k+t-1}_3}[x,\vec{y}]\bigg)\\
& + \lambda^2 \bigg(\bigg(  L^\op_{\mathcal{B}^{k+t-1}_2}[\vec{y}]+      L^\op_{\mathcal{C}^{k+t-1}_2}[\vec{y}]\bigg)x\bigg)\bigg)]
  +        L^\op_{\mathcal{A}^{k+t+m-2-i}_{m+1}}[\vec{x},\sum\limits^l_{t=0}\lambda^{(t)}\bigg( R \bigg(R\bigg(  L^\op_{\mathcal{B}^{k+t-1}_2}[\vec{y}]+              L^\op_{\mathcal{A}^{n-1}_2}[\vec{y}]\bigg)x\bigg)\\
& +\lambda R\bigg( \bigg(  L^\op_{\mathcal{B}^{k+t-1}_2}[\vec{y}]+              L^\op_{\mathcal{A}^{n-1}_2}[\vec{y}]\bigg)x\bigg)
 +\lambda \bigg( R\bigg(  L^\op_{\mathcal{B}^{k+t-1}_2}[\vec{y}]
 +              L^\op_{\mathcal{A}^{n-1}_2}[\vec{y}]\bigg)x\bigg)\\
 &+ R\bigg(   {L^*}^\op_{\mathcal{B}^{k+t-1}_3}[x,\vec{y}]+    {L^*}^\op_{\mathcal{A}^{n-1}_3}[x,\vec{y}]\bigg)
 + \lambda \bigg(   {L^*}^\op_{\mathcal{B}^{k+t-1}_3}[x,\vec{y}]+   {L^*}^\op_{\mathcal{C}^{k+t-1}_3}[x,\vec{y}]\bigg)\\
& + \lambda^2 \bigg(\bigg(  L^\op_{\mathcal{B}^{k+t-1}_2}[\vec{y}]+      L^\op_{\mathcal{C}^{k+t-1}_2}[\vec{y}]\bigg)x\bigg)\bigg)]\bigg)x\bigg)
 + \lambda \tbinom{m-1}{i}R\bigg( R\bigg(  L^\op_{\mathcal{B}^{k+t+m-2-i}_{m+1}}[\vec{x},\sum\limits^l_{t=0}\lambda^{(t)}\bigg( R \bigg(R\bigg(  L^\op_{\mathcal{B}^{k+t-1}_2}[\vec{y}]+              L^\op_{\mathcal{A}^{n-1}_2}[\vec{y}]\bigg)x\bigg)\\
& +\lambda R\bigg( \bigg(  L^\op_{\mathcal{B}^{k+t-1}_2}[\vec{y}]+              L^\op_{\mathcal{A}^{n-1}_2}[\vec{y}]\bigg)x\bigg)
 +\lambda \bigg( R\bigg(  L^\op_{\mathcal{B}^{k+t-1}_2}[\vec{y}]
 +              L^\op_{\mathcal{A}^{n-1}_2}[\vec{y}]\bigg)x\bigg)\\
 &+ R\bigg(   {L^*}^\op_{\mathcal{B}^{k+t-1}_3}[x,\vec{y}]+    {L^*}^\op_{\mathcal{A}^{n-1}_3}[x,\vec{y}]\bigg)
 + \lambda \bigg(   {L^*}^\op_{\mathcal{B}^{k+t-1}_3}[x,\vec{y}]+   {L^*}^\op_{\mathcal{C}^{k+t-1}_3}[x,\vec{y}]\bigg)
 + \lambda^2 \bigg(\bigg(  L^\op_{\mathcal{B}^{k+t-1}_2}[\vec{y}]+      L^\op_{\mathcal{C}^{k+t-1}_2}[\vec{y}]\bigg)x\bigg)\bigg)]
 \\
& +        L^\op_{\mathcal{A}^{k+t+m-2-i}_{m+1}}[\vec{x},\sum\limits^l_{t=0}\lambda^{(t)}\bigg( R \bigg(R\bigg(  L^\op_{\mathcal{B}^{k+t-1}_2}[\vec{y}]+              L^\op_{\mathcal{A}^{n-1}_2}[\vec{y}]\bigg)x\bigg) +\lambda R\bigg( \bigg(  L^\op_{\mathcal{B}^{k+t-1}_2}[\vec{y}]+              L^\op_{\mathcal{A}^{n-1}_2}[\vec{y}]\bigg)x\bigg)\\
 &+\lambda \bigg( R\bigg(  L^\op_{\mathcal{B}^{k+t-1}_2}[\vec{y}]
 +              L^\op_{\mathcal{A}^{n-1}_2}[\vec{y}]\bigg)x\bigg)
 + R\bigg(   {L^*}^\op_{\mathcal{B}^{k+t-1}_3}[x,\vec{y}]+    {L^*}^\op_{\mathcal{A}^{n-1}_3}[x,\vec{y}]\bigg)
 + \lambda \bigg(   {L^*}^\op_{\mathcal{B}^{k+t-1}_3}[x,\vec{y}]+   {L^*}^\op_{\mathcal{C}^{k+t-1}_3}[x,\vec{y}]\bigg)\\
& + \lambda^2 \bigg(\bigg(  L^\op_{\mathcal{B}^{k+t-1}_2}[\vec{y}]+      L^\op_{\mathcal{C}^{k+t-1}_2}[\vec{y}]\bigg)x\bigg)\bigg)]\bigg)x\bigg)
+ \lambda^2 \tbinom{m-2}{i} R\bigg( \bigg(  L^\op_{\mathcal{B}^{k+t+m-2-i}_{m+1}}[\vec{x},\sum\limits^l_{t=0}\lambda^{(t)}\bigg( R \bigg(R\bigg(  L^\op_{\mathcal{B}^{k+t-1}_2}[\vec{y}]+              L^\op_{\mathcal{A}^{n-1}_2}[\vec{y}]\bigg)x\bigg)\\
& +\lambda R\bigg( \bigg(  L^\op_{\mathcal{B}^{k+t-1}_2}[\vec{y}]+              L^\op_{\mathcal{A}^{n-1}_2}[\vec{y}]\bigg)x\bigg)
 +\lambda \bigg( R\bigg(  L^\op_{\mathcal{B}^{k+t-1}_2}[\vec{y}]
 +              L^\op_{\mathcal{A}^{n-1}_2}[\vec{y}]\bigg)x\bigg)+ R\bigg(   {L^*}^\op_{\mathcal{B}^{k+t-1}_3}[x,\vec{y}]+    {L^*}^\op_{\mathcal{A}^{n-1}_3}[x,\vec{y}]\bigg)\\
 \end{aligned}
 \]
 \[
 \begin{aligned}
& + \lambda^2 \bigg(\bigg(  L^\op_{\mathcal{B}^{k+t-1}_2}[\vec{y}]+      L^\op_{\mathcal{C}^{k+t-1}_2}[\vec{y}]\bigg)x\bigg)\bigg)]\bigg)x\bigg) 
 + \lambda^2 \tbinom{m-1}{i}R\bigg( \bigg(  L^\op_{\mathcal{C}^{k+t+m-2-i}_{m+1}}[\vec{x},\sum\limits^l_{t=0}\lambda^{(t)}\bigg( R \bigg(R\bigg(  L^\op_{\mathcal{B}^{k+t-1}_2}[\vec{y}]+              L^\op_{\mathcal{A}^{n-1}_2}[\vec{y}]\bigg)x\bigg)\\
& +\lambda R\bigg( \bigg(  L^\op_{\mathcal{B}^{k+t-1}_2}[\vec{y}]+              L^\op_{\mathcal{A}^{n-1}_2}[\vec{y}]\bigg)x\bigg)
 +\lambda \bigg( R\bigg(  L^\op_{\mathcal{B}^{k+t-1}_2}[\vec{y}]
 +              L^\op_{\mathcal{A}^{n-1}_2}[\vec{y}]\bigg)x\bigg)+ \lambda \bigg(   {L^*}^\op_{\mathcal{B}^{k+t-1}_3}[x,\vec{y}]+   {L^*}^\op_{\mathcal{C}^{k+t-1}_3}[x,\vec{y}]\bigg)\\
 &+ R\bigg(   {L^*}^\op_{\mathcal{B}^{k+t-1}_3}[x,\vec{y}]+    {L^*}^\op_{\mathcal{A}^{n-1}_3}[x,\vec{y}]\bigg)
 + \lambda \bigg(   {L^*}^\op_{\mathcal{B}^{k+t-1}_3}[x,\vec{y}]+   {L^*}^\op_{\mathcal{C}^{k+t-1}_3}[x,\vec{y}]\bigg)\\
& + \lambda^2 \bigg(\bigg(  L^\op_{\mathcal{B}^{k+t-1}_2}[\vec{y}]+      L^\op_{\mathcal{C}^{k+t-1}_2}[\vec{y}]\bigg)x\bigg)\bigg)]\bigg)x\bigg)
+\lambda^3 \tbinom{m-2}{i}R\bigg( \bigg(  L^\op_{\mathcal{D}^{k+t+m-2-i}_{m+1}}[\vec{x},\sum\limits^l_{t=0}\lambda^{(t)}\bigg( R \bigg(R\bigg(  L^\op_{\mathcal{B}^{k+t-1}_2}[\vec{y}]+              L^\op_{\mathcal{A}^{n-1}_2}[\vec{y}]\bigg)x\bigg)\\
& +\lambda R\bigg( \bigg(  L^\op_{\mathcal{B}^{k+t-1}_2}[\vec{y}]+              L^\op_{\mathcal{A}^{n-1}_2}[\vec{y}]\bigg)x\bigg)
 +\lambda \bigg( R\bigg(  L^\op_{\mathcal{B}^{k+t-1}_2}[\vec{y}]
 +              L^\op_{\mathcal{A}^{n-1}_2}[\vec{y}]\bigg)x\bigg)\\
 &+ R\bigg(   {L^*}^\op_{\mathcal{B}^{k+t-1}_3}[x,\vec{y}]+    {L^*}^\op_{\mathcal{A}^{n-1}_3}[x,\vec{y}]\bigg)
 + \lambda \bigg(   {L^*}^\op_{\mathcal{B}^{k+t-1}_3}[x,\vec{y}]+   {L^*}^\op_{\mathcal{C}^{k+t-1}_3}[x,\vec{y}]\bigg)\\
& + \lambda^2 \bigg(\bigg(  L^\op_{\mathcal{B}^{k+t-1}_2}[\vec{y}]+      L^\op_{\mathcal{C}^{k+t-1}_2}[\vec{y}]\bigg)x\bigg)\bigg)]\bigg)x\bigg) ) {=} 0;\ m = |\vec{x}|{\geq}3,\ |\vec{y}|=2.
\end{aligned}
\]
Composition between relations $(R10)$ and $(R10)$, we get a sum of the expressions, which are consequences of~$(R10)$ itself.
Below, we show this example for  $\sum\limits_{t=0}^l \sum\limits^{m-1}_{i=0}\lambda^{(t)}\lambda^i
M_{t,i}$, where
\[
\begin{aligned}
M_{t,i}&=
\tbinom{m-1}{i}R^2\bigg( {L^*}^\op_{\mathcal{B}^{k+t+m-2-i}_{m+1}}[x,\vec{x},\sum\limits^l_{t=0}\sum\limits^{m-2}_{i=0}\lambda^{(t)}\lambda^i \left(
\tbinom{m-2}{i}R\bigg(  {L^*}^\op_{\mathcal{B}^{k+t+m-3-i}_m}[x,\vec{x}]
+   {L^*}^\op_{\mathcal{A}^{k+t+m-3-i}_m}[x,\vec{x}]\bigg) \right. 
\\
&+\lambda \tbinom{m-3}{i}\bigg(  {L^*}^\op_{\mathcal{B}^{k+t+m-3-i}_m}[x,\vec{x}]\bigg)
+\lambda \tbinom{m-2}{i} \bigg(    {L^*}^\op_{\mathcal{C}^{k+t+m-3-i}_m}[x,\vec{x}]\bigg)\bigg)+\lambda^2 \tbinom{m-3}{i}\bigg(   {L^*}^\op_{\mathcal{D}^{k+t+m-4-i}_m}[x,\vec{x}]\bigg)\\
&
 + \tbinom{m-2}{i}R\bigg(R\bigg(          L^\op_{\mathcal{B}^{k+t+m-3-i}_{m}}[\vec{x}]
 +            L^\op_{\mathcal{A}^{k+t+m-3-i}_{m}}[\vec{x}]\bigg)x\bigg) + \lambda \tbinom{m-2}{i}R\bigg(\bigg(          L^\op_{\mathcal{B}^{k+t+m-3-i}_{m}}[\vec{x}]
 +            L^\op_{\mathcal{A}^{k+t+m-3-i}_{m}}[\vec{x}]\bigg)x\bigg)\\
& + \lambda \tbinom{m-2}{i}\bigg( R\bigg(          L^\op_{\mathcal{B}^{k+t+m-3-i}_{m}}[\vec{x}]
 +            L^\op_{\mathcal{A}^{k+t+m-3-i}_{m}}[\vec{x}]\bigg)x\bigg)
 + \lambda^2 \tbinom{m-3}{i} \bigg( \bigg(          L^\op_{\mathcal{B}^{k+t+m-3-i}_{m}}[\vec{x}]\bigg)x\bigg) \\
& + \lambda^2 \tbinom{m-2}{i}\bigg( \bigg(            L^\op_{\mathcal{C}^{k+t+m-3-i}_{m}}[\vec{x}]\bigg)x\bigg)
 \left. +\lambda^3 \tbinom{m-3}{i}\bigg( \bigg(              L^\op_{\mathcal{D}^{k+t+m-3-i}_{m}}[\vec{x}]\bigg)x\bigg) \right)]
\\
&+ {L^*}^\op_{\mathcal{A}^{k+t+m-2-i}_{m+1}}[x,\vec{x},\sum\limits^l_{t=0}\sum\limits^{m-2}_{i=0}\lambda^{(t)}\lambda^i \left(
\tbinom{m-2}{i}R\bigg(  {L^*}^\op_{\mathcal{B}^{k+t+m-3-i}_m}[x,\vec{x}]
+   {L^*}^\op_{\mathcal{A}^{k+t+m-3-i}_m}[x,\vec{x}]\bigg) \right. \\
&+\lambda \tbinom{m-3}{i}\bigg(  {L^*}^\op_{\mathcal{B}^{k+t+m-3-i}_m}[x,\vec{x}]\bigg)
+\lambda \tbinom{m-2}{i} \bigg(    {L^*}^\op_{\mathcal{C}^{k+t+m-3-i}_m}[x,\vec{x}]\bigg)\bigg)+\lambda^2 \tbinom{m-3}{i}\bigg(   {L^*}^\op_{\mathcal{D}^{k+t+m-4-i}_m}[x,\vec{x}]\bigg)\\
 &+ \tbinom{m-2}{i}R\bigg(R\bigg(          L^\op_{\mathcal{B}^{k+t+m-3-i}_{m}}[\vec{x}]
 +            L^\op_{\mathcal{A}^{k+t+m-3-i}_{m}}[\vec{x}]\bigg)x\bigg)+ \lambda \tbinom{m-2}{i}R\bigg(\bigg(          L^\op_{\mathcal{B}^{k+t+m-3-i}_{m}}[\vec{x}]
 +            L^\op_{\mathcal{A}^{k+t+m-3-i}_{m}}[\vec{x}]\bigg)x\bigg)\\
& + \lambda \tbinom{m-2}{i}\bigg( R\bigg(          L^\op_{\mathcal{B}^{k+t+m-3-i}_{m}}[\vec{x}]
 +            L^\op_{\mathcal{A}^{k+t+m-3-i}_{m}}[\vec{x}]\bigg)x\bigg)
 + \lambda^2 \tbinom{m-3}{i} \bigg( \bigg(          L^\op_{\mathcal{B}^{k+t+m-3-i}_{m}}[\vec{x}]\bigg)x\bigg) \\
& + \lambda^2 \tbinom{m-2}{i}\bigg( \bigg(            L^\op_{\mathcal{C}^{k+t+m-3-i}_{m}}[\vec{x}]\bigg)x\bigg)
 \left. +\lambda^3 \tbinom{m-3}{i}\bigg( \bigg(              L^\op_{\mathcal{D}^{k+t+m-3-i}_{m}}[\vec{x}]\bigg)x\bigg) \right)]\bigg) 
\\
&+\lambda \tbinom{m-2}{i}R\bigg( {L^*}^\op_{\mathcal{B}^{k+t+m-2-i}_{m+1}}[x,\vec{x},\sum\limits^l_{t=0}\sum\limits^{m-2}_{i=0}\lambda^{(t)}\lambda^i \left(
\tbinom{m-2}{i}R\bigg(  {L^*}^\op_{\mathcal{B}^{k+t+m-3-i}_m}[x,\vec{x}]
+   {L^*}^\op_{\mathcal{A}^{k+t+m-3-i}_m}[x,\vec{x}]\bigg) \right.\\ 
&+\lambda \tbinom{m-3}{i}\bigg(  {L^*}^\op_{\mathcal{B}^{k+t+m-3-i}_m}[x,\vec{x}]\bigg)
+\lambda \tbinom{m-2}{i} \bigg(    {L^*}^\op_{\mathcal{C}^{k+t+m-3-i}_m}[x,\vec{x}]\bigg)\bigg)+\lambda^2 \tbinom{m-3}{i}\bigg(   {L^*}^\op_{\mathcal{D}^{k+t+m-4-i}_m}[x,\vec{x}]\bigg)
\end{aligned}
 \]
 \[
 \begin{aligned}
&
 + \tbinom{m-2}{i}R\bigg(R\bigg(          L^\op_{\mathcal{B}^{k+t+m-3-i}_{m}}[\vec{x}]
 +            L^\op_{\mathcal{A}^{k+t+m-3-i}_{m}}[\vec{x}]\bigg)x\bigg)+ \lambda \tbinom{m-2}{i}R\bigg(\bigg(          L^\op_{\mathcal{B}^{k+t+m-3-i}_{m}}[\vec{x}]
 +            L^\op_{\mathcal{A}^{k+t+m-3-i}_{m}}[\vec{x}]\bigg)x\bigg)\\
& + \lambda \tbinom{m-2}{i}\bigg( R\bigg(          L^\op_{\mathcal{B}^{k+t+m-3-i}_{m}}[\vec{x}]
 +            L^\op_{\mathcal{A}^{k+t+m-3-i}_{m}}[\vec{x}]\bigg)x\bigg)
 + \lambda^2 \tbinom{m-3}{i} \bigg( \bigg(          L^\op_{\mathcal{B}^{k+t+m-3-i}_{m}}[\vec{x}]\bigg)x\bigg) \\
& + \lambda^2 \tbinom{m-2}{i}\bigg( \bigg(            L^\op_{\mathcal{C}^{k+t+m-3-i}_{m}}[\vec{x}]\bigg)x\bigg)
 \left. +\lambda^3 \tbinom{m-3}{i}\bigg( \bigg(              L^\op_{\mathcal{D}^{k+t+m-3-i}_{m}}[\vec{x}]\bigg)x\bigg) \right)]\bigg)
\\
&+\lambda \tbinom{m-1}{i} R\bigg( {L^*}^\op_{\mathcal{C}^{k+t+m-2-i}_{m+1}}[x,\vec{x},\sum\limits^l_{t=0}\sum\limits^{m-2}_{i=0}\lambda^{(t)}\lambda^i \left(
\tbinom{m-2}{i}R\bigg(  {L^*}^\op_{\mathcal{B}^{k+t+m-3-i}_m}[x,\vec{x}]
+   {L^*}^\op_{\mathcal{A}^{k+t+m-3-i}_m}[x,\vec{x}]\bigg) \right. \\
&+\lambda \tbinom{m-3}{i}\bigg(  {L^*}^\op_{\mathcal{B}^{k+t+m-3-i}_m}[x,\vec{x}]\bigg)
+\lambda \tbinom{m-2}{i} \bigg(    {L^*}^\op_{\mathcal{C}^{k+t+m-3-i}_m}[x,\vec{x}]\bigg)\bigg)+\lambda^2 \tbinom{m-3}{i}\bigg(   {L^*}^\op_{\mathcal{D}^{k+t+m-4-i}_m}[x,\vec{x}]\bigg)\\
& + \tbinom{m-2}{i}R\bigg(R\bigg(          L^\op_{\mathcal{B}^{k+t+m-3-i}_{m}}[\vec{x}]
 +            L^\op_{\mathcal{A}^{k+t+m-3-i}_{m}}[\vec{x}]\bigg)x\bigg)
 + \lambda \tbinom{m-2}{i}R\bigg(\bigg(          L^\op_{\mathcal{B}^{k+t+m-3-i}_{m}}[\vec{x}]
 +            L^\op_{\mathcal{A}^{k+t+m-3-i}_{m}}[\vec{x}]\bigg)x\bigg)\\
& + \lambda \tbinom{m-2}{i}\bigg( R\bigg(          L^\op_{\mathcal{B}^{k+t+m-3-i}_{m}}[\vec{x}]
 +            L^\op_{\mathcal{A}^{k+t+m-3-i}_{m}}[\vec{x}]\bigg)x\bigg)
 + \lambda^2 \tbinom{m-3}{i} \bigg( \bigg(          L^\op_{\mathcal{B}^{k+t+m-3-i}_{m}}[\vec{x}]\bigg)x\bigg) \\
& + \lambda^2 \tbinom{m-2}{i}\bigg( \bigg(            L^\op_{\mathcal{C}^{k+t+m-3-i}_{m}}[\vec{x}]\bigg)x\bigg)
 \left. +\lambda^3 \tbinom{m-3}{i}\bigg( \bigg(              L^\op_{\mathcal{D}^{k+t+m-3-i}_{m}}[\vec{x}]\bigg)x\bigg) \right)]\bigg)\bigg)
\\
&+\lambda^2 \tbinom{m-2}{i}R\bigg(  {L^*}^\op_{\mathcal{D}^{k+t+m-3-i}_{m+1}}[x,\vec{x},\sum\limits^l_{t=0}\sum\limits^{m-2}_{i=0}\lambda^{(t)}\lambda^i \left(
\tbinom{m-2}{i}R\bigg(  {L^*}^\op_{\mathcal{B}^{k+t+m-3-i}_m}[x,\vec{x}]
+   {L^*}^\op_{\mathcal{A}^{k+t+m-3-i}_m}[x,\vec{x}]\bigg) \right. \\
&+\lambda \tbinom{m-3}{i}\bigg(  {L^*}^\op_{\mathcal{B}^{k+t+m-3-i}_m}[x,\vec{x}]\bigg)
+\lambda \tbinom{m-2}{i} \bigg(    {L^*}^\op_{\mathcal{C}^{k+t+m-3-i}_m}[x,\vec{x}]\bigg)\bigg)+\lambda^2 \tbinom{m-3}{i}\bigg(   {L^*}^\op_{\mathcal{D}^{k+t+m-4-i}_m}[x,\vec{x}]\bigg)\\
& + \tbinom{m-2}{i}R\bigg(R\bigg(          L^\op_{\mathcal{B}^{k+t+m-3-i}_{m}}[\vec{x}]
 +            L^\op_{\mathcal{A}^{k+t+m-3-i}_{m}}[\vec{x}]\bigg)x\bigg) + \lambda \tbinom{m-2}{i}R\bigg(\bigg(          L^\op_{\mathcal{B}^{k+t+m-3-i}_{m}}[\vec{x}]
 +            L^\op_{\mathcal{A}^{k+t+m-3-i}_{m}}[\vec{x}]\bigg)x\bigg)\\
& + \lambda \tbinom{m-2}{i}\bigg( R\bigg(          L^\op_{\mathcal{B}^{k+t+m-3-i}_{m}}[\vec{x}]
 +            L^\op_{\mathcal{A}^{k+t+m-3-i}_{m}}[\vec{x}]\bigg)x\bigg)
 + \lambda^2 \tbinom{m-3}{i} \bigg( \bigg(          L^\op_{\mathcal{B}^{k+t+m-3-i}_{m}}[\vec{x}]\bigg)x\bigg) \\
& + \lambda^2 \tbinom{m-2}{i}\bigg( \bigg(            L^\op_{\mathcal{C}^{k+t+m-3-i}_{m}}[\vec{x}]\bigg)x\bigg)
 \left. +\lambda^3 \tbinom{m-3}{i}\bigg( \bigg(              L^\op_{\mathcal{D}^{k+t+m-3-i}_{m}}[\vec{x}]\bigg)x\bigg) \right)]\bigg)
\\
 &+ \tbinom{m-1}{i}R^2\bigg(R\bigg(  L^\op_{\mathcal{B}^{k+t+m-2-i}_{m+1}}[\vec{x},\sum\limits^l_{t=0}\sum\limits^{m-2}_{i=0}\lambda^{(t)}\lambda^i \left(
\tbinom{m-2}{i}R\bigg(  {L^*}^\op_{\mathcal{B}^{k+t+m-3-i}_m}[x,\vec{x}]
+   {L^*}^\op_{\mathcal{A}^{k+t+m-3-i}_m}[x,\vec{x}]\bigg) \right. \\
&+\lambda \tbinom{m-3}{i}\bigg(  {L^*}^\op_{\mathcal{B}^{k+t+m-3-i}_m}[x,\vec{x}]\bigg)
+\lambda \tbinom{m-2}{i} \bigg(    {L^*}^\op_{\mathcal{C}^{k+t+m-3-i}_m}[x,\vec{x}]\bigg)\bigg)+\lambda^2 \tbinom{m-3}{i}\bigg(   {L^*}^\op_{\mathcal{D}^{k+t+m-4-i}_m}[x,\vec{x}]\bigg)\\
 &+ \tbinom{m-2}{i}R\bigg(R\bigg(          L^\op_{\mathcal{B}^{k+t+m-3-i}_{m}}[\vec{x}]
 +            L^\op_{\mathcal{A}^{k+t+m-3-i}_{m}}[\vec{x}]\bigg)x\bigg) + \lambda \tbinom{m-2}{i}R\bigg(\bigg(          L^\op_{\mathcal{B}^{k+t+m-3-i}_{m}}[\vec{x}]
 +            L^\op_{\mathcal{A}^{k+t+m-3-i}_{m}}[\vec{x}]\bigg)x\bigg)\\
& + \lambda \tbinom{m-2}{i}\bigg( R\bigg(          L^\op_{\mathcal{B}^{k+t+m-3-i}_{m}}[\vec{x}]
 +            L^\op_{\mathcal{A}^{k+t+m-3-i}_{m}}[\vec{x}]\bigg)x\bigg)
 + \lambda^2 \tbinom{m-3}{i} \bigg( \bigg(          L^\op_{\mathcal{B}^{k+t+m-3-i}_{m}}[\vec{x}]\bigg)x\bigg) \\
& + \lambda^2 \tbinom{m-2}{i}\bigg( \bigg(            L^\op_{\mathcal{C}^{k+t+m-3-i}_{m}}[\vec{x}]\bigg)x\bigg)
 \left. +\lambda^3 \tbinom{m-3}{i}\bigg( \bigg(              L^\op_{\mathcal{D}^{k+t+m-3-i}_{m}}[\vec{x}]\bigg)x\bigg) \right)]
 \\
& +        L^\op_{\mathcal{A}^{k+t+m-2-i}_{m+1}}[\vec{x},\sum\limits^l_{t=0}\sum\limits^{m-2}_{i=0}\lambda^{(t)}\lambda^i \left(
\tbinom{m-2}{i}R\bigg(  {L^*}^\op_{\mathcal{B}^{k+t+m-3-i}_m}[x,\vec{x}]
+   {L^*}^\op_{\mathcal{A}^{k+t+m-3-i}_m}[x,\vec{x}]\bigg) \right. \\
&+\lambda \tbinom{m-3}{i}\bigg(  {L^*}^\op_{\mathcal{B}^{k+t+m-3-i}_m}[x,\vec{x}]\bigg)
+\lambda \tbinom{m-2}{i} \bigg(    {L^*}^\op_{\mathcal{C}^{k+t+m-3-i}_m}[x,\vec{x}]\bigg)\bigg)+\lambda^2 \tbinom{m-3}{i}\bigg(   {L^*}^\op_{\mathcal{D}^{k+t+m-4-i}_m}[x,\vec{x}]\bigg)\\
 &+ \tbinom{m-2}{i}R\bigg(R\bigg(          L^\op_{\mathcal{B}^{k+t+m-3-i}_{m}}[\vec{x}]
 +            L^\op_{\mathcal{A}^{k+t+m-3-i}_{m}}[\vec{x}]\bigg)x\bigg)
 + \lambda \tbinom{m-2}{i}R\bigg(\bigg(          L^\op_{\mathcal{B}^{k+t+m-3-i}_{m}}[\vec{x}]
 +            L^\op_{\mathcal{A}^{k+t+m-3-i}_{m}}[\vec{x}]\bigg)x\bigg)\\
& + \lambda \tbinom{m-2}{i}\bigg( R\bigg(          L^\op_{\mathcal{B}^{k+t+m-3-i}_{m}}[\vec{x}]
 +            L^\op_{\mathcal{A}^{k+t+m-3-i}_{m}}[\vec{x}]\bigg)x\bigg)
 + \lambda^2 \tbinom{m-3}{i} \bigg( \bigg(          L^\op_{\mathcal{B}^{k+t+m-3-i}_{m}}[\vec{x}]\bigg)x\bigg) \\
 \end{aligned}
 \]
 \[
 \begin{aligned}
& + \lambda^2 \tbinom{m-2}{i}\bigg( \bigg(            L^\op_{\mathcal{C}^{k+t+m-3-i}_{m}}[\vec{x}]\bigg)x\bigg)
 \left. +\lambda^3 \tbinom{m-3}{i}\bigg( \bigg(              L^\op_{\mathcal{D}^{k+t+m-3-i}_{m}}[\vec{x}]\bigg)x\bigg) \right)]\bigg)x\bigg)
\\
 &+ \lambda \tbinom{m-1}{i}R^2\bigg(\bigg(  L^\op_{\mathcal{B}^{k+t+m-2-i}_{m+1}}[\vec{x},\sum\limits^l_{t=0}\sum\limits^{m-2}_{i=0}\lambda^{(t)}\lambda^i \left(
\tbinom{m-2}{i}R\bigg(  {L^*}^\op_{\mathcal{B}^{k+t+m-3-i}_m}[x,\vec{x}]
+   {L^*}^\op_{\mathcal{A}^{k+t+m-3-i}_m}[x,\vec{x}]\bigg) \right. 
\\
&+\lambda \tbinom{m-3}{i}\bigg(  {L^*}^\op_{\mathcal{B}^{k+t+m-3-i}_m}[x,\vec{x}]\bigg)
+\lambda \tbinom{m-2}{i} \bigg(    {L^*}^\op_{\mathcal{C}^{k+t+m-3-i}_m}[x,\vec{x}]\bigg)\bigg)
+\lambda^2 \tbinom{m-3}{i}\bigg(   {L^*}^\op_{\mathcal{D}^{k+t+m-4-i}_m}[x,\vec{x}]\bigg)\\
 &+ \tbinom{m-2}{i}R\bigg(R\bigg(          L^\op_{\mathcal{B}^{k+t+m-3-i}_{m}}[\vec{x}]
 +            L^\op_{\mathcal{A}^{k+t+m-3-i}_{m}}[\vec{x}]\bigg)x\bigg)
 + \lambda \tbinom{m-2}{i}R\bigg(\bigg(          L^\op_{\mathcal{B}^{k+t+m-3-i}_{m}}[\vec{x}]
 +            L^\op_{\mathcal{A}^{k+t+m-3-i}_{m}}[\vec{x}]\bigg)x\bigg)\\
& + \lambda \tbinom{m-2}{i}\bigg( R\bigg(          L^\op_{\mathcal{B}^{k+t+m-3-i}_{m}}[\vec{x}]
 +            L^\op_{\mathcal{A}^{k+t+m-3-i}_{m}}[\vec{x}]\bigg)x\bigg)
 + \lambda^2 \tbinom{m-3}{i} \bigg( \bigg(          L^\op_{\mathcal{B}^{k+t+m-3-i}_{m}}[\vec{x}]\bigg)x\bigg) \\
& + \lambda^2 \tbinom{m-2}{i}\bigg( \bigg(            L^\op_{\mathcal{C}^{k+t+m-3-i}_{m}}[\vec{x}]\bigg)x\bigg)
 \left. +\lambda^3 \tbinom{m-3}{i}\bigg( \bigg(              L^\op_{\mathcal{D}^{k+t+m-3-i}_{m}}[\vec{x}]\bigg)x\bigg) \right)]
\\
& +        L^\op_{\mathcal{A}^{k+t+m-2-i}_{m+1}}[\vec{x},\sum\limits^l_{t=0}\sum\limits^{m-2}_{i=0}\lambda^{(t)}\lambda^i \left(
\tbinom{m-2}{i}R\bigg(  {L^*}^\op_{\mathcal{B}^{k+t+m-3-i}_m}[x,\vec{x}]
+   {L^*}^\op_{\mathcal{A}^{k+t+m-3-i}_m}[x,\vec{x}]\bigg) \right. \\
&+\lambda \tbinom{m-3}{i}\bigg(  {L^*}^\op_{\mathcal{B}^{k+t+m-3-i}_m}[x,\vec{x}]\bigg)
+\lambda \tbinom{m-2}{i} \bigg(    {L^*}^\op_{\mathcal{C}^{k+t+m-3-i}_m}[x,\vec{x}]\bigg)\bigg)+\lambda^2 \tbinom{m-3}{i}\bigg(   {L^*}^\op_{\mathcal{D}^{k+t+m-4-i}_m}[x,\vec{x}]\bigg)\\
 &+ \tbinom{m-2}{i}R\bigg(R\bigg(          L^\op_{\mathcal{B}^{k+t+m-3-i}_{m}}[\vec{x}]
 +            L^\op_{\mathcal{A}^{k+t+m-3-i}_{m}}[\vec{x}]\bigg)x\bigg) + \lambda \tbinom{m-2}{i}R\bigg(\bigg(          L^\op_{\mathcal{B}^{k+t+m-3-i}_{m}}[\vec{x}]
 +            L^\op_{\mathcal{A}^{k+t+m-3-i}_{m}}[\vec{x}]\bigg)x\bigg)\\
& + \lambda \tbinom{m-2}{i}\bigg( R\bigg(          L^\op_{\mathcal{B}^{k+t+m-3-i}_{m}}[\vec{x}]
 +            L^\op_{\mathcal{A}^{k+t+m-3-i}_{m}}[\vec{x}]\bigg)x\bigg)
 + \lambda^2 \tbinom{m-3}{i} \bigg( \bigg(          L^\op_{\mathcal{B}^{k+t+m-3-i}_{m}}[\vec{x}]\bigg)x\bigg) \\
& + \lambda^2 \tbinom{m-2}{i}\bigg( \bigg(            L^\op_{\mathcal{C}^{k+t+m-3-i}_{m}}[\vec{x}]\bigg)x\bigg)
 \left. +\lambda^3 \tbinom{m-3}{i}\bigg( \bigg(              L^\op_{\mathcal{D}^{k+t+m-3-i}_{m}}[\vec{x}]\bigg)x\bigg) \right)]\bigg)x\bigg)
\\
& + \lambda \tbinom{m-1}{i}R\bigg( R\bigg(  L^\op_{\mathcal{B}^{k+t+m-2-i}_{m+1}}[\vec{x},\sum\limits^l_{t=0}\sum\limits^{m-2}_{i=0}\lambda^{(t)}\lambda^i \left(
\tbinom{m-2}{i}R\bigg(  {L^*}^\op_{\mathcal{B}^{k+t+m-3-i}_m}[x,\vec{x}]
+   {L^*}^\op_{\mathcal{A}^{k+t+m-3-i}_m}[x,\vec{x}]\bigg) \right. \\
&+\lambda \tbinom{m-3}{i}\bigg(  {L^*}^\op_{\mathcal{B}^{k+t+m-3-i}_m}[x,\vec{x}]\bigg)
+\lambda \tbinom{m-2}{i} \bigg(    {L^*}^\op_{\mathcal{C}^{k+t+m-3-i}_m}[x,\vec{x}]\bigg)\bigg)+\lambda^2 \tbinom{m-3}{i}\bigg(   {L^*}^\op_{\mathcal{D}^{k+t+m-4-i}_m}[x,\vec{x}]\bigg)\\
& + \tbinom{m-2}{i}R\bigg(R\bigg(          L^\op_{\mathcal{B}^{k+t+m-3-i}_{m}}[\vec{x}]
 +            L^\op_{\mathcal{A}^{k+t+m-3-i}_{m}}[\vec{x}]\bigg)x\bigg) + \lambda \tbinom{m-2}{i}R\bigg(\bigg(          L^\op_{\mathcal{B}^{k+t+m-3-i}_{m}}[\vec{x}]
 +            L^\op_{\mathcal{A}^{k+t+m-3-i}_{m}}[\vec{x}]\bigg)x\bigg)\\
& + \lambda \tbinom{m-2}{i}\bigg( R\bigg(          L^\op_{\mathcal{B}^{k+t+m-3-i}_{m}}[\vec{x}]
 +            L^\op_{\mathcal{A}^{k+t+m-3-i}_{m}}[\vec{x}]\bigg)x\bigg)
 + \lambda^2 \tbinom{m-3}{i} \bigg( \bigg(          L^\op_{\mathcal{B}^{k+t+m-3-i}_{m}}[\vec{x}]\bigg)x\bigg) \\
& + \lambda^2 \tbinom{m-2}{i}\bigg( \bigg(            L^\op_{\mathcal{C}^{k+t+m-3-i}_{m}}[\vec{x}]\bigg)x\bigg)
 \left. +\lambda^3 \tbinom{m-3}{i}\bigg( \bigg(              L^\op_{\mathcal{D}^{k+t+m-3-i}_{m}}[\vec{x}]\bigg)x\bigg) \right)]
 \\
& +        L^\op_{\mathcal{A}^{k+t+m-2-i}_{m+1}}[\vec{x},\sum\limits^l_{t=0}\sum\limits^{m-2}_{i=0}\lambda^{(t)}\lambda^i \left(
\tbinom{m-2}{i}R\bigg(  {L^*}^\op_{\mathcal{B}^{k+t+m-3-i}_m}[x,\vec{x}]
+   {L^*}^\op_{\mathcal{A}^{k+t+m-3-i}_m}[x,\vec{x}]\bigg) \right. \\
&+\lambda \tbinom{m-3}{i}\bigg(  {L^*}^\op_{\mathcal{B}^{k+t+m-3-i}_m}[x,\vec{x}]\bigg)
+\lambda \tbinom{m-2}{i} \bigg(    {L^*}^\op_{\mathcal{C}^{k+t+m-3-i}_m}[x,\vec{x}]\bigg)\bigg)+\lambda^2 \tbinom{m-3}{i}\bigg(   {L^*}^\op_{\mathcal{D}^{k+t+m-4-i}_m}[x,\vec{x}]\bigg)\\
& + \tbinom{m-2}{i}R\bigg(R\bigg(          L^\op_{\mathcal{B}^{k+t+m-3-i}_{m}}[\vec{x}]
 +            L^\op_{\mathcal{A}^{k+t+m-3-i}_{m}}[\vec{x}]\bigg)x\bigg)
 + \lambda \tbinom{m-2}{i}R\bigg(\bigg(          L^\op_{\mathcal{B}^{k+t+m-3-i}_{m}}[\vec{x}]
 +            L^\op_{\mathcal{A}^{k+t+m-3-i}_{m}}[\vec{x}]\bigg)x\bigg)\\
& + \lambda \tbinom{m-2}{i}\bigg( R\bigg(          L^\op_{\mathcal{B}^{k+t+m-3-i}_{m}}[\vec{x}]
 +            L^\op_{\mathcal{A}^{k+t+m-3-i}_{m}}[\vec{x}]\bigg)x\bigg)
 + \lambda^2 \tbinom{m-3}{i} \bigg( \bigg(          L^\op_{\mathcal{B}^{k+t+m-3-i}_{m}}[\vec{x}]\bigg)x\bigg) \\
& + \lambda^2 \tbinom{m-2}{i}\bigg( \bigg(            L^\op_{\mathcal{C}^{k+t+m-3-i}_{m}}[\vec{x}]\bigg)x\bigg)
 \left. +\lambda^3 \tbinom{m-3}{i}\bigg( \bigg(              L^\op_{\mathcal{D}^{k+t+m-3-i}_{m}}[\vec{x}]\bigg)x\bigg) \right)]\bigg)x\bigg)
\end{aligned}
 \]
 \[
 \begin{aligned}
 &+ \lambda^2 \tbinom{m-2}{i} R\bigg( \bigg(  L^\op_{\mathcal{B}^{k+t+m-2-i}_{m+1}}[\vec{x},\sum\limits^l_{t=0}\sum\limits^{m-2}_{i=0}\lambda^{(t)}\lambda^i \left(
\tbinom{m-2}{i}R\bigg(  {L^*}^\op_{\mathcal{B}^{k+t+m-3-i}_m}[x,\vec{x}]
+   {L^*}^\op_{\mathcal{A}^{k+t+m-3-i}_m}[x,\vec{x}]\bigg) \right. \\
&+\lambda \tbinom{m-3}{i}\bigg(  {L^*}^\op_{\mathcal{B}^{k+t+m-3-i}_m}[x,\vec{x}]\bigg)
+\lambda \tbinom{m-2}{i} \bigg(    {L^*}^\op_{\mathcal{C}^{k+t+m-3-i}_m}[x,\vec{x}]\bigg)\bigg)+\lambda^2 \tbinom{m-3}{i}\bigg(   {L^*}^\op_{\mathcal{D}^{k+t+m-4-i}_m}[x,\vec{x}]\bigg)\\
 &+ \tbinom{m-2}{i}R\bigg(R\bigg(          L^\op_{\mathcal{B}^{k+t+m-3-i}_{m}}[\vec{x}]
 +            L^\op_{\mathcal{A}^{k+t+m-3-i}_{m}}[\vec{x}]\bigg)x\bigg) + \lambda \tbinom{m-2}{i}R\bigg(\bigg(          L^\op_{\mathcal{B}^{k+t+m-3-i}_{m}}[\vec{x}]
 +            L^\op_{\mathcal{A}^{k+t+m-3-i}_{m}}[\vec{x}]\bigg)x\bigg)\\
& + \lambda \tbinom{m-2}{i}\bigg( R\bigg(          L^\op_{\mathcal{B}^{k+t+m-3-i}_{m}}[\vec{x}]
 +            L^\op_{\mathcal{A}^{k+t+m-3-i}_{m}}[\vec{x}]\bigg)x\bigg)
 + \lambda^2 \tbinom{m-3}{i} \bigg( \bigg(          L^\op_{\mathcal{B}^{k+t+m-3-i}_{m}}[\vec{x}]\bigg)x\bigg) \\
& + \lambda^2 \tbinom{m-2}{i}\bigg( \bigg(            L^\op_{\mathcal{C}^{k+t+m-3-i}_{m}}[\vec{x}]\bigg)x\bigg)
 \left. +\lambda^3 \tbinom{m-3}{i}\bigg( \bigg(              L^\op_{\mathcal{D}^{k+t+m-3-i}_{m}}[\vec{x}]\bigg)x\bigg) \right)]\bigg)x\bigg) 
 \\
& + \lambda^2 \tbinom{m-1}{i}R\bigg( \bigg(  L^\op_{\mathcal{C}^{k+t+m-2-i}_{m+1}}[\vec{x},\sum\limits^l_{t=0}\sum\limits^{m-2}_{i=0}\lambda^{(t)}\lambda^i \left(
\tbinom{m-2}{i}R\bigg(  {L^*}^\op_{\mathcal{B}^{k+t+m-3-i}_m}[x,\vec{x}]
+   {L^*}^\op_{\mathcal{A}^{k+t+m-3-i}_m}[x,\vec{x}]\bigg) \right. \\
&+\lambda \tbinom{m-3}{i}\bigg(  {L^*}^\op_{\mathcal{B}^{k+t+m-3-i}_m}[x,\vec{x}]\bigg)
+\lambda \tbinom{m-2}{i} \bigg(    {L^*}^\op_{\mathcal{C}^{k+t+m-3-i}_m}[x,\vec{x}]\bigg)\bigg)+\lambda^2 \tbinom{m-3}{i}\bigg(   {L^*}^\op_{\mathcal{D}^{k+t+m-4-i}_m}[x,\vec{x}]\bigg)\\
& + \tbinom{m-2}{i}R\bigg(R\bigg(          L^\op_{\mathcal{B}^{k+t+m-3-i}_{m}}[\vec{x}]
 +            L^\op_{\mathcal{A}^{k+t+m-3-i}_{m}}[\vec{x}]\bigg)x\bigg)
 + \lambda \tbinom{m-2}{i}R\bigg(\bigg(          L^\op_{\mathcal{B}^{k+t+m-3-i}_{m}}[\vec{x}]
 +            L^\op_{\mathcal{A}^{k+t+m-3-i}_{m}}[\vec{x}]\bigg)x\bigg)\\
& + \lambda \tbinom{m-2}{i}\bigg( R\bigg(          L^\op_{\mathcal{B}^{k+t+m-3-i}_{m}}[\vec{x}]
 +            L^\op_{\mathcal{A}^{k+t+m-3-i}_{m}}[\vec{x}]\bigg)x\bigg)
 + \lambda^2 \tbinom{m-3}{i} \bigg( \bigg(          L^\op_{\mathcal{B}^{k+t+m-3-i}_{m}}[\vec{x}]\bigg)x\bigg) \\
& + \lambda^2 \tbinom{m-2}{i}\bigg( \bigg(            L^\op_{\mathcal{C}^{k+t+m-3-i}_{m}}[\vec{x}]\bigg)x\bigg)
 \left. +\lambda^3 \tbinom{m-3}{i}\bigg( \bigg(              L^\op_{\mathcal{D}^{k+t+m-3-i}_{m}}[\vec{x}]\bigg)x\bigg) \right)]\bigg)x\bigg)
\\
  &+\lambda^3 \tbinom{m-2}{i}R\bigg( \bigg(  L^\op_{\mathcal{D}^{k+t+m-2-i}_{m+1}}[\vec{x},\sum\limits^l_{t=0}\sum\limits^{m-2}_{i=0}\lambda^{(t)}\lambda^i \left(
\tbinom{m-2}{i}R\bigg(  {L^*}^\op_{\mathcal{B}^{k+t+m-3-i}_m}[x,\vec{x}]
+   {L^*}^\op_{\mathcal{A}^{k+t+m-3-i}_m}[x,\vec{x}]\bigg) \right. 
\\
&+\lambda \tbinom{m-3}{i}\bigg(  {L^*}^\op_{\mathcal{B}^{k+t+m-3-i}_m}[x,\vec{x}]\bigg)
+\lambda \tbinom{m-2}{i} \bigg(    {L^*}^\op_{\mathcal{C}^{k+t+m-3-i}_m}[x,\vec{x}]\bigg)\bigg)+\lambda^2 \tbinom{m-3}{i}\bigg(   {L^*}^\op_{\mathcal{D}^{k+t+m-4-i}_m}[x,\vec{x}]\bigg)\\
 &+ \tbinom{m-2}{i}R\bigg(R\bigg(          L^\op_{\mathcal{B}^{k+t+m-3-i}_{m}}[\vec{x}]
 +            L^\op_{\mathcal{A}^{k+t+m-3-i}_{m}}[\vec{x}]\bigg)x\bigg) + \lambda \tbinom{m-2}{i}R\bigg(\bigg(          L^\op_{\mathcal{B}^{k+t+m-3-i}_{m}}[\vec{x}]
 +            L^\op_{\mathcal{A}^{k+t+m-3-i}_{m}}[\vec{x}]\bigg)x\bigg)\\
& + \lambda \tbinom{m-2}{i}\bigg( R\bigg(          L^\op_{\mathcal{B}^{k+t+m-3-i}_{m}}[\vec{x}]
 +            L^\op_{\mathcal{A}^{k+t+m-3-i}_{m}}[\vec{x}]\bigg)x\bigg)
 + \lambda^2 \tbinom{m-3}{i} \bigg( \bigg(          L^\op_{\mathcal{B}^{k+t+m-3-i}_{m}}[\vec{x}]\bigg)x\bigg) \\
& + \lambda^2 \tbinom{m-2}{i}\bigg( \bigg(            L^\op_{\mathcal{C}^{k+t+m-3-i}_{m}}[\vec{x}]\bigg)x\bigg)
 \left. +\lambda^3 \tbinom{m-3}{i}\bigg( \bigg(              L^\op_{\mathcal{D}^{k+t+m-3-i}_{m}}[\vec{x}]\bigg)x\bigg) \right)]\bigg)x\bigg) ) {=} 0;\ m = |\vec{x}|{\geq}3.
\end{aligned}
\]

\subsection*{Acknowledgments}
The study was supported by a grant from the Russian Science Foundation No.~23-71-10005, \url{https://rscf.ru/project/23-71-10005/}. The author is grateful 
to V.~Yu. Gubarev for discussions and useful comments.

\end{document}